\pdfoutput=1
\documentclass[pdftex, noinfoline]{imsart}
 
\usepackage[english]{babel}
\usepackage{parskip}
\usepackage{amsmath,amsfonts,amsthm, amssymb}
\usepackage{graphicx} 
\usepackage{enumitem} 
\usepackage{bbm}
\usepackage{xcolor, latexsym} 
\usepackage{subfigure} 
\usepackage{bm}

\usepackage{titletoc}  

\usepackage[pdftex, 
colorlinks,citecolor=blue,urlcolor=blue,filecolor=blue,linkcolor=blue,backref=page,
bookmarks=true, breaklinks, hypertexnames=false, pdfborder={0 0 0 0}]{hyperref} 
\usepackage[normalem]{ulem} 
\usepackage{float}
\usepackage{placeins}
\usepackage{soul} 
\soulregister\cite7 
\soulregister\ref7
\soulregister\pageref7

\newcommand{\given}{\,|\,}

\newcommand{\Beta}{\text{Beta}}
\newcommand{\FDP}{\text{FDP}}
\newcommand{\FNP}{\text{FNP}}
\DeclareMathOperator*{\FDR}{\text{FDR}}
\DeclareMathOperator*{\FNR}{\text{FNR}}
\DeclareMathOperator*{\argmax}{arg\,max}

\DeclareMathOperator*{\B}{\text{\bf B}}
\DeclareMathOperator*{\G}{\text{\bf G}}
\newcommand{\T}{\mathsf{T}}
\newcommand{\fF}{\text{FD}}
\newcommand{\fT}{\text{TD}}
\newcommand{\adj}{\mathsf{adj}}

\allowdisplaybreaks

\theoremstyle{plain} 
\newtheorem{theorem}{Theorem}
\newtheorem{corollary}{Corollary}
\newtheorem{lemma}{Lemma}

\newtheorem{prop}{Proposition}
\newtheorem{remark}{Remark}

\endlocaldefs
\usepackage{natbib}
\usepackage{multirow}
\usepackage{chngpage}
\usepackage{siunitx}
\usepackage{booktabs}
\usepackage{makecell}
\usepackage[T1]{fontenc}
\usepackage{float}
\usepackage{grffile}
\usepackage{url} 
\usepackage{graphicx,psfrag,epsf}

\usepackage{bbm, dsfont}
\usepackage{mathrsfs} 
\usepackage[ruled,vlined]{algorithm2e}
\SetKwInput{KwInput}{Input}                
\SetKwInput{KwOutput}{Output}

\usepackage{xr-hyper}
\usepackage{titletoc}
\defcitealias{cast20MT}{CR20}

\begin{document}

\begin{frontmatter}

\title{Empirical Bayes large-scale multiple testing for high-dimensional binary outcome data}

\runtitle{Empirical Bayes multiple testing for sparse binary data}

\author{\fnms{Yu-Chien Bo} \snm{Ning}\ead[label=e1]{bycning@hsph.harvard.edu}}
 
 
\affiliation{Harvard University}

\address{
Harvard T.H. Chan School of Public Health\\
677 Huntington Ave\\
Boston, MA 02155\\
\printead{e1}}

\runauthor{Ning}

\begin{abstract}
This paper explores the multiple testing problem for sparse high-dimensional data with binary outcomes. 
We propose novel empirical Bayes multiple testing procedures based on a spike-and-slab posterior and then evaluate their performance in controlling the false discovery rate (FDR). A surprising finding is that the procedure using the default conjugate prior (namely, the $\ell$-value procedure) can be overly conservative in estimating the FDR. To address this, we introduce two new procedures that provide accurate FDR control.
Sharp frequentist theoretical results are established for these procedures, and numerical experiments are conducted to validate our theory in finite samples. To the best of our knowledge, we obtain the first {\it uniform} FDR control result in multiple testing for high-dimensional data with binary outcomes under the sparsity assumption. 
\end{abstract}

\begin{keyword}[class=MSC]
\kwd[Primary ]{62G10, 62G20}
\end{keyword}

\begin{keyword}
\kwd{Binomial distribution}
\kwd{empirical Bayes}
\kwd{false discovery rate}
\kwd{multiple testing}
\kwd{sparse binary data}
\kwd{spike-and-slab posterior}
\end{keyword}

\end{frontmatter}


\section{Introduction} \label{sec:intro}

Large-scale multiple testing problems frequently arise in modern statistical applications, such as in astronomy, biology, medical studies, etc. 
Despite their broad popularity in practice, previous theoretical research has focused almost exclusively on Gaussian models.
In this paper, we explore the multiple testing problem for high-dimensional data with binary outcomes. 
The dataset consists of $m \times n$ binary outcomes, denoted as
$\mathcal{D} = \{Z_{ij}, i = 1, \dots, m, \ j = 1, \dots, n\}$, $m, n \geq 1$, and the model is given by 
\begin{align}
\label{model}
Z_{ij} \overset{ind}{\sim} \text{Ber}(\theta_j), \quad i = 1, \dots, m, \ j = 1, \dots, n,
\end{align}
where $m$ is the number of samples, and $\theta = (\theta_1, \dots, \theta_n)'$ is an $n$-dimensional unknown $n$ vector with each $\theta_j \in [0, 1]$ for all $j \in \{1, \dots, n\}$. 
Let $X_j = \sum_{i = 1}^m Z_{ij}$, then $X_j \overset{ind}{\sim} \text{Bin}(m, \theta_j)$. 

High-dimensional binary data are prevalent in statistical and machine learning applications.
For example, in genomics, DNA and RNA data are often translated to binary format before conducting data analysis.
In the crowdsourcing problem from machine learning, $m$ workers are asked to provide binary label assignments for $n$ (typically $n \gg m$) objects to determine the true labels \citep{gao16, butucea18}. For instance, the Space Warps project \citep{marshal15}, which initiated in 2013, asked 37,000 citizen scientists to participate in classifying 11 million images to identify gravitational lenses over an eight-month period. Data collected from this project contain a large number of binary outcomes, each indicating the existence of gravitational lenses in an image identified by a citizen scientist. 

Given the high-dimensional nature of these examples, we adopt the following sparsity assumption: denote $\theta_0$ as the true value of $\theta$, 
assuming $\theta_0 \in l_0[s_n]$ for $s_n \leq n$ such that
\begin{align}
\label{ell0}
l_0[s_n] = \{\theta \in {[0, 1]^n},\ \# \{j: \theta_j \neq 1/2\} \leq s_n\}.
\end{align}
The value $1/2$ chosen in \eqref{ell0} corresponds to a scenario where data are generated from a random stochastic Bernoulli process where the parameter is 1/2, indicating the absence of signals. This is analogous to the sparse Gaussian sequence model when the mean is zero. 

The multiple testing problem considered in this paper involves simultaneously testing the following hypotheses:
\begin{align}
\label{testing}
H_{0j}: \theta_{0,j} = 1/2 \quad \text{versus} \quad H_{1j}: \theta_{0,j} \neq 1/2, \ \  j = 1, \dots, n.
\end{align}
We construct multiple testing procedures using an empirical Bayes posterior, where sparsity is imposed through a spike-and-slab prior. In this paper, we will provide an in-depth frequentist theoretical analysis of these procedures.

The literature on Bayesian high-dimensional analysis has developed rapidly in recent years. However, theoretical works on these methods have predominantly focused on Gaussian settings \citep[e.g.,][]{johnstone04, cast12, cast15, martin17, ning20, rockova18, bai22, ray22}. In the realm of multiple testing, \citet[][hereafter, CR20]{cast20MT} studied the Gaussian sequence model and provided sharp theoretical results for multiple testing procedures based on empirical Bayes posteriors, demonstrating their superior performance in controlling the false discovery rate (FDR). More recently, \citet{abraham21} extended these findings to accommodate models with sub-Gaussian noise. These studies offer theoretical insights into the FDR control of empirical Bayes multiple testing procedures.

Besides Gaussian models, progress in high-dimensional models dealing with discrete data has lagged. \citet{mukherjee15} initiated the study of the testing problem for the sparse binary regression model which involves testing a global null hypothesis versus sparse alternatives. Their study revealed a key difference in the signal detection boundary between the binary regression model and the Gaussian sequence model for a small $m$.
Indeed, we also discover several differences between the analysis of the sparse binomial model and the sparse Gaussian model using the empirical Bayes method; see Section \ref{sec:contribution} for comments on this point. 
In the following sections, let us introduce the prior and the empirical Bayes approach in detail.

\subsection{The spike-and-slab posterior for the binomial model} 
The spike-and-slab prior given as follows:
\begin{align}
\label{prior}
\theta \given w \sim \bigotimes_{j=1}^n \{ (1-w) \delta_{1/2} + w \gamma\},
\end{align}
where $w \in (0, 1)$ is a weight, $\delta_{1/2}$ is the Dirac measure at $1/2$, which is the value at the null hypothesis in \eqref{testing}, and $\gamma \sim \text{Beta}(\alpha, \alpha)$, where $\alpha \in \mathbb{R}^+$. 
If choosing $\alpha = 1$, then $\gamma \sim \text{Unif}(0, 1)$, the uniform distribution on $[0, 1]$.
We hence refer the prior in \eqref{prior} with the uniform distribution as the {\it spike-and-uniform slab} prior. 
By combining the model in \eqref{model} and the spike-and-uniform slab prior using Bayes' theorem, the posterior is given by 
\begin{align}
\label{posterior}
 P^\pi(\theta \given X, w)  = \bigotimes_{j=1}^n \{\ell(X_j) \delta_{1/2} + (1-\ell(X_j)) \mathcal{G}_{X_j} \},
\end{align}
where $\mathcal{G}_{X_j} = \Beta(\theta_j; X_j + 1, m - X_j + 1)$ and
\begin{align}
\label{lval}
\ell(x) = \ell(x; w) = P^{\pi}(\theta = 1/2 \given X = x, w) = \frac{(1-w) b(x)}{(1-w) b(x) + wg(x)},
\end{align}
where $b = b_{1/2}$, $b_\theta = \text{Bin}(m, \theta)$ is the binomial distribution with parameters $m$ and $\theta$, and 
$g(x) = \int b_\theta(x) \gamma(\theta) d\theta = (m + 1)^{-1}$ is a constant. 

We choose $\alpha = 1$ in the prior of $\gamma$ because, for other values, $g(x)$ is no longer constant. This makes the posterior, and hence the multiple testing procedure, not only more difficult to analyze but also offers no additional benefit, as the resulting procedure does not guarantee correct FDR control; see Section \ref{sec:disc} for further discussion. 

\subsection{The empirical Bayes approach} 
\label{sec:eb}

The empirical Bayes approach first estimates $w$ and then plugs in its estimated value into the posterior. 
Let $L(w)$ be the logarithm of the marginal density of $X$ given $w$ given by 
\begin{align}
\label{eqn:marginal-post}
L(w) = \sum_{j = 1}^n \log b(X_j) + \sum_{j=1}^n \log (1 +w \beta(X_j)),
\end{align}
where $\beta(u) = (g/b)(u) - 1$.
We estimate $w$ by solving
\begin{align}
\label{w-hat}
\hat w = \argmax_{w \in [1/n, 1]}  L(w).
\end{align}
Then $\hat w$ is the marginal maximum likelihood estimator (MMLE). 
The lower bound $1/n$ in this optimization is imposed to prevent $\hat w$ from being too small, which is crucial to effectively control the FDR when all signals are very close to 1/2.
The same constraint is imposed by \citetalias{cast20MT} for the empirical Bayes multiple testing approach based on the Gaussian sequence model for a similar reason.
The solution of $\hat w$ is determined by the score function given by
\begin{align}
\label{score-fn}
{S}(w) = \frac{\partial}{\partial w} L(w)
= \sum_{j=1}^n \beta(X_j, w), \quad
\beta(x, w) = \frac{\beta(x)}{1+w\beta(x)}.
\end{align}
If ${S}(1) < 0$ and ${S}(1/n) > 0$, then a solution for $S(w) = 0$ exists and is unique, as $\beta(x, w)$ is monotone decreasing on $x \in [m/2, m]$ due to $b(x)$ is monotone increasing on $[m/2, m]$ and $g(x)$ is a constant. Otherwise, the solution of \eqref{w-hat} will be at the boundary of the interval $[1/n, 1]$.
In Sections \ref{sec:bound-w-hat} and \ref{sec:bound-m1-mtilde} of the supplementary material \citep{ning23supp}, we will conduct an in-depth analysis for the score function in both scenarios.

\subsection{Our contribution}  
\label{sec:contribution}

This paper presents novel methodologies as well as sharp theoretical results for the multiple testing problem stated in \eqref{testing}.
Specifically, three multiple testing procedures are proposed: 
\begin{itemize}
\item The $\ell$-value procedure, a.k.a. the local FDR introduced by Bradley Efron \citep[see Chapter 5 of][]{efron10}. However, to avoid confusion between FDR and the procedure itself, we follow \citetalias{cast20MT} and refer to it as the $\ell$-value;
\item The adjusted $\ell$-value ($\adj\ell$-value) procedure, which is less conservative comparing to the $\ell$-value procedure (see Section \ref{sec:cl-q});
\item The $q$-value procedure, which alters the significance region used in the $\ell$-value procedure, so that it provides a global measure for FDR control, instead of the local FDR in the $\ell$-value procedure (see Section \ref{sec:cl-q}). 
\end{itemize}
Our theoretical results include:
\begin{itemize}
\item we show that the $\ell$-value procedure using the spike-and-uniform prior is too conservative but the $\adj\ell$-value procedure as well as the $q$-value procedure allow a correct {\it uniform} FDR control for arbitrary sparse signals. In particular, the $q$-value procedure can attain the exact target level of FDR control for large signals (Theorem \ref{thm:unif-fdr-q-Cl}).

\item we demonstrate that both the $\adj\ell$-value and the $q$-value procedures can effectively control the multiple testing risk (the sum of FDR and FNR (false negative rate) as defined in \eqref{risk}) for large signals, while the $\ell$-value procedure cannot (Theorems \ref{thm:fdr-qval} \& \ref{thm:fnr-Clval-qval} and Lemma \ref{lem:fnr-lval}). 

\item we show all three procedures are thresholding-based procedures and analyze their corresponding thresholds. We also obtain the lower bound of the testing boundary for all thresholding-based multiple testing procedures, given by $\sqrt{\log(n/s_n)/(2m)}$, for $m \gg (\log n)^2$ (see Proposition \ref{prop:lowerbound}).
\end{itemize}

Our analysis for the above results is inspired by the work of \citet{johnstone04} and \citetalias{cast20MT} on the Gaussian sequence model.
However, a notable difficulty in our setting is the need for tight bounds shaper than the Chernoff bound for both centered and non-centered binomial distributions, in order to control the MMLE $\hat w$ in a neighborhood of $s_n/n$ (up to some constant), which is an essential step for obtaining uniform FDR control.
To this end, we derive tight bounds by first approximating the binomial distribution with a Gaussian and then carefully controlling the approximation error. This task is especially delicate for non-centered binomial distributions, which are asymmetric and thus can be poorly approximated by a Gaussian. Relevant inequalities for various binomial distributions are collected in Section \ref{sec:bound-binom} of \citet{ning23supp}, and they could be of independent interest for other studies.

Another challenge in our analysis is to keep the assumption on $m$ as minimal as possible. To achieve this, we carefully study the behavior of certain quantities involving small signals that are close to, but not exactly, $1/2$. This allows us to control the error contributions from these signals so that they do not dominate the upper bound of the FDR without imposing a stronger assumption on $m$. Relevant proofs can be found in Sections \ref{sec:pf-unif-fdr}, \ref{sec:pf-fdr-fnr}, and \ref{sec:control-small-signals} of \citet{ning23supp}.

Before concluding this section, it is worth mentioning that the aim of this paper is to conduct a frequentist analysis for FDR control of Bayesian multiple testing procedures. Our interest lies in the FDR rather than the Bayes FDR (BFDR). The BFDR is defined as the FDR integrated over the prior distribution, given by
$$
\text{BFDR}(\T; w, \gamma) = \int_{\theta \in [0, 1]^n} \text{FDR}(\theta, \T; w, \gamma) d\Pi(\theta).
$$
The last display suggests that while controlling the FDR ensures control of the BFDR, the reverse is not necessarily true. Readers can refer to Proposition 1 in \citetalias{cast20MT} for a formal proof of this point. More importantly, controlling the BFDR does not provide information about how the FDR behaves under arbitrary sparsity patterns for $\theta_0$.

\subsection{Outline of this paper}
 
The rest of the paper proceeds as follows: Section \ref{sec:lval} studies the $\ell$-value procedure and its thresholding rule. 
Section \ref{sec:cl-q} introduces the $\adj\ell$-value and $q$-value procedures and compares their thresholds with that of the $\ell$-value procedure. 
Sections \ref{sec:unif-fdr} presents the uniform FDR control result for the three empirical Bayes multiple testing procedures.  The lower bound of the testing boundary is provided in Section \ref{sec:fdr-fnr-large}; this section also studies the FNR control and the multiple testing risk for the three procedures.
Numerical experiments are conducted in Section \ref{sec:sim}.
The conclusion and discussion are given in Section \ref{sec:disc}.
All the proofs are included in the supplementary material in \citet{ning23supp}. 

\subsection{Notation}
Let $b_\theta(x) = \text{Bin}(m, \theta)$ be the density function of the binomial distribution with parameters $m$ and $\theta$ and $b(x) = b_{1/2}(x)$.
Denote the upper tail probability of $b_\theta(x)$ as $\bar \B_\theta(u) = \sum_{x = u}^m b_\theta(x)$ and similarly, $\bar \B(u) = \sum_{x= u}^m b(x)$.
Let  $\B_\theta(u) = 1 - \bar \B_\theta(u)$ and $\B(u) = 1 - \bar \B(u)$. 
The symbol $\phi(x)$ stands for the standard normal distribution, and $\Phi(x)$ is the cdf.
For any cdf function, say $F(x)$, let $\bar F(x) = 1 - F(x)$. For any two real numbers $a_1$ and $a_2$, let $a_1 \vee a_2 = \max\{a_1, a_2\}$, $a_1 \wedge a_2 = \min\{a_1, a_2\}$, and $a_1 \lesssim a_2$ as $a_1 \leq C a_2$ for some constant $C$. For two sequences $c_n$ and $d_n$ depending on $n$, $c_n \ll d_n$ stands for $c_n/d_n \to 0$ as $n \to \infty$, 
$c_n \asymp d_n$ stands for that there exists constants $a, a' > 0$ such that $a c_n \leq d_n \leq a'c_n$,
and $c_n \sim d_n$ stands for $c_n - d_n = o(c_n)$, where $o(1)$ is a deterministic sequence going to 0 with $n$. 
The indicator function is denoted by $\mathbbm{1}\{\cdot\}$.
For a parameter $\theta$ and a set $\T$ containing $n$ test functions such as $\T = (\T_1, \dots, \T_n)$, 
the false discovery rate and the false negative rate are defined as
$
\FDR(\theta,\T) = \mathbb{E}_{\theta_0}\FDP(\theta, \T) 
\ \text{and} \ 
\FNR(\theta, \T) = \mathbb{E}_{\theta_0}\FNP(\theta, \T),
$
where $\FDP(\theta,\T)$ and $\FNP(\theta, \T)$ are the false discovery proportion and the false negative proportion respectively given by 
\begin{align}
\label{fdp-fnp}
\FDP(\theta, \T) = \frac{\sum_{j=1}^n \mathbbm{1}\{\theta_j = 1/2\} \T_j}{1 \vee \sum_{j=1}^n \T_j},
\quad
\FNP(\theta, \T) = \frac{\sum_{j=1}^n \mathbbm{1}\{\theta_j \neq 1/2\}(1-\T_j)}{1 \vee \sum_{j=1}^n \mathbbm{1}\{\theta_j = 1/2\}}.
\end{align}

\section{The $\ell$-value procedure}
\label{sec:lval}

Before introducing the procedure, we first define the $\ell$-value. An $\ell$-value, also known as the local false discovery rate \citep{efron04}, is the probability that the null hypothesis is true conditional on the test statistics equals to the observed value.
Based on this definition, $\ell(x) = \ell(x; w)$ in \eqref{lval} is the $\ell$-value.

The $\ell$-value procedure is constructed as follows: first, estimate the MMLE $\hat w$ by solving \eqref{w-hat}. 
Second, compute $\hat \ell(x) = \hat \ell(x; \hat w)$ by substituting $\hat w$.
Last, determine a cutoff value $t \in (0, 1)$ and choose to reject or accept a null hypothesis based on whether $\hat \ell(x) \leq t$ or $\hat \ell(x) > t$.
A summary of this procedure is given in Algorithm \ref{algo:l-val}.

\begin{center}
\begin{minipage}{.48\linewidth}
\begin{algorithm}[H]
\label{algo:l-val}
\caption{The $\ell$-value procedure}\label{algo:lval}
\KwData{$\mathcal{D} = \{Z_{ij}, i = 1, \dots, m, j = 1, \dots, n\}$}
{\bf Input:} A pre-specified value $t \in (0, 1)$
\begin{enumerate}
\item[0.] Compute $X_j = \sum_{i=1}^m Z_{ij}$;
\item[1.] Compute $\hat w$ in \eqref{w-hat};
\item[2.] Evaluate $\hat \ell_j = \ell(X_j; \hat w)$ in \eqref{lval};
\item[3.] Obtain $\T_j^{\ell} = \mathbbm{1}\{\hat \ell_j \leq t\}$;
\end{enumerate}
{\bf Output:} $\T_1^{\ell}, \dots, \T_n^{\ell}$.\\[.2cm]
\end{algorithm}
\end{minipage}
\end{center}

\subsection{Analyzing the threshold of the $\ell$-value procedure}
\label{sec:th-lval}
In the following lemma, we show that the $\ell$-value procedure is a thresholding-based procedure.
\begin{lemma}
\label{lem:th-l}
For a fixed $t \in (0, 1)$ and $w \in (0, 1)$,  
let $r(w,t) = \frac{wt}{(1-w)(1-t)}$, 
consider the test function $\T^\ell = \mathbbm{1}\{\ell(x;  w, g) \leq t\}$ where $\ell(\cdot)$ is given in \eqref{lval}. Then,
$\T^\ell = \mathbbm{1}\{|x - m/2| \geq mt^\ell_n\}$, 
where for 
\begin{align}
\label{eqn:etal}
\eta^\ell(\cdot) = \frac{1}{m}(b/g)^{-1}(\cdot),
\end{align}
we have $t^{\ell}_n := t^{\ell}_n ( w, t) = \eta^\ell ( r( w, t)) - 1/2$.
\end{lemma}

In the next lemma, we derive an asymptotic bound for $\eta_\ell(\cdot)$. The non-asymptotic upper and lower bounds of this quantity are also obtained in Lemma \ref{lem:bound-etal} of Section \ref{sec:threshold} in \citet{ning23supp}.

\begin{lemma}
\label{lem:th-lval}
For $t_n^\ell(w, t)$ defined in Lemma \ref{lem:th-l} with a fixed $t\in (0, 1)$ and some $w \in (0, 1)$, 
let $r(w,t) = \frac{wt}{(1-w)(1-t)}$,
if $\log^2 (1/(r(w, t)))/m \to 0$ as $m \to \infty$, then
\begin{align}
\label{eqn:th-lval}
t_n^\ell(w, t) \sim \sqrt{\frac{1}{2m} 
\left(\log \left(\frac{1}{r(w,t)}\right) + \log \left( \frac{\sqrt{2}(1+m)}{\sqrt{\pi m}} \right)\right)}.
\end{align}
\end{lemma}

\begin{remark}
\label{rmk-1}
Consider $w \asymp s_n/n$, which holds for the MMLE $\hat w$ when a number of large signals is present (see Lemma \ref{lem:bound-what} in Section \ref{sec:bound-w-hat} of \citet{ning23supp} for details), then 
$$
t_n^\ell(w, t) \sim \sqrt{\frac{1}{2m} \left(\log\left(\frac{n}{s_n}\right) + \log \left( \frac{\sqrt{2}(1+m)}{\sqrt{\pi m}} \right)\right)},
$$ 
which is larger than $\sqrt{ \log(n/s_n)/(2m)}$, the lower bound among all thresholding-based multiple testing procedures as established in Proposition \ref{prop:lowerbound}. If $n/s_n = m^\alpha$ for a fixed $\alpha < 1/2$, then $t_n^\ell(w, t) \sim \sqrt{(\alpha+1/2) \log m/(2m)}$, which misses the optimal constant when $\alpha$ is smaller than $1/2$. In general, the threshold of the $\ell$-value can be sub-optimal, which leads the $\ell$-value procedure to be overly conservative for controlling the FDR, as shown in Lemma \ref{lem:unif-fdr-l}.
\end{remark}

\section{The $\adj\ell$-value and $q$-value procedures}
\label{sec:cl-q}

In this section, we introduce the $\adj\ell$-value and $q$-value procedures and then analyze their thresholds respectively. 
\subsection{Introducing the $\adj\ell$-value and $q$-value procedures}

Given that the thresholding rule of the $\ell$-value procedure can be sub-optimal as discussed in Remark \ref{rmk-1}. We introduce the adjusted $\ell$-value ($\adj \ell$-value) that can improve the threshold of the $\ell$-value through replacing $g(x)$ in \eqref{lval} with $\sqrt{2/(\pi m)}(1+m) g(x)$. The $\adj\ell$-value is defined as follows:
\begin{align}
\label{Clval}
\adj\ell(x; w, g) = \frac{(1-w) b(x)}{(1-w)b(x) + w \sqrt{\frac{2}{\pi m}}(1+m) g(x)}.
\end{align}
In fact, one can also view the $\adj\ell$-value as the $\ell$-value with the slab density of the spike-and-slab prior chosen as $\gamma \sim C_m \times \text{Unif}[0, 1]$ with $C_m = \sqrt{2/(\pi m)}(m+1)$.

Next, we introduce the $q$-value \citet{storey03}, which is the probability that the null hypothesis is true conditionally on the test statistics being larger than the observed value. This allows the spike-and-uniform slab prior to be maintained. Let $Y = X - m/2$ and $y =  x - m/2$, the $q$-value is defined as
\begin{align}
\label{qval}
q(x; w) = P^{\pi}(\theta = 1/2 \given |Y| \geq |y|, w) 
= \frac{(1-w) \bar{\mathbf{B}}(m/2 + |y|)}{(1-w) \bar{\mathbf{B}}(m/2 + |y|) + w \bar\G(m/2 + |y|)},
\end{align}
where $\bar \G(u) = \sum_{x = u}^m g(u)$.
Since $g(\cdot) = (m+1)^{-1}$, $\bar \G(m/2 + |y|) = (m/2 - |y| +1)/(1+m)$.

The $\adj\ell$-value procedure is constructed as follows: first, obtain $\hat w$ from \eqref{w-hat} (note that this step is the same as the $\ell$-value procedure.) Next, evaluate $\adj\ell(x; \hat w, g)$. Finally, choose the cutoff value $t \in (0, 1)$ to reject or accept each null hypothesis.
The $q$-value procedure is constructed in a similar way, except that we replace $\adj\ell(x; \hat w, g)$ in the second step with $q(x; \hat w)$. Details of the two procedures are given in Algorithm \ref{algo:Cl-q-val}.

\begin{algorithm}
\label{algo:Cl-q-val}
\vspace{.05cm}
\caption{The $\adj\ell$- and $q$-value procedures}\label{algo:lval}
\begin{minipage}{0.45\textwidth}
{\bf \normalsize The $\adj\ell$-value procedure:} \\[0.2cm]
{\bf Data:} $\mathcal{D} = \{Z_{ij}, i = 1, \dots, m, j = 1, \dots, n\}$\\
{\bf Input:} A pre-specified value $t \in (0, 1)$,
\begin{enumerate}
\item[a0.] Compute $X_j = \sum_{i=1}^m Z_{ij}$;
\item[a1.] Compute $\hat w$ as in \eqref{w-hat};
\item[a2.] Evaluate $\widehat{\adj\ell}_j = \adj\ell(X_j; \hat w, g)$ using \eqref{Clval};
\item[a3.] Obtain $\T_j^{\adj\ell} = \mathbbm{1}\{\widehat {\adj\ell}_j \leq t\}$;
\end{enumerate}
{\bf Output:} $\T_1^{\adj\ell}, \dots, \T_n^{\adj\ell}$.
\end{minipage}%
\begin{minipage}{0.45\textwidth}
\vspace{-.07cm}
{\bf \normalsize The $q$-value procedure:}\\[0.19cm]
{\bf Data:} $\mathcal{D} = \{Z_{ij}, i = 1, \dots, m, j = 1, \dots, n\}$\\
{\bf Input:} A pre-specified value $t \in (0, 1)$
\begin{enumerate}
\item[q0.] Compute $X_j = \sum_{i=1}^m Z_{ij}$;
\item[q1.] Compute $\hat w$ as in \eqref{w-hat};
\item[q2.] Evaluate $\hat q_j = q(X_j; \hat w)$ using \eqref{qval};
\item[q3.] Obtain $\T_j^{q} = \mathbbm{1}\{\hat q_j \leq t\}$;
\end{enumerate}
{\bf Output:} $\T_1^{q}, \dots, \T_n^{q}$.
\end{minipage}
\end{algorithm}

\subsection{Analyzing the threshold of the $\adj\ell$- and $q$-value procedures}

We now analyze the threshold of the two procedures. 
Lemma \ref{lem:th-q-cl} provides the asymptotic bounds for each procedure, and their non-asymptotic bounds are deferred to Section \ref{sec:threshold} of \citet{ning23supp}.

\begin{lemma}
\label{lem:th-q-cl}
For a fixed $t \in (0, 1)$ and $w \in (0, 1)$, define $r(w, t) = \frac{wt}{(1-w)(1-t)}$,
\begin{enumerate}
\item[(a)] let $\T^{\adj\ell} = \mathbbm{1}\{\adj\ell(x; w, g) \leq t\}$ be the test function with $\adj\ell(\cdot)$ in \eqref{Clval}, then $\T^{\adj\ell} = \mathbbm{1}\{|x - m/2| \geq m t_n^{\adj\ell}\}$ and 
$t_{m}^{\adj\ell} = t_{m}^{\adj\ell}(w, t) = \eta^{\adj\ell}(r(w, t)) - 1/2$, where
\begin{align}
\label{eqn:etaCl}
\eta^{\adj\ell}(u) = \frac{1}{m} (b/g)^{-1} \left(\frac{\sqrt{2}(1+m) u}{\sqrt{\pi m}}\right),
\end{align}
if $m \to \infty$ and $\log^2(1/r(w,t))/m \to 0$, then
\begin{align}
\label{eqn:tCl}
t_{m}^{\adj\ell}(w,t) \sim \sqrt{\frac{\log (1/r(w,t))}{2m}}.
\end{align}
\item[(b)] let $\T^q = \mathbbm{1}\{q(x; w, g) \leq t\}$ be the test function with $q(\cdot)$ in \eqref{qval}, then $\T^q = \mathbbm{1}\{|x - m/2| \geq mt^{q}_n\}$ and $t^q_{m} := t^q_{m}(w, t) = \eta^q (r(w, t)) - 1/2$, where
\begin{align}
\label{eqn:etaq}
\eta^q(\cdot) = \frac{1}{m}\left(\bar{\mathbf{B}}/\bar\G\right)^{-1}(\cdot),
\end{align}
if $m \to \infty$, then
\begin{align}
\label{eqn:tq}
t^q_{m}(w,t) \sim \sqrt{\frac{\log (1/r(w,t))}{2m}}.
\end{align}
\end{enumerate}
\end{lemma}

\begin{lemma}
\label{lem:compare-thresholds}
Consider the three thresholds $\eta^\ell(u)$, $\eta^{\adj\ell}(u)$, and $\eta^{q}(u)$ given in \eqref{eqn:etal}, \eqref{eqn:etaCl}, and \eqref{eqn:etaq} respectively, for any $u \in (0, 1)$ and $m > 0$, we have $\eta^{\adj\ell}(u) \leq \eta^{\ell}(u)$ and $\eta^{q}(u) \leq \eta^{\ell}(u)$. 
\end{lemma}

\begin{figure}[!h]
  \subfigure[]{\includegraphics[width=0.45\textwidth]{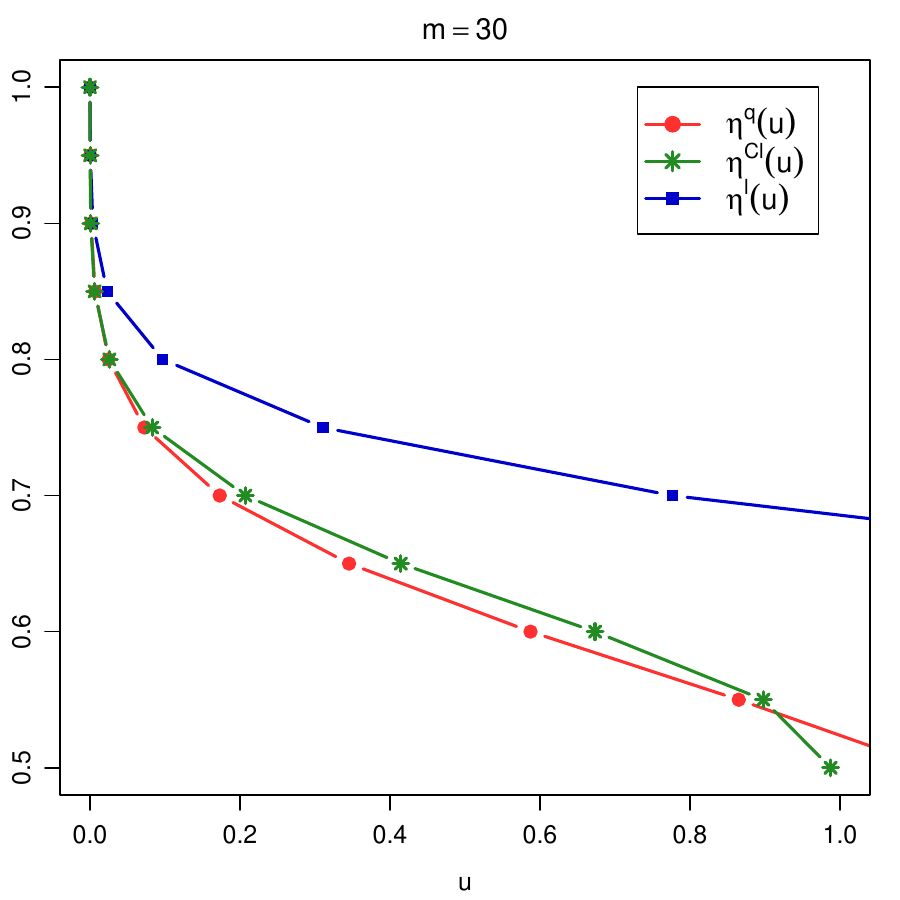}} \quad
  \subfigure[]{\includegraphics[width=0.45\textwidth]{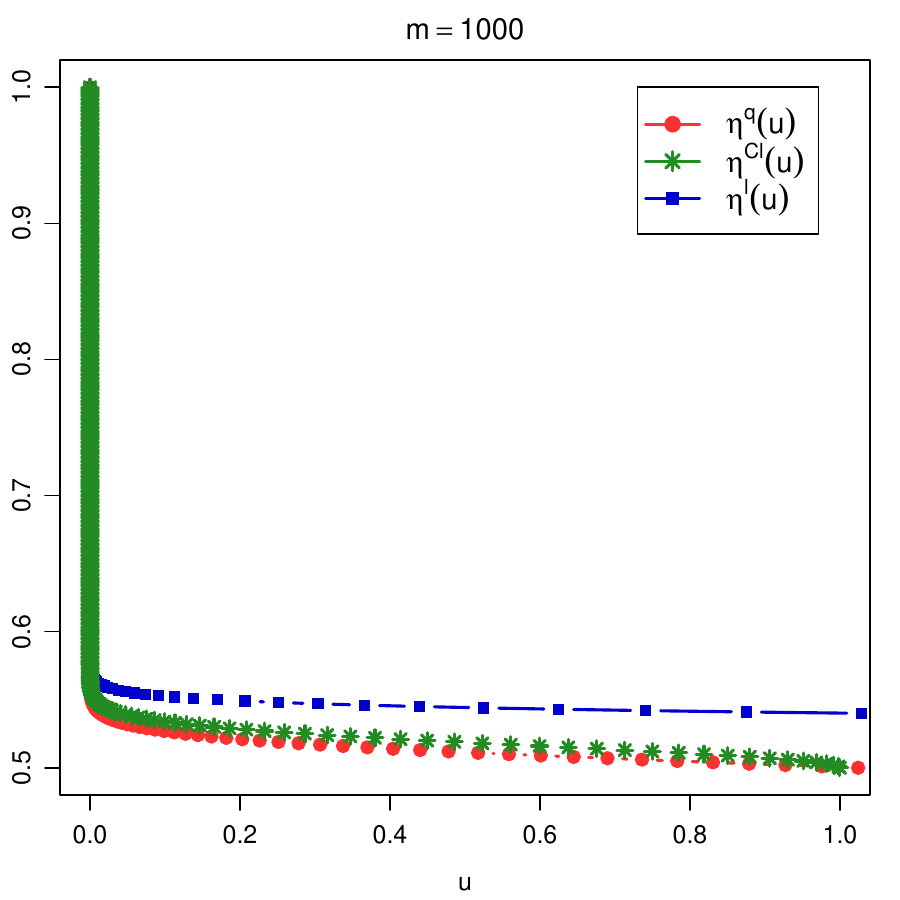}} 
\caption{Plots of $\eta^q(u), \eta^{\adj\ell}(u), \eta^\ell(u)$ for $X \sim \text{Bin}(m, 1/2)$ with (a) $m = 30$ and (b) $m = 1,000$.}
\label{fig:three-th}
\end{figure}

In Figure \ref{fig:three-th}, we compare the three thresholds, $\eta^\ell(u)$, $\eta^{\adj\ell}(u)$, and $\eta^q(u)$, with a small $m$ ($m = 30$) (left) and a large $m$ ($m = 1000$) (right). 
We observe that as $m$ increases, $\eta^q(u)$ and $\eta^{\adj\ell}(u)$ are getting closer, while 
$\eta^\ell(u)$ can be significantly larger than the two.

\section{Uniform FDR control for the $\ell$-, $\adj\ell$-, and $q$-value procedures}
\label{sec:unif-fdr}

In this section, we examine the FDR control of the three procedures introduced in the previous two sections. We provide affirmative results showing that all three procedures ensure valid FDR control while adapting to a broad class of sparsity patterns.

We first study the $\adj\ell$-value and the $q$-value procedures.
In the next theorem, we show both procedures could allow a {\it uniform} control for the FDR under the sparsity assumption $\theta_0 \in l_0[s_n]$ for any $s_n \leq n^{v_1}$ with $v_1 \in (0, 1)$. 
\begin{theorem}
\label{thm:unif-fdr-q-Cl}
Let $w = \hat w$ be the MMLE in \eqref{w-hat},
consider the parameter space $l_0[s_n]$ in \eqref{ell0} with $s_n \leq n^{v_1}$ for some $v_1 \in (0, 1)$, if $m \gg (\log n)^2$, then for the $\adj\ell$-value procedure, there exists a constant $K_1 > 0$ depends on $v_1$ such that for any $t \leq 4/5$ and a sufficiently large $n$,
$$
\sup_{\theta_0 \in l_0[s_n]} \FDR(\theta_0, \T^{\adj\ell}) \leq \frac{K_1 t \log \log n}{\sqrt{\log n}}.
$$
For the $q$-value procedure, there exist a constant $K_2 > 0$ depends on $v_1$ such that for any $t \leq 4/5$ and a sufficiently large $n$,
$$
\sup_{\theta_0 \in l_0[s_n]} \FDR(\theta_0, \T^{q}) \leq K_2 t \log (1/t).
$$
\end{theorem}

\begin{remark}
\label{rmk:m}
Theorem \ref{thm:unif-fdr-q-Cl} requires $m \gg (\log n)^2$. However, we speculate that the conclusion could still hold under a milder assumption on $m$, as \citet{mukherjee15} proved that the necessary condition for constructing powerful two-sided testing procedures under any sparsity assumptions is $m \gg \log n$ (instead of $(\log n)^2$). Indeed, our simulations in Section \ref{sec:sim} show that even when $m = (\log n)^2$, both procedures perform well in FDR control. Nonetheless, we believe our assumption on $m$ cannot be further relaxed with the current proof technique, as the tail bounds for binomial distributions employed in the analysis are already sharp.
 \end{remark}

\begin{remark}
By adopting a similar argument as that of Theorem 2 in \citetalias{cast20MT},
the $log(1/t)$ term in the upper bound of the FDR result for the $q$-value procedure in Theorem \ref{thm:unif-fdr-q-Cl} can be removed if replacing each $\T_j^q$ in Step q3 of Algorithm \ref{algo:Cl-q-val} with $\T_j^q\mathbbm{1}\{\hat w > w_n\}$ for $w_n = {\log n}/{n}$. 
\end{remark}

The proof of Theorem \ref{thm:unif-fdr-q-Cl} is left to Section \ref{sec:pf-unif-fdr} of \citet{ning23supp}. 
Our proof strategy is inspired by that used in \citetalias{cast20MT} for the Gaussian sequence model. The essential part of our proof---also, the most difficult part---is to establish a tight concentration bound for $\hat w$ to control it within a close neighborhood near $s_n/n$ (up to some constant). To do so, we study two different scenarios for $\hat w$ depending on whether \eqref{w-hat} has a unique solution or not. We also obtain sharp bounds for several quantities related to the score function (see Section \ref{sec:bound-m1-mtilde} in \citet{ning23supp} for more details).
In particular, we pay close attention to small signals near $1/2$ to ensure they do not dominate in the process of bounding those quantities.

There are three major differences in adapting the proof strategy from \citetalias{cast20MT} to study our methods. First, we need to obtain sharp lower and upper bounds for a binomial distribution, which are much more challenging than those for a Gaussian distribution. Our strategy involves first controlling the approximation error between the binomial distribution and a corresponding Gaussian distribution, and then using existing bounds for Gaussian distributions. However, the approximation error has a complicated expression and is nontrivial to handle.

Second, when dealing with the ratio of binomial distributions under the null and alternative hypotheses. In the Gaussian sequence model, the two distributions differ only in their means. However, for binomial distributions, both the means and variances are different.  
This distinction increases the complexity of analyzing the score function for bounding $\hat w$ (see Remark \ref{rmk:diff-Gaussian-binomial} in \citet{ning23supp}). 

Last, there is a difference in obtaining the bound for the FDR in our model than that in the Gaussian sequence model with the spike-and-Laplace (or Cauchy)-slab prior, as the suboptimality of the threshold of the posterior given in \eqref{posterior} is not present in their procedures. Consequently, in addition to controlling the errors caused by estimating those signals below the threshold of the $q$-value or the $\adj\ell$-value, one must also analyze those slightly larger signals in between the threshold and $\eta^\ell(u)$ from the posterior, ensuring that the accumulative error from estimating these signals do not cause trouble in bounding the FDR.

Now let's return to the $\ell$-value procedure. In the next lemma, we show that this procedure can also achieve uniform FDR control, but with a much smaller upper bound for the FDR compared to those of the previous two procedures, suggesting that the $\ell$-value procedure can be too conservative.

\begin{lemma}
\label{lem:unif-fdr-l}
For the $\ell$-value in \eqref{lval} with $w = \hat w$ be the MMLE in \eqref{w-hat}, under the same condition as in Theorem \ref{thm:unif-fdr-q-Cl}, there exists a constant $K_3$ depends on $v_1$ such that for any $t \leq 4/5$ and a sufficiently large n,
$$
\sup_{\theta_0 \in l_0[s_n]} \FDR(\theta_0, \T^{\ell}) \leq \frac{K_3 t (\log m + \log\log n)}{\sqrt{m \log n}}.
$$
\end{lemma}

The proof of Lemma \ref{lem:unif-fdr-l} is similar to that of the $\adj\ell$-value procedure in Theorem \ref{thm:unif-fdr-q-Cl} and can be found in Section \ref{sec:pf-unif-fdr} of \citet{ning23supp}.

\section{FDR and FNR control for large signals}
\label{sec:fdr-fnr-large}

In this section, we focus on the set containing `large' signals given by 
\begin{align}
\label{L0}
\Theta_0[s_n, a] = \left\{
\theta \in l_0[s_n]: |\theta_j - 1/2| \geq a \zeta(s_n/n), \ j \in \mathcal{S}_\theta, \ |\mathcal{S}_\theta| = s_n, a > 0
\right\},
\end{align}
where, for any $\omega \in (0, 1)$, 
\begin{align}
\label{eqn:zeta}
\zeta(\omega) = \sqrt{\frac{1}{2m} \log \left(\frac{1}{\omega}\right)}.
\end{align}
We first establish the lower bound for a large class of thresholding-based multiple testing procedures.
We then study the $q$-value procedure and show it can effectively control the FDR at any target level $t \in (0, 1)$.
Last, we examine the FNR control and the multiple testing risk for the three procedures.

\subsection{The lower bound for the testing boundary}
\label{sec:lb-mt-risk}

Let $\mathcal{T}$ stand for a class of thresholding-based multiple testing procedures. Then
for any test $\T \in \mathcal{T}$, $\T = \{\T_1, \dots, \T_n\}$, let
\begin{align}
\label{eqn:test}
\T_j(X) = \mathbbm{1}\{X_j - m/2 \geq m \tau_1(X) \ \text{or} \ m/2 - X_j \geq m\tau_2(X)\},
\quad 1\leq j \leq n,
\end{align}
for some positive measurable functions $\tau_1(X)$ and $\tau_2(X)$. 
The multiple testing risk is given by
\begin{align}
\label{risk}
\mathfrak{R}(\theta, \T) = \FDR(\theta, \T) + \FNR(\theta, \T).
\end{align}
The next proposition establishes the lower bound for all thresholding-type tests in $\mathcal{T}$.
\begin{prop}
\label{prop:lowerbound}
Let $\mathcal{T}$ be a class of thresholding-based multiple testing procedures, suppose that for some $v_1 \in (0, 1)$ for which $s_n \leq n^{v_1}$ and $\theta_0 \in \Theta_0[s_n, a]$ in \eqref{L0}, if $m \gg (\log n)^2$,  then for any $\T \in \mathcal{T}$ in \eqref{eqn:test} and any positive $a < 1$, we have
\begin{align}
\liminf_{n \to \infty} \inf_{\T \in \mathcal{T}} \sup_{\theta_0 \in \Theta_0[s_n, a]} 
\mathfrak{R}(\theta_0, \T) \geq 1.
\end{align}
\end{prop}

The proof of Proposition \ref{prop:lowerbound} is left to Section \ref{sec:pf-lowerbound} of \citet{ning23supp}. In the proof, we obtain both upper and lower bounds of the inverse of 
$\bar \B(u)$,
the upper tail of the distribution $\text{Bin}(m, 1/2)$. These bounds are new and can be found in Lemma \ref{lem:inv-bin-bound} of \citet{ning23supp}.

\subsection{FDR and FNR control for the $\ell$-value, $\adj\ell$-value, and $q$-value procedures}

In the next theorem, we show that the $q$-value procedure can successfully control FDR at an arbitrary targeted level $t \in (0, 1)$ for $\theta_0 \in \Theta_0[s_n, a]$ for any $a \geq 1$. 

\begin{theorem}
\label{thm:fdr-qval}
For the $q$-value given in \eqref{qval} and let $w = \hat w$ be the MMLE in \eqref{w-hat}, suppose $s_n \leq n^{v_1}$ for some $v_1 \in (0, 1)$ and $m \gg (\log n)^2$, then, for any fixed $t \in (0, 1)$ and $a \geq 1$,
$$
\lim_{n \to \infty} \sup_{\theta_0 \in \Theta_0[s_n, a]} \FDR(\theta_0, \T^q) 
= \lim_{n \to \infty} \inf_{\theta_0 \in \Theta_0[s_n, a]} \FDR(\theta_0, \T^q) 
= t.
$$
\end{theorem}

Next, we obtain an upper bound for the FNR of the $q$-value and the $\adj\ell$-value procedures respectively. Our result shows that both procedures can effectively control the FNR for large signals as sample size increases.
 
\begin{theorem}
\label{thm:fnr-Clval-qval}
Let $\hat w$ be the MMLE given in \eqref{w-hat}, if $s_n \leq n^{v_1}$ for some $v_1 \in (0, 1)$ and $m \gg (\log n)^2$, then, for a fixed $t \in (0, 1)$ and any $a \geq 1$, as $n\to \infty$,
\begin{itemize}
\item[(i)] for the $\adj\ell$-value given in \eqref{Clval} and $w = \hat w$, 
$\sup_{\theta_0 \in \Theta_0[s_n, a]} \FNR(\theta_0, \T^{\adj\ell}) \to 0,$
\item[(ii)] for the $q$-value given in \eqref{qval} and $w = \hat w$, 
$\sup_{\theta_0 \in \Theta_0[s_n, a]} \FNR(\theta_0, \T^{q}) \to 0.$
\end{itemize}
\end{theorem}
By combining Theorem \ref{thm:fnr-Clval-qval} and the uniform FDR control result in Section \ref{sec:unif-fdr}, we obtain the upper bound for the multiple testing risk for the $\adj\ell$-value and $q$-value procedures in the next corollary.

\begin{corollary}
Let $\hat w$ be the MMLE given in \eqref{w-hat}, for a fixed $t \in (0, 1)$ and the risk function $\mathfrak{R}(\theta_0, \cdot)$ given in \eqref{risk}, if $s_n \leq n^{v_1}$ for some $v_1 \in (0, 1)$ and $m \gg (\log n)^2$, then, for $\theta_0 \in \Theta_0[s_n]$ and any $a \geq 1$, as $n \to \infty$, 
$ \sup_{\theta_0 \in \Theta_0[s_n, a]} \mathfrak{R}(\theta_0, \T^{\adj\ell}) \to 0$ and
$\sup_{\theta_0 \in \Theta_0[s_n, a]} \mathfrak{R}(\theta_0, \T^q) \to t$.
\end{corollary}

Last, we show a negative result on FNR control for the $\ell$-value procedure.
\begin{lemma}
\label{lem:fnr-lval}
Let $\hat w$ be the MMLE in \eqref{w-hat} and $t \in (0,1)$ be a fixed value, if $s_n \leq n^{v_1}$ and $m \gg (\log n)^2$ for some $v_1 \in (0, 1)$, then for $\theta_0 \in \Theta_0[s_n, a]$ with $a = 1$ and the $\ell$-value in \eqref{lval} with $w = \hat w$, as $n \to \infty$,
$\sup_{\theta_0 \in \Theta_0[s_n, 1]} \FNR(\theta_0, \T^{\ell}) \to 1. $
\end{lemma}
 
\begin{remark}
Note that one can apply the same argument as the proof of Theorem \ref{thm:fnr-Clval-qval} to show that the $\ell$-value procedure can effectively control the FNR for those signals that are bounded away from the boundary $1/2 \pm \sqrt{ \left(\log(n/s_n) + \log (\sqrt{m})\right) /(2m)}$. 
The $\ell$-value procedure fails to control the FNR for signals in between this boundary and $1/2 \pm \sqrt{\log(n/s_n)/(2m)}$.
\end{remark}
 
\section{Numerical experiments}
\label{sec:sim}
In this section, we conduct numerical experiments to compare the three procedures based on the $\ell$-value, $q$-value, and $\adj\ell$-value. We also compare them with the Benjamini-Hochberg (BH) procedure, which is a benchmark procedure that is routinely used in practice. 
To run the BH procedure, we use the standard package `p.adjust' in $\mathsf{R}$.

Two sets of simulations are conducted: the first set aims to validate our theoretical results in Sections \ref{sec:unif-fdr} and \ref{sec:fdr-fnr-large}. For this purpose, we generate data from the binomial distribution $\text{Bin}(m, \vartheta_0)$ with a fixed $m$.
In the second study, we allow $m$ to vary with $j$, hence denoted it by $m_j$. This scenario is more realistic in certain practical problems, such as the crowdsourcing problem mentioned in \nameref{sec:intro}, where the number of workers assigned to each task can vary. 

Besides the two studies, additional simulations are conducted to compare the performance of the three procedures by choosing two different priors for $\gamma$. Their results can be found in Section \ref{sec:ad-sim} of the supplementary material \citep{ning23supp}.

\subsection{Comparing the $\ell$-value, $q$-value, $\adj\ell$-value, and BH procedures}
\label{sec:sim-1}

Data are generated as follows: for each dataset, we choose $\vartheta_0$, $s_n$, $n$, and $m$, and generate $n - s_n$ independent samples from $\text{Bin}(m, 1/2)$ and $s_n$ independent samples from $\text{Bin}(m, \vartheta_0)$, then the true value $\theta_0 = (\underbrace{\vartheta_0, \dots, \vartheta_0}_{s_n}, \underbrace{1/2, \dots, 1/2}_{n - s_n})$. The dimension of each dataset is $m \times n$. 
We estimate $\hat w$ in \eqref{w-hat} using the `optim' function in $\mathsf{R}$ and apply the four procedures ($\ell$-value, $q$-value, $\adj\ell$-value, and BH procedures) to multiple test at three different target FDR levels: $t = 0.05$, $0.1$, and $0.2$.  
Last, we repeat each experiment 10,000 times and calculate the average value for the FDR. 
In each experiment, we set $n = 10,000$ and choose $\vartheta_0$ to be one of 45 equally spaced values between 0.5 and 0.95.
Three different values for the ratio $s_n/n$ are considered: $0.001$, $0.1$, and $0.5$, representing super-sparse, sparse, and dense scenarios and three different values for $m$ are chosen: $(\log n)^2 \approx 85, 200$, and $1000$, representing small, medium, and large $m$ cases.

\begin{figure}[h!]
\centering
	  \subfigure[$m = 85$, $s_n/n = 0.001$]{\includegraphics[width=0.33\textwidth]{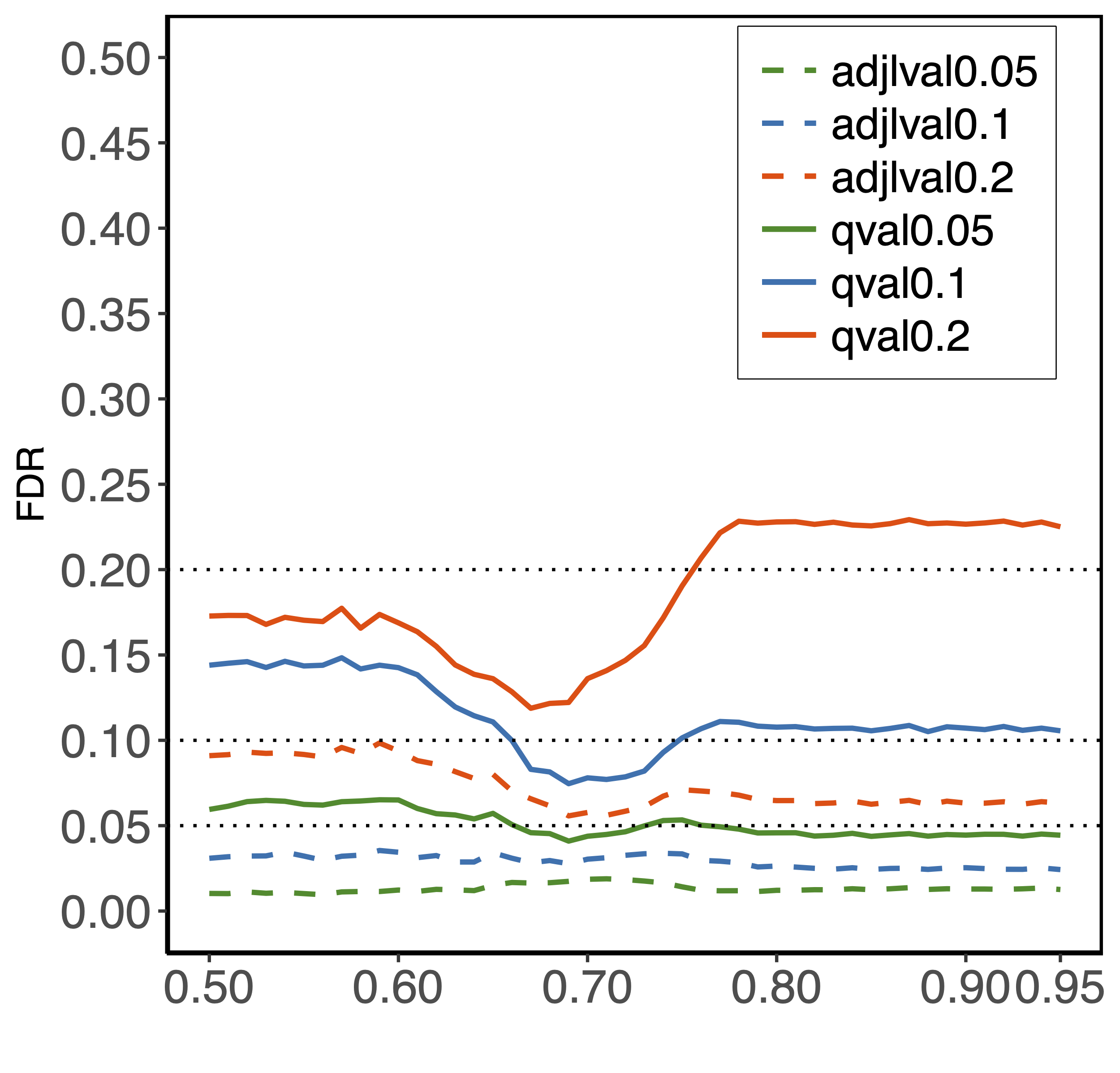}}%
     \subfigure[$m = 85$, $s_n/n = 0.1$]{\includegraphics[width=0.33\textwidth]{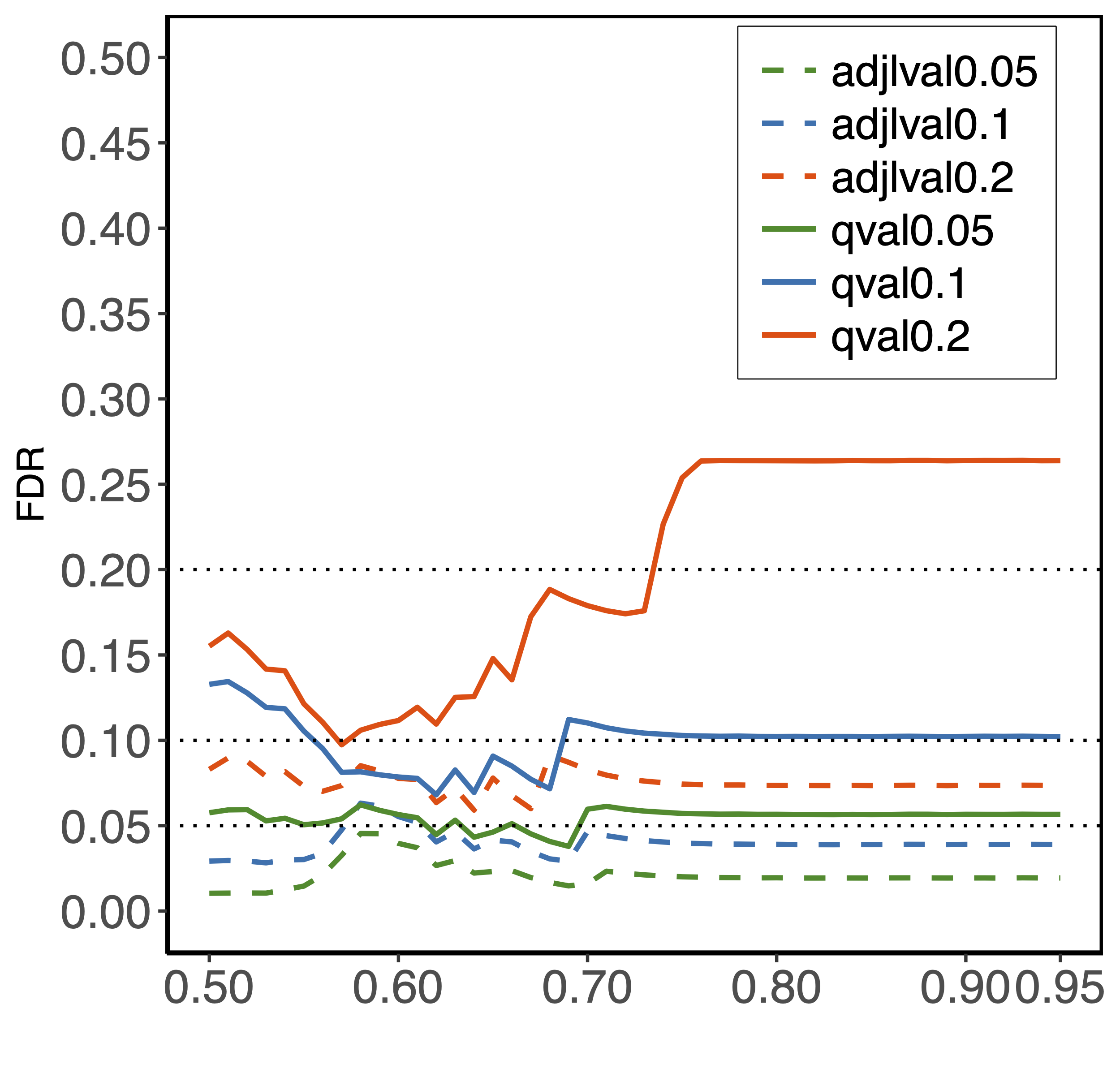}}%
    \subfigure[$m = 85$, $s_n/n = 0.5$]{\includegraphics[width=0.33\textwidth]{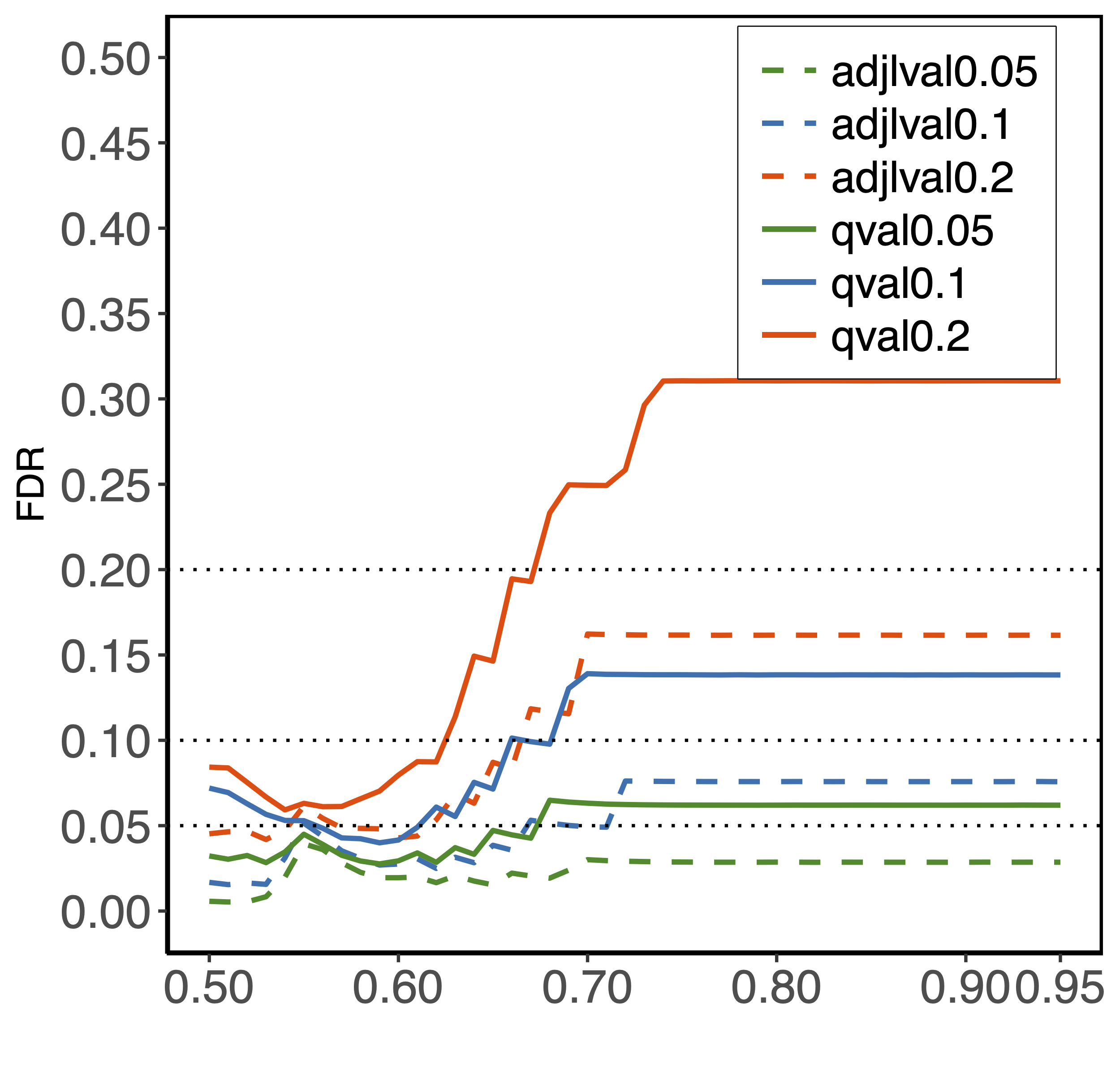}} 
    \subfigure[$m = 200$, $s_n/n = 0.001$]{\includegraphics[width=0.33\textwidth]{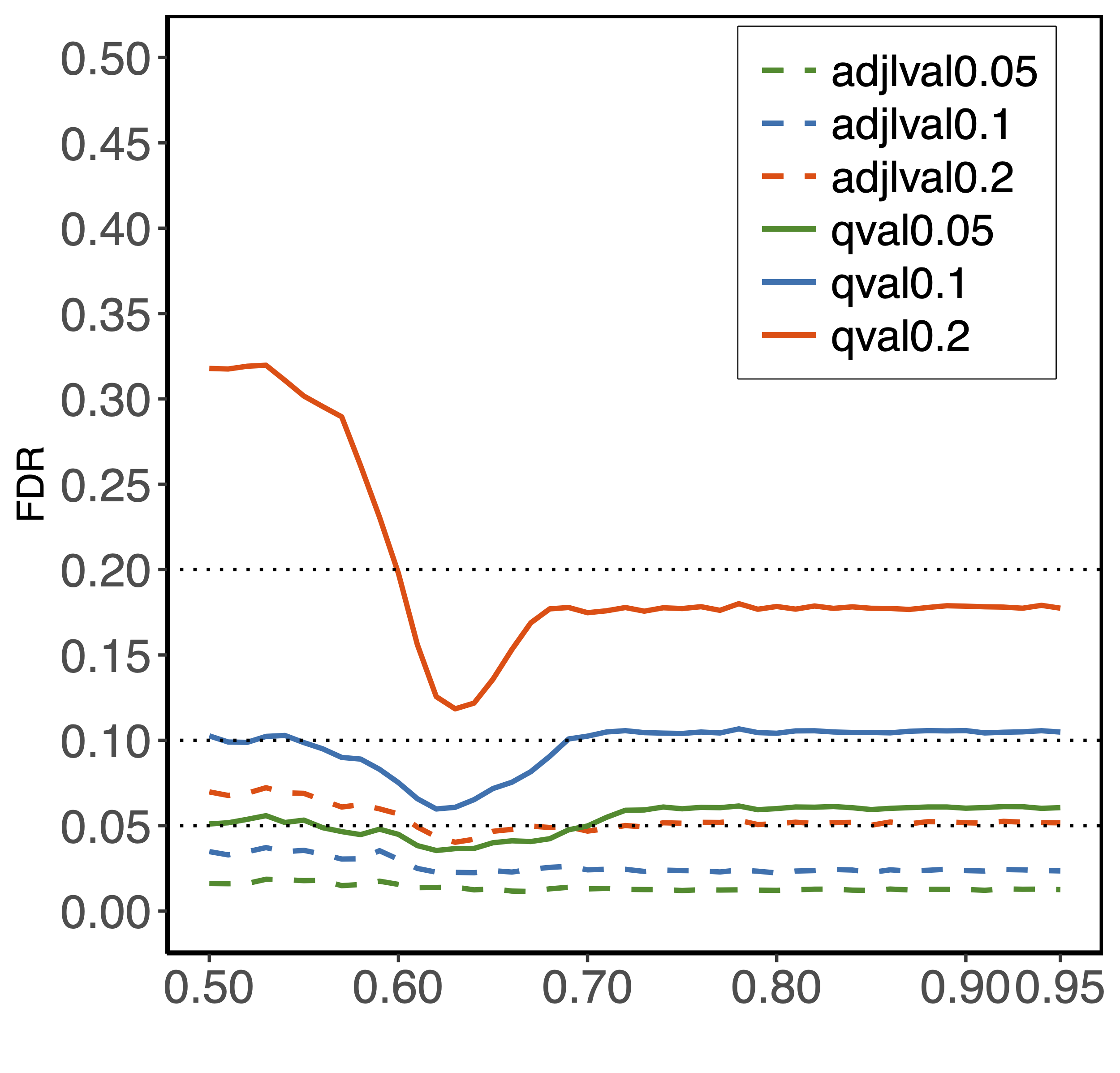}}%
     \subfigure[$m = 200$, $s_n/n = 0.1$]{\includegraphics[width=0.33\textwidth]{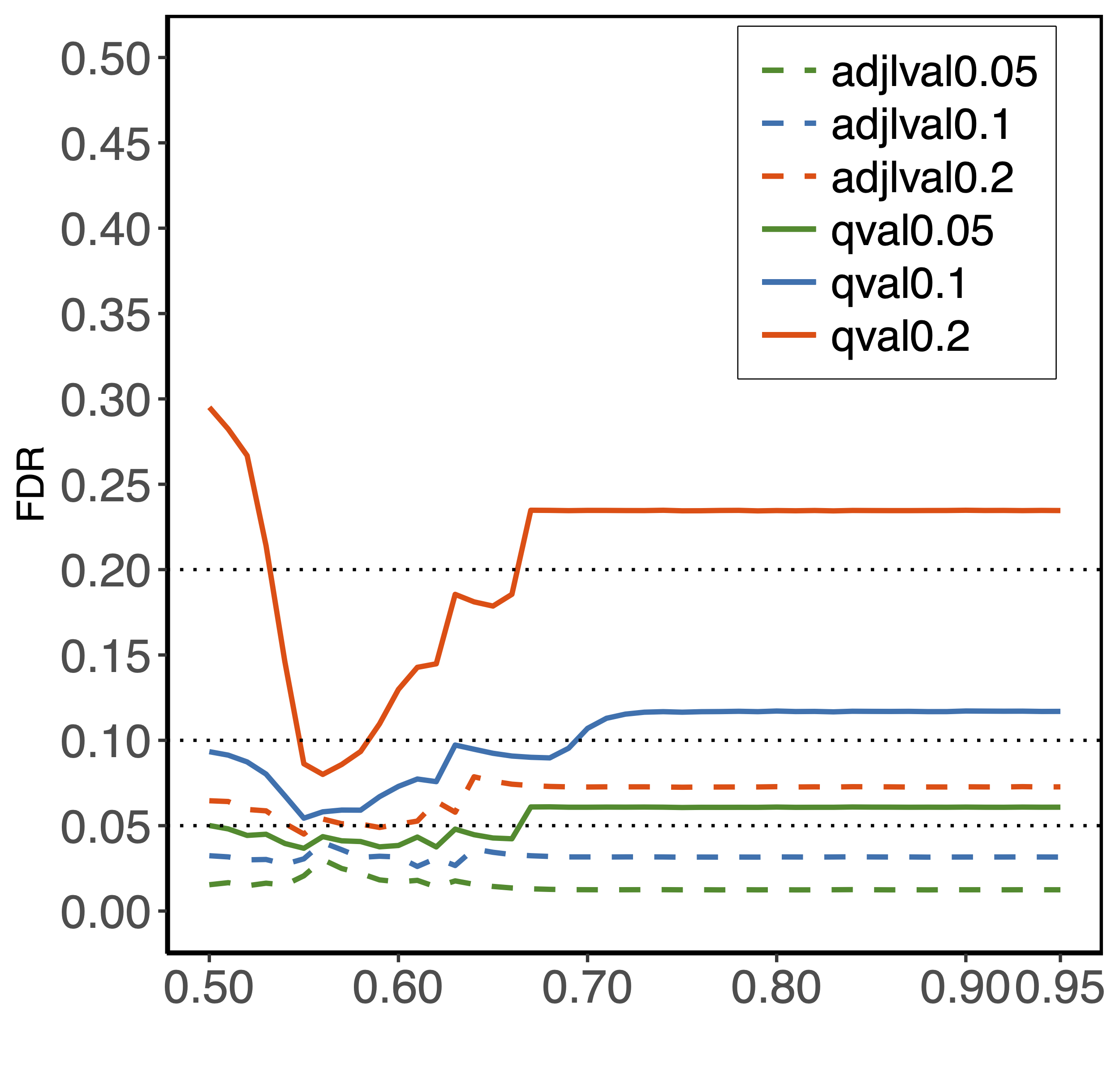}}%
    \subfigure[$m = 200$, $s_n/n = 0.5$]{\includegraphics[width=0.33\textwidth]{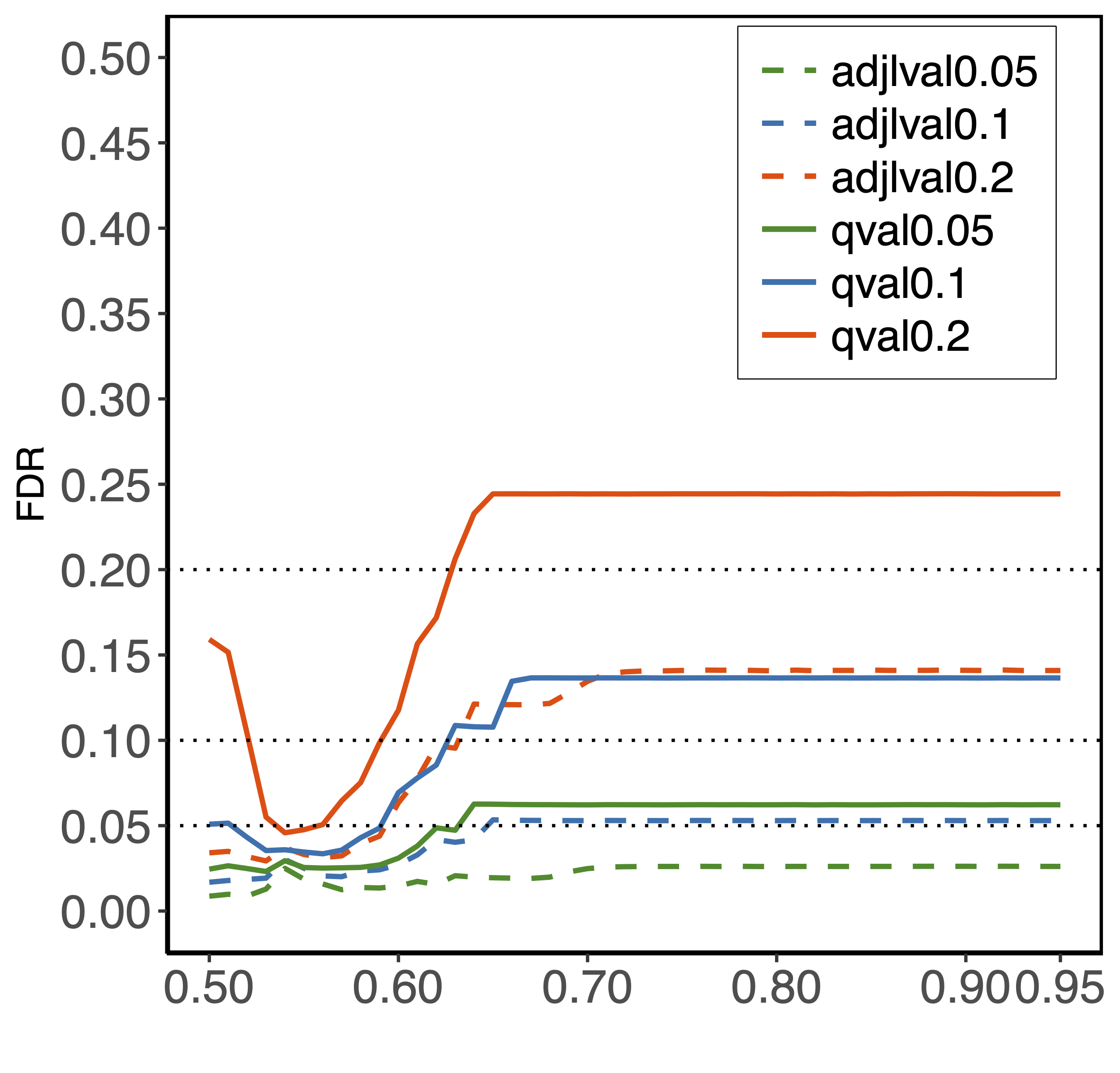}} 
    \subfigure[$m = 1,000$, $s_n/n = 0.001$]{\includegraphics[width=0.33\textwidth]{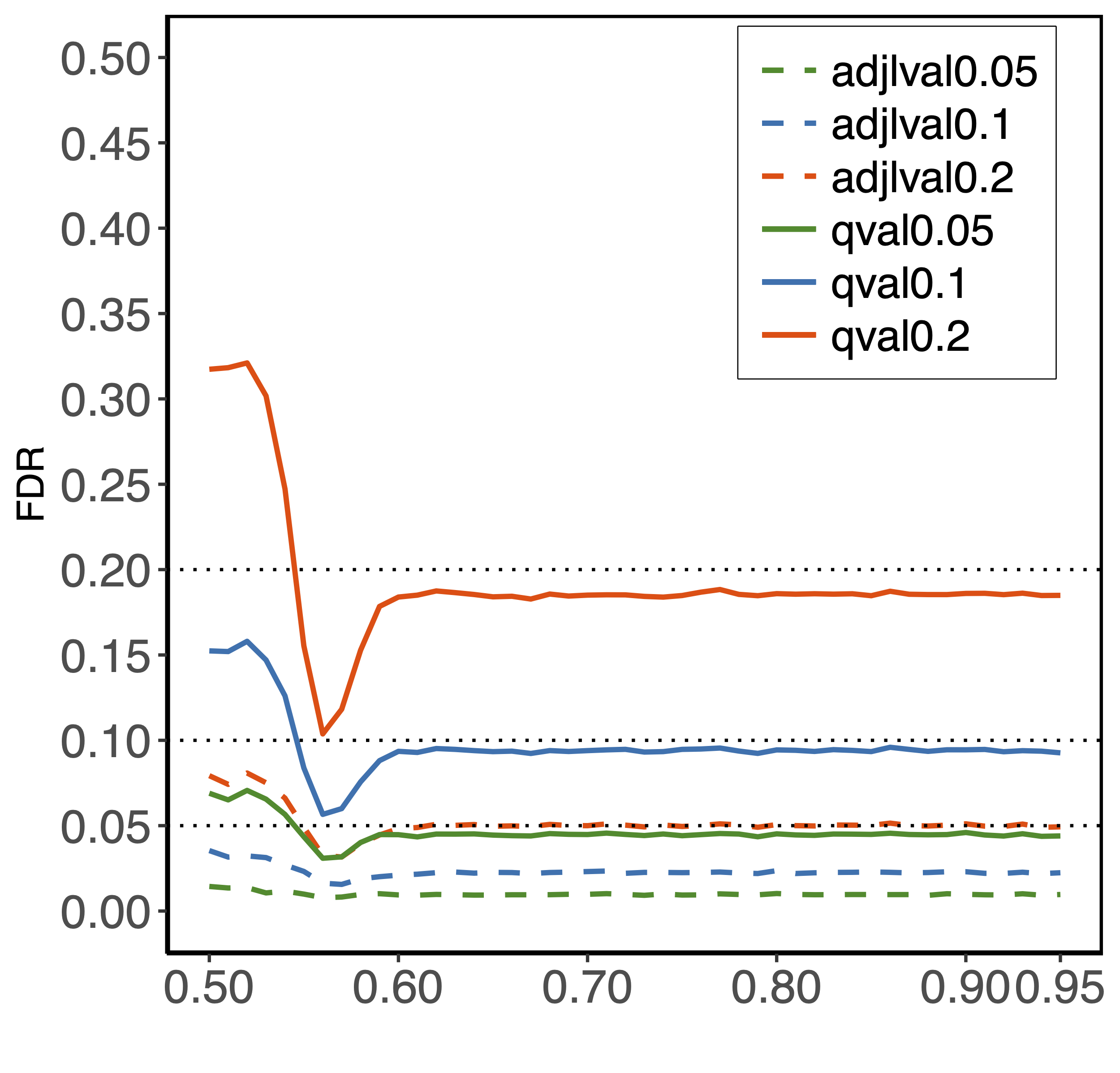}}%
     \subfigure[$m = 1,000$, $s_n/n = 0.1$]{\includegraphics[width=0.33\textwidth]{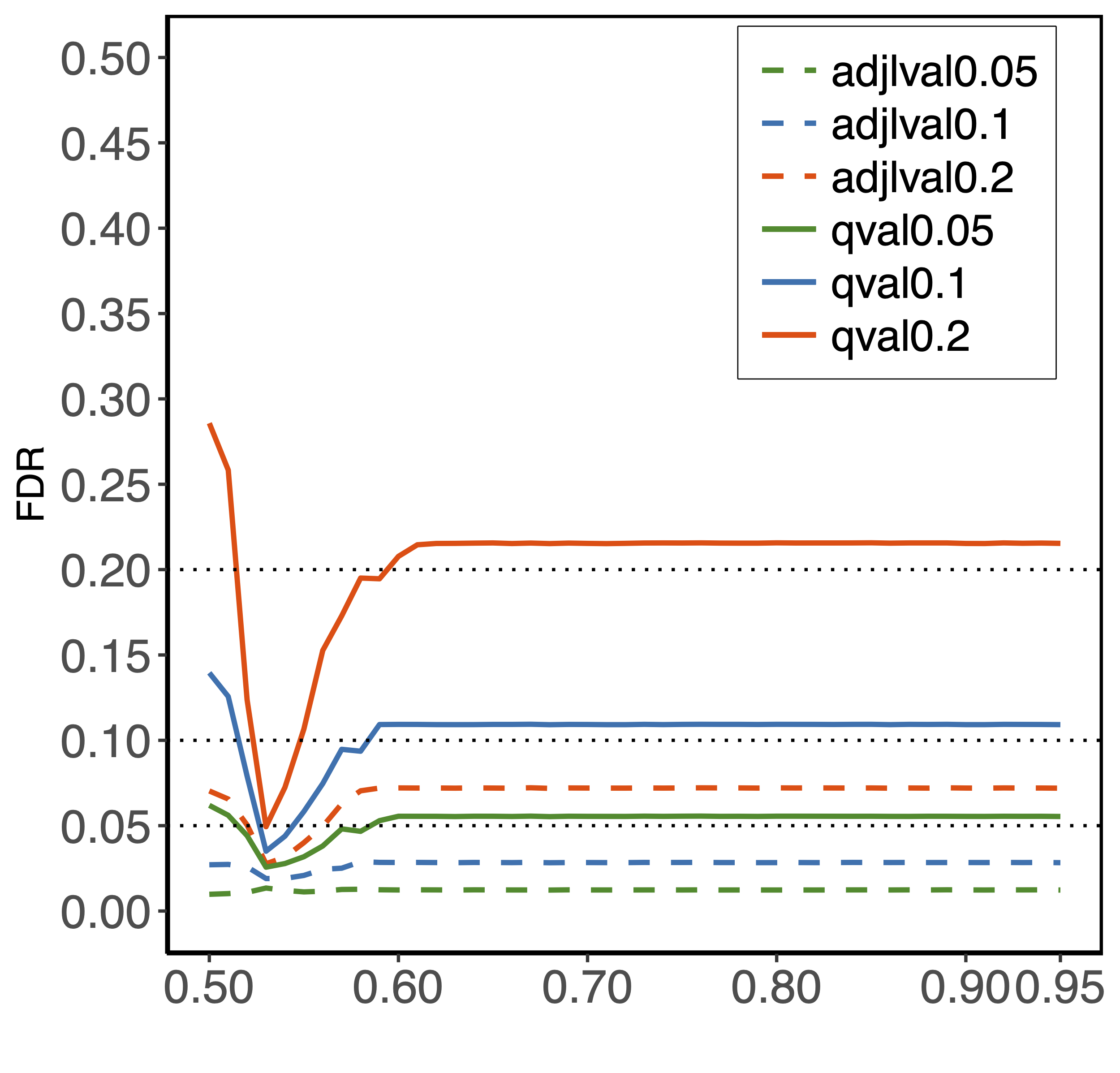}}%
    \subfigure[$m = 1,000$, $s_n/n = 0.5$]{\includegraphics[width=0.33\textwidth]{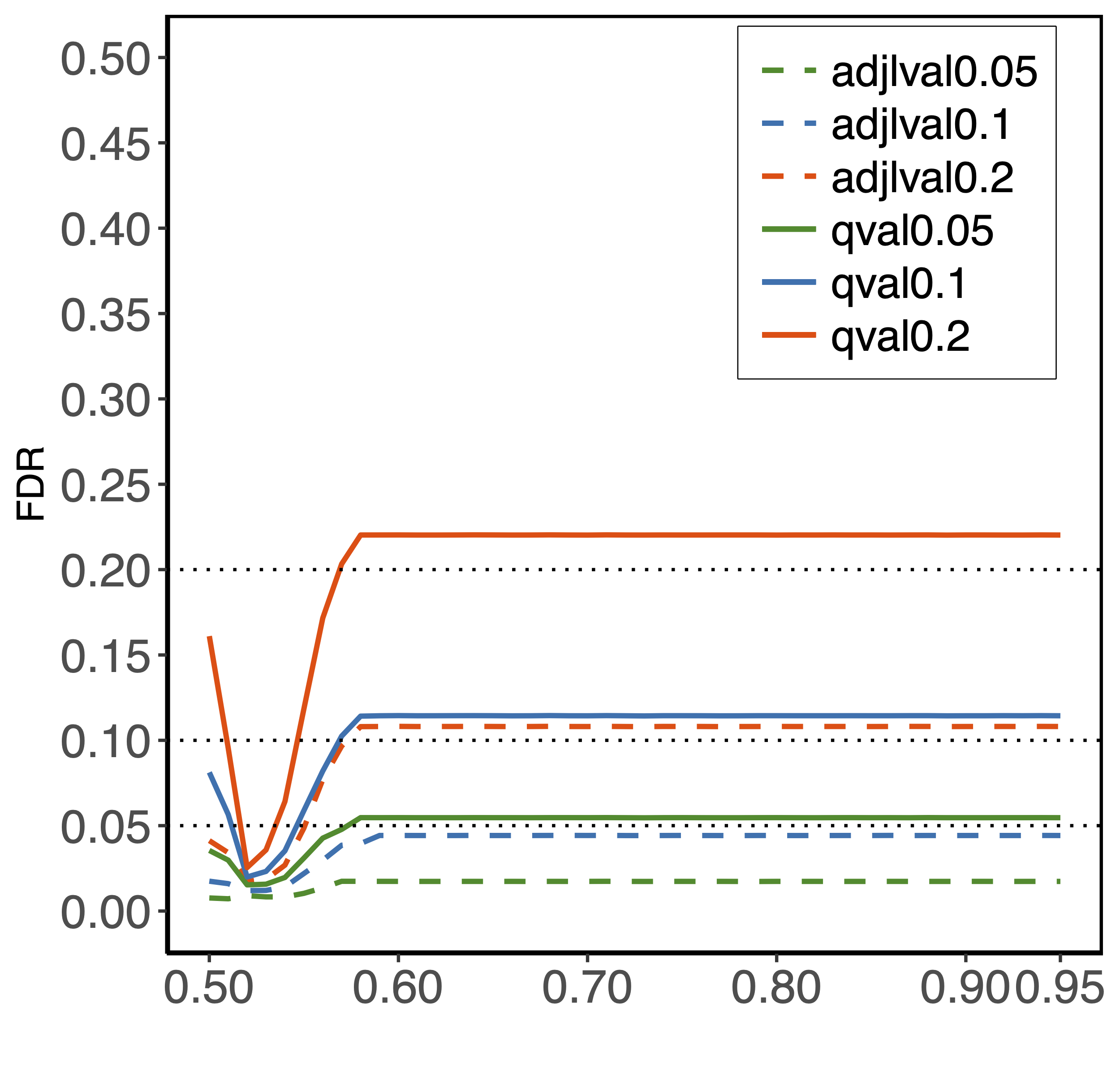}} 
	\caption{The estimated FDR of the $\adj\ell$-value (dash) and the $q$-value (solid) procedures at $t = 0.05$ (blue), $t = 0.1$ (green), and $t = 0.2$ (red) respectively with $m = (\log n)^2 \approx 85, 200$, and $1000$ and $s_n/n = 0.001, 0.1$, and $0.5$.}
	\label{fig:compare-q-Cl}
\end{figure}

\begin{figure}[h!]
\centering
	  \subfigure[$m = 85$, $s_n/n = 0.001$]{\includegraphics[width=0.33\textwidth]{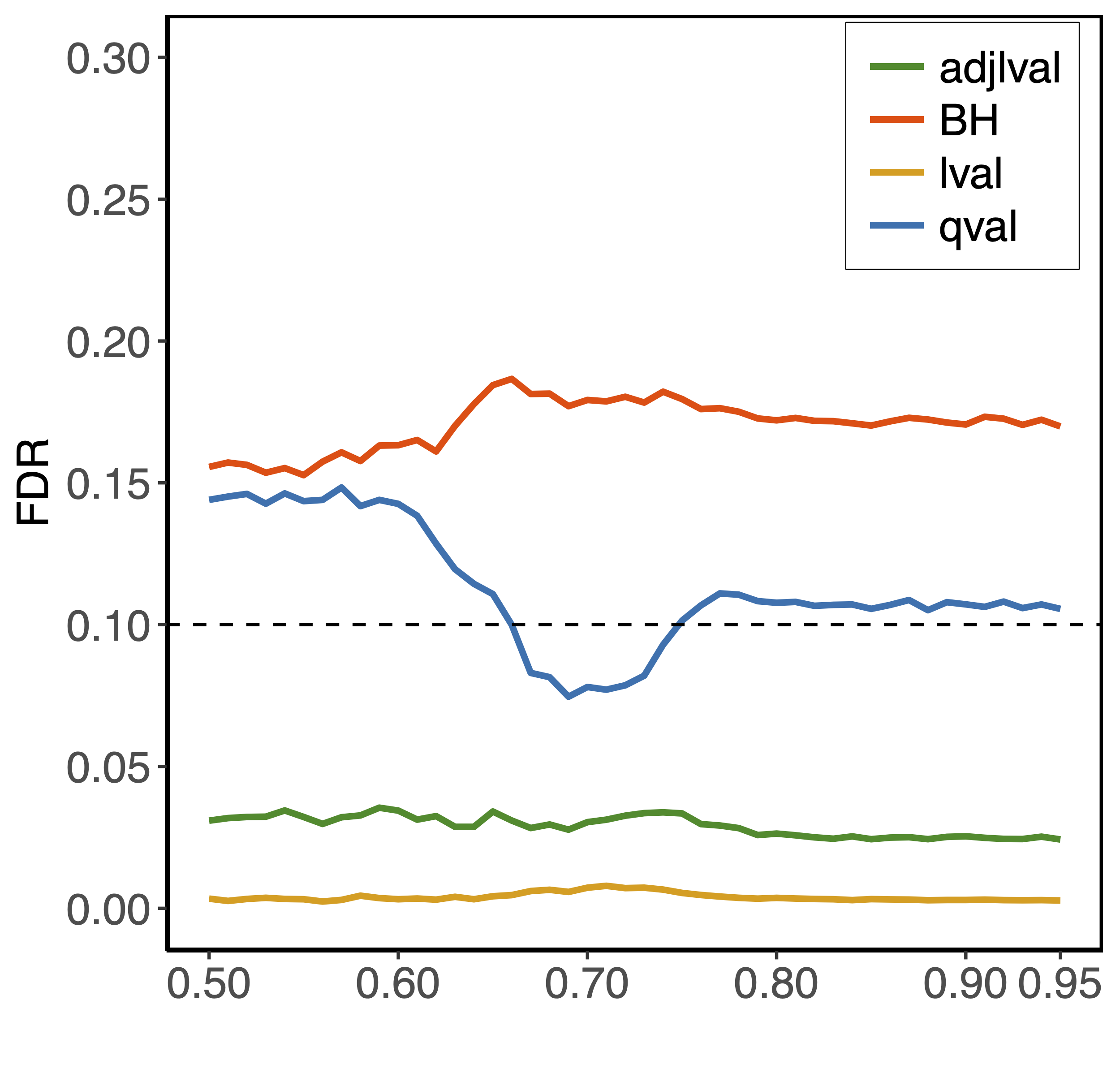}}%
     \subfigure[$m = 85$, $s_n/n = 0.1$]{\includegraphics[width=0.33\textwidth]{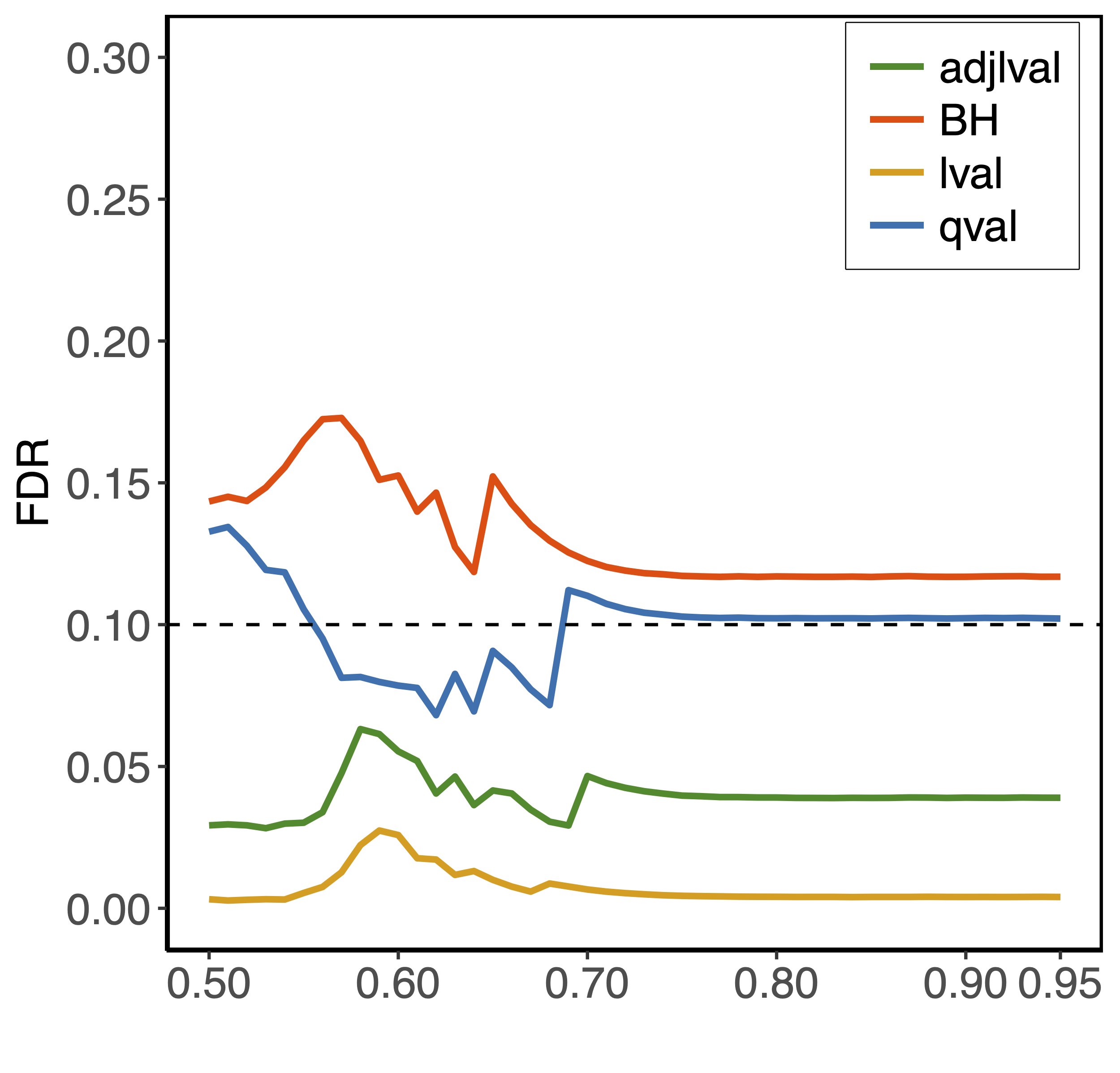}}%
    \subfigure[$m = 85$, $s_n/n = 0.5$]{\includegraphics[width=0.33\textwidth]{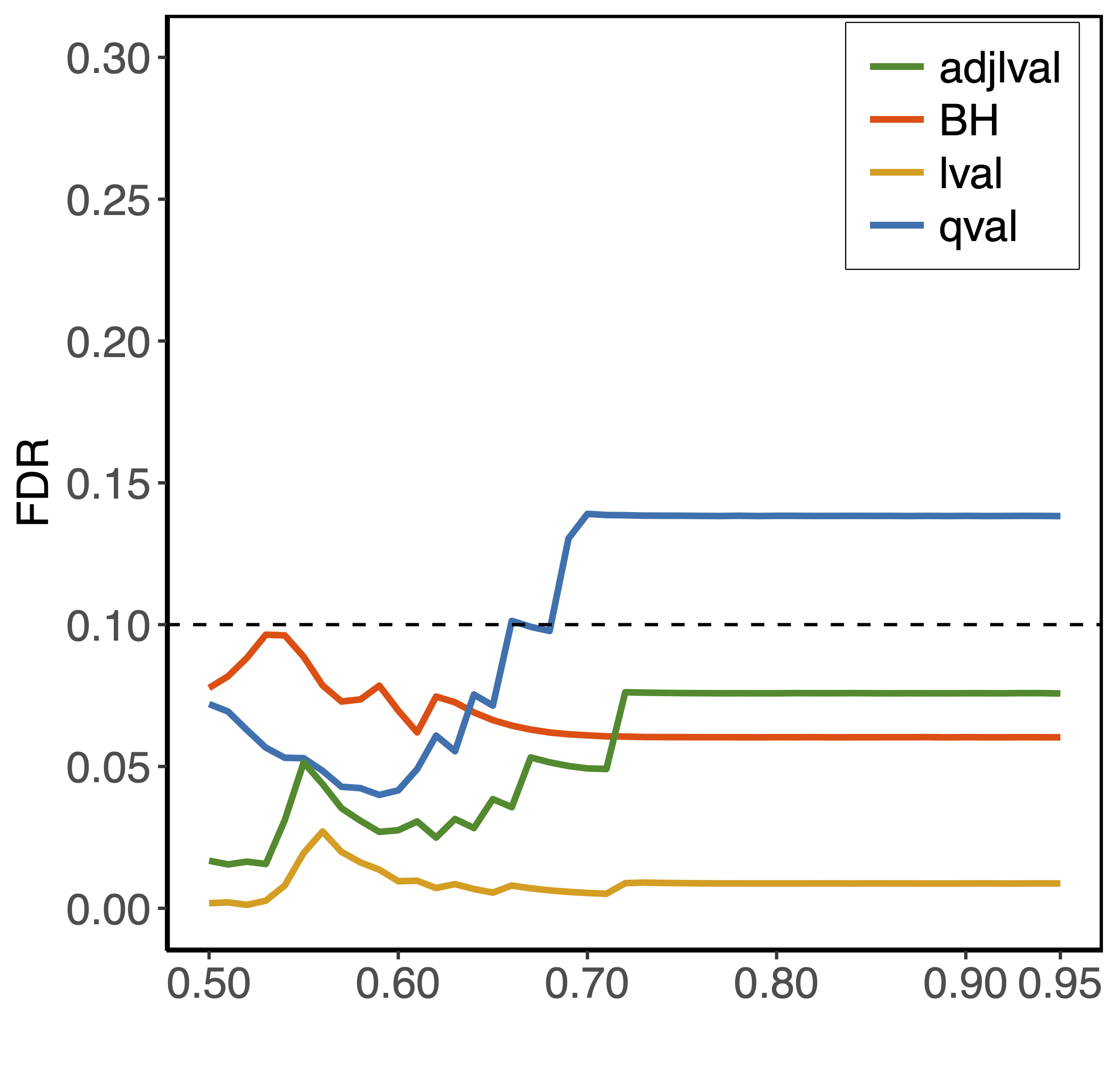}} 
    \subfigure[$m = 200$, $s_n/n = 0.001$]{\includegraphics[width=0.33\textwidth]{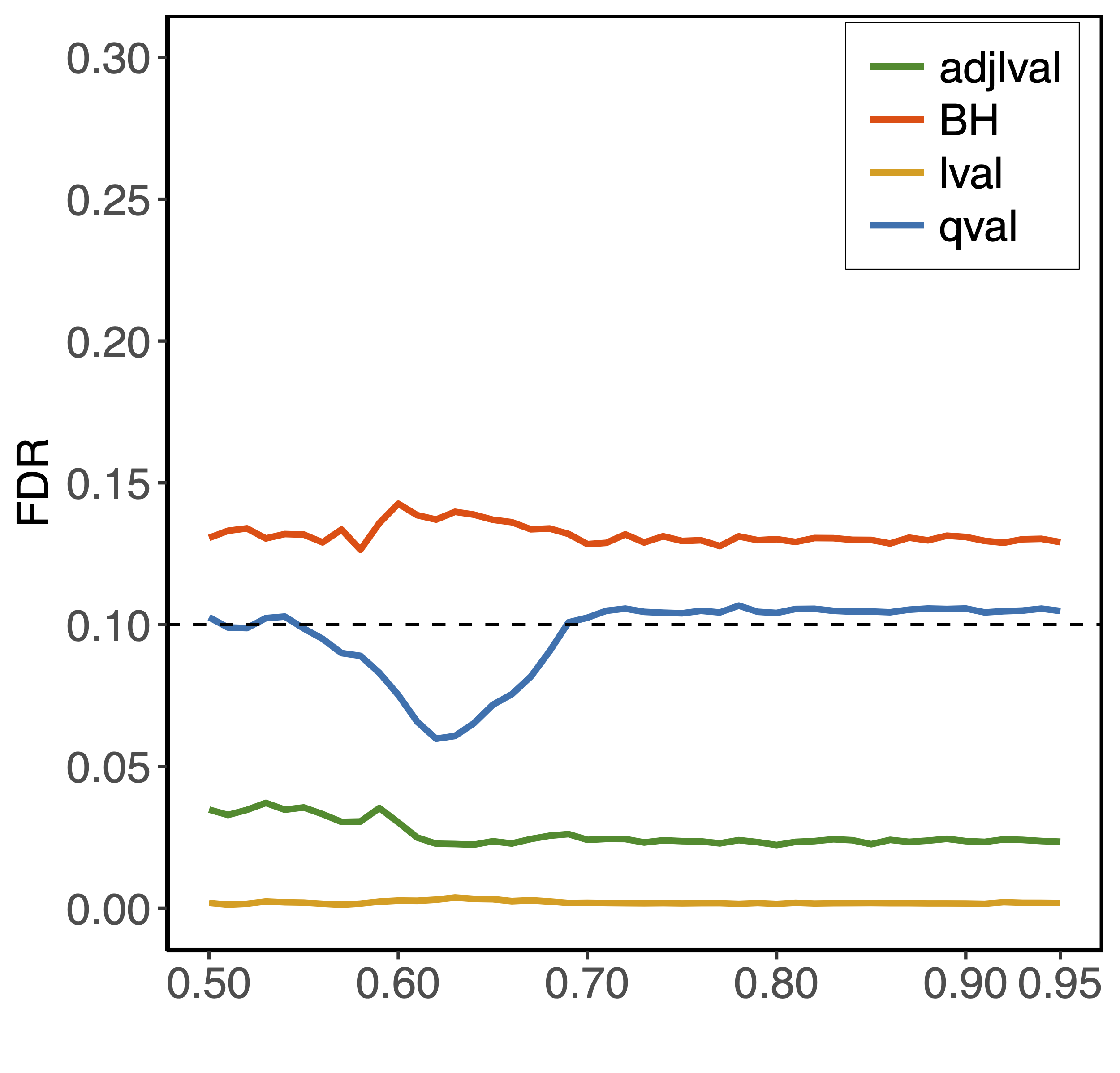}}%
     \subfigure[$m = 200$, $s_n/n = 0.1$]{\includegraphics[width=0.33\textwidth]{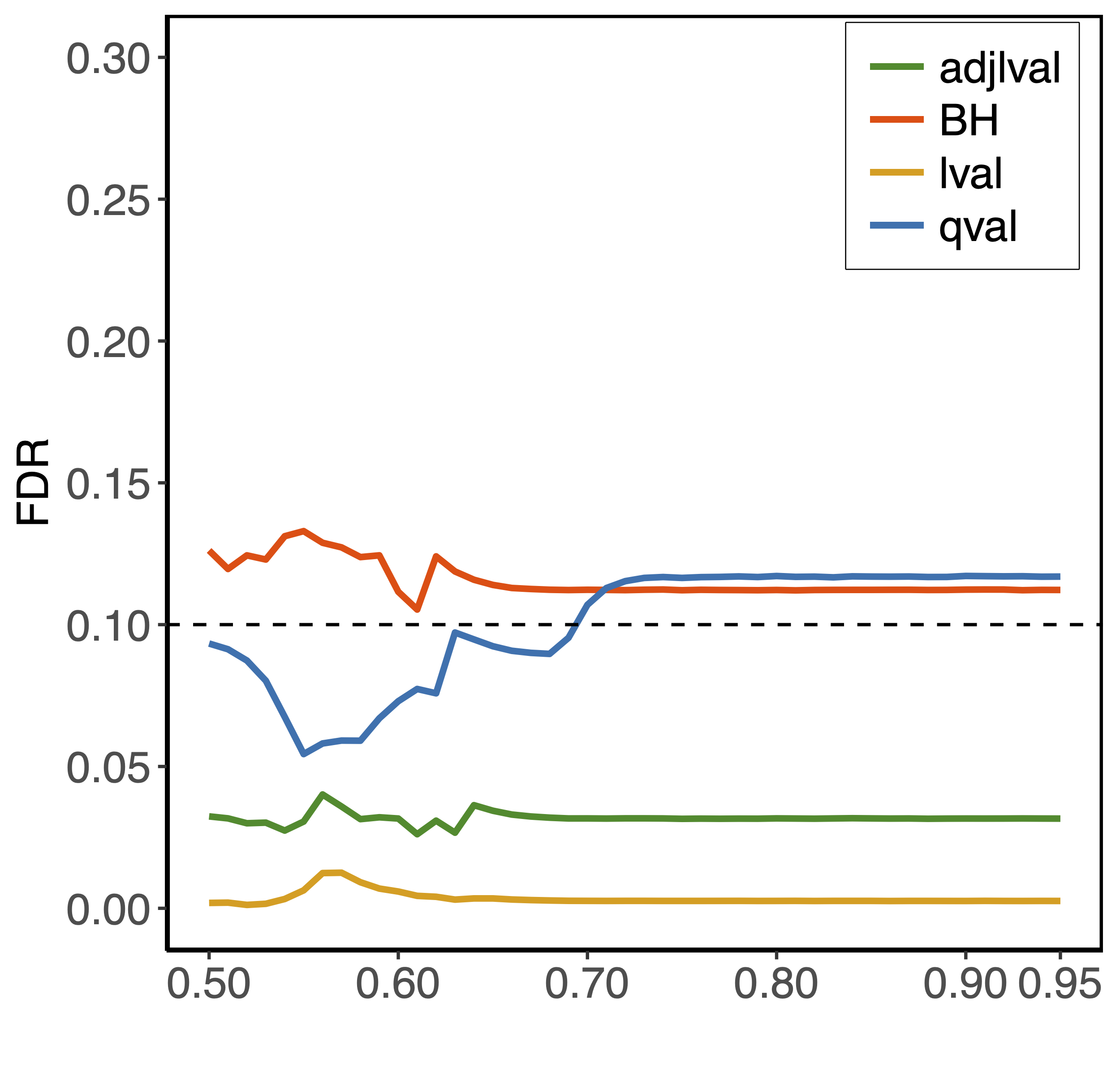}}%
    \subfigure[$m = 200$, $s_n/n = 0.5$]{\includegraphics[width=0.33\textwidth]{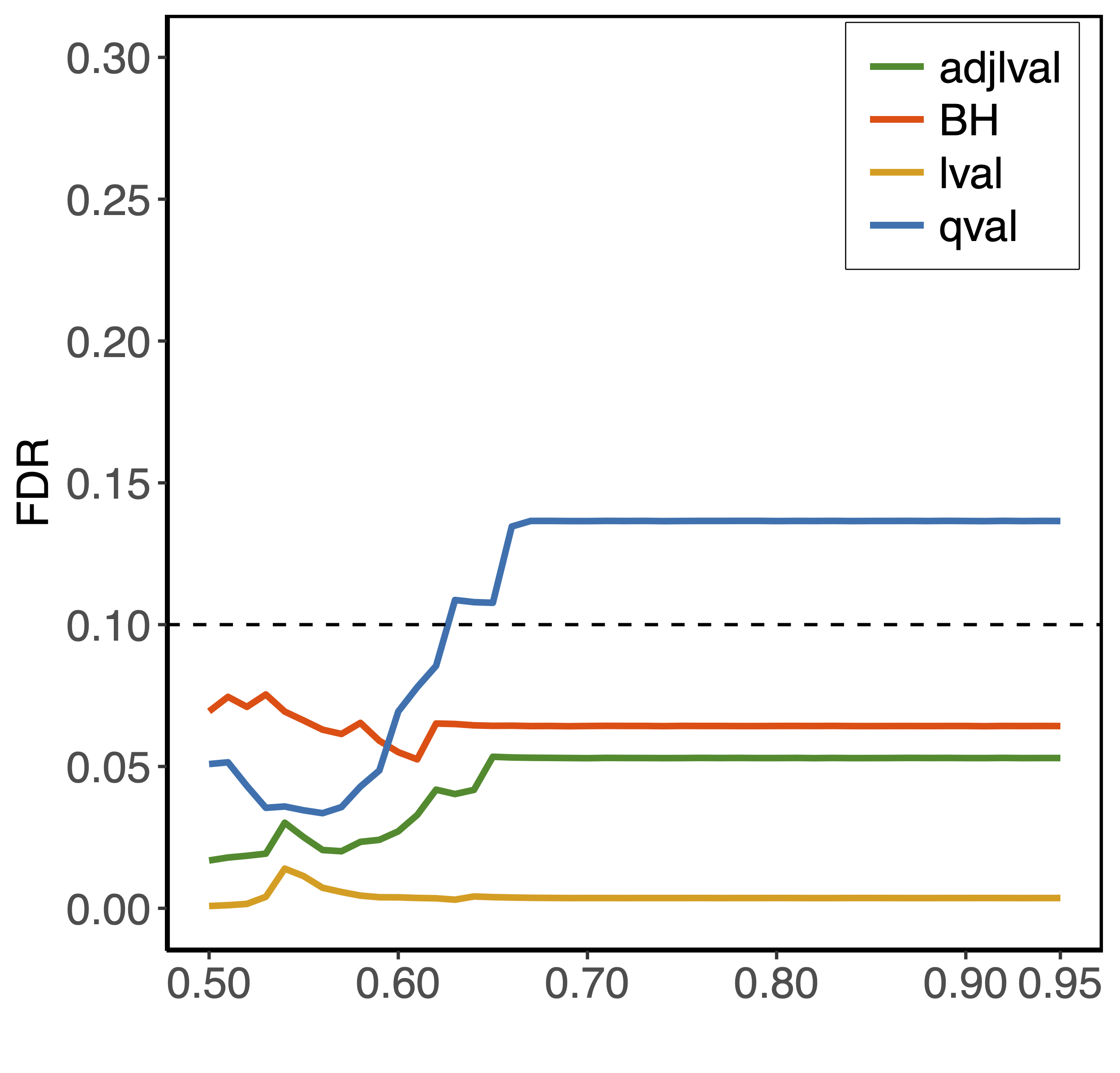}} 
    \subfigure[$m = 1,000$, $s_n/n = 0.001$]{\includegraphics[width=0.33\textwidth]{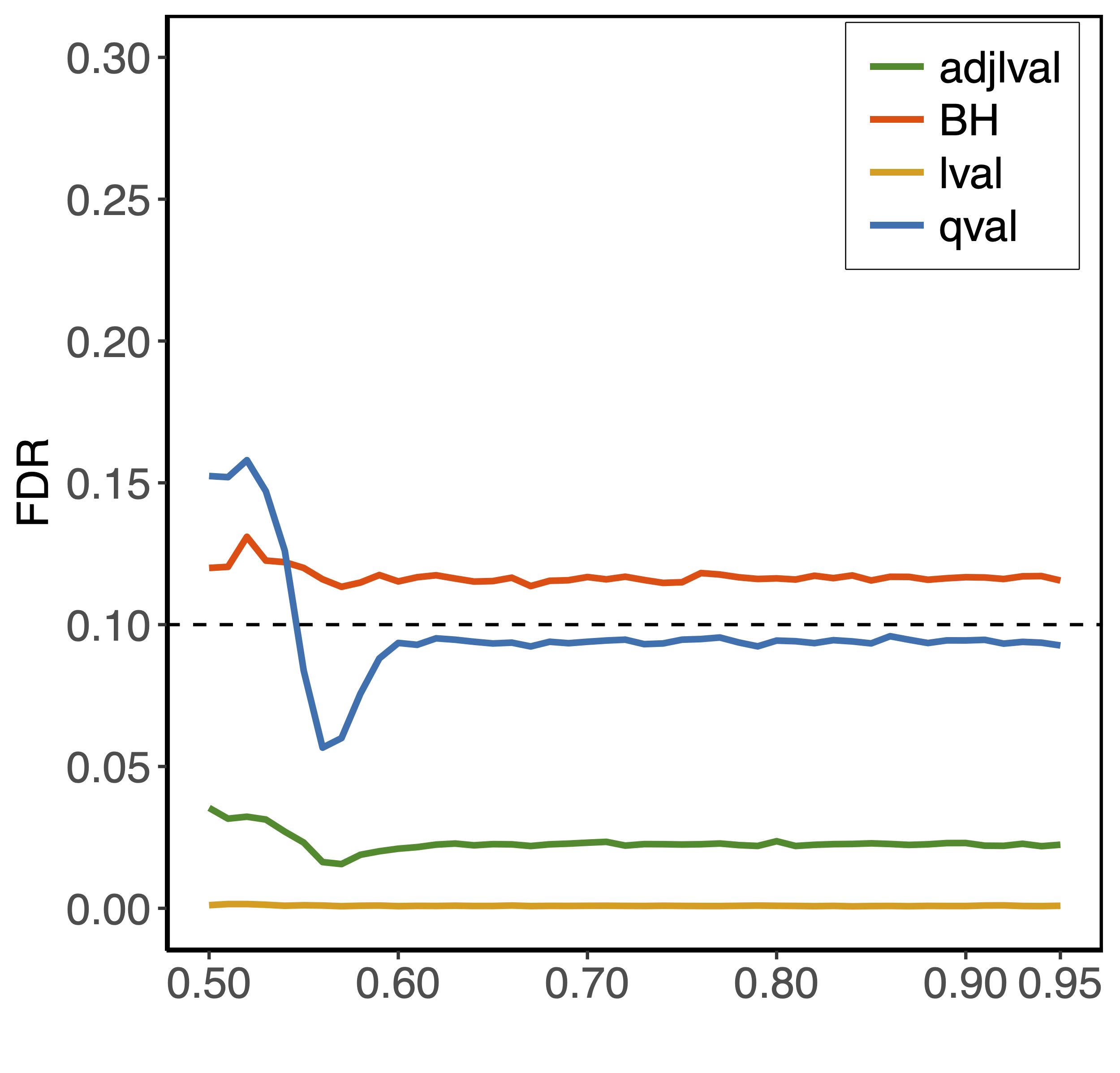}}%
     \subfigure[$m = 1,000$, $s_n/n = 0.1$]{\includegraphics[width=0.33\textwidth]{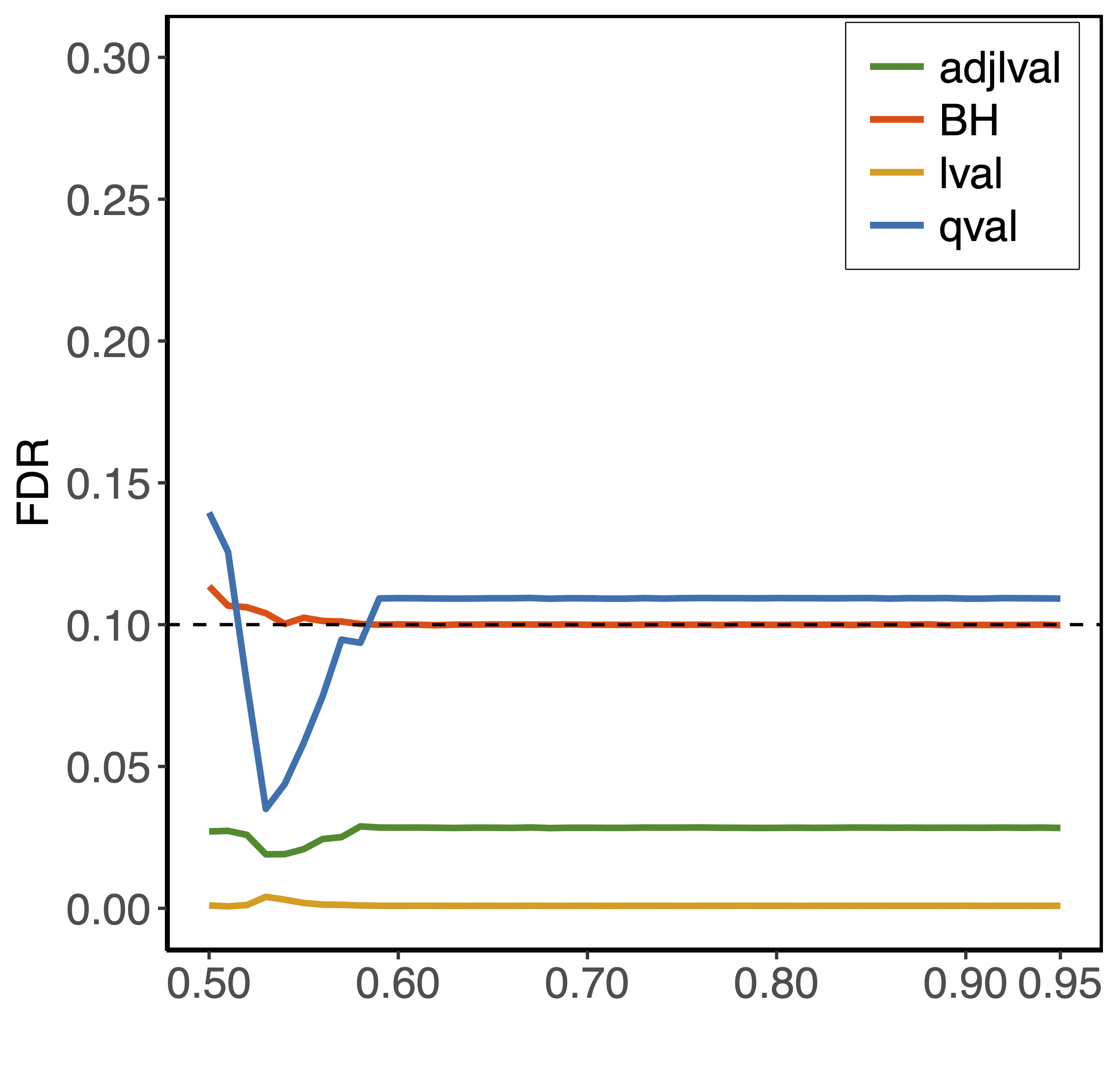}}%
    \subfigure[$m = 1,000$, $s_n/n = 0.5$]{\includegraphics[width=0.33\textwidth]{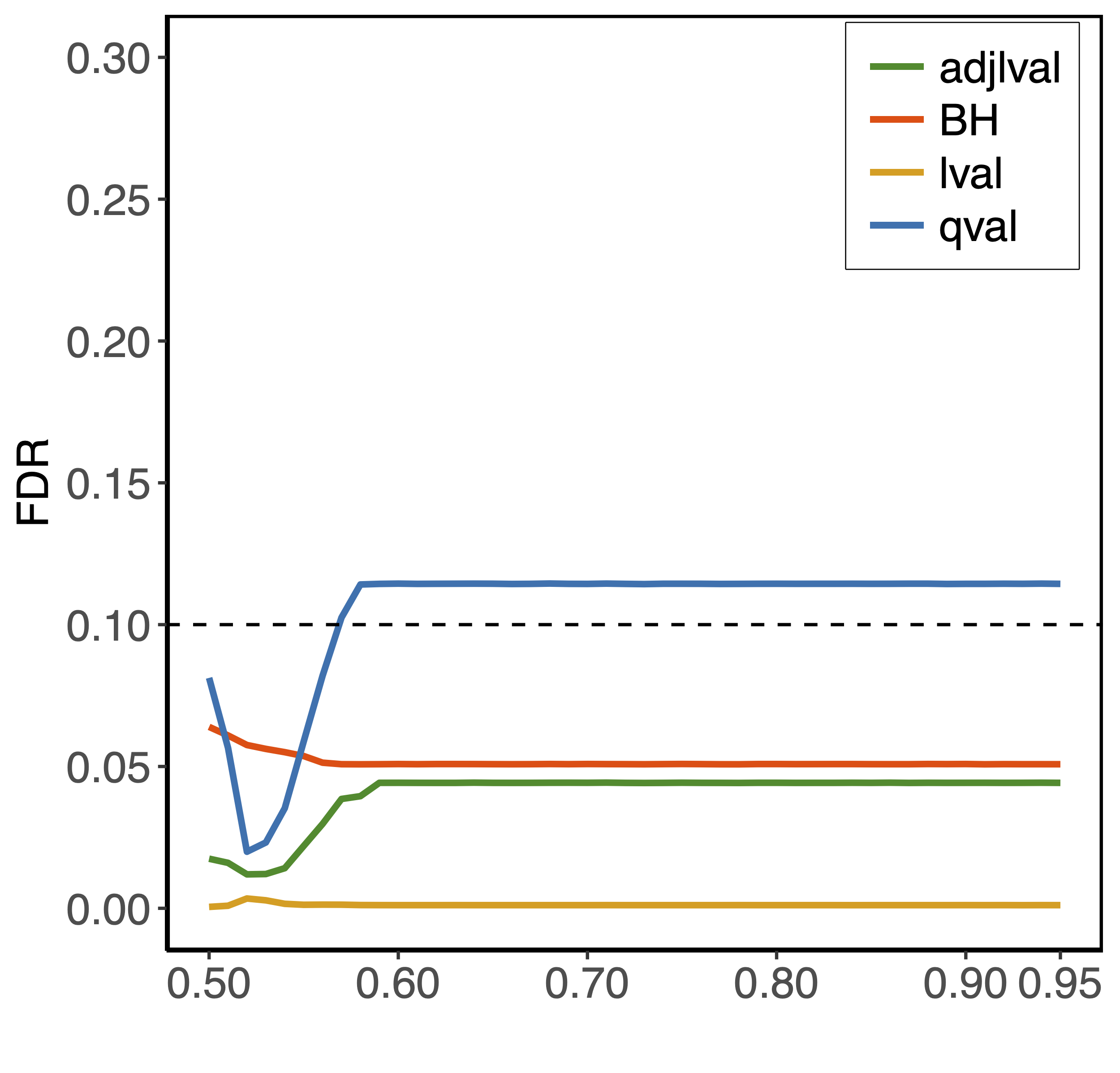}} 
	\caption{The estimated FDR of the $\ell$-value (yellow), the $\adj\ell$-value (green), the $q$-value (blue) and the BH procedures (red) with $t = 0.1$, $m = (\log n)^2 \approx 85, 200$, and $1000$, and $s_n/n = 0.001, 0.1$, $0.5$. }
	\label{fig:compare-all-proc}
\end{figure}

Simulation results are presented in Figures \ref{fig:compare-q-Cl} and \ref{fig:compare-all-proc}. Each of nine subplots in both figures represents a pair of $m$ and $s_n$. 
Each line indicates the average values of 10,000 estimated FDRs from 45 different values of $\vartheta_0$.
Figure \ref{fig:compare-q-Cl} plots the FDR results of $q$-value and $\adj\ell$-value procedures.
The three solid lines in each subplot represents the FDR of the $q$-value procedure at $t=0.2$ (red), $t = 0.1$ (blue), and $t = 0.05$ (green). Similarly, the three dashed lines represent those of the $\adj\ell$-value procedure. Figure \ref{fig:compare-all-proc} compares the FDR of the BH (red), $q$- (blue), $\adj\ell$- (green), and $\ell$-value (yellow) procedures with the significant level set to $t=0.1$.

First, let us compare the $q$-value procedure with the $\adj\ell$-value procedure. 
Theorem \ref{thm:unif-fdr-q-Cl} suggests that the FDR of the $\adj\ell$-value procedure is smaller than that of the $q$-value procedure regardless the value of $\vartheta_0$. Indeed, we also observe this phenomenon in Figure \ref{fig:compare-q-Cl}. Additionally, we found that the $q$-value procedure often overshoots the FDR when signals are dense and $m$ is small (see (b) and (c) in Figure \ref{fig:compare-q-Cl}), while the $\adj\ell$-value procedure can successfully control the FDR below the targeted level regardless the sparsity level and the value of $m$.
Second, we found that when $\vartheta_0$ is slightly above $0.5$, the $q$-value procedure significantly overestimates the FDR, especially in super-sparse and sparse scenarios (see (d) and (g)). This suggests that the constant $K_2$ in Theorem \ref{thm:unif-fdr-q-Cl} can be large especially when $s_n/n$ is small. 
Nonetheless, the FDR of the $q$-value procedure quickly converges to the targeted level as $\vartheta_0$ moves away from $1/2$ (see (f) and (i)).
Third, we observe a bump of the estimated FDR when transiting from a small value of
$\vartheta_0$ to a larger value in both procedures. The depth of the bump for the $q$-value procedure is significant and has not been observed in simulations of the Gaussian sequence model with a similar multiple testing procedure using the empirical Bayes approach (see Figure 1 of \citetalias{cast20MT}).
This bump does not vanish even when $m$ is large, which highlights the major difference between the same procedure between the sparse Binomial model and the sparse Gaussian model. It would be worthwhile to study the phase transition of the multiple testing risk between the small and the large signal regions to understand the differences between the two models. However, given the already substantial length of our analysis, we leave this extension for future work.
Last, despite that our theorems require $m \gg (\log n)^2$. Results in Figure \ref{fig:compare-q-Cl} suggests that this assumption can be further weakened, as we discussed in Remark \ref{rmk:m}.

Next, we compare the two procedures with the $\ell$-value and BH procedures. In Figure \ref{fig:compare-all-proc}, we observe that the estimated FDR of the $\ell$-value procedure is indeed smaller than those of the other two empirical Bayes multiple testing procedures, which aligns with the conclusion in Lemma \ref{lem:unif-fdr-l}. In fact, the FDR of the $\ell$-value procedure is excessively small---almost 0---across all nine scenarios, confirming that the $\ell$-value procedure is overly conservative. 
In the same figure, we also observe that the threshold from the BH procedure (at least the default one used in the $\mathsf{R}$ package) is sensitive to the change of the sparsity level, as its performance varies between different sparsity levels. For example, it overestimates the FDR in the super-sparse case (see (a), (d), (g)), but underestimates the FDR in dense case (see (c), (f), (i)). 
In contrast, the $q$-value and $\adj\ell$-value procedures are more stable in controlling the FDR across various sparsity regions.
We speculate that one could improve the BH procedure by adjusting its rejection region (e.g., through boosting as introduced by \citet{wang22}), and we leave this for future investigation.

\subsection{Simulation when $X_j \overset{ind}{\sim} \text{Bin}(m_j, \theta_j)$ with $m_j$ varies with $j$}
\label{sec:cs}

Consider $X_j \overset{ind}{\sim} \text{Bin}(m_j, \theta_j)$ with $m_1, \dots, m_j \leq m$, the logarithm of the marginal posterior distribution in \eqref{eqn:marginal-post} then becomes
$$
L^C(w) = \sum_{j=1}^n \log b_j(X_j) + \sum_{j=1}^n \log (1+w \beta_j(X_j)),
$$
where $\beta_j(u) = (g_j/b_j)(u) - 1$, $g_j = (1+m_j)^{-1}$ and $b_j(u) = \text{Bin}(u; m_j, 1/2)$, and the score function changes to
$$
S^C(w) = \sum_{j = 1}^n \beta_j(X_j, w), \quad \beta_j(x, w) = \frac{\beta_j(x)}{1 + w\beta_j(x)}.
$$
The MMLE can be estimated by solving
\begin{align}
\label{w-hat-C}
\hat w^C = \argmax_{w \in [1/n, 1]}  L^C(w),
\end{align}
If $\displaystyle \min_{1 \leq j \leq n} m_j \gg (\log n)^2$, then one can easily check that our theorems and lemmas in Sections \ref{sec:unif-fdr} and \ref{sec:fdr-fnr-large} still hold.

\begin{figure}[h!]
\centering
	  \subfigure[$\lambda = 85$, $s_n/n = 0.001$]{\includegraphics[width=0.33\textwidth]{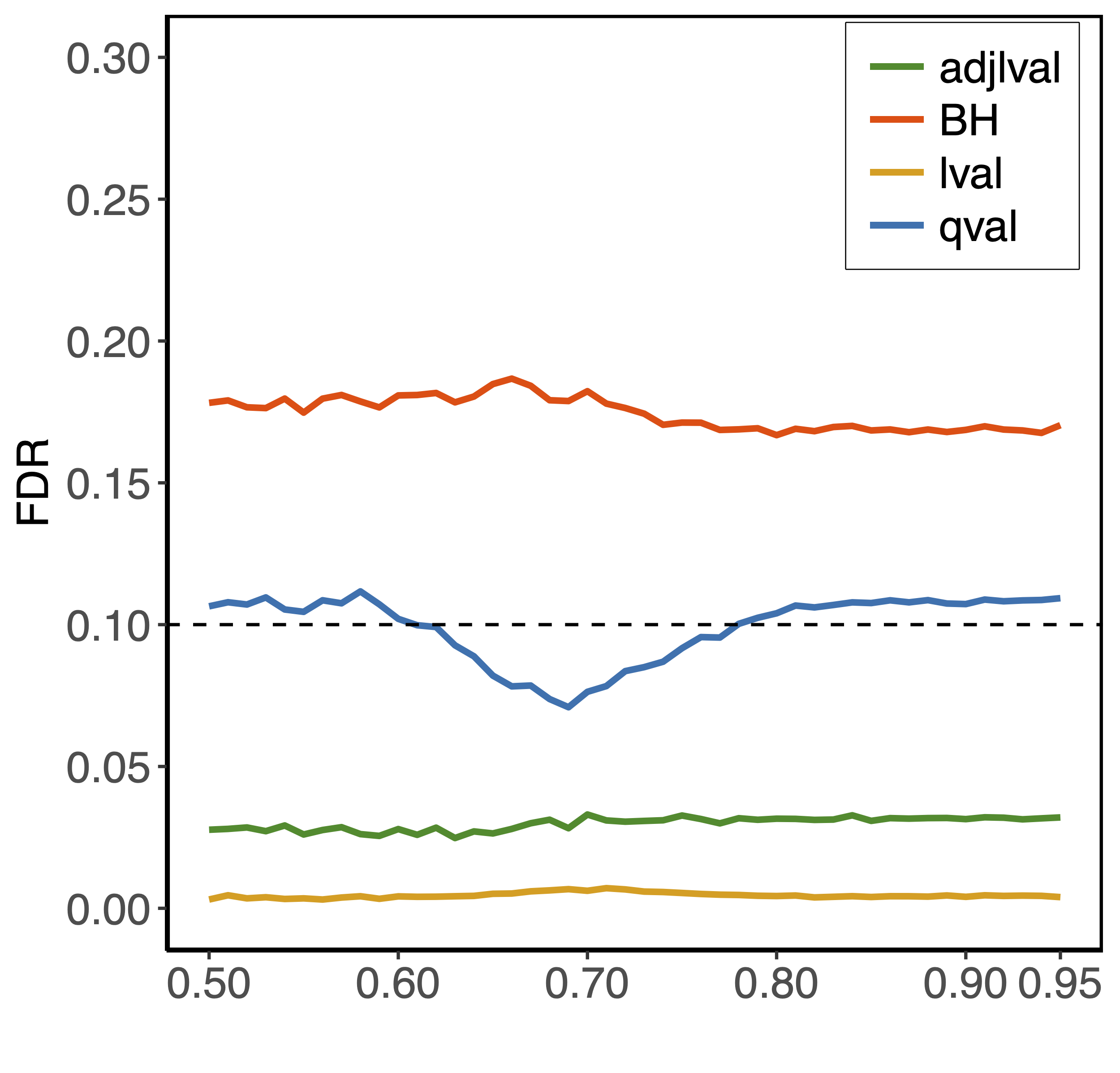}}%
     \subfigure[$\lambda = 85$, $s_n/n = 0.1$]{\includegraphics[width=0.33\textwidth]{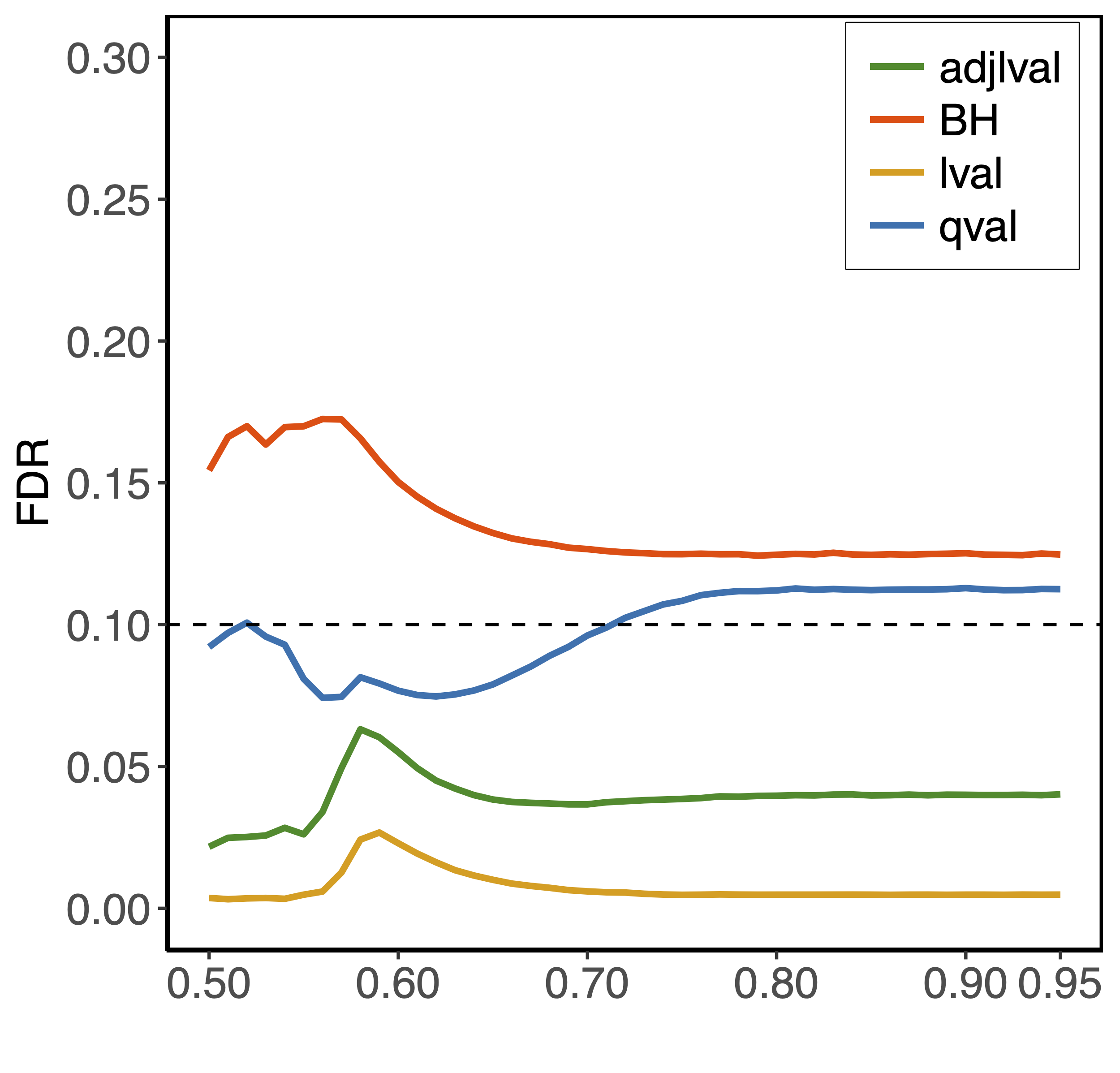}}%
    \subfigure[$\lambda = 85$, $s_n/n = 0.5$]{\includegraphics[width=0.33\textwidth]{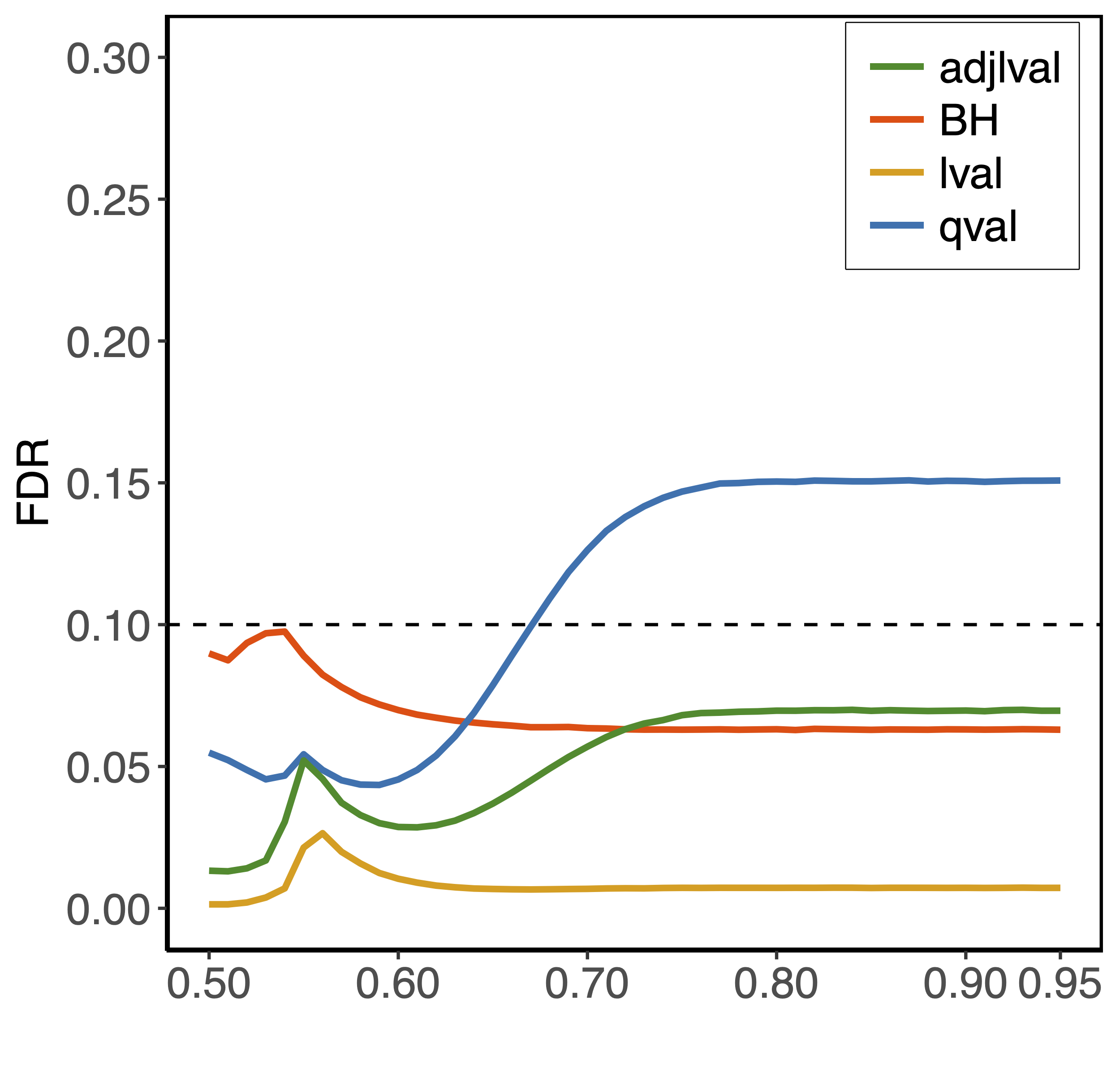}} 
    \subfigure[$\lambda = 200$, $s_n/n = 0.001$]{\includegraphics[width=0.33\textwidth]{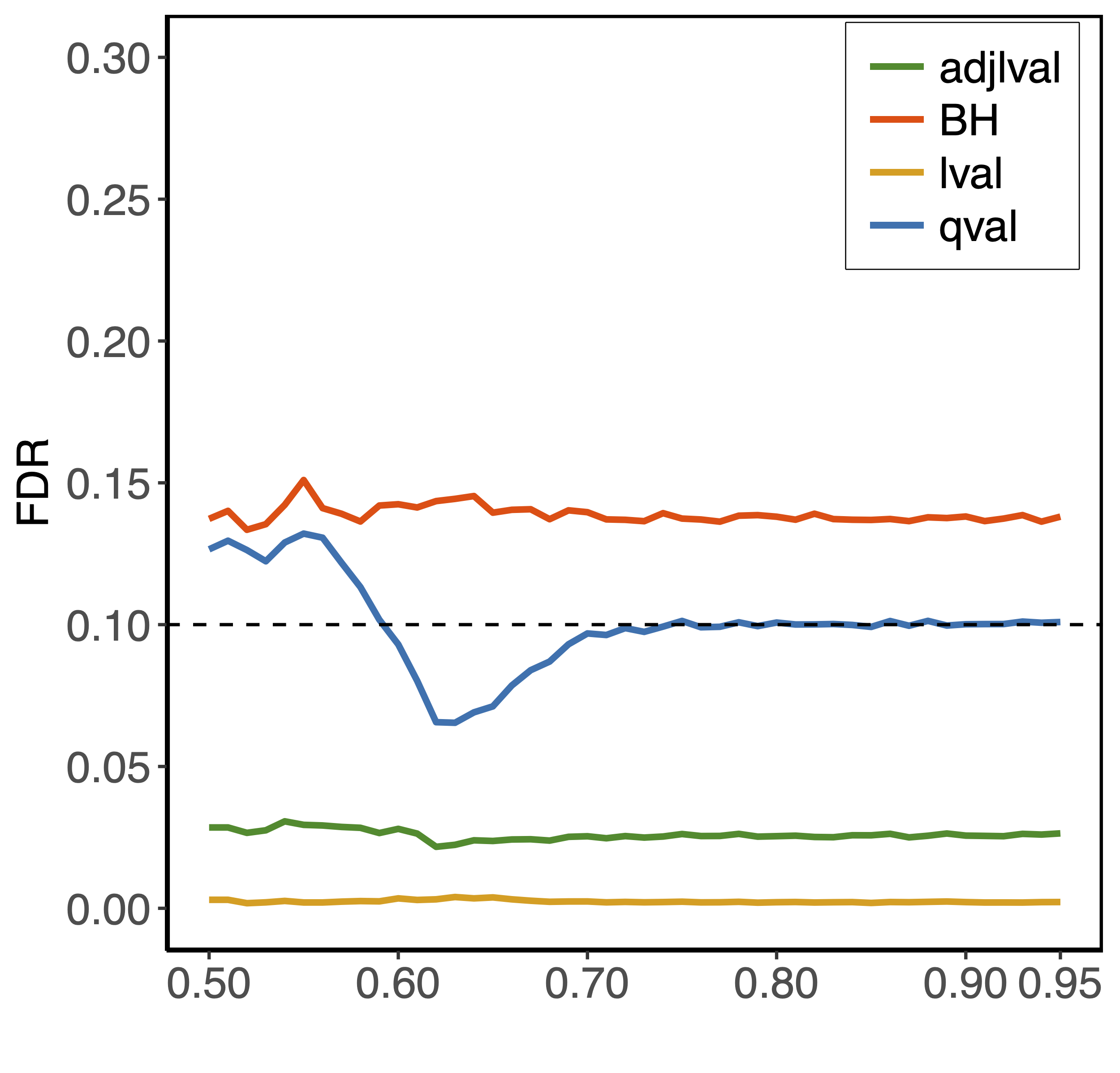}}%
     \subfigure[$\lambda = 200$, $s_n/n = 0.1$]{\includegraphics[width=0.33\textwidth]{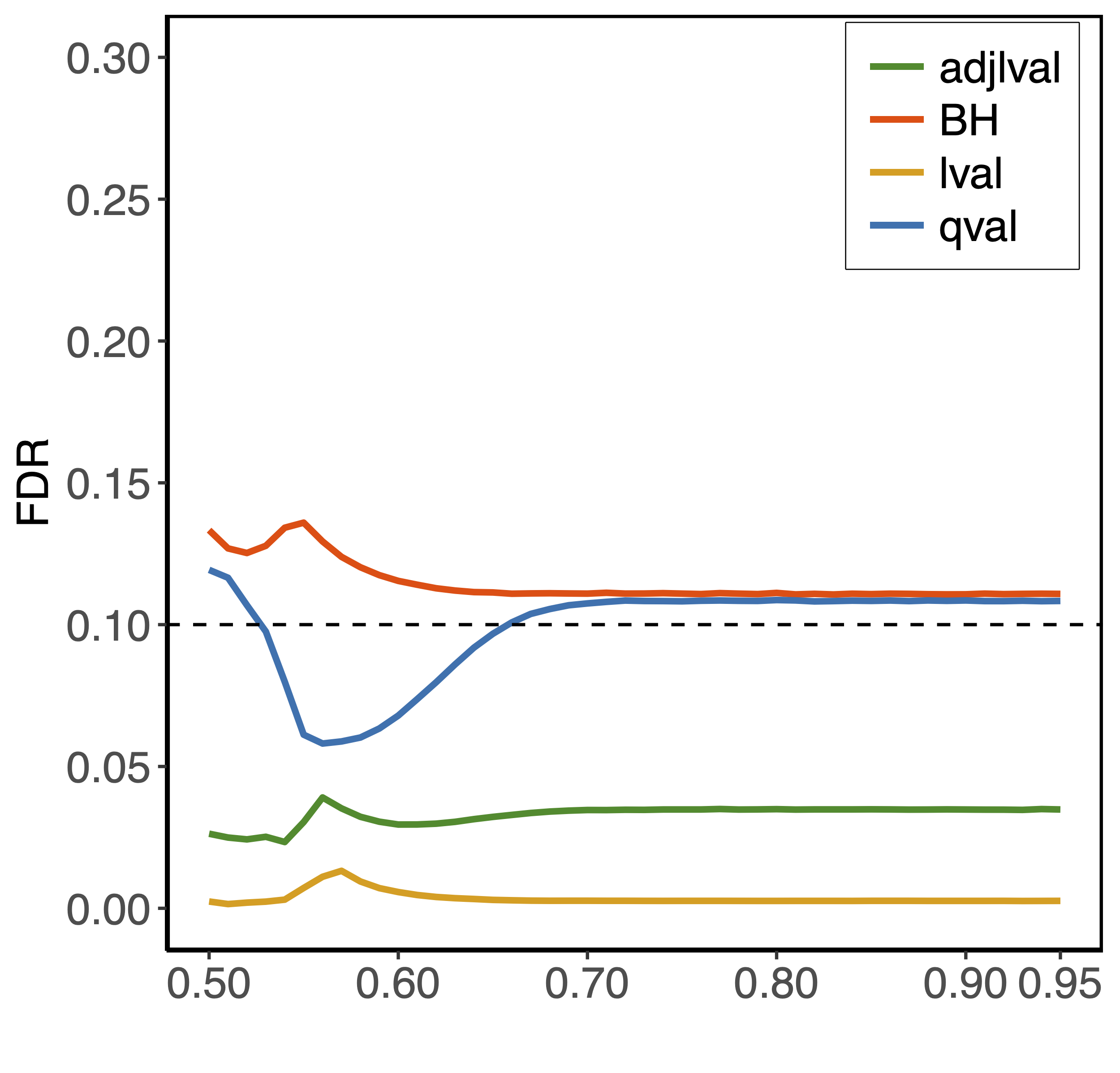}}%
    \subfigure[$\lambda = 200$, $s_n/n = 0.5$]{\includegraphics[width=0.33\textwidth]{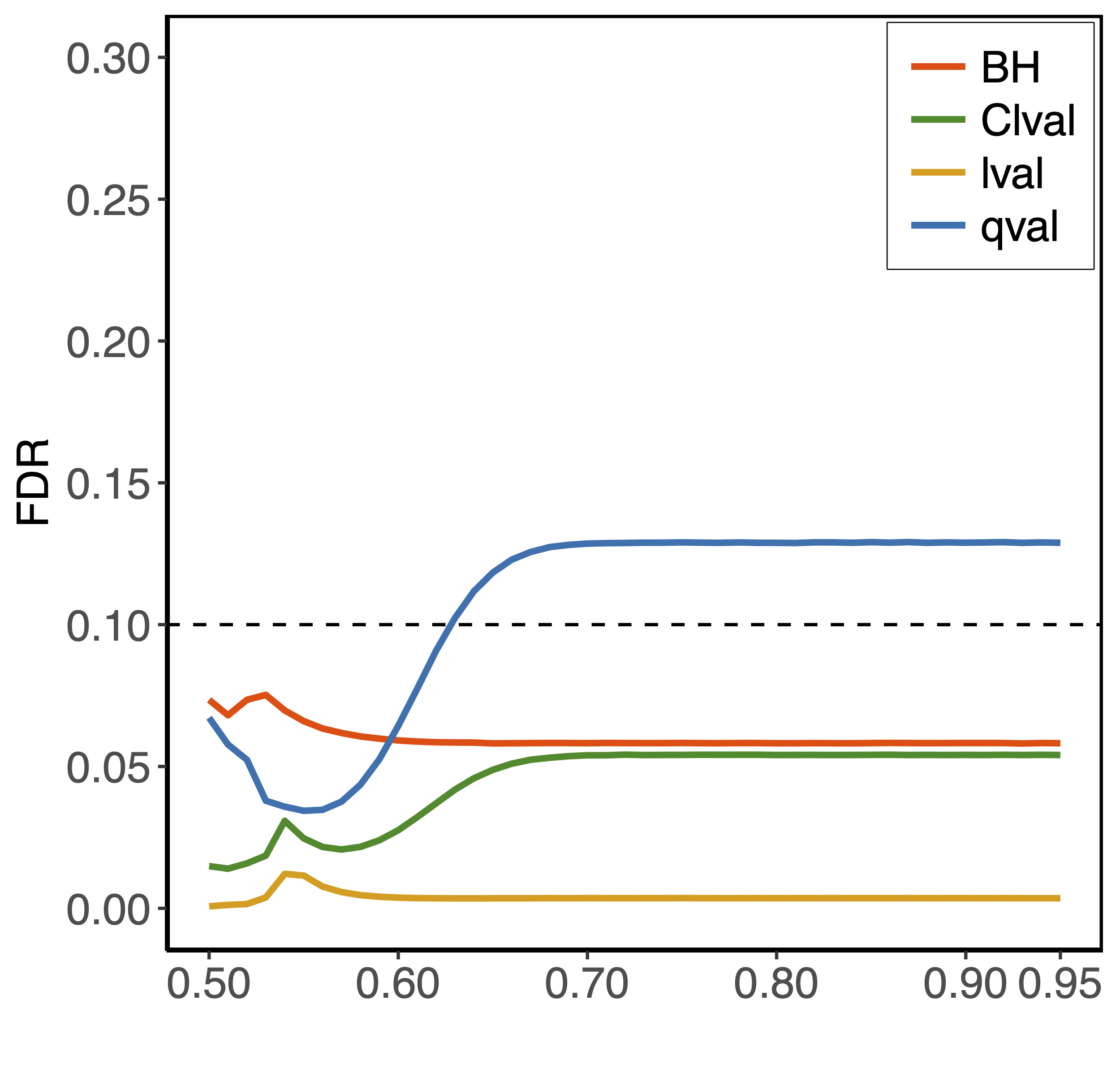}} 
    \subfigure[$\lambda = 1,000$, $s_n/n = 0.001$]{\includegraphics[width=0.33\textwidth]{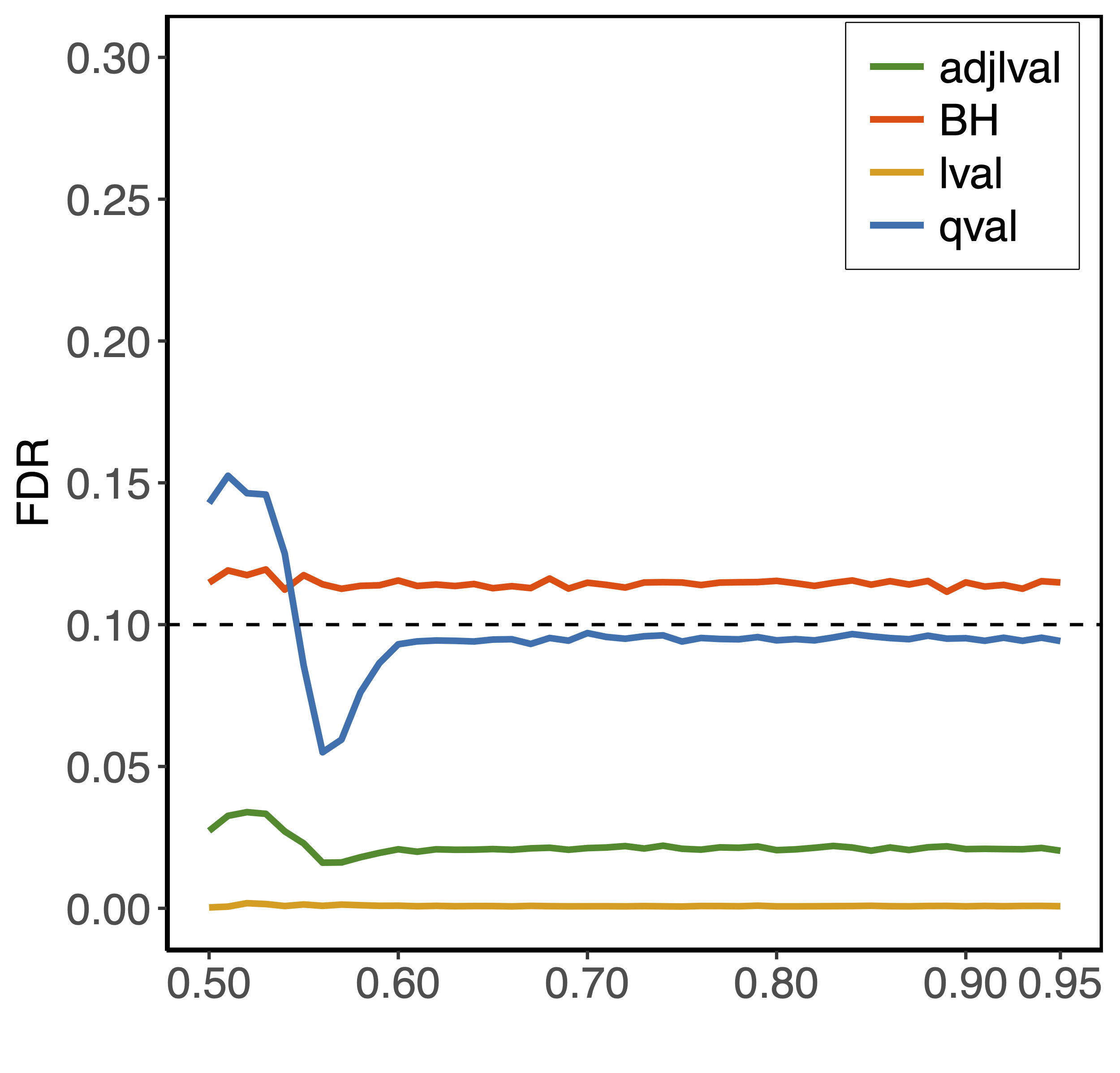}}%
     \subfigure[$\lambda = 1,000$, $s_n/n = 0.1$]{\includegraphics[width=0.33\textwidth]{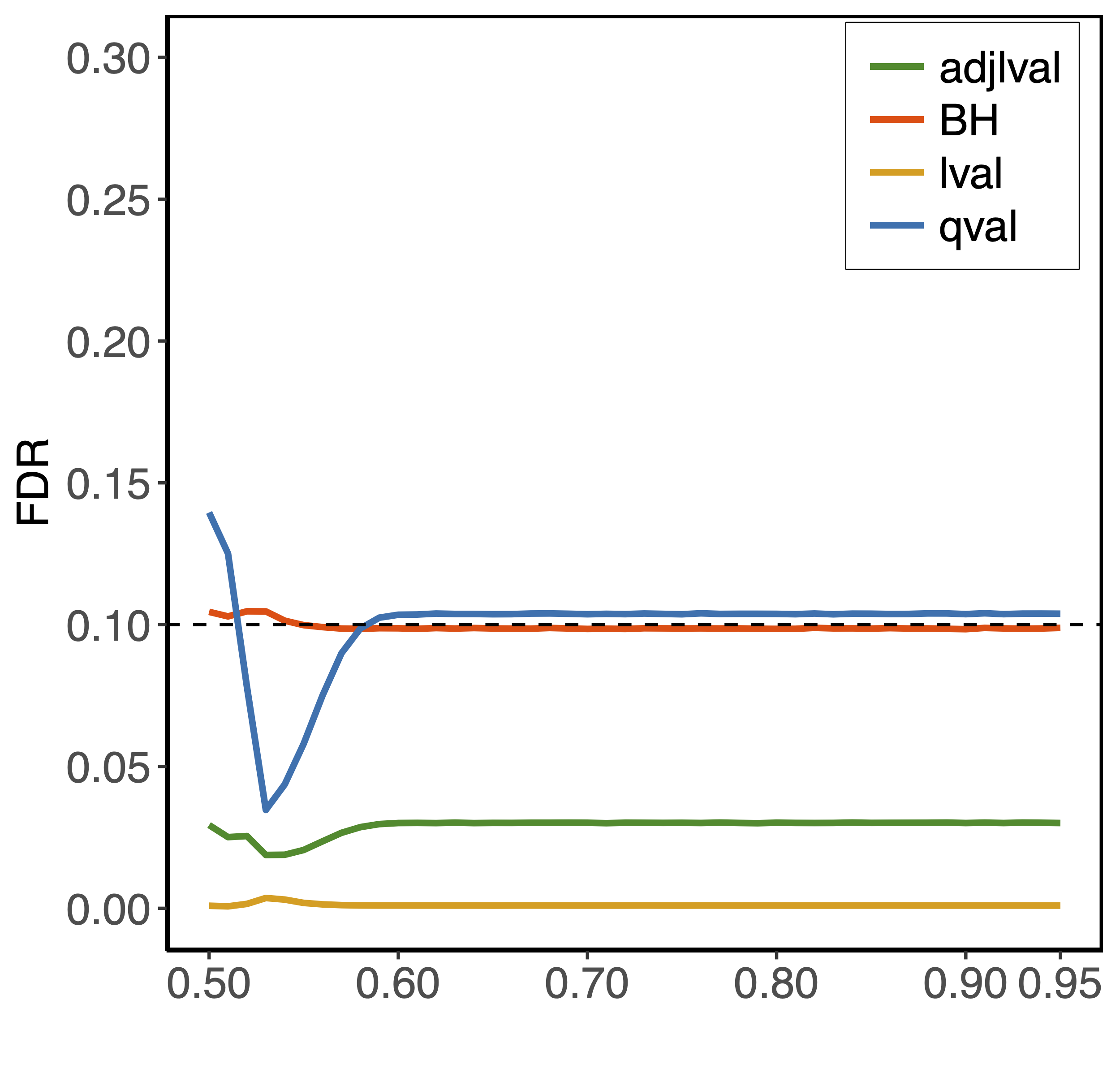}}%
    \subfigure[$\lambda = 1,000$, $s_n/n = 0.5$]{\includegraphics[width=0.33\textwidth]{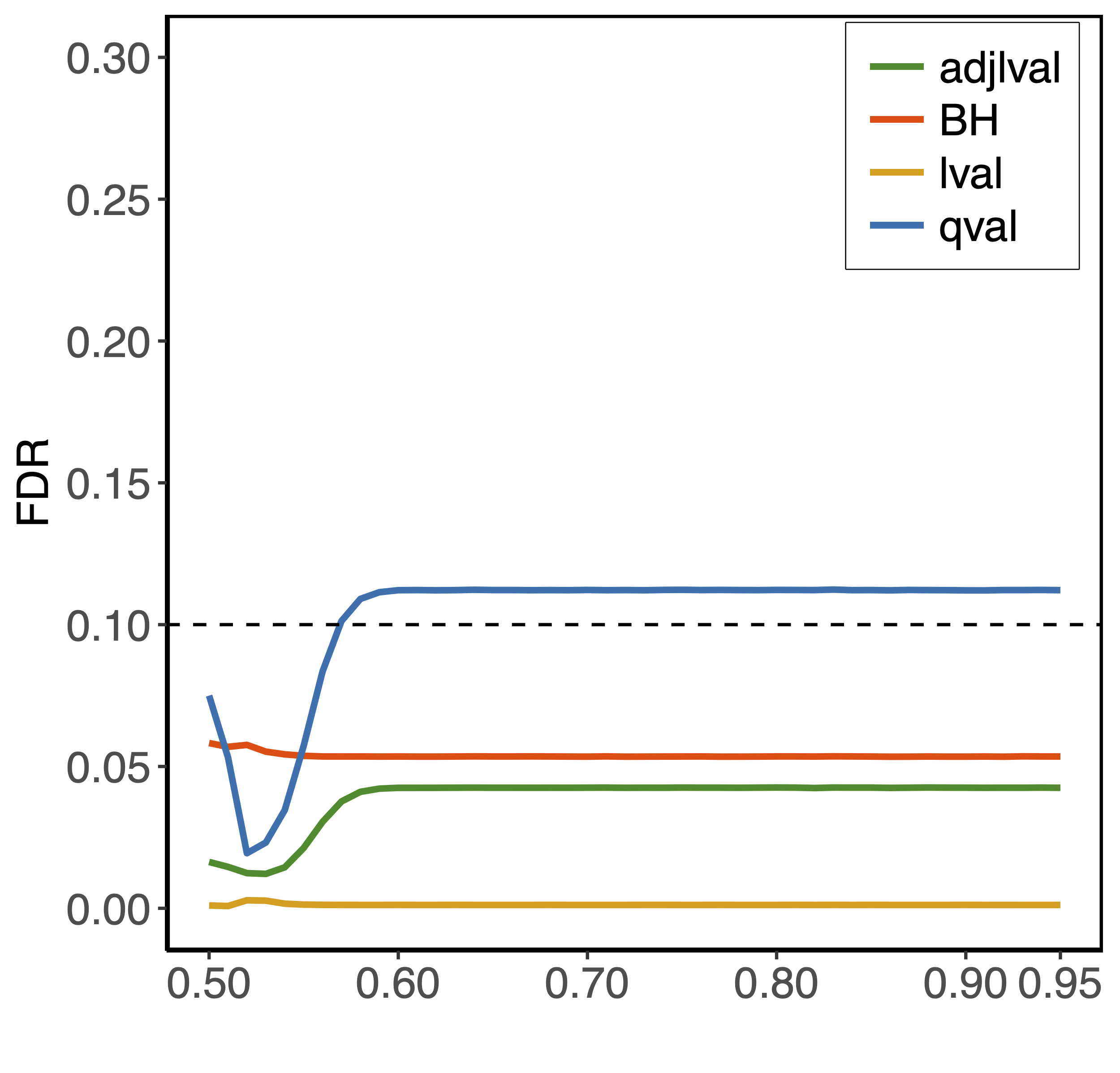}} 
	\caption{The estimated FDR of the $\ell$-value (yellow), the $\adj\ell$-value (green), the $q$-value (blue) and the BH procedures (red) with $m_j \sim \text{Poisson}(\lambda)$ with $\lambda = (\log n)^2 \approx 85, 200, 1000$, $t = 0.1$, and $s_n/n = 0.001, 0.1, 0.5$ respectively. }
	\label{fig:compare-all-proc-cs}
\end{figure}

In Figure \ref{fig:compare-all-proc-cs}, we conduct simulations for this setting.
The simulation is conducted similarly to the previous section, but with two notable differences: 1) instead of choosing a fixed $m$, $m_j$ is randomly drawn from a $poisson(\lambda)$ independently with three different choice of $\lambda = $ 85, $200$, and 1000, and 2) the constant $C_m$ appearing in the $\adj\ell$-value procedure is replaced with $C_{m_j} = \sqrt{2/(\pi m_j)}(m_j + 1)$. 
Despite the variability in $m_j$ across distributions, 
the results from each subplot in Figure \ref{fig:compare-all-proc-cs} are similar to those in Figure \ref{fig:compare-all-proc}. This similarity suggests that the theoretical results obtained in Sections \ref{sec:unif-fdr} and \ref{sec:fdr-fnr-large} remain valid, at least for the current approach for generating $m_j$. 

\section{Discussion}
\label{sec:disc}

We've introduced three empirical Bayes multiple testing procedures for sparse binary sequences, the $\ell$-value procedure, the $\adj\ell$-value procedure, and the $q$-value procedure. In depth frequentist theoretical analysis for these procedures were conducted. Our results suggest that the $q$-value procedure and the $\adj\ell$-value procedure can achieve excellent FDR control for sparse signals, while the $\ell$-value procedure is overly conservative. Our theoretical results were verified through simulation studies. 

It is worth noting that the $\ell$-value uses the spike-and-uniform-slab prior. One might ask whether the threshold of the $\ell$-value can improve by choosing other conjugate priors for the slab density, such as $\gamma \sim \text{Beta}(\alpha, \alpha)$ for some $\alpha > 1$.   
The answer is negative. Let us recalculate $g(x)$ (we use $g(x, \alpha)$ to indicate its dependence on $\alpha$), the quantity that determines the threshold of the $\ell$-value, 
$$
g(x, \alpha) = \int_\theta b_\theta(x) \text{Beta}(\theta; \alpha, \alpha) d\theta = 
{m \choose x} \frac{\Gamma(2\alpha)\Gamma(x + \alpha) \Gamma(m-x+\alpha)}{(\Gamma(\alpha))^2\Gamma(m + 2\alpha)}.
$$
Clearly, $g(x, \alpha)$ is a nonlinear function of $x$ for any $\alpha > 1$. Using the well known approximation for the gamma function $\Gamma(z+a) \sim \Gamma(z) z^{a}$ for any fixed $a$ as $z \to \infty$, then
$$
g(x, \alpha) \sim \left( \frac{x(m-x)}{m^2} \right)^{\alpha - 1} \frac{\Gamma(2\alpha)}{m(\Gamma(\alpha))^2}.
$$ 
In general, this expression is not close to $\sqrt{{2}/(\pi m)}(m+1)$, the multiplying factor for calibrating the threshold of the $\ell$-value procedure. 
Therefore, choosing other values for $\alpha$ in the prior won't resolve the issue encountered with the $\ell$-value procedure. Simulation studies that confirm this point is conducted in Section \ref{sec:ad-sim} of \citet{ning23supp}.

In sum, the present work provides the first theoretical results for empirical Bayes multiple testing on high-dimensional binary outcomes. 
Our results serve as an important initial step for exploring the multiple testing problem for a broader class of models involving discrete outcomes data. Several exciting directions are worth pursuing. First, it is of interest to design new approaches in a similar vein for other discrete high-dimensional models such as the sparse binary regression model in \citet{mukherjee15} and the Ising model in \citet{mukherjee18}.
It would be also important to develop a similar methodology to handle one-sided tests, instead of the two-sided tests that considered in this paper.
The second direction is to study the frequentist coverage of credible sets from the posterior distribution by following the line of research of \citet{vdp17, belister20, cast20} and the minimax risk for the model under various loss functions, e.g., the expected Hamming loss as studied by \citet{butucea18}. 
Furthermore, one could pursue the testing problem from a decision theoretical perspective by characterizing the gap between the risk of a data-driven multiple testing procedure and the risk of an oracle restricted to permutation equivariant decision rules as in \citet{greenshtein09}.

\section*{Acknowledgement}

The author would like to thank Isma\"el Castillo for offering valuable suggestions and travel support at the early stage of this paper. The authors would also like to thank the Associate Editor and three referees for insightful comments.

\begin{supplement}
The supplement \citet{ning23supp} includes the proofs of the results stated in this paper. 
\end{supplement}

\bibliographystyle{chicago}
\bibliography{citation}


\clearpage
\setcounter{equation}{0}
\setcounter{figure}{0}
\setcounter{table}{0}
\setcounter{page}{1}
\makeatletter
\renewcommand{\theequation}{S\arabic{equation}}
\renewcommand{\thefigure}{S\arabic{figure}}
\setcounter{lemma}{0}
\renewcommand{\thelemma}{S\arabic{lemma}}
\setcounter{theorem}{0}
\renewcommand{\thetheorem}{S\arabic{theorem}}
\setcounter{remark}{0}
\renewcommand{\theremark}{S\arabic{remark}}
\setcounter{corollary}{0}
\renewcommand{\thecorollary}{S\arabic{corollary}}
\setcounter{prop}{0}
\renewcommand{\theprop}{S\arabic{prop}}
\renewcommand{\bibnumfmt}[1]{[S#1]}
\renewcommand{\citenumfont}[1]{S#1}
\setcounter{section}{0}
\renewcommand{\thesection}{S\arabic{section}}

\begin{frontmatter}

\title{Supplement to ``empirical Bayes large-scale multiple testing for high-dimensional binary outcome data''}



\author{\fnms{Yu-Chien Bo} \snm{Ning}\ead[label=e2]{bycning@hsph.harvard.edu}}

\affiliation{Harvard University\thanksmark{m1}}

\address{Harvard T.H. Chan School of Public Health \\
677 Huntington Ave\\
Boston, MA 02155\\
\printead{e2}}

\runauthor{Ning}

\begin{abstract}
This supplemental material includes the proofs for the results stated in the main paper as well as additional simulation results.
\end{abstract}

\end{frontmatter}



\section*{Table of Contents}
{\linespread{-1}\selectfont
\startcontents
\printcontents{ }{1}{}
}

\section{Summary of contents}

This document contains the proofs of all the theorems and lemmas presented in the main article. It is organized as follows: 
Section \ref{sec:mt-quantities} summarizes several useful quantities related the multiple testing procedures that are frequently used in our proofs. 
Section \ref{sec:proof-sketch} provides high-level sketch of the proof for the uniform FDR control result for the $q$-value procedure in Theorem \ref{thm:unif-fdr-q-Cl}.
Both Sections \ref{sec:relation-th} and \ref{sec:bound-p-mt} study the quantities related to the threshold in the $\ell$-value, the $\adj\ell$-value, and the $q$-value procedures.
In particular, Section \ref{sec:relation-th} establishes precise upper and lower bounds for $\eta^\ell$, $\eta^q$, and $\eta^{\adj\ell}$, which generalized the asymptotic bounds given in Lemma \ref{lem:th-q-cl}. As the two ratios $(b/g)(\cdot)$ and $(\bar \B/\bar \G)(\cdot)$ govern the behavior of the corresponding thresholds, they are thoroughly analyzed in this section. We also derive non-asymptotic bounds for the difference between $\eta^\ell$ and $\zeta$ and the differences between $\xi$ and $\eta^q$, $\eta^{\adj\ell}$ respectively. 

Our main theorems are proved in Sections \ref{sec:pf-unif-fdr} and \ref{sec:pf-fdr-fnr}. 
The derivation of the tight concentration bound for the MMLE $\hat w$ is given in Section \ref{sec:bound-w-hat}. The three quantities $\tilde m(w)$, $m_1(t, w)$, $m_2(t, w)$ in \eqref{tildem}--\eqref{m2} respectively, which are closely related the score function are thoroughly studied in Section \ref{sec:bound-m1-mtilde}.

In Section \ref{sec:beta}, we exam the function $\beta(u)$ in \eqref{score-fn} closely. The bounds we obtained for this quantity play a crucial role in our analysis of $\tilde m(w)$, $m_1(t, w)$, and $m_2(t, w)$ in Section \ref{sec:bound-m1-mtilde}.
Section \ref{sec:pf-lowerbound} gives the proof of Proposition \ref{prop:lowerbound}. Our proof hinges on obtaining suitable bounds for the inverse tail of the Binomial distribution $\text{Bin}(m, 1/2)$. These bounds are derived in Lemma \ref{lem:inv-bin-bound} in Section \ref{sec:bound-binom}. Additionally, we have gathered several useful lemmas related to the binomial distribution in the same section. Some of these lemmas are crucial in obtaining sharp bounds for a binomial distribution.
We believe those lemmas can be of an independent interest for other studies related to the binomial distribution as well.
Last, auxiliary lemmas are given in Section \ref{sec:aux-lemma} and additional simulation studies are conducted in Section \ref{sec:ad-sim}, where we choose the beta distribution $\text{Beta}(\alpha, \alpha)$ for the spike-and-slab prior with $\alpha = 5$ and $\alpha = 10$.

\section{Introducing several useful multiple-testing related quantities}
\label{sec:mt-quantities}

In this section, we introduce several useful quantities that will be frequently used in this supplemental material.

\begin{enumerate}
\item The posterior distribution is given by 
\begin{align*}
	P^\pi(\theta \given X, w) & = 
	\bigotimes_{j=1}^n
	\left\{
		\ell_j(X_j; w) \delta_{1/2} + (1-\ell_j(X_j; w)) \mathcal{G}_{X_j}
	\right\},
\end{align*}
where $\mathcal{G}_{x} = \Beta(\theta; x + 1, m-x + 1)$, $b(x) = \text{Bin}(x, m/2)$,
$g(x) = (m+1)^{-1}$, and
$$
\ell(x; w) = P^\pi(\theta = 1/2 \given X, w)
= \frac{
(1 - w) b(x)
}{
(1 - w) b(x) + w g(x)
}.
$$

\item The logarithm of the marginal posterior of $w$ is given by 
$$
L(w \given X) 
= \sum_{j=1}^n \left(\log b(X_{j}) + \log (1 + w\beta(X_{j}))\right),
$$
where $\beta(\cdot) = (g/b)(\cdot) - 1$,
and the score function is given by 
\begin{align}
\label{scorefn}
S(w) = \sum_{j=1}^n \beta(X_j, w), 
\quad \beta(x, w) = \frac{\beta(x)}{1 + w \beta(x)}.
\end{align}

\item The thresholds of the $\ell$-, $q$-, and $\adj\ell$-value are given by $t_n^\ell(w, t) = \eta^\ell(r(w,t)) - 1/2$,
$t_n^{\adj\ell}(w, t) = \eta^{\adj\ell}(r(w,t)) - 1/2$, and
$t_n^q(w, t) = \eta^q(r(w,t)) - 1/2$, respectively,
where 
\begin{align*}
\eta^\ell(u) & = \frac{1}{m}\left(\frac{b}{g}\right)^{-1}(u), \\
\eta^{\adj\ell}(u) &= \frac{1}{m}\left(\frac{b}{g}\right)^{-1}\left(\frac{\sqrt{2}(1+m) u}{\sqrt{\pi m}}\right),\\
\eta^q(u) &= \frac{1}{m}\left(\frac{\bar \B}{\bar \G}\right)^{-1}(u).
\end{align*}
These quantities will be studied in Section \ref{sec:threshold}.

\item Let $w^\star$ be the solution of 
$\mathbb{E}_0 S(w^\star) = 0$, where $S(\cdot)$ is the score function in \eqref{scorefn}, 
we introduce $\tilde m(w^\star)$ and $m_1(\theta_{0,j}, w^\star)$ such that
\begin{align}
\label{sol-wstar}
(n - s_n) \tilde m(w^\star) = \sum_{j: \theta_{0,j} \neq 1/2} m_1(\theta_{0,j}, w^\star),
\end{align}
where 
\begin{align}
& \tilde m(w) = - \mathbb{E}_0 \beta(x, w) 
= - \sum_{u=0}^m \beta(x, w) b(x), 
\label{tildem}\\  
& m_1(\theta, w) = \mathbb{E}_\theta \beta(x, w) 
= \sum_{u=0}^m \beta(x, w) b_\theta(x),
\label{m1}
\end{align}
We also define
\begin{align}
m_2(\theta, w) = \mathbb{E}_\theta \beta(x, w)^2 
= \sum_{u=0}^m \beta(x, w)^2 b_\theta(x).
\label{m2}
\end{align}
Note that the above three quantities play an important role in bounding the MMLE $\hat w$.

\item Let $\zeta_n(w) = \sqrt{\frac{1}{2m} \log(1/w)}$,
$\xi_n(w)$ is the solution of $\beta(m/2 + m\xi(w)) = 1/w$,
and $\nu_n(w)$ is the solution of $\beta(m/2 + m\nu_n(w)) = 0$.
The latter two $\xi_n(w)$ and $\nu_w$ will be analyzed in Section \ref{sec:beta}.

\item Let's define the number of false and true discoveries for a vector of tests $\T = (\T_1, \dots, \T_n)$ as
\begin{align}
\label{eqn:FD-TD}
\text{FD}_\T(t, w) = \sum_{j: \theta_{0,j} = 1/2} \T_j(t,  w),
\quad 
\text{TD}_\T(t, w) = \sum_{j: \theta_{0,j} \neq 1/2} \T_j(t,  w).
\end{align}
Also, define the false discovery proportion (FDP) and the false negative proportion (FNP) are given by
\begin{align*}
\FDP(\theta, \T) = \frac{\sum_{j=1}^n \mathbbm{1}\{\theta_j = 1/2\} \T_j}{1 \vee \sum_{j=1}^n \T_j}, 
\quad 
\FNP(\theta, \T) = \frac{\sum_{j=1}^n \mathbbm{1}\{\theta_j \neq 1/2\}(1-\T_j)}{1 \vee \sum_{j=1}^n \mathbbm{1}\{\theta_j = 1/2\}}.
\end{align*}
Last, let $\FDR(\theta_0, \T) = \mathbb{E}_{\theta_0} \FDP(\theta_0, \T)$ and $\FNR(\theta_0, \T) = \mathbb{E}_{\theta_0} \FNP(\theta_0, \T)$ be the false discovery rate (FDR) and the false negative rate (FNR) respectively.
\end{enumerate}

\section{Proof sketch for Theorem \ref{thm:unif-fdr-q-Cl}}
\label{sec:proof-sketch}

Before presenting the detailed proof, we provide a high-level sketch of the proof for the uniform FDR control result for the $q$-value procedure stated in Theorem~\ref{thm:unif-fdr-q-Cl}. The proof of the other two procedures only requires subtle modification: replacing the corresponding threshold with the threshold of the $q$-value.  
The proof consists of two main parts. First, we derive a tight concentration bound for $\hat w$.
By using concentration arguments, we show $\hat w$ is close to $w^\star$, which allows us to replace the empirical quantity with its expectation. 
We then introduce $w_1$ and $w_2$ to be the solutions of  
\begin{align}
\label{eqn:w1w2}
\sum_{j \in S_0} m_1(\theta_{0,j}, w) = (1-\kappa)(n-s_0)\tilde m(w), \quad \
\sum_{j \in S_0} m_1(\theta_{0,j}, w) = (1+\kappa)(n-s_0)\tilde m(w).
\end{align}
respectively. 
If their solutions exist, 
we show both $w_1 \asymp s_n/n$ and $w_2 \asymp s_n/n$ in Lemma \ref{lem:bound-w1-w2}. 
By the monotonicity of $m_1$ and $\tilde m$ and choosing $\kappa$ to be sufficiently small, then $w^\star \asymp s_n/n$.
The details of this part are given in Section \ref{sec:bound-w-hat}.

Next, we obtain an upper bound for the FDR uniformly over all $\theta \in l_0[s_n]$. By monotonicity and the relation between $w_1, w_2$ and $\hat w$, we then show that
\begin{align*}
\sup_{\theta_0 \in l_0[s_n]} \FDR(\theta_0, \T^q(t, \hat w)) & \lesssim \frac{\mathbb{E}_{\theta_0} \text{FD}_{\T^q}(t, \hat w)}{\mathbb{E}_{\theta_0} \text{TD}_{\T^q}(t, \hat w)} + o_p(1) 
\lesssim \frac{\mathbb{E}_{\theta_0} \text{FD}_{\T^q}(t, w_1)}{\mathbb{E}_{\theta_0} \text{TD}_{\T^q}(t, w_2)} + o_p(1).
\end{align*}
To bound the first term of the upper bound in the last display, we use the assumption $w_2 \geq 1/n$ (and hence $\hat w > 1/n$) to show that the denominator is bounded below by $C(n - s_0) w_2$. We then show that the ratio is bounded by $\displaystyle \frac{(n - s_0) w_1 t}{(n - s_0) w_2} \lesssim t$.
A detailed derivation and a sharper upper bound for the FDR are provided in Section \ref{sec:pf-thm1}. 

So far, we assume that solutions for the two equations in \eqref{eqn:w1w2} exist. In Section \ref{sec:pf-thm1}, we also need to analyze the case where the solutions do not exist. In Lemma \ref{lem:bound-w1-w2}, we show that the MMLE $\hat w$ is bounded by $w_0$ with a large probability for $w_0 \in [1/n, \rho_n/n]$, where $\rho_n = o(n)$. By plugging-in this bound, we then show that the FDR is bounded by $Ct \log (1/t)$ for some positive constant $C > 0$. By combining the upper bounds from both cases and noting that $\log(1/t) > 1$, we obtain the final result $\sup_{\theta_0 \in l_0[s_n]} \FDR(\theta_0, \T^q(t, \hat w)) \lesssim t \log (1/t)$.

\section{Relations between $\eta^\ell$, $\eta^q$, $\eta^{\adj\ell}$, $\zeta$, and $\xi$}
\label{sec:relation-th}

In this section, we establish the relation between $\eta^\ell$, $\eta^q$, $\eta^{\adj\ell}$, $\zeta$, and $\xi$. 
First, we examine the monotonicity property of the two functions $(b/g)(\cdot)$ and $(\bar\B /\bar\G)(\cdot)$. 
We then obtain non-asymptotic bounds for $\eta^\ell$, $\eta^q$, $\eta^{\adj\ell}$ in Section \ref{sec:threshold}.
Next, we establish the relation between the three quantities and $\zeta$, $\xi$ in Section \ref{sec:compare-th}. 
Last, we prove Lemmas \ref{lem:th-l}, \ref{lem:th-lval} and \ref{lem:th-q-cl} in Section \ref{sec:pf-three-lemmas} and Lemma \ref{lem:compare-thresholds} in Section \ref{sec:pf-compare-th}.
Below and throughout this paper, we assume that $m$ is even for simplicity. Our results can be easily extended to the case where $m$ is odd with minimal modifications.

\subsection{Monotonicity for $(b/g)(\cdot)$ and $(\bar\B /\bar\G)(\cdot)$}
\label{sec:mono-ratios}

\begin{lemma}
\label{lem:ratio-phi-g}
Let $b(x) = \text{Bin}(m, 1/2)$ and $g(x) = \int b_\theta(x) \gamma(\theta) d\theta = (m + 1)^{-1}$,
then the function $(b/g)(m/2 + |y|)$ is symmetric at $y = 0$ and is monotone increasing on $y \in [-m/2, 0)$ and monotone decreasing on $y \in [0, m/2)$.
\end{lemma}
 
\begin{proof}
Since $|y|$ is symmetric around 0 and note that $g = (1+m)^{-1}$ is a constant, $b/g(m/2 + |y|)$ is also symmetric around 0. Due to $b(\cdot)$ is the binomial distribution centered at $m/2$, we immediately obtain that $b/g(m/2 + |y|)$ is monotone decreasing on $y \in [0, m/2)$. Its symmetry about $m/2$ implies that the function is monotone increasing on $y \in [-m/2, 0)$. 
\end{proof}

\begin{lemma}
\label{lem:ratio-barPhi-barG}
Let $\bar \B(x)$ be the upper tail probability of $\text{Bin}(m, 1/2)$ and $\bar \G(x) = (m - x)/(m+1)$,
then the function $(\bar\B /\bar\G)(m/2 + |y|)$ is symmetric at $y = 0$ and is monotone increasing on $y \in [-m/2, 0)$ and monotone decreasing on $y \in [0, m/2)$.
\end{lemma}
 
 \begin{proof}
 It is trivial to verify $(\bar{\mathbf{B}}/\bar\G)(m/2 + |y|)$ is symmetric at $y = 0$, we thus omit the proof. By plugging-in the expressions of $\bar{\mathbf{B}}(\cdot)$ and $\bar \G(\cdot)$ respectively, we have 
 \begin{align*}
 	& \frac{\bar{\mathbf{B}}}{\bar\G}(m/2 + y) - \frac{\bar{\mathbf{B}}}{\bar\G}(m/2 + y + 1) \\
	& \quad = (m+1) \left(
	\frac{\sum_{z= m/2 + y}^m b(z)}{m/2 - y +1} 
	- \frac{\sum_{z= m/2 + y +1}^m b(z)}{m/2 - y }
	\right) \\
	& \quad = (m+1) \left(
	\frac{(m/2 - y)b(m/2 + y) - \sum_{z= m/2 + y +1}^m b(z)}{(m/2 - y +1)(m/2 - y )} 
	\right).
 \end{align*}
The last line is always positive as $b(z) < b(m/2 + y)$ for all $z \geq m/2 + y+1$. 
Thus, $(\bar{\mathbf{B}}/\bar\G)(m/2 + |y|)$ is monotone decreasing on $y \in [0, m/2]$.
The proof of the ratio is monotone increasing is similar.
 \end{proof}

\subsection{Bounding $\eta^\ell$, $\eta^q$, and $\eta^{\adj\ell}$}
\label{sec:threshold}

\begin{lemma}
\label{lem:bound-etal}
Let $\displaystyle \eta^{\ell}(u) = \frac{1}{m}(b/g)^{-1}(u)$ for $u \in (0, 1)$, where 
$b(u) = \text{Bin}(m, 1/2)$ and $g(u) = (1+m)^{-1}$,
if $|\eta^\ell(u) - 1/2| \leq \eta_\circ$ for some $u \in (0, 1)$ and a fixed constant $\eta_\circ < 1/2$, 
then 
\begin{align*}
\eta^\ell(u) 
& \leq \frac{1}{2} + \sqrt{\frac{1}{2m} \left(
\log (1/u) +  \log \left(\frac{\sqrt{2}(1+m)}{\sqrt{\pi m (1-4 \eta_\circ^2)}}\right) + \frac{1}{12m}
\right)}, \\
\eta^\ell(u) 
& \geq \frac{1}{2} + \sqrt{\frac{1}{2m} \left(
\log (1/u) +  \log \left(\frac{\sqrt{2}(1+m)}{\sqrt{\pi m (1-4 \eta_\circ^2)}}\right) 
\right) 
\left(
1 + \frac{8\eta_\circ^2}{3(1-4\eta_\circ^2)^2}
\right)^{-1}}.
\end{align*}
Moreover, if $m \gg \log^2(1/u)$ for some $u \in (0, 1)$, then 
\begin{align}
\label{lemma-S3}
\eta^\ell(u)  - \frac{1}{2} \sim  \sqrt{\frac{1}{2m} \left(
\log \left(\frac{1}{u} \right) 
+ \log \left(\frac{\sqrt{2}(1+m)}{\sqrt{\pi m}}\right)
\right)}.
\end{align}
\end{lemma}

\begin{proof}
We first write $\eta^\ell(u) = \frac{1}{m} \left(b/g\right)^{-1}(u)$ 
as $1/u = (g/b)(m \eta^\ell(u))$, 
which implies 
\begin{align}
\label{pf:bound-etal}
{m \choose {m \eta^\ell(u)}} = \frac{2^m u}{m + 1}.
\end{align}
Let $\tilde \eta^\ell := \tilde \eta^\ell(u) = \eta^\ell(u) - 1/2$, 
then by \eqref{binomial-bound} in Lemma \ref{lem:bin-coef-bound},
we have
\begin{align}
\label{pf:bound-etal2}
{m \choose {m \eta^\ell(u)}} =
\frac{\sqrt{2} e^{- m T(1/2 + \tilde \eta^\ell, 1/2) + m\log 2 + \omega(\tilde \eta^\ell)}}
{\sqrt{\pi m (1-4(\tilde \eta^\ell)^2)}},
\end{align}
where $T(a, p) = a\log(a/p) + (1-a)\log((1-a)/(1-p))$ and $\omega(\tilde \eta^\ell) \sim 1/(12m)$.
By the second point in $(d)$ in Lemma \ref{lem:entropy-bound}, one can bound $T(\cdot, \cdot)$ in the last display by 
$$
2 (\tilde \eta^\ell)^2 \leq T(1/2 + \tilde \eta^\ell, 1/2) \leq 2(\tilde \eta^\ell)^2 + \frac{8 (\tilde \eta^\ell)^4}{3 (1-4(\tilde \eta^\ell)^2)}.
$$
Using the last display, we obtain both upper and lower bounds for \eqref{pf:bound-etal2}. Using \eqref{pf:bound-etal}, we thus have
$$
2m(\tilde \eta^\ell)^2 \leq \log (1/u) +  \log \left(\frac{\sqrt{2}(1+m)}{\sqrt{\pi m (1-4(\tilde \eta^\ell)^2)}}\right) + \frac{1}{12m},
$$
and 
$$
2m(\tilde\eta^\ell)^2\left(1 + \frac{8(\tilde \eta^\ell)^2}{3(1-4(\tilde\eta^\ell)^2)^2}\right) \geq \log (1/u) + \log \left(\frac{\sqrt{2}(1+m)}{\sqrt{\pi m(1-4(\tilde \eta^\ell)^2)}}\right).
$$ 
Using the assumption $\tilde \eta^\ell < \eta_\circ$, we obtain the upper and lower bounds for $\eta^\ell(u)$.
If $m \gg \log^2(1/u) \to 0$, then $m(\tilde\eta^\ell)^4 \to 0$. We thus have
$1 - 4(\tilde \eta^\ell)^2 \sim 1$ and $\frac{m(\tilde \eta^\ell)^4}{(1-4(\tilde\eta^\ell)^2)^2} \to 0$. 
The last two displays imply \eqref{lemma-S3}.
\end{proof}

\begin{lemma}
\label{lem:bound-etaq}
Let $\displaystyle \eta^q(u) = \frac{1}{m} (\bar\G/\bar \B)^{-1}(u)$ for $u \in (0, 1)$, where 
$\bar \G(k) = \sum_{x =k}^m g(x)$ and $\bar \B(k) = \sum_{x=k}^m b(x)$, if 
$\eta^q(u) -1/2 \leq \eta_\circ$, $0 < \eta_\circ < 1/2$ is a fixed constant, 
then we have
\begin{align*}
\eta^q(u) & \leq \frac{1}{2} + \sqrt{\frac{\log (1/u) + A_1(m, \eta_\circ) - \log (\sqrt{\log(1/u) + A_2(m, \eta_\circ)})}{2m}},\\
\eta^q(u) & \geq \frac{1}{2} + \sqrt{\frac{\log (1/u) + A_2(m, \eta_\circ) - \log (\sqrt{\log(1/u) + A_1(m,\eta_\circ)})}{2m}},
\end{align*}
where 
$A_1(m, \eta_\circ) = - \log \left( \sqrt{\pi}(1/2 - \eta_\circ) \sqrt{1-4\eta_\circ^2}\right)$ and 
$A_2(m, \eta_\circ) =  \log \left( \sqrt{2/\pi} \eta_\circ\right) - (12m)^{-1}$.
Moreover, if $m \to \infty$, then 
$$
\eta^q(u) - \frac{1}{2} 
\sim  \sqrt{\frac{1}{2m} \left(\log \left(\frac{1}{u}\right) \right)}.
$$
\end{lemma}
\begin{proof}
By the definition of $\eta^q(u)$, we obtain $\bar{\mathbf{B}}(m\eta^q(u) ) = u \bar\G(m\eta^q(u))$.
By plugging-in the expression of $\bar{\mathbf{B}}(\cdot)$, we have $\bar\G(m\eta^q(u)) = 1 - m\eta^q(u) /(m+1)$.
Let $\tilde \eta^q: = \tilde \eta^q(u) = \eta^q(u) - 1/2$,
using Lemma \ref{lem:binom-tail} and then the second point in $(d)$ in Lemma \ref{lem:entropy-bound}, 
one obtains
\begin{align}
\label{pf:bound-etaq-1}
\frac{\sqrt{2}e^{-2m(\tilde \eta^q)^2- (12m)^{-1}}}{\sqrt{\pi m (1- 4(\tilde \eta^q)^2)}} 
\leq u \bar \G(m  \eta^q) \leq
\frac{\eta^q \sqrt{2}e^{-2m(\tilde \eta^q)^2}}{2 \tilde \eta^q \sqrt{\pi m (1- 4(\tilde \eta^q)^2)}}.
\end{align}
By plugging-in the expression of $\bar\G(m \eta^q)$, 
the upper bound in \eqref{pf:bound-etaq-1} implies 
\begin{align*}
2m(\tilde \eta^q)^2 + \log \left( \sqrt{2m} \tilde \eta^q \right)
& \leq \log (1/u) + \log \left(\frac{1/2 + \tilde \eta^q}{1/2-\tilde \eta^q}\right) - \log \left(\sqrt{\pi (1-4(\tilde \eta^q)^2)}\right) \\
& \leq 
 \log (1/u) - \log\left(\sqrt{\pi} (1/2 - \eta_\circ) \sqrt{1-4\eta_\circ^2}\right),
\end{align*}
which we used $ - \log (1/2-\tilde \eta^q(u)) <  - \log (1/2 -\eta_\circ)$, $\log(1/2 + \tilde \eta^q) \leq 0$, 
and $-\log(1-4(\tilde \eta^q)^2) \leq - \log(1 - 4\eta_\circ^2)$.
Similarly, 
the lower bound in \eqref{pf:bound-etaq-1} can be bounded as follows:
\begin{align*}
2m(\tilde \eta^q)^2 + \log \left(\sqrt{2m} \tilde \eta^q \right)
& \geq \log (1/u) + \log(\sqrt{2/\pi} \tilde \eta^q) -\frac{1}{12m}\\
& \geq \log (1/u) + \log(\sqrt{2/\pi} \eta_\circ) -\frac{1}{12m}.
\end{align*}
The previous two displays implies the upper and the lower bounds for $\tilde \eta^q$.
If $m \to \infty$, then we obtain the second result by noting that $1/(12m^2) \to 0$, and $\log(1/u) \gg \log\sqrt{\log(1/u)}$ for any $u \in (0, 1)$.
\end{proof}

\begin{lemma}
\label{lem:bound-etaCl}
Let $\eta^{\adj\ell}(u) = \frac{1}{m}\left(\frac{b}{g}\right)^{-1}\left(u \sqrt{\frac{2(1+m)^2}{\pi m}}\right)$, $u \in (0, 1)$, where $g(u) = (1+m)^{-1}$ and $b(u) = \text{Bin}(m, 1/2)$,
if $\eta^{\adj\ell}(u) - 1/2 \leq \eta_\circ$ and $\eta_\circ < 1/2$ is a fixed constant, 
then we have
\begin{align*}
\eta^{\adj\ell}(u)
& \leq \frac{1}{2} + \sqrt{\frac{1}{2m} \left(\log (1/u) -  \log (1-4 \eta_\circ^2) + \frac{1}{12m}\right)}, \\
\eta^{\adj\ell}(u) 
& \geq \frac{1}{2} + \sqrt{\frac{1}{2m} 
\left(\log (1/u) - \log (1-4 \eta_\circ^2)\right) 
\left(1 + \frac{8\eta_\circ^2}{3(1-4\eta_\circ^2)^2} \right)^{-1}}.
\end{align*}
Moreover, if $m \gg \log^2(1/u)$, then 
$$
\eta^{\adj\ell}(u) - \frac{1}{2} \sim \sqrt{\frac{1}{2m} \log \left(\frac{1}{u}\right)}.
$$
\end{lemma}

\begin{proof}
The proof is similar to that of Lemma \ref{lem:bound-etal}. 
Results follows by replacing $u$ in $(b/g)^{-1}(u)$ in Lemma \ref{lem:bound-etal} with 
$\sqrt{\frac{2}{\pi m}} (1+m) u$.
\end{proof}

\subsection{Comparing $\eta^\ell(r(w,t))$, $\eta^q(r(w,t))$, $\eta^{\adj\ell}(r(w,t))$ with $\xi(w)$ and $\zeta(w)$}
\label{sec:compare-th}

Let $\tilde \eta^\ell(u) = \eta^\ell(u) - 1/2$, $\tilde \eta^q(u) = \eta^q(u) - 1/2$, and $\tilde \eta^{\adj\ell}(u) = \eta^{\adj\ell}(u) - 1/2$, where $\eta^\ell$, $\eta^q$, and $\eta^{\adj\ell}$ are given in \eqref{eqn:etal}, \eqref{eqn:etaq}.
In this section, we obtain upper bounds for the absolute value of the differences between 
$\tilde \eta^\ell(\cdot)$, $\tilde \eta^q(\cdot)$, $\tilde \eta^{\adj\ell}(\cdot)$ and $\xi(\cdot)$ nad $\zeta(\cdot)$ respectively. 

\begin{lemma}
\label{lem:etaq-zeta}
For $\eta^{q}(u)$ and $\zeta(w)$ given in \eqref{eqn:etaq} and \eqref{eqn:zeta} respectively, 
let $u = r(w,t) = \frac{wt}{(1-w)(1-t)}$, 
then for any $w \leq w_0(t)$,
$w_0(t)$ is sufficiently small, and a fixed $t \in (0, 1)$, there exists a constant $\eta_\circ$ such that $\tilde \eta^q  \leq \eta_\circ$, $\eta_\circ \leq 1/2$ and $C = C(w_0, t, m, \eta_\circ)$ such that
$$
\big|\tilde \eta^{q}(r(w,t)) - \zeta(w)\big| \leq \frac{\log((1-t)/t) + C}{\sqrt{2m\log(1/w)}}.
$$
\end{lemma}
\begin{proof}
Using that $|\sqrt{a} - \sqrt{b}| = |a-b|/(\sqrt{a} + \sqrt{b})$ for any $a, b > 0$, we have
\begin{align}
\label{pf:lemma-S6}
|\tilde \eta^{q}(r(w,t))  - \zeta(w)|
& = \frac{2m|(\tilde \eta^{q}(r(w,t)))^2 - \zeta^2(w)| }{2m(\tilde \eta^{q}(r(w,t)) + \zeta(w) )}.
\end{align}
By Lemma \ref{lem:bound-etaq},
for the same $A_1, A_2$ given in that lemma, let $R = (1-w)(1-t)/t \leq (1-t)/t$, 
if $(\tilde \eta^{q}(r(w,t)))^2 \geq \zeta^2(w)$, then
$$
2m(\tilde \eta^{q}(r(w,t)))^2 - 2m\zeta^2(w) \leq \log ((1-t)/t) + A_1 - \log(\sqrt{\log (wt/((1-w)(1-t)) + A_2)}).
$$
If $(\tilde \eta^{q}(r(w,t)))^2 < \zeta^2(w)$, then 
$$
2m\zeta^2(w) - 2m(\tilde \eta^{q}(r(w,t)))^2 \leq - \log((1-t)/t) - \log(1-w) - A_2 + \log(\sqrt{\log (wt/((1-w)(1-t)) + A_1)}).
$$
Using the assumption $0 < w\leq w_0$ and then
let $C= C_1 \vee C_2$, where 
$$
C_1 = A_1 - \log(\sqrt{\log (t/(1-t)) - \log(1-w_0) + A_2})
$$ 
and 
$$C_2 = A_2 - \log(\sqrt{\log (t/(1-t)) - \log(1-w_0) + A_1}) + \log(1-w_0),$$
we thus obtain the upper bound for 
the numerator of \eqref{pf:lemma-S6}.
The denominator of \eqref{pf:lemma-S6} can be bounded from below by $2m \zeta(w) = \sqrt{2m\log(1/w)}$. By plugging-in the upper bound of the numerator and the lower bound of the denominator, we obtain the result. 
\end{proof}

\begin{lemma}
\label{lem:etaCl-zeta}
For $\eta^{\adj\ell}(u)$ and $\zeta(w)$ given in \eqref{eqn:etaCl} and \eqref{eqn:zeta} respectively, let $u = r(w,t)$, if $w \leq w_0(t)$ for a sufficiently small $w_0$, then for a fixed $t\in (0, 1)$ and a constant $C(w_0)$ depending on $w_0$, we have
$$
\big|\tilde \eta^{\adj\ell}(r(w,t)) - \zeta(w)\big| \leq \frac{\log(t(1-t)^{-1}) + K(\eta_\circ) \log(1/w) + C(w_0)}{\sqrt{2m\log(1/w)}},
$$
where $\displaystyle K(\eta_\circ) = \frac{8\eta_\circ^2}{3(1-4\eta_\circ^2)}$ for some $\eta_\circ < 1/2$.
\end{lemma}
\begin{proof}
Let us denote $R = (1-t)(1-w)/t$. By Lemma \ref{lem:bound-etaCl}, 
if $2m(\tilde \eta^{\adj\ell}(r(w,t)))^2 \geq 2m\zeta^2(w)$, then
\begin{align*}
& 2m(\tilde \eta^{\adj\ell}(r(w,t)))^2 - 2 m\zeta^2(w) 
\leq \log R - \log (1-4\eta_\circ^2) + 1/(12m) := U_1.
\end{align*}
If $2m(\tilde \eta^{\adj\ell}(r(w,t)))^2 < 2m\zeta^2(w)$, 
then 
\begin{align*}
& 2m\zeta^2(w) - 2m(\tilde \eta^{\adj\ell}(r(w,t)) )^2 
\geq \frac{
- \log R + \log(1-4\eta_\circ^2) - \log (w) K(\eta_\circ)
}{1+K(\eta_\circ)} := -U_2
\end{align*}
Using that $|\sqrt{a} - \sqrt{b}| = |a-b|/(\sqrt{a} + \sqrt{b})$ for any $a, b > 0$, 
if $U_1 \geq -U_2$, then 
\begin{align*}
|\tilde \eta^{\adj\ell}(r(w,t))  - \zeta(w) |
& = \frac{|\left(\tilde \eta^{\adj\ell}(r(w,t))\right)^2  - \zeta^2(w) |}{\tilde \eta^{\adj\ell}(r(w,t))  + \zeta(w))} 
\leq \frac{ \log R - \log (1-4\eta_\circ^2) + 1/(12m)}{2m(\tilde \eta^{\adj\ell}(r(w,t))  + \zeta(w))}  \\
& \leq \frac{ \log (t(1-t)^{-1}) + \log 2 + 1/(12m)}{2m(\tilde \eta^{\adj\ell}(r(w,t))  + \zeta(w))}.
\end{align*}
If $U_1 < -U_2$, then 
\begin{align*}
\big|\tilde \eta^{\adj\ell}(r(w,t))  - \zeta(w) \big|
& \leq \frac{\log R + K(\eta_\circ)\log (1/w)}
{2m(1+K(\eta_\circ))\left(\tilde \eta^{\adj\ell}(r(w,t))  + \zeta(w)\right)}  \\
& \leq \frac{ - \log (1-w_0) + \log(t(1-t)^{-1}) + K(\eta_\circ) \log(1/w)}{2m\tilde \eta^{\adj\ell}(r(w,t))  + \zeta(w))}
\end{align*}
By combining the above two cases, we obtain
\begin{align*}
\big|\tilde \eta^{\adj\ell}(r(w,t))  - \zeta(w) \big| 
& \leq \frac{\log(t(1-t)^{-1}) + \log 2 + K(\eta_\circ) \log(1/w) - \log(1-w_0)}{2m\tilde \eta^{\adj\ell}(r(w,t))  + \zeta(w))}\\
& \leq \frac{\log(t(1-t)^{-1}) + K(\eta_\circ) \log(1/w) + C(w_0)}{\sqrt{2m\log(1/w)}},
\end{align*}
where $C(w_0) = \log2 - \log(1-w_0)$.
\end{proof}

\begin{lemma}
\label{lem:etaq-etaCl-xi}
For $\eta^{q}(u)$ and $\eta^{\adj\ell}(u)$ given in \eqref{eqn:etaq} and \eqref{eqn:etaCl} respectively, let $\xi(w)$ be the solution of $\beta(u) = 1/w$ as given in Lemma \ref{lem:beta}, then for a sufficiently large $m$,
$$
\tilde \eta^{q}(r(w,t)) \leq \tilde \eta^{\adj\ell}(r(w,t)) \leq \xi(w).
$$
\end{lemma}
\begin{proof}
From Lemma \ref{lem:compare-thresholds}, we already have $\tilde \eta^q \leq \tilde \eta^{\adj\ell}$, thus it is sufficient to show $\tilde \eta^{\adj\ell}(u) \leq \xi(w)$.
Let $D = \left(1 + \frac{8\xi_\circ^2}{3(1-4\xi_\circ^2)^2}\right)^{-1}$ for some fixed $\xi_\circ < 1/2$ and $m > M_0$ for a sufficiently large $M_0$, then by Lemma \ref{lem:bound-etaCl}, we have
\begin{align*}
2m \xi^2(w) - 2m (\tilde \eta^{\adj\ell}(r(t, w)))^2 
& \geq D \left[
\log(1+1/w) + \log \left(\frac{\sqrt{2}(1+m)}{\sqrt{\pi m}}\right) 
\right]  \\
& \quad  - \log (1/r(t,w)) + \log(1-4\eta_\circ^2) -1/(12m)\\
& \geq D \log (\sqrt{2 M_0/\pi}) 
+ \log \left( 
(1 + 1/w)^D(1-4\eta_\circ^2)r(t, w) 
\right) - \frac{1}{12M_0}.
\end{align*}
Choosing $M_0$ such that 
$\log (\sqrt{2 M_0/\pi}) - (12DM_0)^{-1} > - \log \left( (1 + 1/w) [ (1-4\eta_\circ^2)r(t, w)]^{1/D}\right)$
(such $M_0$ always exists for any fixed $D$, as both $w, t \neq 0$ or $1$ and $\eta_\circ$ bounded away from $1/2$), then the last line in the last display is positive, which implies $\tilde \eta^{\adj\ell} \leq \xi(w)$ (as $\tilde \eta^{\adj\ell}$ and $\xi(w)$ are both positive.)
\end{proof}

\begin{lemma}
\label{lem:etal-xi}
Let $\xi(w)$ be the solution for $\beta(u) = 1/w$ as in Lemma \ref{lem:beta}.
For $\eta^{\ell}(u)$ defined in Lemma \ref{lem:bound-etal},
suppose $\tilde \eta^{\ell}(r(w,t)) \leq \eta_\circ$, $\eta_\circ <1/2$ is fixed, and $w \leq w_0(t)$, there exists some constant $C > 0$ depending on $t, \eta_\circ, w_0$ such that for all $t \in (0, 1)$, 
$$
|\tilde \eta^{\ell}(r(w,t)) - \xi(w)| \leq \frac{|\log(t(1-t)^{-1})| + C}{2m(\xi(w) + \tilde \eta^{\ell}(r(w,t)) )}
$$
\end{lemma}
\begin{proof}
Let us denote $\tilde \eta^{\ell}(r(w,t)) = \eta^{\ell}(r(w,t)) - 1/2$.
Using the upper bound of $\eta^{\ell}(\cdot)$ in Lemma \ref{lem:bound-etal}, we obtain
\begin{align*}
2m(\tilde \eta^{\ell}(r(w,t)))^2 - 2 m\xi^2(w) 
& \leq -\log(1-w) - \log (1-4\eta_\circ^2) + 1/(12m) + \log (t(1-t)^{-1})\\
& \leq \log (t(1-t)^{-1}) + D_1,
\end{align*}
where $D_1$ is a fixed constant.
The second inequality in the last display holds because $w \leq w_0$ and $\eta_\circ$ is smaller than 1/2.
On the other hand, using the lower bound of $\eta^{\ell}(\cdot)$ in Lemma \ref{lem:bound-etal}, 
let $D_2 = \left(1+ \frac{8\eta_\circ^2}{3(1-4\eta_\circ^2)^2}\right)^{-1}$, then 
\begin{align*}
2m \xi^2 - 2 m(\tilde \eta^\ell(r(w,t)))^2 
& \leq |\log(t(1-t)^{-1})| + (D_2 - 1) \log(wt(1-t)^{-1}) \\
& \quad + D_2 \left(\log(1-w) + \log(1-4\eta_\circ^2) 
- \log\left(\frac{\sqrt{2}(1+m)}{\sqrt{\pi m}}\right)\right)\\
& \leq |\log(t(1-t)^{-1})| - (D_2 - 1) \log(1-t).
\end{align*}
Since $t$ is a fixed constant between $0$ and $1$, let $D_3 =  - (D_2 - 1) \log(1-t) > 0$.
By combining the two upper bounds and letting $C = \max\{D_1, D_3\}$, using the fact
$|a - b| = \frac{|a^2 - b^2|}{a + b}$ for any $a, b >0$,
we obtain the result.
\end{proof}

\subsection{Proof of Lemmas \ref{lem:th-l}, \ref{lem:th-lval} and \ref{lem:th-q-cl}}
\label{sec:pf-three-lemmas}

First, we prove Lemmas \ref{lem:th-l} \& \ref{lem:th-lval}. By the definition of the $\ell$-value,
$$
\ell(x; w, g) \leq t \quad 
\Longleftrightarrow \quad
\frac{b}{g}(x) \leq r(w, t).
$$ 
Since $b/g(x)$ is symmetric at $x = m/2$ and is monotone decreasing on $x \in [m/2, m]$ and monotone increasing on $x \in [0, m/2)$ by Lemma \ref{lem:ratio-phi-g},
the last display implies $x - m/2 \geq m\eta^\ell(r(w, t)) - m/2$ if $x \geq m/2$
and $m/2 - x \geq m\eta^\ell(r(w, t)) - m/2$ if $x \leq m/2$.
By combining the two cases, by invoking Lemma \ref{lem:bound-etal} leads to the result. 

Next, we prove Lemma \ref{lem:th-q-cl}. We omit the proof of part $(a)$, as it is similar to the proof of Lemma \ref{lem:th-l}. What left is to prove part $(b)$. By the definition of the $q$-value,
$$
q(x; w, g) \leq t \quad 
\Longleftrightarrow \quad
\frac{\bar{\mathbf{B}}}{\bar\G}(m/2 + |y|) \leq r(w, t).
$$ 
By Lemma \ref{lem:ratio-barPhi-barG}, $(\bar{\mathbf{B}}/\bar\G)(\cdot)$ is symmetric at $y = 0$ and is monotone decreasing on $x \in [m/2, m]$ and monotone increasing on $x \in [0, m/2)$; therefore, 
$x - m/2 \geq m t_n^q$ when $x \in [m/2, m]$ and
$m/2 - x \geq mt_n^q$ when $x \in [0, m/2)$.
The result follows by invoking Lemma \ref{lem:bound-etaq}.

\subsection{Proof of Lemma \ref{lem:compare-thresholds}}
\label{sec:pf-compare-th}

Recall that $b(u)$ is symmetric at $m/2$ and is monotone decreasing on $[m/2, m]$, $g(u)$ is a constant, hence $(b/g)^{-1}(u)$ is also a monotone decreasing function on $0 < u<1/2$ and symmetric at $1/2$. Thus, we obtain $\eta^{\adj\ell}(u) \leq \eta^\ell(u)$ for all $u \in (0, 1)$.

Next, we show $\eta^{q}(x) < \eta^{\ell}(x)$. Consider the function 
$$
f(u) = \bar \G \left(\bar \B^{-1}(u)\right) = 1 - \frac{\bar\B^{-1}(u)}{1+m}, \ u \in (0, 1/2).
$$
By calculation, $f'(u) = (g/b)(\bar \B^{-1}(u))$, which is a decreasing function on $(0, 1/2)$. Thus, $f(u)$ is strictly concave on $(0,1/2)$. Also note that $f(0) = 0$, by the mean-value theorem,
$\bar\G(\bar \B^{-1}(u)) \geq u (g/b) (\bar \B^{-1}(u))$.
Since for any integer $x$, $m > x > m/2$, there exists one-to-one mapping to $\bar \B^{-1}(u)$ for $u \in (0,1/2)$, so for such $x$, we have $(\bar \B/\bar \G)(x) \leq (b/g)(x)$, which implies $\eta^{q}(x) \leq \eta^{\ell}(x)$.

\section{Bounding $P_{\theta_0 = 1/2}(\ell \leq t)$, $P_{\theta_0 = 1/2}(q \leq t)$, and $P_{\theta_0 = 1/2}(\adj\ell \leq t)$}
\label{sec:bound-p-mt}

\begin{lemma}
\label{lem:Plval}
For $\ell(\cdot)$ defined in \eqref{lval}, let $r(w, t) = \frac{wt}{(1-w)(1-t)}$ for any fixed $t \in (0, 1)$ and $w \leq w_0 \in (0, 1)$, suppose $\log^2 (1/r(w,t))/m \to 0$ as $m \to \infty$, let $M = m - 1$,
then as $m \to \infty$,
\begin{align}
\label{eqn:bound-lval}
 \frac{\phi(\sqrt{M} \varepsilon)}{\sqrt{M}\varepsilon}
\leq
P_{\theta_0 = 1/2}(\ell(X) \leq t) \leq  \frac{2(1+o(1)) \phi(\sqrt{M} \varepsilon)}{\sqrt{M}\varepsilon},
\end{align}
where $\displaystyle \varepsilon = \frac{2\eta^\ell(r(w,t)) -  M}{M}$.
Moreover, for any constant $C > 2\sqrt{2/\pi}$, we have
\begin{align}
\label{eqn:ub-lval}
P_{\theta_0 = 1/2}(\ell(X) \leq t) \leq \frac{C r(w,t)}{\sqrt{m}}.
\end{align}
\end{lemma}

\begin{proof} 
We use $P(\cdot)$ as the shorthand notation for $P_{\theta_0 = 1/2}(\cdot)$.
By the definition of $\ell(\cdot)$, 
\begin{align}
\label{pf:lqcl-1}
P(\ell(X) \leq t) = P((b/g)(X) \leq r(w, t)).
\end{align}
Let $\tilde U = X - m/2$, $\eta^\ell(u) = \frac{1}{m}(b/g)^{-1}(u)$, 
and $\tilde \eta^\ell(u) = \eta^\ell(u) - 1/2$, 
then \eqref{pf:lqcl-1} implies
\begin{align}
\label{pf:lqcl-2}
P((b/g)(\tilde U + m/2) \leq r(w, t)) = 
P\left(|\tilde U| \geq m\tilde \eta^\ell (r(w,t))\right)
= 2\bar \B(m \eta^{\ell} (r(w,t))),
\end{align}
as $(b/g)(\cdot)$ is symmetric at $\tilde U = 0$.
Let $K = m\eta^\ell(r(w,t)) - 1$, $M = m - 1$,
and $\varepsilon = (2K - M)/M$, from Lemma \ref{lem:binom-tail-carter}, we have
$$
\bar \B(m \eta^{\ell} (r(w,t))) = 
\bar \Phi \left(\varepsilon \sqrt{M} \right)\exp(A_m(\varepsilon)),
$$
where $A_m(\varepsilon) = - M\varepsilon^4 \gamma(\varepsilon) - \log (1-\varepsilon^2)/2 - \lambda_{m-K+1} + r_{K+1}$ with $\gamma(\varepsilon) \sim 1/12$, $\lambda_{m-K+1} =O(1/m)$ and $r_{K+1} = O(1/m)$.
By invoking Lemma \ref{lem:gaussian-bound}, the last display can be bounded by 
$$
\bar \B(m \eta^{\ell} (r(w,t))) 
\leq \frac{\phi \left(\varepsilon \sqrt{M} \right)\exp(A_m(\varepsilon))}{\varepsilon \sqrt{M}}.
$$
To obtain the upper bound in \eqref{eqn:bound-lval}, what left is show $A_m(\varepsilon) = o(1)$. 
By Lemma \ref{lem:bound-etal}, 
$K \sim m/2 + \sqrt{m(\log (1/r(w,t)) + \log \sqrt{2m/\pi})/2}$ for a sufficiently large $m$.
Since $\varepsilon = (2K- M)/M$,
$\varepsilon \sim \sqrt{2(\log (1/r(w,t)) + \log \sqrt{2m/\pi})/m}$.
Using the assumption $m \gg \log^2 (1/r(w, t)))$, we thus have $A_m(\varepsilon) \to 0$.
Therefore, the last display implies 
$$
P((b/g)(\tilde U + m/2) \leq r(w, t)) 
\leq \frac{2 (1+o(1))\phi( \varepsilon \sqrt{M})}{\varepsilon \sqrt{M}}.
$$
To prove the lower bound in \eqref{eqn:bound-lval}, we use the lower bound of the Gaussian tail in Lemma \ref{lem:gaussian-bound} to obtain
$$
\bar \Phi(\varepsilon\sqrt{M}) 
\geq \frac{M\varepsilon^2}{1+M\varepsilon^2} \frac{\phi(\varepsilon \sqrt{M})}{\varepsilon \sqrt{M}}
\geq \frac{\phi(\varepsilon \sqrt{M})}{2\varepsilon \sqrt{M}}, 
$$
as long as $\sqrt{M} \varepsilon > 1$.
Since $A_m(\varepsilon)) = o(1)$, we have
$$
P((b/g)(\tilde U + m/2) \leq r(w, t)) 
= 2 \bar\B(m \eta^\ell(r(w,t)))
\geq \frac{\phi(\varepsilon \sqrt{M})}{\varepsilon \sqrt{M}}.
$$
We thus complete the proof of \eqref{eqn:bound-lval}.

To prove \eqref{eqn:ub-lval},
by invoking the Bernstein inequality as given in Lemma \ref{bernstein-inq} (choosing $D = 1$, $V = m/4$, and $A = m\tilde \eta^\ell(r(w,t))$), we have
\begin{align*}
P( |\tilde U| \geq m\tilde \eta^\ell(r(w,t))) \leq 2\exp\left(
- \frac{m^2(\tilde \eta^\ell(r(w,t)))^2}{m/2 + m\tilde \eta^\ell(r(w,t)) /3}
\right).
\end{align*}
Since $-\log (r(w,t))/m \to 0$ and $\log (\sqrt{m})/m \to 0$, we have $m/2 \gg m\tilde \eta^\ell(r(w,t)) /3$.
Using that $2m \left(\tilde \eta^\ell(r(w,t))\right)^2 = - \log (r(w,t)) + \log (C_m\sqrt{m})+ o(1)$ for some positive constant $C_m < \sqrt{2/\pi}(1+1/m) < 2\sqrt{2/\pi}$, the last display is bounded by 
$2\exp(- 2m (1- o(1)) (\tilde \eta^\ell(r(w,t)))^2) \leq C {r(w,t)}/{\sqrt{m}}$ for a sufficiently large $m$ for any $C > C_m$. 
\end{proof}

\begin{lemma}
\label{lem:PClval}
For $\adj\ell(\cdot)$ given in \eqref{Clval}, let $r(w, t) = \frac{wt}{(1-w)(1-t)}$
for any fixed $t \in (0, 1)$ and $w \leq w_0 \in (0, 1)$, 
suppose $(\log (1/r(w,t)))^2/m \to 0$ as $m \to \infty$,
then
\begin{align}
\label{eqn:bound-Clval}
 \frac{\phi(\sqrt{M}\tilde \varepsilon)}{\sqrt{M}\tilde\varepsilon}
\leq
P_{\theta_0 = 1/2}(\adj\ell(X) \leq t) \leq  \frac{2(1+o(1)) \phi(\sqrt{M} \tilde\varepsilon)}{\sqrt{M}
\tilde\varepsilon},
\end{align}
where $\displaystyle \tilde \varepsilon = \frac{2\eta^{\adj\ell}(r(w,t)) -  m + 1}{m-1}$.
Moreover, we have
\begin{align}
\label{eqn:ub-Clval}
P_{\theta_0 = 1/2}(\adj\ell(X) \leq t) \leq C' r(w,t),
\end{align}
for some constant $C' > 2$.
\end{lemma}

\begin{proof}
The proof is similar to that of Lemma \ref{lem:Plval}.
Again, we use $P(\cdot)$ as the shorthand notation for $P_{\theta_0 = 1/2}(\cdot)$.
By the definition of $\adj\ell(\cdot)$ in \eqref{Clval}, we have
$$
P(\adj\ell(X) \leq t) = P\left((b/g)(X) \leq  \sqrt{2/(\pi m)} r(w,t) (1+m)\right).
$$
By replacing the upper bound in \eqref{pf:lqcl-1} with $\sqrt{2/(\pi m)} r(w,t) (1+m)$
and let $\tilde U = X - m/2$, we obtain
$$
P((b/g)(\tilde U + m/2) \leq r(w, t)) = 
P\left(|\tilde U| \geq \tilde \eta^{\adj\ell}(r(w,t))\right), 
$$
where $\tilde \eta^{\adj\ell}(r(w,t)) = \eta^{\adj\ell}(r(w,t))- 1/2$.
Using Lemma \ref{lem:bound-etaCl} and the assumption $\log^2 (1/r(w,t))/m \to 0$,
then
$2m \left(\eta^{\adj\ell}(r(w,t))\right)^2 = - \log (r(w,t)) + o(1)$.
The remaining proof is exact the same as that for Lemma \ref{lem:Plval} (with $\eta^{\adj\ell}(\cdot)$ is replaced by $\eta^{\ell}(\cdot)$).
\end{proof}

\begin{lemma}
\label{lem:Pqval}
For $q(\cdot)$ defined in \eqref{qval}, define $r(w, t) = \frac{wt}{(1-w)(1-t)}$, 
for any fixed $t \in (0, 1)$ and $w \in (0, 1)$,  then
$P_{\theta_0 = 1/2}(q(X) \leq t) = 2 r(w, t) \bar \G(\eta^q(r(w,t))) \leq 2 r(w, t)$.
\end{lemma}

\begin{proof}
Since $\eta^q(\cdot) = \frac{1}{m}(\bar{\mathbf{B}}/\bar\G)^{-1}(\cdot)$,
let $\tilde U = X - m/2$, then 
\begin{align*}
P_{\theta_0 = 1/2}(q(X) \leq t) 
& = P_{\theta_0 = 1/2}((\bar{\mathbf{B}}/\bar\G)(X) \leq r(w,t))
= P_{\theta_0 = 1/2}\left(|\tilde U| \geq m(\eta^q(r(w,t)) - 1/2) \right)\\
& = 2 \bar \B(m \eta^q(r(w,t))) 
= 2 r(w,t)\bar \G(m\eta^q(r(w,t)))
\leq 2r(w, t).
\end{align*}
\end{proof}


\section{Proof of results in Section \ref{sec:unif-fdr}}
\label{sec:pf-unif-fdr}

In this section, we prove Theorem \ref{thm:unif-fdr-q-Cl} and Lemma \ref{lem:unif-fdr-l}. 
In our proof, we will use results from Lemma \ref{lem:bound-what}, which gives a concentration bound for the MMLE $\hat w$, and Lemma \ref{lem:m1-cntl-small-signals}, which controls those `small' signals (i.e., between $1/2 \pm \xi(w)$ for $\xi(w)$ given in \eqref{eqn:xi}.) Proofs of the two lemmas are given in Sections \ref{sec:bound-w-hat} and \ref{sec:control-small-signals} respectively.

\subsection{Proof of Theorem \ref{thm:unif-fdr-q-Cl}}
\label{sec:pf-thm1}
We proof the results for the $\adj\ell$-value and the $q$-value procedures together. 
We will refer $\T$ as either $\T^{\adj\ell}$ or $\T^q$ when their proof are the same. 
The proof is divided into two parts based on whether a solution for the equation \eqref{eqn:w1} exists or not. 

{\bf Case 1. \eqref{eqn:w1} has a solution.} Using $(i)$ of Lemma \ref{lem:bound-what}, 
there exists a constant $C > 0$ and a fixed $\kappa \in (0, 1)$ such that
$$
P_{\theta_0}(\hat w \not\in [w_2, w_1]) 
\leq e^{-C\kappa^2 nw_1\tilde m(w_1)} + e^{-C\kappa^2 nw_2\tilde m(w_2)}
\leq 2e^{- 0.4 C\kappa^2 s_n}.
$$
The second upper bound in the last display is obtained since $\tilde m(w_1) \geq 0.4$ and $\tilde m(w_2) \geq 0.4$ for a sufficiently large $m$ by Lemma \ref{lem:bound-mtilde}
and $w_2 \leq w_1 \lesssim s_n/n$ by Lemma \ref{lem:bound-w1-w2}.
Therefore,
\begin{align}
\FDR(& \theta_0, \T(t, \hat w))
= 
\mathbb{E}_{\theta_0} 
\left(
\frac{\fF_\T(t, \hat w)}{\max\{1, \fF_\T(t, \hat w) + \fT_\T (t, \hat w)\}}
\right)
\nonumber \\
& \leq 
\mathbb{E}_{\theta_0} 
\left(
\frac{\fF_\T(t, \hat w)}{\max\{1, \fF_\T(t, \hat w) + \fT_\T (t, \hat w)\}} \mathbbm{1}\{w_2 \leq \hat w \leq w_1\}
\right) + P_{\theta_0}(\hat w \not\in [w_2, w_1]) 
\nonumber \\
& \leq 
\mathbb{E}_{\theta_0} 
\left(
\frac{\fF_\T(t, \hat w)}{\max\{1, \fF_\T(t, \hat w) + \fT_\T (t, \hat w)\}} \mathbbm{1}\{w_2 \leq \hat w \leq w_1\}
\right) 
+ 2 e^{- 0.4 C\kappa^2 s_n}.
\label{pf:unif-fdr-0}
\end{align}
Since both $\fF_\T(t, w)$ and $\fT_\T(t, w)$ are monotone functions of $w$ (one can check this from their definitions in \eqref{eqn:FD-TD}), by the monotonicity of the function $x \to x/(1+x)$ and Lemma \ref{lem:boundFDR}, the first term in the summation of \eqref{pf:unif-fdr-0} can be bounded by
\begin{align}
\label{pf:unif-fdr-1}
\mathbb{E}_{\theta_0} 
\left(
\frac{\fF_\T(t, w_1)}{\max\{1, \fF_\T(t, w_1) + \fT_\T (t, w_2)\}} 
\right)
\leq 
\exp\left(-\mathbb{E}_{\theta_0} \fT_\T(t, w_2)\right)
+ \frac{12\mathbb{E}_{\theta_0} \fF_\T(t, w_1)}{\mathbb{E}_{\theta_0} \fT_\T(t, w_2)}.
\end{align} 
What remains is to obtain a lower bound for $\mathbb{E}_{\theta_0} \fT_\T(t, w_2)$ and an upper bound for $\mathbb{E}_{\theta_0} \fF_\T(t, w_1)$.
 
{\it Lower bound for} $\mathbb{E}_{\theta_0} \fT_{\T}(w_2)$.
From Lemma \ref{lem:etaq-etaCl-xi}, we have $\tilde \eta^{\adj\ell}(w_2) \leq \xi(w_2)$ and $\tilde \eta^{q}(w_2) \leq \xi(w_2)$ for a sufficiently large $m$, thus for either $\T = \T^{\adj\ell}$ or $\T = \T^q$, 
\begin{align}
\mathbb{E}_{\theta_0} \fT_{\T}(w_2)
&\geq 
\sum_{j:\theta_{0,j} \neq 1/2} P_{\theta_{0,j}}\left(|X_j - m/2| \geq m\xi(w_2)\right)
\nonumber \\
& = \sum_{j:\theta_{0,j} \neq 1/2} 
\left(\bar\B_{\theta_{0,j}}(m/2 + m\xi(w_2)) 
+ {\B}_{\theta_{0,j}}(m/2 - m\xi(w_2) + 1)\right). 
\label{pf:unif-fdr-2}
\end{align}

Recall the set $\mathcal{J}_0$ defined in \eqref{J0}, let's further define 
$$
\mathcal{J}_0^+ = \{1 \leq j \leq n: \theta_{0,j} - 1/2 \geq \xi(w)/K\}
\ \text{and} \
\mathcal{J}_0^- = \{1 \leq j \leq n: 1/2 - \theta_{0,j} > \xi(w)/K\},
$$
then \eqref{pf:unif-fdr-2} leads to
\begin{align}
\mathbb{E}_{\theta_0} \fT_{\T}(w_2) \geq 
& \sum_{j \in \mathcal{J}_0^+} 
\left(\bar\B_{\theta_{0,j}}(m/2 + m\xi(w_2)) 
+ {\B}_{\theta_{0,j}}(m/2 - m\xi(w_2) + 1)\right)
\label{pf:unif-fdr-2-1}
 \\
&  + 
\sum_{j \in \mathcal{J}_0^-} 
\left(\bar\B_{\theta_{0,j}}(m/2 + m\xi(w_2)) 
+ {\B}_{\theta_{0,j}}(m/2 - m\xi(w_2) + 1)\right).
\label{pf:unif-fdr-2-2}
\end{align}

We first obtain a lower bound for \eqref{pf:unif-fdr-2-1}. 
If $\theta_{0,j} - 1/2 > \xi(w_2)$, then by Lemma \ref{lem:binom-tail-slud},
\begin{align*}
\eqref{pf:unif-fdr-2-1}
& \geq 
\frac{1}{2} \sum_{j:\theta_{0,j} \geq 1/2} 
\bar\B_{\theta_{0,j}}\left(m/2 + m\xi(w_2)\right)
+\frac{1}{2} \sum_{j: \theta_{0,j} \geq 1/2} \bar \Phi\left(
\frac{\sqrt{m}(\xi(w_2) - (\theta_{0,j} - 1/2))}{\sqrt{\theta_{0,j}}}
\right) \\
& \geq 
\frac{1}{2} \sum_{j:\theta_{0,j} \geq 1/2} 
\bar\B_{\theta_{0,j}}\left(m/2 + m\xi(w_2)\right)
+ \frac{1}{2} 
\sum_{j: \theta_{0,j} \geq 1/2} 
\bar \Phi\left(2\sqrt{m}|\xi(w_2) - (\theta_{0,j} - 1/2)|\right)\\
& \geq \frac{w_2}{2} \sum_{j \in \mathcal{J}_0^+} m_1(\theta_{0,j}, w_2)
\end{align*}
by Corollary \ref{cor:m1} and using that 
$T_m(\theta - 1/2, \xi(w_2)) \leq \frac{1}{\sqrt{1-4\xi^2(w_2)}} = 1 + o(1)$.

If $\xi(w_2) \geq \theta_{0,j}  - 1/2$, by the assumption $(\log n)^2 \ll m$ (and thus $m \xi^4 \to 0$) and Lemma \ref{lem:lb-bin-cdf}, we have 
\begin{align*}
\bar \B_{\theta_{0,j}}(m/2 + m \xi(w_2)) 
\geq 
\frac{1}{2}\sqrt{\frac{1-2(\theta_{0,j} - 1/2)}{1-2\xi(w_2)}}
\bar\Phi\left(
\frac{2\sqrt{m}(\xi(w_2) - (\theta_{0,j} - 1/2))}{\sqrt{1-4(\theta_{0,j} - 1/2)^2}}
\right).
\end{align*}
Using the lower bound in Lemma \ref{lem:gaussian-bound}, we have
\begin{align*}
& \frac{\bar\Phi\left(
\frac{2\sqrt{m}(\xi(w_2) - (\theta_{0,j} - 1/2))}{\sqrt{1-4(\theta_{0,j} - 1/2)^2}}
\right)}{\bar\Phi \left( 2\sqrt{m}(\xi(w_2) - (\theta_{0,j} - 1/2)) \right)} \\
& \geq 
\frac{\sqrt{1-4\xi^2(w_2)}}{2}
\exp\left(- \frac{8m(\theta_{0,j}-1/2)^2 (\xi(w_2)-(\theta_{0,j}-1/2))^2}{1-4(\theta_{0,j}-1/2)^2}\right).
\end{align*}
Since $\theta_{0,j} - 1/2$ is at most $\xi(w_2)$, by the assumption $m \xi^4(w_2) \to 0$, the exponential term in the last display is $e^{-o(1)}$ and thus is bounded from below by $1/2$ for a sufficiently large $m$.
Therefore, by plugging-in the lower bound in the last display, we obtain
\begin{align*}
\bar \B_\theta(m/2 + m \xi(w_2)) 
& \geq \frac{1}{4} \sqrt{\frac{1-2(\theta_{0,j}-1/2)}{1-2\xi(w_2)}} \bar\Phi\left(2\sqrt{m}(\xi(w_2) - (\theta_{0,j}-1/2))\right) \\
& \geq \frac{1}{4}\bar\Phi\left(2\sqrt{m}(\xi(w_2) - (\theta_{0,j}-1/2))\right).
\end{align*}
We thus obtain
$$
\eqref{pf:unif-fdr-2-1} \geq \frac{1}{8} \sum_{j \in \mathcal{J}_0^+} 
\left(\bar \B_{\theta_0,j}(m/2 + m\xi(w_2)) 
+ \bar \Phi\left(2\sqrt{m}(\xi(w_2) - (\theta_{0,j} - 1/2))\right)\right) 
\geq \frac{w_2}{8} \sum_{j \in \mathcal{J}_0^+} m_1(\theta_{0,j}, w_2).
$$
The lower bound for \eqref{pf:unif-fdr-2-2} can be obtained using a similar argument using the fact that $\B_\theta(m/2 - m\xi_2 + 1) = \bar \B_{1/2 + |\theta-1/2|}(m/2+m\xi_2)$.

By combining both cases, we have
$\eqref{pf:unif-fdr-2} \geq C' {w_2} \sum_{j \in \mathcal{J}_0} m_1(\theta_{0,j}, \mu)$ for some constant $C' \leq 1/8$.
By \eqref{eqn:w2} and Lemma \ref{lem:m1-cntl-small-signals}, 
there exists constants $C, D > 0$ such that 
$$
 \sum_{j \in \mathcal{J}_0} m_1(\theta_{0,j}, \mu) = 
(1+\kappa)(n - s_0) \tilde m(w_2) - Cn^{1-D} \tilde m(w_2).
$$
Note that the second term in the last display is of a smaller order of the first one, as $0.4 < \tilde m(w_2) \leq 1.1$ for a sufficiently large $m$ from Lemma \ref{lem:bound-mtilde}; therefore, 
\begin{align}
\label{pf:lb-TD}
\mathbb{E}_{\theta_0} \text{TD}_{\T(w_2)} \geq \eqref{pf:unif-fdr-2}
\geq C' (n-s_0) w_2 \tilde m(w_2).
\end{align}

{\it Upper bound for} $\mathbb{E}_{\theta_0} \fF_\T(t, w_1)$.
We will consider the $q$-value procedure and the $\adj\ell$-value procedure separately.
For the $q$-value procedure, by Lemma \ref{lem:Pqval} and Lemma \ref{lem:bound-w1-w2}, we have
\begin{align}
\mathbb{E}_{\theta_0} \fF_{\T^{q}}(w_1)
& \leq \sum_{j: \theta_{0,j} = 1/2} \mathbb{E}_{\theta_0} \T_j^{q}(w_1)
= \sum_{j: \theta_{0,j} = 1/2} P_{\theta_0} (q(X_j; w_1) \leq t)
\nonumber \\
& \lesssim \frac{w_1 t (n - s_0)}{(1-w_1)(1-t)}
= C_1 (n-s_0) w_1 t,
\label{pf:ub-FD-q}
\end{align}
where $C_1 = 1/((1-w_1)(1-t))$.

For the $\adj\ell$-value procedure, by Lemma \ref{lem:PClval} and Lemma \ref{lem:bound-w1-w2}, we obtain
\begin{align}
\mathbb{E}_{\theta_0} \fF_{\T^{\adj\ell}}(w_1)
& \leq \sum_{j: \theta_{0,j} = 1/2} \mathbb{E}_{\theta_0} \T_j^{\adj\ell}(w_1)
= \sum_{j: \theta_{0,j} = 1/2} P_{\theta_0}(\adj\ell(X_j; w_1) \leq t) 
\nonumber \\
& \leq (n - s_0) \frac{ \phi( \varepsilon \sqrt{M}) e^{A_m}}{\varepsilon \sqrt{M}},
\label{pf:unif-fdr-3}
\end{align}
where $\varepsilon = \frac{2K-M}{M}$, $K = m \eta^{\adj\ell}(r(w_1, t)) - 1$, $M = m-1$, and $A_m \lesssim e^{-C M \varepsilon^4}$. 
By Lemma \ref{lem:etaCl-zeta}, 
\begin{align*}
\eta^{\adj\ell}(r(w_1, t)) - 1/2 
& \leq \zeta(w_1) + \frac{\log (t(1-t)^{-1})}{\sqrt{2m\log(1/w_1)}} - K(\eta_\circ) \frac{\log(1/w_1)}{\sqrt{2m}} + \frac{C}{\sqrt{2m \log(1/w_1)}} \\
& \leq \sqrt{\frac{\log(1/w_1)}{2m}} + o(1).
\end{align*}
Since $w_1 \lesssim s_n/n = n^{v_1 - 1}$ and $\log (1/w_1) \asymp \log n \ll m$,
Therefore,
\begin{align}
\eqref{pf:unif-fdr-3} 
\leq 
\frac{(n - s_0) \phi( 2\sqrt{M} (\eta^{\adj\ell}(r(w_1, t)) - 1/2)) e^{A_m}}
{2\sqrt{M}(\eta^{\adj\ell}(r(w_1, t)) - 1/2)}
\lesssim 
\frac{(n - s_0) w_1}{\sqrt{2\log (1/w_1)}}
\asymp \frac{(n - s_0) w_1}{\sqrt{2 \log n}}.
\label{pf:ub-FD-Cl}
\end{align}

{\it Upper bound for the FDR.}
We now combine the previous results to obtain the upper bound for \eqref{pf:unif-fdr-1} and \eqref{pf:unif-fdr-0}.
For the $q$-value procedure, by combining the lower bound in \eqref{pf:lb-TD} and the upper bound in \eqref{pf:ub-FD-q}, we have
$$
\eqref{pf:unif-fdr-1}
\leq e^{-C'(n-s_0) w_2\tilde m(w_2)} + C_2 t.
$$
Using this inequality to bound \eqref{pf:unif-fdr-0}, then for a sufficiently large $n$, 
\begin{align}
\label{pf:fdr-q-case1}
\sup_{\theta_0 \in l_0[s_n]} 
\FDR(\theta_0, & \T^q(t, \hat w)) 
\lesssim e^{-C'(n-s_0) w_2\tilde m(w_2)} + C_2 t.
\end{align}
For the FDR of the $\adj\ell$-value procedure, using \eqref{pf:lb-TD} and \eqref{pf:ub-FD-Cl} lead to
\begin{align}
\label{pf:fdr-Cl-case1}
\sup_{\theta_0 \in l_0[s_n]} 
\FDR(\theta_0, \T^{\adj\ell}(t, \hat w)) \lesssim e^{- C'(n - s_0)w_2\tilde m(w_2)} + \frac{C_3}{\sqrt{\log n} }.
\end{align}

{\bf Case 2: \eqref{eqn:w1} does not have a solution}. 
If \eqref{eqn:w1} does not have a solution, then we must have 
$$
\sum_{j \in S_0} m_1(\theta_{0,j}, w) < (1-\kappa) (n-s_0) \tilde m(w)
$$
due to $\tilde m(w)$ is continuous and monotone increasing and $m_1(u, w)$ is continuous and monotone decreasing (see Lemma \ref{lem:m1-tildem-mnt}). By $(ii)$ in Lemma \ref{lem:bound-what}, there exists a constant $C$ and $\Delta_n = nw_0 \tilde m(w_0)$ such that 
$P_{\theta_0}(\hat w\geq w_0) \leq e^{-C\kappa^2 \Delta_n}$.
Then, for either $\T = \T^{\adj\ell}$ or $\T = \T^q$, we have 
\begin{align}
\FDR(\theta_0, \T(t, \hat w))
& \leq P_{\theta_0}(j: \theta_{0,j} = 1/2, \T_j(t, \hat w) = 1) \nonumber\\
& \leq (n-s_0) P_{\theta_0 = 1/2}(\T_j(t, w_0) = 1)  + P_{\theta_0}(\hat w\geq w_0) \nonumber\\
& \leq (n-s_0) P_{\theta_0 = 1/2}(\T_j(t, w_0) = 1)  + e^{-C\kappa^2 \Delta_n}.
\end{align}
For the $q$-value procedure, by invoking Lemma \ref{lem:Pqval}, we have
$P_{\theta_0 = 1/2}(\T_j^q(t, w_0) = 1) = P_{\theta_0 = 1/2}(q(w_0) \leq t) \leq 2r(w_0, t) = \frac{2 w_0 t}{(1-w_0)(1-t)}$, then
\begin{align}
\label{pf:fdr-q-case2}
\FDR(\theta_0, \T^q(t, \hat w))
\leq \frac{2 (n-s_0) w_0 t}{(1-w_0) (1-t)} + e^{-C\kappa^2 \Delta_n}
\leq \frac{2 (n-s_0) w_0 \tilde m(w_0) t}{0.4 (1-w_0) (1-t)} + e^{-C\kappa^2 \Delta_n},
\end{align}
which we used $0.4 \leq \tilde m(w) \leq 1.1$ from Lemma \ref{lem:bound-mtilde}.
Note that $w_0$ is the solution of \eqref{eqn:w0}, thus for a sufficiently large $n$ and $t \leq 4/5$, the first term in the upper bound in \eqref{pf:fdr-q-case2} is bounded by 
$$
\frac{2 t \Delta_n}{0.4(1- n^{-1})(1-t)} \leq 25 (1+o(1)) t \Delta_n.
$$
Next, consider the $\adj\ell$-value procedure, 
by Lemma \ref{lem:PClval} and the upper bound in \eqref{eqn:bound-Clval},
let $\tilde \varepsilon = \sqrt{{2\log(1/r(w_0,t))}/{m}}$, then
$$
P_{\theta_0 = 1/2}(\adj\ell(x) \leq t) \lesssim 
\frac{2\phi(\sqrt{M} \tilde\varepsilon)}{\sqrt{M} \tilde\varepsilon}
\lesssim \frac{2e^{- M^2  \tilde\varepsilon^2/2}}{\sqrt{2M \log (1/r(w_0,t))} }
\lesssim \frac{\sqrt{2}r(w_0, t)}{\sqrt{M \log (1/r(w_0,t))}}.
$$
Therefore, for any $t \leq 4/5$,
\begin{align}
\label{pf:fdr-Cl-case2}
\FDR(\theta_0, \T^{\adj\ell}(t, \hat w))
& \lesssim \frac{2(n-s_0) w_0 t}{(1-w_0)(1-t) \sqrt{\log n}} + e^{-C\kappa^2 \Delta_n} 
 \nonumber \\
& \leq \frac{25(1+o(1))\Delta_n t}{ \sqrt{\log n}} + e^{-C\kappa^2 \Delta_n}. 
\end{align}

{\bf Combining Case 1 and Case 2.}
Combining \eqref{pf:fdr-q-case1} and \eqref{pf:fdr-q-case2} leads to 
$$
\sup_{\theta_0 \in l_0[s_n]} \FDR(\theta_0, \T^{q}(t, \hat w)) 
\leq \max\{e^{-C'(n-s_0)w_2\tilde m(w_2)} + C_2 t, \ 25(1+o(1)) t\Delta_n + e^{-C\kappa^2 \Delta_n} \}.
$$
Using that $\tilde m(\cdot) \in [0.4, 1.1]$ for a sufficiently large $m$ from Lemma \ref{lem:bound-mtilde},
we obtain  
\begin{align}
\label{pf:fdr-q}
(n-s_0) w_2 \tilde m(w_2) \gtrsim (n-s_0) w_0 \tilde m(w_0) \geq C'' \Delta_n,
\end{align}
for some $C'' > 0$.
We then have
$$
\sup_{\theta_0 \in l_0[s_n]} \FDR(\theta_0, \T^{q}(t, \hat w)) 
\lesssim \max\{e^{-C'' \Delta_n}+ C_2 t, 23 t \Delta_n + e^{-C\kappa^2 \Delta_n}\}.
$$
Choosing $\Delta_n = \max\{\frac{1}{C''}, \frac{\log (1/t)}{C\kappa^2}\}$, 
then
$\sup_{\theta_0 \in l_0[s_n]} \FDR(\theta_0, \T^{q}(t, \hat w))  \lesssim t\log(1/t)$.

The result for the $\adj\ell$-value procedure can be obtained similarly. By combining \eqref{pf:fdr-Cl-case1} and \eqref{pf:fdr-Cl-case2}, we have
$$
\sup_{\theta_0 \in l_0[s_n]} \FDR(\theta_0, \T^{\adj\ell}(t, \hat w)) 
\leq \max\left\{ e^{-C'' \Delta_n} + \frac{C_3}{\sqrt{\log n}}, \ \frac{K t \Delta_n}{\sqrt{\log n}} + e^{-C\kappa^2 \Delta_n} \right\}.
$$
By choosing $\Delta_n = \max\{\frac{\log \log n}{2C\kappa^2}, \frac{\log \log n}{2C''}\}$, we obtain the upper bound
$$
\sup_{\theta_0 \in l_0[s_n]} \FDR(\theta_0, \T^{\adj\ell}(t, \hat w)) 
\lesssim \frac{t\log \log n}{\sqrt{\log n}}.
$$
We thus complete the proof.

\subsection{Proof of Lemma \ref{lem:unif-fdr-l}}
\label{sec:pf-unif-fdr-l}

Similar to the proof in Section \ref{sec:pf-thm1}, we also consider two cases depending on whether \eqref{eqn:w1} has a solution or not. 
If a solution exists, then by following the same proof of the $\adj\ell$-value procedures in Section \ref{sec:pf-thm1}, we have
\begin{align}
\label{pf:unif-fdr-l-1}
\FDR(\theta_0, \T^\ell(t, \hat w)) \leq 
\exp\left(-\mathbb{E}_{\theta_0} \fT_{\T^\ell} (t, w_2)\right)
+ \frac{12\mathbb{E}_{\theta_0} \fF_{\T^\ell}(t, w_1)}{\mathbb{E}_{\theta_0} \fT_{\T^\ell}(t, w_2)}
+ P_{\theta_0} (\hat w \not\in [w_2, w_1]).
\end{align}

{\it Lower bound for $\mathbb{E}_{\theta_0} \fT_{\T^\ell} (w_1)$.} By the definition of $\fT_{\T^\ell} (t, w_1)$,
\begin{align}
\mathbb{E}_{\theta_0} & \fT_{\T^\ell} (t, w_1)
= \sum_{j: \theta_{0,j} \neq 1/2} \mathbb{E}_{\theta_0} \T_j^\ell(t, w_1)
= \sum_{j: \theta_{0,j} \neq 1/2} P_{\theta_0}(\ell(X_j; w_1) \leq t) 
\nonumber \\
& \geq \sum_{j: \theta_{0,j} > 1/2} \bar\B_{\theta_{0,j}} (m/2 + m \tilde \eta^\ell(t, w_2))
+ \sum_{j: \theta_{0,j} < 1/2} \bar\B_{\theta_{0,j}} (m/2 - m \tilde \eta^\ell(t, w_2)).
\label{pf:unif-fdr-l-2}
\end{align}
From Lemma \ref{lem:etal-xi}, we have
$
|2m(\tilde \eta^{\ell}(r(w,t)))^2 -2m\xi^2 (w)| \leq |\log(t(1-t)^{-1}))^2| + C
$ for some constant $C$.
Therefore, by Lemma \ref{lem:ratio-bin-cdf}, we obtain
\begin{align*}
\eqref{pf:unif-fdr-l-2} \geq C_t
\left(
\sum_{j: \theta_{0,j} > 1/2} \bar\B_{\theta_{0,j}} (m/2 + m \xi(w_2))
+ \sum_{j: \theta_{0,j} < 1/2} \bar\B_{\theta_{0,j}} (m/2 - m \xi(w_2))
\right),
\end{align*}
for some positive fixed $C_t$ depending on $t$ and $\theta_{0}$.
Then, for a large enough $m$ and a constant $C'$ depends on $t$ and $\theta_0$, we have
\begin{align}
\label{pf:unif-fdr-l-3}
\mathbb{E}_{\theta_0} \fT_{\T^\ell} (w_1)
\geq C' (n-s_0) w_2 \tilde m(w_2).
\end{align}

{\it Upper bound for $\mathbb{E}_{\theta_0} \fF_{\T^\ell} (w_2)$.}
We have
\begin{align}
\mathbb{E}_{\theta_0} \fF_{\T^\ell} (w_1)
& = \sum_{j: \theta_{0,j} = 1/2} \mathbb{E}_{\theta_0} \T_j^\ell(w_1)
= \sum_{j: \theta_{0,j} = 1/2} P_{\theta_0}(\ell(X_j; w_1) \leq t) 
\nonumber \\
& = (n - s_0) \bar \B(m \eta^\ell(r(w_1, t))).
\label{pf:unif-fdr-l-2}
\end{align}
By Lemma \ref{lem:binom-tail-carter}, 
let $M = m-1$ and $K = m \eta^{\ell}(r(w_1, t)) - 1$, for $\varepsilon = \frac{2K-M}{M}$, 
and then by Lemma \ref{lem:gaussian-bound},
$$
\bar \B(m \eta^{\ell}(r(w_1, t)))
= \bar \Phi( \varepsilon \sqrt{M}) e^{A_m}
\leq \frac{ \phi( \varepsilon \sqrt{M}) e^{A_m}}{\varepsilon \sqrt{M}}.
$$
Therefore, \eqref{pf:unif-fdr-l-2} can be bounded by 
\begin{align*}
\frac{(n - s_0) \phi( 2\sqrt{M} (\eta^{\ell}(r(w_1, t)) - 1/2)) e^{A_m}}
{2\sqrt{M}(\eta^{\ell}(r(w_1, t)) - 1/2)}
\lesssim 
\frac{(n - s_0) w_1}{\sqrt{2m (\log (1/w_1) + \log \sqrt{m})}}.
\end{align*}
Since $w_1 \lesssim s_n/n = n^{v_1-1}$ from Lemma \ref{lem:bound-w1-w2}, by Lemma \ref{lem:bound-mtilde}, we obtain 
\begin{align}
\label{pf:unif-fdr-l-4}
\mathbb{E}_{\theta_0} \fF_{\T^\ell} (w_1)
\lesssim
\frac{(n-s_0) w_1}{\sqrt{2(1-v_1)m \log n}}
\lesssim 
\frac{(n-s_0) w_2 \tilde m(w_2)}{\sqrt{m \log n}}.
\end{align}
Now we combine \eqref{pf:unif-fdr-l-2} and \eqref{pf:unif-fdr-l-4} and then use \eqref{pf:fdr-q}, 
\begin{align*}
\FDR(\theta_0, \T^\ell(t, \hat w))
& \lesssim e^{-C'(n-s_0) w_2 \tilde m(w_2) + \frac{2Ks_n}{\sqrt{m} \xi(w_2)}}
+ \frac{12 (n-s_0) w_2 \tilde m(w_2)/\sqrt{m\log n}}{(n-s_0) w_2\tilde m(w_2) - 2Ks_n/\sqrt{m\xi^2(w_2)}} \\
& \leq 
e^{-C'' \Delta_n + \frac{2Ks_n}{\sqrt{m} \xi(w_2)}} 
+ \frac{12}{\sqrt{m \log n} }.
\end{align*}
Now consider the case if a solution of \eqref{eqn:w1} does not exist, then 
$$
\FDR(\theta_0, \T^\ell(t, \hat w))
\leq (n-s_0) P_{\theta_0 = 1/2}(\T_j^\ell(t, w_0) = 1)  + e^{-C\kappa^2 \Delta_n}.
$$
By the upper bound in Lemma \ref{lem:Plval}, we have
$$
P_{\theta_0 = 1/2}(\ell(x) \leq t) \lesssim 
\frac{\phi(\sqrt{M} \tilde \varepsilon')}{\sqrt{M} \tilde\varepsilon'},
$$
where $\tilde \varepsilon' \sim \sqrt{2(\log(1/r(w_0,t)) + \log(\sqrt{m}))/m}$.
Furthermore,
$$
\frac{\phi(\sqrt{M} \tilde\varepsilon')}{\sqrt{M} \tilde\varepsilon'}
\lesssim \frac{e^{- M^2  (\tilde\varepsilon')^2/2}}{\sqrt{2M (\log (1/r(w_0,t)) + \log\sqrt{m})}}
\lesssim \frac{r(w_0, t)/\sqrt{m}}{\sqrt{2M (\log (1/r(w_0,t)) + \log \sqrt{m})}}.
$$
Therefore, for any $t \leq 4/5$, 
$$
\FDR(\theta_0, \T^{\ell}(t, \hat w))
\lesssim \frac{2(n-s_0) w_0 t}{(1-w_0) (1-t) \sqrt{m \log n}} + e^{-C\kappa^2 \Delta_n} 
\leq \frac{25(1+o(1))\Delta_n t}{\sqrt{m\log n}} + e^{-C\kappa^2 \Delta_n}.
$$

By combining the two cases considered above, we arrive at 
\begin{align*}
\FDR(\theta_0, \T^{\ell}(t, \hat w)) 
\leq 
& \max \Big\{e^{-C'' \Delta_n + \frac{2Ks_n}{\sqrt{m} \xi(w_2)}}  
+ \frac{12}{\sqrt{m \log n}}, \ \frac{25(1+o(1))\Delta_n t}{\sqrt{m\log n}} + e^{-C\kappa^2 \Delta_n}
\Big\}.
\end{align*}
The upper bound in the lemma follows by choosing $\Delta_n = \max\{\frac{1}{2C''} \log (m \log n), \frac{1}{2C\kappa^2} \log(m \log n)\}$.

\section{Proof of results in Section \ref{sec:fdr-fnr-large}}
\label{sec:pf-fdr-fnr}

In this section, we prove results in Section \ref{sec:fdr-fnr-large}. The proof of Theorem \ref{thm:fdr-qval} is given in Section \ref{sec:pf-fdr-q} and proofs of Theorem \ref{thm:fnr-Clval-qval} and Lemma \ref{lem:fnr-lval} are provided in \ref{sec:pf-fnr-Cl-q-l}.

\subsection{Proof of Theorem \ref{thm:fdr-qval}}
\label{sec:pf-fdr-q}

We first obtain the upper bound.
By definition, 
\begin{align*}
\sup_{\theta_0 \in \Theta_0[s_n]} \FDR (\theta_0, \T^q(t, \hat w))
= \sup_{\theta_0 \in \Theta_0[s_n]} 
\mathbb{E}_{\theta_0} \left[
\frac{\sum_{j=1}^n \mathbbm{1}\{\theta_{0,j} = 1/2\} \T_j^q(t, \hat w)}
{1 \vee \sum_{j=1}^n \T_j^q(t, \hat w)}
\right],
\end{align*}
where recall that $\T_j^q(t, \hat w) = \mathbbm{1}\{q(X_j; \hat w, g) \leq t\}$ 
and $q(x; w, g) = \left(1 + \frac{w}{1-w} (g/b)(x)\right)^{-1}$.
Let $\Omega_n = \Omega_0 \cap P(\hat w \in [w_1, w_2])$,
where $\Omega_0 = \{\# \{j \in \mathcal{S}_0, |X_j - m/2|  > bm\zeta_n\} \geq s_n - K_n\}$, 
and $w_1$ and $w_2$ are the solutions for \eqref{eqn:w1} and \eqref{eqn:w2} respectively,
then the last display can be bounded by 
\begin{align} 
\sup_{\theta \in \Theta_0[s_n]}
\mathbb{E}_{\theta_0} 
\left[
\frac{\sum_{j=1}^n \mathbbm{1}\{\theta_{0,j} = 1/2\} \T^q_j(t, \hat w)}{1 \vee \sum_{j=1}^n \T^q_j(t, \hat w)} 
\mathbbm{1}\{\Omega_n\}
\right] + P(\Omega_0^c) + P(\hat w \not\in [w_1, w_2]).
\label{ebfdr-1}
\end{align}
By Lemma \ref{bound-Omega0} and $(i)$ in Lemma \ref{lem:bound-what}, we obtain $P(\Omega_0^c) + P(\hat w \not\in [w_1, w_2]) = o(1)$.
Also, on $\Omega_n$, the denominator of the first term in \eqref{ebfdr-1} can be bounded from below by 
$$
\sum_{j=1}^n \T^q_j(t, \hat w)
\geq
\sum_{j=1}^n \mathbbm{1}\{\theta_{0,j} = 1/2\}\T_j^q(t, \hat w) +
s_n - K_n.
$$
Thus, by concavity and monotonicity of the function $x \in (0, \infty) \to x/(x+1)$, \eqref{ebfdr-1} can be bounded by 
\begin{align}
\label{ebfdr-2}
& \sup_{\theta_0 \in \Theta_0[s_n]} \FDR(\theta_0, \T^q(t, \hat w))
\leq \sup_{\theta \in \Theta_0[s_n]}
\mathbb{E}_{\theta_0} 
\left[
\frac{\sum_{j=1}^n \mathbbm{1}\{\theta_{0,j} = 1/2\} \T^q_j(t, w_1)} 
{\sum_{j=1}^n \mathbbm{1}\{\theta_{0,j} = 1/2\} \T^q_j(t, w_1) + s_n -K_n} 
\right] + o(1) \nonumber \\
& \qquad \quad
\leq \frac{\sup_{\theta_0 \in \Theta_0[s_n]}
\mathbb{E}_{\theta_0} 
\left(\sum_{j=1}^n \mathbbm{1}\{\theta_{0,j} = 1/2\} \T^q_j (t, w_1)\right)}
{\sup_{\theta_0 \in \Theta_0[s_n]}
\mathbb{E}_{\theta_0} \left(\sum_{j=1}^n \mathbbm{1}\{\theta_{0,j} = 1/2\} \T^q_j (t, w_1)\right) + s_n - K_n} + o(1).
\end{align}
To bound the first term of \eqref{ebfdr-2}, for $\eta^q$ given in \eqref{lem:bound-etaq} and by Lemma \ref{lem:Pqval} and then Lemma \ref{lem:bound-w1-w2},
\begin{align*}
\sup_{\theta_0 \in \Theta_0[s_n]} 
& \mathbb{E}_{\theta_0} \left(\sum_{j=1}^n \mathbbm{1}\{\theta_{0,j} = 1/2\} \T_j^q(t; w_1) \right)
\leq C(n - s_n) r(w_1, t)\\
&= C(n - s_n) w_1 t (1-t)^{-1} (1-w_1)^{-1} 
\lesssim s_n (1- C's_n/n) t(1-t)^{-1} (1 + \epsilon),
\end{align*}
as $(1-w_1)^{-1} \leq 1+\epsilon$ for some $\epsilon = o(1)$. 
Therefore, as long as $s_n/n \to 0$, \eqref{ebfdr-2} can be bounded by 
$$
\frac{(1+\epsilon)(1-\epsilon') s_n t(1-t)^{-1}}{(1-\epsilon')(1+\epsilon) s_n t(1-t)^{-1} + s_n - K_n}
\to \frac{t(1-t)^{-1}}{t(1-t)^{-1} + 1 - o(1)} \to t,
$$
as $K_n = o(s_n)$ and $\epsilon, \epsilon' = o(1)$.

Next, we prove the result $\lim_{n\to \infty}\inf_{\theta_0 \in \Theta_0[s_n]} \FDR(\theta_0, \T^q(t, \hat w)) \to t$.
By definition, 
we have
\begin{align*}
& \inf_{\theta_0 \in \Theta_0[s_n]} \mathbb{E}_{\theta_0} 
\left(\frac{\fF_{\T^q}(t, \hat w)}{\fF_{\T^q}(t, w) + \fT_{\T^q}(t, \hat w)} \mathbbm{1}\{\hat w \in [w_1, w_2] \}\right) \\
& \quad \geq 
\inf_{\theta_0 \in \Theta_0[s_n]} \mathbb{E}_{\theta_0} 
\Big(
\frac{\mathbb{E}_{\theta_0} (\fF_{\T^q}(t, \hat w)) (1-\delta)}
{\mathbb{E}_{\theta_0} (\fF_{\T^q}(t, \hat w)) (1-\delta) + s_n} 
\mathbbm{1}\{\hat w \in [w_1, w_2] \} \\
& \qquad \qquad \times 
\mathbbm{1}\{|\fF_{\T^q}(t, \hat w) - \mathbb{E}_{\theta_0} (\fF_{\T^q}(t, \hat w))|  
\leq \delta \mathbb{E}_{\theta_0} (\fF_{\T^q}(t, \hat w))\} 
\Big),
\end{align*}
for some small $\delta$ to be specified later.
On the event $\hat w \in [w_1, w_2]$, we have
\begin{align*}
\mathbb{E}_{\theta_0}(\fF_{\T^q}(t, \hat w)) 
& \geq \mathbb{E}_{\theta_0}(\fF_{\T^q}(t, w_2))
= \sum_{j: \theta_{0,j} = 1/2} \mathbb{E}_{\theta_0} \T_j^q(t, w_2)\\
& = (n - s_n) P_{\theta_{0}} \left(|X_j - m/2| \geq m\eta^q(r(w,t)) - m/2\right)\\
& = 2 (n - s_n) r(w_2, t) \bar\G(m\eta^q(r(w_2, t)) - m/2)
\end{align*}
For any $w \in (0, 1)$ and fixed $t \in (0, 1)$, and by Lemma \ref{lem:etaq-zeta}, we have
$$
\bar\G(m\eta^q(r(w_2, t)) - m/2)
= \frac{m/2 - m\eta^q(r(w_2, t)) }{m+1}
= \frac{m/2 - m\zeta(w_2) + o(m\zeta(w_2))}{m+1}.
$$
Then, there exist an $\varepsilon \in (0, 1)$ such that 
\begin{align*}
\mathbb{E}_{\theta_0}(\fF_{\T^q}(t, \hat w)) 
& \geq (n - s_n) r(w_2, t)(1-\varepsilon)
= (n-s_n) w_2 t (1-w_2)^{-1} (1-t)^{-1} (1-\varepsilon)\\
& = w_2 t (1-w_2)^{-1} (1-t)^{-1} (1-\varepsilon) \sum_{j \in \mathcal{S}_0} m_1(\theta_{0,j}, w_2) (1+\kappa)^{-1} (\tilde m(w_2))^{-1}\\
& \geq (1-\varepsilon')^2 s_n (1+\kappa)^{-1} t(1-t)^{-1}.
\end{align*}
The last inequality is obtained by invoking Lemma \ref{lem:lb-for-m1}.

On the other hand, using Chebychev's inequality, 
\begin{align*}
\sup_{\theta_0 \in \Theta[s_n]} 
& P_{\theta_0} 
\left(|\fF_{\T^q}(t, \hat w) - \mathbb{E}_{\theta_0} (\fF_{\T^q}(t, \hat w))| 
> \delta \mathbb{E}_{\theta_0} (\fF_{\T^q}(t, \hat w))\right) \\
& \leq \frac{\text{Var}_{\theta_0} (\fF_{\T^q}(t, \hat w))}
{\delta^2 (\mathbb{E}_{\theta_0}(\fF_{\T^q}(t, \hat w)))^2}
\leq \frac{1}{\delta^2 \mathbb{E}_{\theta_0}(\fF_{\T^q}(t, \hat w))} \to 0,
\end{align*}
for any fixed $\delta \in (0, 1)$, as $s_n \to \infty$. 

By combining the relevant lower bounds obtained above, we obtain
\begin{align*}
\lim_{n\to \infty}\inf_{\theta_0 \in \Theta_0[s_n]} \FDR(\theta_0, \T^q(t, \hat w))
& \geq \frac{(1-\varepsilon')^{-1}(1-\delta)(1+\kappa)^{-1} t(1-t)^{-1} s_n}{(1-\varepsilon')^{-1}(1-\delta)(1+\nu)^{-1} t(1-t)^{-1}s_n + s_n} +o(1) \\
& \to t + o(1) \to t, \quad \text{as} \ n\to \infty,
\end{align*}
if letting $\kappa \to 0$ and $\delta \to 0$ (but not too fast as long as $\delta s_n \to \infty$; e.g., choosing $\delta = 1/\sqrt{s_n}$).

\subsection{Proof of the FNR results for the $\ell$-value, $q$-value, and C$\ell$-value procedures}
\label{sec:pf-fnr-Cl-q-l}

\subsubsection{Proof of Theorem \ref{thm:fnr-Clval-qval}}

Since the $q$-value procedure is less conservative than the $\adj\ell$-value procedure,
it is sufficient to prove the result for the $\adj\ell$-value procedure only. 
By the definition of $\FNR$, we have
\begin{align}
\label{pf-fnr-Cl}
\sup_{\theta_0 \in \Theta_0[s_n]} \FNR(\theta_0, \T^{\adj\ell}(t, \hat w)) 
= 
\sup_{\theta_0 \in \Theta_0[s_n]} 
\mathbb{E}_{\theta_0} \left(
\frac{s_n - 
\sum_{j=1}^n \mathbbm{1}\{\theta_{0,j} \neq 1/2\} \T_j^{\adj\ell}(t, \hat w))
}{s_n \vee 1} \right).
\end{align}
Let $\tilde \eta^{\adj\ell}(r(\hat w, t)) = \eta^{\adj\ell}(r(\hat w, t)) - 1/2$, then
$$
P_{\theta_0} \T_j^{\adj\ell}(t, \hat w)
\geq P_{\theta_0}(|X_j - m/2| \geq m\tilde \eta^{\adj\ell}(r(\hat w, t))).
$$
From Lemma \ref{lem:etaCl-zeta}, we have
$$
\tilde \eta^{\adj\ell}(r(\hat w, t)) 
\leq \zeta(\hat w) + \frac{C(t, w_0)}{\sqrt{2m \log(1/\hat w)}} + K(\eta_0) \sqrt{\frac{\log(1/\hat w)}{2m}},
$$
where $C(t,w_0) = \log(t(1-t)^{-1}) + C(w_0)$, $C(w_0)$ is a constant.

Consider the event $\mathcal{W} = \{\hat w \in [w_2, w_1]\}$, note that $P_{\theta_0} (\mathcal{W}^c) = o(1)$ by Lemma \ref{lem:bound-what}. On the event $\mathcal{W}$, by Lemma \ref{lem:bound-w1-w2}, we obtain
$\zeta(w_1) \leq \zeta(\hat w) \leq \zeta(w_2) \leq \zeta(C's_n/n)$ and 
$\log(n/(C' s_n)) \leq \log(1/w_1) \leq \log(1/\hat w) \leq \log(1/w_2) < \log n$.
Since $m \gg (\log n)^2$, for a sufficiently large $n$, $\tilde \eta^{\adj\ell}(r(\hat w, t)) \leq 2\zeta(C's_n/n)$. 
We thus can bound $P_{\theta_0} \T_j^{\adj\ell}(t, \hat w)$ from below by 
$P_{\theta_0}(|X_j - m/2| \geq 2m\zeta(C's_n/n))$.
Applying Lemma \ref{bound-Omega0} for some $K_n$ to be specify later (choosing $b = 2$ in Lemma \ref{bound-Omega0}), we obtain 
{\small 
\begin{align*}
\eqref{pf-fnr-Cl}
& \leq 
\sup_{\theta_0 \in \Theta_0[s_n]} 
\left\{
\mathbb{E}_{\theta_0} \left(
\frac{s_n - 
\sum_{j=1}^n \mathbbm{1}\{\theta_{0,j} \neq 1/2\} \T_j^{\adj\ell}(t, \hat w))
}{s_n \vee 1} \mathbbm{1}\{\Omega_n \cap \mathcal{W}\}\right) +  P_{\theta_0}(\mathcal{W}^c) + P_{\theta_0} (\Omega_n^c)
\right\} \\
&\leq  \frac{K_n}{s_n \vee 1} + o(1).
\end{align*}}
\vspace{-0.6cm}

Now choosing $K_n = \max (2s_n p_n, s_n/\log s_n)$ with 
$p_n = 2\bar\B(m/2 + k \sqrt{m \log(n/s_n)/2})$, then $K_n = o(s_n)$.
Thus, the last display goes to $0$ as $s_n \to \infty$.

\subsubsection{Proof of Lemma \ref{lem:fnr-lval}}

Introducing the event $\mathcal{W} = \{\hat w \in [w_2, w_1]\}$, we have
\begin{align}
\lim_n \inf_{\theta_0 \in \Theta_0[s_n]} \FNR(\theta_0, \T^{\ell})
& \geq 
\inf_{\theta_0 \in \Theta_0[s_n]} 
\mathbb{E}_{\theta_0} 
\left(
\frac{\sum_{j =1 }^n \mathbbm{1}\{\theta_{0,j} \neq 1/2\} (1-\T_j^\ell (t; \hat w))}{s_n \vee 1}
\mathbbm{1}\{\mathcal{W}\}
\right) 
\nonumber \\
& \geq 1 - \sup_{\theta_0 \in \Theta_0[s_n]} 
\mathbbm{E}_{\theta_0} \left(
\frac{\sum_{j=1}^n \mathbbm{1}\{\theta_{0,j} \neq 1/2\} \T_j^\ell (t, w_2)}{s_n}
\right)
\nonumber \\
& \geq 1 - \sup_{\theta_0 \in \Theta_0[s_n]} \left(\max_j P_{\theta_0} (\ell(X_j) \leq t)\right),
\label{pf:fnr-lval}
\end{align}
which we used the fact that $\T_j^\ell(t, w)$ is a decreasing function as $w$ increases for each $j$. 
By the definition of the $\ell$-value, 
\begin{align*}
\sup_{\theta_0 \in \Theta_0[s_n]}\left( \max_j P_{\theta_0}(\ell(X_j) \leq t) \right)
& = \sup_{\theta_0 \in \Theta_0[s_n]} 
\left(\max_j P_{\theta_0}\left(\left(b/{g}\right)(X_j) \leq r(w_2,t)\right)\right) \\
& \leq \sup_{\theta_0 \in \Theta_0[s_n]} 
\left(\max_j P_{\theta_0} \left(|\tilde u_j| \geq m\tilde \eta^\ell(r(w_2, t)) - m\zeta(s_n/n)\right)\right),
\end{align*}
which we used $|\theta_{0,j}-1/2| \geq \zeta(s_n/n)$, 
where $\tilde U_j = X_j - m \theta_{0,j}$ is a centered variable and $\tilde \eta^\ell(\cdot) = (b/g)^{-1}(\cdot) - 1/2$.
Applying the Bernstein's inequality in Lemma \ref{bernstein-inq} and letting $A = m\tilde \eta^\ell(r(w_2, t)) - m\zeta_n(n/s_n)$, 
$M = 1$, and $V = \sum_{j=1}^m \theta_{0,j}(1-\theta_{0,j}) \leq m/4$, we obtain
$$
 \sup_{\theta_0 \in \Theta_0[s_n]} 
\left(\max_j P_{\theta_0} \left(|\tilde U_j| \geq m\tilde \eta^\ell(r(w_2, t)) - m\zeta(s_n/n)\right)\right)
\leq 2\exp\left(- \frac{A^2}{m/2 + 2A/3}\right).
$$
By Lemma \ref{lem:etal-xi}, 
$2m|\tilde \eta^\ell(r(w_2, t))) - (\xi(w_2))| \leq C_t/ \xi(w_2)$, and $C_t$ is a fixed constant depending on $t$. 
Since $\xi(w_2) \sim \sqrt{(\log(1/w_2) + \log(\sqrt{m}))/(2m)}$ by \eqref{eqn:xi} and $w_2 \leq w_1 \lesssim s_n/n$ by Lemma \ref{lem:bound-w1-w2}, $2A/3 \leq m/2$ for a sufficiently large $m$. Therefore, the last display can be further bounded by 
\begin{align*}
& 2\exp\left(- \frac{1}{m} 
\left(
\sqrt{\frac{m}{2} \log \left(\left(\frac{C' n }{s_n}\right) + \log(\sqrt{m}) \right)} -  \sqrt{\frac{2m C_t^2}{\log (\sqrt{m} n)}}
- \sqrt{\frac{m}{2} \log \left(\frac{n}{s_n}\right)}
\right)^2
\right) \\
& \qquad 
\leq 2 \exp\left(
- \frac{C' \log(\sqrt{m})}{2} +o(1)
\right) \lesssim 2 m^{- C'/4} \to 0, \ \text{as} \ m \to \infty,
\end{align*}
which we used the inequality $(a - b)^2 > (a^2 +b^2)/2$ for any $a, b > 0$ to obtain the first upper bound. 
Therefore, we have $\eqref{pf:fnr-lval} \geq 1 - 2 m^{- C'/4}$, which goes to $1$ as $m \to \infty$. 

\section{A tight concentration bound for the MMLE $\hat w$}
\label{sec:bound-w-hat}

In this section, 
we obtain a concentration bound for $\hat w$, which is essential for obtaining the uniform FDR control result for our multiple testing procedures, as already shown in the proof of Theorem \ref{thm:unif-fdr-q-Cl}. 
This bound needs to be sharp, and ideally, shall control $\hat w$ around the value $\frac{s_n}{n}$ up to some constant. However, since $\hat w$ is a random quantity, we study the solution for $\mathbb{E}_{\theta_0} S(w) = 0$ instead, denoted by $w^\star$.
Let us consider the following equation and let $w_1$ be its solution:
\begin{align}
& \sum_{j \in \mathcal{S}_0} m_1(\theta_{0,j}, w) = (1-\kappa)(n - s_0) \tilde m (w), 
\label{eqn:w1}
\end{align}
where $\kappa \in (0, 1)$ is a fixed constant, $w \in [w_0, 1)$, $\theta_0 \in l_0[s_n]$,
$\mathcal{S}_0 = \{1 \leq j \leq n: \theta_{0,j} \neq 1/2\}$,
and $\tilde m(w)$ and $m_1(\theta_{0,j}, w)$ are quantities related to the score function given in \eqref{tildem} and \eqref{m1} respectively.
Depending on $\kappa, n, m, s_n$ and the true value $\theta_0$, a solution for \eqref{eqn:w1} may or may not exist. 
If a solution exists, it must be unique, as $\tilde m(w)$ is monotone increasing and $\tilde m(u, w)$ is monotone decreasing (see Lemma \ref{lem:m1-tildem-mnt}).
In $(i)$ of Lemma \ref{lem:bound-what}, we will show that $\hat w \in [w_2, w_1]$ in high probability, where
$w_2$ is the solution of 
\begin{align}
& \sum_{j \in \mathcal{S}_0} m_1(\theta_{0,j}, w) = (1+\kappa)(n - s_0) \tilde m (w).
\label{eqn:w2}
\end{align}
By Lemma \ref{lem:bound-w1-w2}, we have $w_2 \leq w_1 \lesssim s_n/n$ and $w_1/K \leq w_2$.
Hence, $\hat w$ shall concentrate on a neighborhood around $s_n/n$ up to some constant. 

If a solution for \eqref{eqn:w1} does not exist,
we will show in the second point in Lemma \ref{lem:bound-what} that $\hat w \leq w_0$ with a high probability if $1/n \leq w_0 \leq \rho_n/n$ for some $\rho_n \ll n$.
The lower bound for $w_0$ is the lower bound we imposed for estimating $\hat w$; see \eqref{w-hat}.

In the next lemma, we show that both $\tilde m(w)$ and $m_1(u, w)$ are monotone, continuous, and nonnegative functions, which are helpful for analyzing $\hat w$ later on.

\begin{lemma}
\label{lem:m1-tildem-mnt}
For $w \in (0, 1)$ and $u \in [0, m]$,
$\tilde m(w)$ is a nonnegative, continuous, monotone increasing function and $m_1(u, w)$ is a continuous, monotone decreasing function. 
\end{lemma}

\begin{proof}
Since $\beta(u, w)$ is a decreasing function with $w$,
by their definitions,
clearly, $\tilde m(\cdot)$ is a monotone increasing function and 
$m_1(u, w)$ is a decreasing function. By Lemma \ref{lem:bound-mtilde}, $\tilde m(w)$ is nonnegative. 
The continuity result for $\tilde m(w)$ (resp. $m_1(u, w)$) follows by noting that $\beta(u,w)$ is continuous and dominates the term $\beta(u, w)b(u)$ (resp. $\beta(u, w)b_\theta(u)$) by $g(u) + b(u)$ (resp. $g(u) + b_\theta(u)$) up to a constant.
\end{proof}

\begin{lemma}
\label{lem:bound-w1-w2}
Let $w_1$ and $w_2$ be the solution of \eqref{eqn:w1} and \eqref{eqn:w2} respectively,
then 
$$
w_1/K < w_2 < w_1 \lesssim s_n/n.
$$
for some constant $K > 4$.
\end{lemma}

\begin{proof}
We first show $w_1 \lesssim s_n$.
By the definition of $m_1(\theta, w)$,
\begin{align}
(1-\kappa)(n - s_0) \tilde m (w_1) \leq \sum_{j \in \mathcal{S}_0} m_1(\theta_{0,j}, w_1)
\leq s_n \max_{j \in \mathcal{S}_0} m_1(\theta_{0,j}, w_1),
\label{pf:bound-w1-w2-1}
\end{align}
where $\mathcal{S}_0 = \{j: \theta_{0,j} \neq 1/2\}$.
The upper bound in \eqref{pf:bound-w1-w2-1} can be further bounded by invoking Lemma \ref{lem:ub-m1}.
Depending on the value of $\theta_{0,j}$, $m_1(\theta_{0,j}, w_1)$ is bounded by \eqref{eqn:ub-m1}, \eqref{eqn:ub-m1-2} or \eqref{eqn:ub-m1-3}.
If $|\theta_{0,j} - 1/2| > \frac{\Lambda}{\sqrt{2m}}$ for some fixed constant $\Lambda$, then
$$ 
m_1(\theta_{0,j}, w_1) \leq \frac{2}{w_1} \bar \B_\theta(m/2 + m\xi(w_1))  
+ \frac{4}{w_1} \bar \Phi\left(2\sqrt{m} (\xi(w_1) - |\mu_{0,j}|)\right) 
\leq \frac{6}{w_1},
$$ 
where $\mu_{0,j} = \theta_{0,j} - 1/2$.
If $\frac{1}{2m\xi(w_1)} \leq |\theta_{0,j} - 1/2| \leq \frac{\Lambda}{\sqrt{2m}}$, then $1-4\mu_{0,j}^2 \approx 1$ and
$$
m_1(\theta_{0,j}, w_1) \lesssim e^{-2m(\xi(w_1) - \mu)^2 + 2m\zeta^2(w_1) + 2m\nu^2}
\leq \frac{1}{w_1} e^{-2m(\xi(w_1) - \mu)^2 + \log \sqrt{m}} \leq \frac{1}{w_1},
$$
which we used that $\xi^2(w) \sim \zeta^2(w) + \nu^2$ in Lemma \ref{lem:relation-xi-zeta-nu}.
If $|\theta_{0,j} - 1/2| < \frac{1}{2m\xi(w_1)}$,  
then
$$
m_1(\theta_{0,j}, w_1) \lesssim \zeta(w_1) + w_1^{C_{\mu_{0,j}}}/\sqrt{m} < 1,
$$
where $C_{\mu_{0,j}} = \frac{4\mu_{0,j}^2}{1-4\mu_{0,j}^2}$.
By combining the three cases, we obtain $\max_{j} m_1(\theta_{0,j}, w_1) \lesssim w_1^{-1}$. 
Therefore, \eqref{pf:bound-w1-w2-1} implies 
$$
w_1 \lesssim \frac{s_n}{(1-\kappa)(n-s_0) \tilde m(w_1)}
\leq \frac{Cs_n}{n}, 
$$
for a fixed constant $C$ depends on the values of $\kappa, v_1, v_2$ and $\mu_{0,j}$,
as $\tilde m(w_1) \in [0.4, 1.1]$ by Lemma \ref{lem:bound-mtilde} and $s_0/n \leq s_n/n \to 0$ as $n \to \infty$. 

Next, we verify $w_2 < w_1$.
From Lemma \ref{lem:m1-tildem-mnt}, 
$m_1(\cdot, w)$ is a continuous monotone decreasing function 
and $\tilde m(w)$ is a continuous monotone increasing function, thus, the ratio
$m_1(\cdot, w)/\tilde m(w)$ is also monotone decreasing. 
\eqref{eqn:w1} and \eqref{eqn:w2} implies that, for any $\kappa \in (0, 1)$,
$$
\frac{\sum_{j \in \mathcal{S}_0} m_1(\theta_{0,j}, w_1)/\tilde m(w_1)}
{\sum_{j \in \mathcal{S}_0} m_1(\theta_{0,j}, w_2)/\tilde m(w_2)}
= \frac{1-\kappa}{1+\kappa}
< 1,
$$
which implies $w_2 < w_1$.

Last, we show $w_2 > w_1/K$. 
Let us introduce the set 
$$
\mathcal{J}_0 := \mathcal{J}(\theta_0, w, K) = \{1 \leq j \leq n: |\theta_{0,j} - 1/2| \geq \zeta(w)/K\},
$$
and define
$$
\mathcal{M}^{\mathcal{S}_0}(w) = \sum_{j\in \mathcal{S}_0} m_1(\theta_{0, j}, w),
\quad 
\mathcal{M}^{\mathcal{J}_0}(w, K') = \sum_{j\in \mathcal{J}_0} m_1(\theta_{0, j}, w).
$$
Since $m_1(\cdot, w)/\tilde m(w)$ is a monotone decreasing function,
it is sufficient to show
\begin{align}
\label{pf:bound-w1-w2-2}
\frac{\mathcal{M}^{\mathcal{S}_0}(w_1/K)}{\tilde m(w_1/K)} 
> 
\frac{\mathcal{M}^{\mathcal{S}_0}(w_2)}{\tilde m(w_2)}
>
\frac{(1+\kappa) \mathcal{M}^{\mathcal{S}_0}(w_1)}{(1-\kappa) \tilde m(w_1)}.
\end{align}
By Lemma \ref{lem:m1-cntl-small-signals} with $w= w_1$ or $w_1/K$,
$$
\sup_{\theta_0 \in l_0[s_n]} 
\sup_{w \in [1/n, 1/\log n]} 
|\mathcal{M}^{\mathcal{S}_0}(w) - \mathcal{M}^{\mathcal{J}_0}(w, K)| \leq C  
n^{1-D}, \quad D \in (0, 1), \ C >0,
$$
Therefore, for some $K > 1$ to be chosen later and using Lemma \ref{lem:m1-ratio}, we have
\begin{align*}
\mathcal{M}^{\mathcal{S}_0}(w_1/K)
\geq 
\mathcal{M}^{\mathcal{J}_0}(w_1/K, K') - Cn^{1-D} 
& \geq K \mathcal{M}^{\mathcal{J}_0}(w_1, K') - Cn^{1-D} \\
& \geq K \mathcal{M}^{\mathcal{S}_0}(w_1) - 2Cn^{1-D}.
\end{align*}
Note that $\mathcal{M}^{\mathcal{S}_0}(w_1) = O(n)$, but the second term in the lower bound of the last display is $o(n)$, then for a sufficiently large $n$, 
$\mathcal{M}^{\mathcal{S}_0}(w_1/K) \geq K\mathcal{M}^{\mathcal{S}_0}(w_1)/2$. 
Furthermore, by Lemma \ref{lem:bound-mtilde}, we obtain
$$
\frac{\mathcal{M}^{\mathcal{S}_0}(w_1/K)}{\tilde m(w_1/K)} 
\geq \frac{K\mathcal{M}^{\mathcal{S}_0}(w_1)}{4\tilde m(w_1)}. 
$$
Then \eqref{pf:bound-w1-w2-2} follows by choosing $K > 4(1+\kappa)/(1-\kappa) > 4$. 
\end{proof}

\begin{lemma}
\label{lem:bound-what}
Let $w_1$ and $w_2$ be solutions of \eqref{eqn:w1} and \eqref{eqn:w2} respectively, suppose $(\log n)^2/m \to 0$, $s_n = n^{v_1}$ and $m= n^{v_2}$ for $v_1 \in (0, 1)$ and $v_2 > \log \log n/\log n$,
\begin{itemize}
\item[(i)] if \eqref{eqn:w1} has a solution,
then for a sufficiently large $n$, there exists some positive constant $C$ such that for 
$\theta_0 \in l_0[s_n]$ and any fixed $\kappa \in (0, 1)$, then 
$$
P_{\theta_0} (\hat w \not\in [w_2, w_1]) 
\leq e^{-C\kappa^2 n w_1 \tilde m(w_1) }+ e^{-C\kappa^2 n w_2 \tilde m(w_2) }.
$$

\item[(ii)] if \eqref{eqn:w1} does not have a solution, let $w_0$ be the solution of 
\begin{align}
\label{eqn:w0}
n w_0 \tilde m(w_0) = \Delta_n, \quad \Delta_n \in [1.1, \rho_n],
\end{align}
for some $1.1 < \rho_n \ll n$,
then for a sufficiently large $n$ and the same $\kappa, C$ as in ${(i)}$,
$$
P_{\theta_0} (\hat w \geq w_0) \leq e^{-C\kappa^2 n w_0 \tilde m(w_0)} = e^{- C\kappa^2 \Delta_n}. 
$$
\end{itemize}
\end{lemma}

\begin{remark}
\label{rmk:w0}
For a sufficiently large $m$, $\tilde m(w_0) \in [0.4, 1.1]$ by Lemma \ref{lem:bound-mtilde}, providing that  $w_0 \to 0$, which is true since we require $\Delta_n/n \to 0$.
Indeed, by rewriting \eqref{eqn:w0}, one immediately obtains $w_0 = \frac{\Delta_n}{n \tilde m(w_0)}$. 
Since $\Delta_n \geq 1.1$, $w_0 \geq 1/n$. Also, since $\Delta_n \leq \rho_n = o(n)$, $w_0 < 1$. Therefore, for $\hat w \in [w_0, 1]$, $\hat w$ still belongs to the range $[1/n, 1]$ given in \eqref{w-hat}.
\end{remark}

\begin{proof}
We first prove $(i)$: since a solution for \eqref{eqn:w1} exists, the event $\{\hat w \geq w_1 \}$ implies $\{{S}(w_1)\geq 0\}$, and then, 
$$
P_{\theta_0} (\hat w \geq w_1) 
= 
P_{\theta_0} ({S}(w_1) \geq 0)
=  
P_{\theta_0} ({S}(w_1)  - \mathbb{E}_{\theta_0} {S}(w_1) \geq - \mathbb{E}_{\theta_0} {S}(w_1))
$$
Using \eqref{eqn:w1}, which implies $\mathbb{E}_{\theta_0} {S}(w_1) = \sum_{j \in \mathcal{S}_0} m_1(\theta_{0, j}, w_1) - (n-s_0) \tilde m(w_1) \leq - \kappa (n-s_0) \tilde m(w_1)$, the last display is bounded by 
\begin{align}
\label{pf:bd-what-1}
P_{\theta_0} ({S}(w_1)  - \mathbb{E}_{\theta_0} {S}(w_1) 
\geq \kappa (n-s_0) \tilde m(w_1)).
\end{align}
We thus need to bound \eqref{pf:bd-what-1}.
Let $W_j = \beta(X_j, w) - \mathbb{E}_{\theta_0} \beta(X_j, w)$, $j \in \{1, \dots, n\}$, then $W_j$ is a centered variable, independently with $W_{j'}$ for $j \neq j'$. 
Also, $|W_j| \leq |\beta(X_j, w_1)| \leq w_1^{-1}$ and 
\begin{align*}
\sum_{j = 1}^n \text{Var}(W_j) \leq \sum_{j=1}^n m_2(\theta_{0,j}, w_1)
& = \sum_{\{j: |\mu_{0,j}| > \frac{\Lambda}{\sqrt{2m}}\}} m_2(\theta_{0,j}, w_1) 
+ \sum_{\{j: |\mu_{0,j}| \leq \frac{\Lambda}{\sqrt{2m}}\}}m_2(\theta_{0,j}, w_1) \\
& = (a) + (b).
\end{align*}
To bound $(a)$, we use that $m_2(\theta_{0,j}, w_1) \lesssim w^{-1} m_1(\theta_{0,j}, w_1)$ from Corollary \ref{cor:m1-m2} and obtain
\begin{align*}
(a) & \lesssim \frac{1}{w} \sum_{\{j: |\mu_{0,j}| > \frac{\Lambda}{\sqrt{2m}}\}} m_1(\theta_{0,j}, w_1)
= \frac{1}{w} (1-\kappa)(n-s_0)\tilde m(w_1) 
- \frac{1}{w} \sum_{\{j: |\mu_{0,j}| \leq \frac{\Lambda}{\sqrt{2m}}\}} m_1(\theta_{0,j}, w_1) \\
& \leq w_1^{-1} (1-\kappa)(n-s_0)\tilde m(w_1)
\leq w_1^{-1} (1-\kappa)n \tilde m(w_1).
\end{align*}

To bound $(b)$, since $\#\{j: |\mu_{0,j}| \leq \frac{\Lambda}{\sqrt{2m}}\} \leq s_0 \leq s_n$ and
by \eqref{eqn:ub-m1-2} and \eqref{eqn:ub-m1-3} in Lemma \ref{lem:ub-m1},
\begin{align*}
(b) 
& \lesssim \frac{1}{w_1} 
\Bigg\{
\sum_{j: \frac{1}{2m\xi} \leq |\mu_{0,j}| \leq \frac{\Lambda}{\sqrt{2m}}}
\left(\frac{1}{\sqrt{m|\xi - \mu_{0,j}|}} + \frac{1}{(m+1)|\mu_{0,j}|} \right)
\exp\left(\frac{-2m\mu_{0,j}^2 + 4m |\mu_{0,j}| \xi}{1-4\mu_{0,j}^2}\right) \\
& \qquad  \quad + 
\sum_{j: |\mu_{0,j}| < \frac{1}{2m\xi}}
\left(\zeta + w^{C_{\mu_{0,j}}}/\sqrt{m}\right)
\Bigg\} \\
& \lesssim \frac{s_n}{w_1} e^{-2m (\xi - \Lambda/\sqrt{2m})^2 + 2m\xi^2}
\leq \frac{s_n}{w_1} e^{-2m (\xi - \Lambda/\sqrt{2m})^2 + 2m \nu^2 + 2m\zeta^2 }
\leq \frac{s_n}{w_1} e^{2m\zeta^2} = \frac{s_n}{w_1^2}.
\end{align*}
where we used $\zeta^2 + \nu^2 \sim \xi^2$ and $2m\zeta(w_1)^2 \asymp \log(1/w_1)$.
Since $w_1 \leq Cs_n/n$ by Lemma \ref{lem:bound-w1-w2},
the last display is bounded by $w_1^{-1} n/C$. 
Therefore, 
$$\sum_{j=1}^n \text{Var}(W_j) \leq (a) + (b) \leq 2w_1^{-1}n \tilde m(w_1).$$ 
Now we are ready to bound \eqref{pf:bd-what-1}. By applying the Bernstein inequality in Lemma \ref{bernstein-inq}, choosing
$A = \kappa(n-s_0)\tilde m(w_1)$, $M \leq w_1^{-1}$, and $V = 2 w_1^{-1} n\tilde m(w_1)$
and noticing that $s_0/n \leq s_n/n = n^{v_1 - 1} \to 0$ as $n\to \infty$, 
we thus obtain
$$
P_{\theta_0}(\hat w \geq w_1)
\leq e^{- \frac{3}{14}\kappa^2 n w_1 \tilde m(w_1)}.
$$

To bound the probability $P_{\theta_0} (\hat w \leq w_2)$ is similar. 
We have
$$
P_{\theta_0} (\hat w \leq w_2) 
= 
P_{\theta_0} ({S}(w_2) \leq 0)
=  
P_{\theta_0} ({S}(w_2)  - \mathbb{E}_{\theta_0} {S}(w_2) \leq - \mathbb{E}_{\theta_0}S(w_2))
$$
Define $\mathbb{E}_{\theta_0} {S}(w_2) = \sum_{j \in {S}_0} m_1(\theta_{0, j}, w_2) - (n-s_0) \tilde m(w_2) \leq - \kappa (n-s_0) \tilde m(w_2)$ by \eqref{eqn:w2}, one then needs to bound the probability
\begin{align*}
P_{\theta_0} ({S}(w_2)  - \mathbb{E}_{\theta_0} {S}(w_2) 
\leq \kappa (n-s_0) \tilde m(w_2)).
\end{align*}
Applying the Bernstein inequality again, one obtains 
$P_{\theta_0}(\hat w \leq w_2) \leq e^{- \frac{3}{14}\kappa^2 n w_2 \tilde m(w_2)}.$
The result in $(i)$ follows by combining the two upper bounds obtained above.

The proof for $(ii)$ is similar. Note that $w_0$ the solution of \eqref{eqn:w0}. Using 
$\sum_{j \in \mathcal{S}_0} m_1(\theta_{0,j}, w_0) < (1-\kappa) (n-s_0) \tilde m(w_0)$ leads to 
\begin{align*}
P_{\theta_0}(\hat w \geq w_0) & = P_{\theta_0}({S}(w_0) \geq 0)
= P_{\theta_0}({S}(w_1) - \mathbb{E}_{\theta_0}{S}(w_0) \geq 
- \mathbb{E}_{\theta_0}{S}(w_0)) \\
& \leq 
 P_{\theta_0}({S}(w_1) - \mathbb{E}_{\theta_0}{S}(w_0) \geq 
- \kappa (n-s_0) \tilde m(w_0)).
\end{align*}
By applying the Bernstein inequality, one then obtains the upper bound. 
\end{proof}

\section{Bounding $\tilde m(w)$, $m_1(\theta, w)$, and $m_2(\theta,w)$}
\label{sec:bound-m1-mtilde}

In this section, we obtain bounds for $\tilde m(w)$, $m_1(u, w)$, and $m_2(u,w)$ in \eqref{tildem}--\eqref{m2} respectively. These are essential quantities for studying the MMLE $\hat w$ in the previous section. 
Although similar quantities arise in the study of the Gaussian sequence model in \citet{johnstone04} and \citetalias{cast20MT}, our bounds are different as we work with the Binomial distribution. We will comment on the difference between those bounds in \citetalias{cast20MT} and ours in Remark \ref{rmk:diff-Gaussian-binomial}.

\subsection{Upper and lower bounds for $\tilde m(w)$}

\begin{lemma}
\label{lem:bound-mtilde}
For $\tilde m(w)$ given in \eqref{tildem}, 
let $\xi = |\xi_n(w)|$ as in Lemma \ref{lem:beta}, 
if $\log (1+1/w)/m \to 0$ as $m \to \infty$,
then
$$
\frac{m/2 - m\xi}{m+1} - \frac{1}{\sqrt{m}(1+w^{-1})}
\leq \tilde m(w) 
\leq 1 + \frac{2w\xi}{1-w},
$$
Furthermore, if $w \to 0$, then for a sufficiently large $m$, $0.4 \leq \tilde m(w) \leq 1.1$.
\end{lemma}

\begin{proof}
Recall that $\tilde m(w) = - \sum_{u=0}^m \beta(u)(1 + w\beta(u))^{-1} b(u)$.
Since $g(u) = (1+m)^{-1}$, we have $\sum_{u=0}^m \beta(u) b(u) = \sum_{u=0}^m g(u) - 1 = 0$.
We can write
\begin{align}
\label{pf:tildem-1}
\tilde m(w) 
= \sum_{u=0}^m \beta(u) b(u) - \sum_{u=0}^m \frac{\beta(u)  b(u)}{1 + w\beta(u)}
= \sum_{u=0}^m \frac{w\beta(u)^2b(u)}{1 + w\beta(u)},
\end{align}
which is always positive as $1 + w\beta(m/2) \geq 0$ for any $w\in (0, 1)$ by Lemma \ref{lem:bound-beta}. Note that $\beta(u)$ is symmetric at $u = m/2$, thus
\begin{align}
\label{pf:bound-tildem-0}
\sum_{u=m/2+1}^m \frac{2w\beta(u)^2b(u)}{1 + w\beta(u)}
\leq 
\tilde m(w) 
\leq
\sum_{u=m/2}^m \frac{2w\beta(u)^2b(u)}{1 + w\beta(u)}.
\end{align}
Note that as we mentioned earlier, we assume $m$ to be an even number throughout the paper; if $m$ is an odd number, the last display needs to be replaced by $\tilde m(w) = \sum_{u = \lceil m/2 \rceil}^m \frac{2w\beta(u)^2b(u)}{1 + w\beta(u)}$. In either case, the bounds for $\tilde m(w)$ stay the same.  

{\it Obtain the lower bound.} The lower bound in \eqref{pf:bound-tildem-0} can be also written as
\begin{align}
\label{pf:bound-tildem-1}
\sum_{u=m/2+1}^m \frac{2w\beta(u)^2b(u)}{1 + w\beta(u)} =
\underbrace{\sum_{u = m/2 + 1}^{m/2 + m\xi} \frac{2w\beta^2(u)  b(u)}{1 + w\beta(u)}}_{(I)}
+ \underbrace{\sum_{u = m/2 + m\xi}^{m} \frac{2w\beta(u) (g - b)(u)}{1 + w\beta(u)}}_{(II)},
\end{align}
which we used the fact that $\beta(u) b(u) = (g-b)(u)$.
Since $1 + w\beta(u) > 0$ for $u \in [m/2+1, m/2 + m\xi]$, $(I) > 0$.
Using that $\beta(u) \geq 0$ when $u > m/2 + m\xi$ by Lemma \ref{lem:beta}, we obtain $g(u) \geq b(u)$ and $w \beta(u) \geq 1$. We thus have
$$
(II) \geq \sum_{u = m/2 + m\xi}^m (g-b)(u)
\geq \frac{m/2 - m\xi}{1+m}  - \bar{\mathbf{B}}(m/2 + m\xi).
$$
By Lemma \ref{lem:binom-tail} and then (2) in Lemma \ref{lem:entropy-bound}, we get 
\begin{align*}
\bar{\mathbf{B}}(m/2 + m\xi)
\leq e^{-m T(1/2+\xi, 1/2)}
\leq e^{-2m\xi^2}.
\end{align*}
By \eqref{eqn:xi}, $2m\xi^2 \sim \log (1+w^{-1}) + \log(\sqrt{m})$, which leads to
$$
(II) \geq \frac{m/2 - m\xi}{m+1} - \frac{K_1}{\sqrt{m}(1+w^{-1})}
$$ 
for some constant $K_1 = 1+o(1)$.
The lower bound for $\tilde m(w)$ follows by combining the lower bounds for $(I)$ (which is 0) and $(II)$.
If $m \to \infty$, then $\tilde m(w) \to 1/2 - \xi(w) \geq 0.4$, providing that $w$ is small.

{\it Obtain the upper bound.} 
Using $\beta(u)b(u) = g(u) - b(u)$, 
the upper bound in \eqref{pf:bound-tildem-0} can be written as
\begin{align}
\label{pf:bound-tildem-1.2}
\tilde m(w) \leq
\underbrace{\sum_{u = m/2}^{m/2 + m\xi} \frac{2w\beta^2(u)  b(u)}{1 + w\beta(u)}}_{(III)}
- \underbrace{\sum_{u = m/2 + m\xi}^{m} \frac{2w\beta(u) b(u)}{1 + w\beta(u)}}_{(IV)}
+ \underbrace{\sum_{u = m/2 + m\xi}^{m} \frac{2w\beta(u) g(u)}{1 + w\beta(u)}}_{(V)}.
\end{align}
We first bound $(III)$. Since $\beta(\cdot)$ is a monotone increasing function on $[m/2, m]$ (see Lemma \ref{lem:monotone-beta}) and $-1< \beta(m/2) < 0$ (see Lemma \ref{lem:solve-beta}), 
we have $1 + w\beta(u) \geq 1 + w\beta(m/2) \geq 1-w$ for any $u \in [m/2, m/2 + m\xi]$.
Thus, using that $\beta(m/2 + m\xi) = 1/w$, we obtain
\begin{align}
(III) & \leq \frac{2 w \beta(m/2+m\xi)}{1 + w\beta(m/2)}
\sum_{u = m/2}^{m/2+m\xi} \beta(u) b(u)
\leq \frac{2 w \beta(m/2+m\xi)}{1 - w} 
\sum_{u = m/2}^{m/2+m\xi} (g(u) - b(u) )\nonumber \\
& \leq \frac{2m \xi}{(1- w)(1+m)} \leq \frac{2\xi}{1-w}.
\label{pf:bound-tildem-2}
\end{align}

Next, we bound $(V)$.
Using that $1 \leq {2w \beta(u)}/(1+ w\beta(u)) \leq 2$ as $w\beta(u) \geq 1$ for $u \geq m/2 + m\xi$,
we obtain
$$
\frac{m/2 - m\xi}{m+1} 
= \sum_{u=m/2 + m\xi}^m g(u)
\leq (V)
\leq \sum_{u=m/2 + m\xi}^m 2 g(u)
= \frac{m - 2m\xi}{m+1}
\leq 1 - 2\xi.
$$
By summing up the upper bounds for $(III)$ and $(V)$, and note that $(IV)$ is smaller than $(V)$ due to $b(u) \leq g(u)$ for $u \in [m/2 + m\xi, m]$, as $\xi > |\nu_n|$ by Lemma \ref{lem:solve-beta}, we obtain the upper bound for $\tilde m(w)$.
\end{proof}

\subsection{Upper and lower bounds for $m_1(\theta, w)$}
\begin{lemma} (Upper bound for $m_1(\theta, w)$)
\label{lem:ub-m1}
For $m_1(\theta, w)$ given in \eqref{m1}, $\theta \in (0, 1)$, and $\xi = |\xi_n(w)|$ as in \eqref{eqn:xi},
suppose $m\xi^4 \to 0$ as $m \to \infty$ and $\mu = \theta - 1/2$,
then
there exist a $w_\circ \in (0,1)$ and a fixed $\mu_0 < 1/2$ such that for any 
$w \leq w_\circ$ and 
$\frac{\Lambda}{\sqrt{2m}} < \mu \leq \mu_0$, $\Lambda > 0$ is a fixed constant, we have
\begin{align}
\label{eqn:ub-m1}
m_1(\theta, w) \leq 
\begin{cases}
\frac{2}{w} \bar \B_\theta(m/2 + m\xi) + 
\frac{1}{w} \bar \Phi \left(2\sqrt{m}|\mu - \xi|\right)T_m(\mu, \xi), & \quad \text{if} \ \theta \geq 1/2\\
\frac{2}{w} {\B}_\theta(m/2 - m\xi + 1) + 
\frac{1}{w} \bar \Phi \left(2\sqrt{m}|\mu -\xi|\right)T_m(\mu, \xi), & \quad \text{if} \ \theta < 1/2,
\end{cases}
\end{align}
where 
\begin{align}
\label{eqn:Tm}
T_m(\mu, \xi) = \frac{|\xi - \mu|}{\mu \sqrt{1-4\xi^2}} .
\end{align}
If $\frac{1}{2m\xi} \leq \mu \leq \frac{\Lambda}{\sqrt{2m}}$, then
\begin{align}
m_1(\theta, w) \lesssim
\left(\frac{\sqrt{2}}{\sqrt{m|\xi - \mu|}}+\frac{1}{(m+1)\mu} \right) 
e^{\frac{-2m\mu^2 + 4m\mu \xi}{1-4\mu^2}}.
\label{eqn:ub-m1-2}
\end{align}
If $0 < \mu < \frac{1}{2m\xi}$, for $C_\mu = \frac{4\mu^2}{1-4\mu^2}$, then
\begin{align}
m_1(\mu, w) \lesssim 2e^{\frac{2}{1-4\mu^2}} \left(\zeta + \frac{w^{C_\mu}}{\sqrt{m}}\right).
\label{eqn:ub-m1-3}
\end{align} 
In addition, if $m$ is sufficiently large such that $\xi \log m \to 0$, then 
\begin{align}
\label{eqn:ub-m1-4}
m_1(\theta, w) \lesssim
\frac{1}{w} \left( \bar \Phi\left(\frac{2\sqrt{m} (\xi - \mu)}{\sqrt{1-4\mu^2}}\right) 
+ \bar \Phi\left(2\sqrt{m} (\xi - \mu)\right) T_m(\mu, \xi)\right).
\end{align}
\end{lemma}

\begin{remark}
As $\xi \sim \sqrt{\frac{\log(1/w) + \log \sqrt{m}}{2m}}$, the assumption $\xi \log m \to 0$ implies $\log(1/w)(\log m)^2/m\to 0$ and $(\log m)^3/m \to 0$, both are stronger conditions than $m \gg (\log n)^2$. Therefore, when we prove the bound for the MMLE $\hat w$ and other relevant results, we use the upper bound in \eqref{eqn:ub-m1} rather than \eqref{eqn:ub-m1-4}.
\end{remark}

\begin{proof} 
We consider only the case $\theta \geq 1/2$, as the proof for $\theta < 1/2$ is similar. Note that $\mu = \theta - 1/2$ in this case, since $\theta \geq 1/2$.
For $m_1(\theta, w)$ given in \eqref{m1}, one can write
\begin{align}
m_1(\theta, w) & = 
\underbrace{\sum_{|u - m/2| \geq m\xi} \beta(u, w) b_\theta(u)}_{(I)} + 
\underbrace{\sum_{|u - m/2| < m\xi} \beta(u, w) b_\theta(u)}_{(II)}.
\label{pf:bound-m1-1}
\end{align}
Since $w \beta(u) \geq 1$ when $|u - m/2| \geq m\xi$ by Lemmas \ref{lem:monotone-beta} and \ref{lem:beta}, 
we have
\begin{align}
(I) \leq \frac{1}{w}\sum_{|u - m/2| \geq m\xi} b_\theta(u)
& = \frac{1}{w}\left(\bar{\mathbf{B}}_\theta(m/2 + m\xi ) + {\B}_\theta (m/2 - m\xi + 1)\right)
\leq \frac{2}{w}\bar{\mathbf{B}}_\theta(m/2 + m\xi).
\label{pf:m1-ub-1}
\end{align}
To bound $(II)$, one writes
\begin{align}
(II) & = \sum_{|u - m/2| < m\xi} 
\left(\beta(u)- \frac{w \beta^2(u)}{1+w\beta(u)}\right)  b_{\theta}(u)  
\leq 
\sum_{|u - m/2| < m\xi} 
\beta(u)  b_{\theta}(u)   \nonumber \\
& \leq
\sum_{m\nu < |u - m/2| < m\xi}  \beta(u)  b_{\theta}(u)
\label{pf:m1-ub-1.2}
\end{align}
for $\nu = |\nu_n(w)|$ given in \eqref{eqn:nu},
as $\beta(u) < 0$ for $|u - m/2| < m\nu$ by Lemma \ref{lem:solve-beta}.
The last display in \eqref{pf:m1-ub-1.2} can be further bounded by 
\begin{align}
\sum_{m \nu < |u - m/2| < m\xi} \frac{g(u) b_{\theta}(u)}{b(u)} 
& = \frac{1}{m+1}  
\sum_{ m\nu < |\tilde u| < m\xi} \frac{b_{\theta}(\tilde u + m/2)}{b(\tilde u+m/2)} \nonumber \\
& < \frac{2}{m+1} \sum_{m\nu < \tilde u < m\xi} 
e^{mT(1/2 + \tilde u/m, 1/2) - m T(1/2+\tilde u/m, 1/2+\mu)},
\label{pf:m1-ub-2}
\end{align}
where
$T(a, p) = a\log (a/p) + (1-a)\log((1-a)/(1-p))$.
By Lemma \ref{lem:entropy-bound}, 
\begin{align}
& mT(1/2 + \tilde u/m, 1/2) - m T(1/2+\tilde u/m, 1/2+\mu) \nonumber \\
& \quad = \frac{2\tilde u^2}{m} - \frac{2m(\mu - \tilde u/m)^2}{1-4\mu^2} 
+ m\varpi_m(\mu, \tilde u/m) \nonumber \\
& \quad 
= \frac{4m}{1-4\mu^2} \left(\frac{\mu \tilde u}{m} - \frac{\mu^2}{2}\right) 
- \frac{8 \mu^2 \tilde u^2}{m(1 - 4\mu^2)} + m\varpi_m(\mu, \tilde u/m),
\label{pf:m1-ub-2-2}
\end{align}
where
\begin{align}
\label{pf:m1-ub-2-3}
\varpi(\mu, \tilde u/m) = 
\frac{(\tilde u/m)^3 \epsilon_1}{3(1/2 + \epsilon_1)^2(1/2-\epsilon_1)^2}
- \frac{(\mu - \tilde u/m)^3 (\mu + \epsilon_2 )}{3(1/2 + \mu + \epsilon_2)^2(1/2 - \mu -\epsilon_2)^2}
\end{align}
for some $\epsilon_1 \in [0, \tilde u/m]$ and 
$\epsilon_2 \in [0, \mu - \tilde u/m]$ if $\mu \geq \tilde u/m$ and $\epsilon_2 \in [\mu - \tilde u/m, 0]$ if 
$\mu < \tilde u/m$.
As $\tilde u$ is at most $m\xi$, 
the first term in \eqref{pf:m1-ub-2-3} is $\lesssim (\tilde u/m)^4 \leq \xi^4$.
We only need to consider the case when the second term is negative (when $\tilde u/(2m) \leq \mu \leq \tilde u/m$), as otherwise, this term can be simply bounded by 0. 
Using $\tilde u \leq m\xi$ again, then $\mu = O(\xi)$ and the second term is bounded by 
$C\xi^4$ for some constant $C > 0$.  
Hence, $m\varpi(\mu, \tilde u/m) \leq C_1 m\xi^4$ for some $C_1 > C$.
We thus obtain
\begin{align}
(II) \leq \eqref{pf:m1-ub-2}
& = 
\frac{2}{m+1} \sum_{m\nu < \tilde u < m\xi}
e^{\frac{4m}{1-4\mu^2} \left(\frac{\mu \tilde u}{m} - \frac{\mu^2}{2} \right)
- \frac{8\mu^2 \tilde u^2}{m(1-4\mu^2)} + m\varpi(\mu, \tilde u/m)}
\nonumber \\
& \leq \frac{2}{m+1} e^{C_1 m\xi^4}
 \sum_{m\nu < \tilde u < m\xi}
e^{\frac{4m}{1-4\mu^2}\left(\frac{\mu \tilde u}{m} - \frac{\mu^2}{2} \right)}
\label{pf:m1-ub-7}
\\
& =
\frac{2}{m+1} e^{- \frac{2m \mu^2}{1-4\mu^2} + C_1 m\xi^4}
\times
\frac{e^{\frac{4\mu m \xi}{1-4\mu^2}} - e^{\frac{4\mu m \nu}{1-4\mu^2}}}{e^{\frac{4\mu}{1-4\mu^2}} - 1}
\nonumber \\
& \leq 
\frac{1 - 4\mu^2}{2(m+1)\mu} 
e^{- \frac{2m \mu^2}{1-4\mu^2}
+ \frac{4m \mu\xi}{1-4\mu^2}
+ C_1 m\xi^4} 
\nonumber\\
& \lesssim
\frac{1}{2(m+1)\mu} e^{\frac{4m\mu \xi - 2m \mu^2}{1-4\mu^2}},
\label{pf:m1-ub-6}
\end{align}
which we used the inequality $e^{x} - 1 \geq x$ for $x > 0$ to obtain the second inequality in the last display and the assumption $m \xi^4 \to 0$ as $m \to \infty$ to obtain the last inequality.
By collecting the relevant bounds obtained above, we arrive at
\begin{align}
m_1(\theta, w) \lesssim \frac{2}{w} \bar \B_\theta(m/2 + m\xi ) 
+ \frac{1}{2(m+1)\mu} e^{\frac{4m\mu \xi - 2m\mu^2}{1-4\mu^2}}
= (a) + (b).
\label{pf:m1-ub-4}
\end{align}
When $\mu \geq \frac{\Lambda}{\sqrt{2m}}$,
using that $g(x) = (m+1)^{-1}$ and $(g/b)(m/2 + m\xi) \asymp w^{-1}$, we have
\begin{align}
(b) = 
& \frac{g(m/2 + m\xi)}{2\mu}  
e^{\frac{4m\mu \xi - 2m\mu^2}{1-4\mu^2}}
= \frac{b(m/2 + m\xi)}{2 w \mu} 
\frac{\phi \left(2\sqrt{m}|\xi - \mu| \right)}
{\phi \left(2\sqrt{m}\xi \right)} 
e^{\frac{8m\mu^2}{1-4\mu^2} (\xi^2 - (\xi - \mu)^2)}
\nonumber \\
& \quad \quad \qquad \qquad 
\leq \frac{\sqrt{m}|\xi - \mu|}{w\mu} 
\bar\Phi \left(2\sqrt{m}|\xi - \mu| \right)
\frac{b(m/2 + m\xi)}{\phi \left(2\sqrt{m}\xi \right)}
e^{\frac{8m\mu^2}{1-4\mu^2} (\xi^2 - (\xi - \mu)^2)}.
\label{pf:m1-ub-5}
\end{align}
The last line is obtained by using the inequality $\phi(x) \leq |x + x^{-1}| \bar\Phi(x) \approx |x| \bar \Phi(x)$ when $1/x \to 0$ (see Lemma \ref{lem:gaussian-bound}).
Note that the exponential term in \eqref{pf:m1-ub-5} is bounded by 
$e^{\frac{8m \xi^4}{1-4\mu^2}} = 1+o(1)$, as $m\xi^4 = o(1)$.
Also, using Lemma \ref{lem:bin-coef-bound} and then Lemma \ref{lem:entropy-bound}, we obtain
\begin{align*}
\frac{b(m/2 + m\xi)}{\phi \left( 2\sqrt{m}\xi \right)}
& \leq 
\frac{1}{\sqrt{m(1-4\xi^2)}}
e^{-m T(1/2+\xi, 1/2) + 2m\xi^2} 
\leq 
\frac{1+o(1)}{\sqrt{m(1-4\xi^2)}}.
\end{align*}
Thus,
\begin{align}
\eqref{pf:m1-ub-5} 
\lesssim
\frac{1}{w} \bar \Phi \left(2\sqrt{m}|\xi - \mu|\right)
\left(
\frac{|\xi - \mu|}{\mu \sqrt{1-4\xi^2}} 
\right).
\label{pf:m1-ub-4-1}
\end{align}
Therefore,
$$
m_1(\theta, w) 
\lesssim
\frac{2}{w} \bar \B_\theta(m/2 + m\xi)
+ 
\frac{1}{w} \bar \Phi \left(2\sqrt{m}|\xi - \mu|\right)
\left(
\frac{|\xi - \mu|}{\mu \sqrt{1-4\xi^2}} 
\right),
$$
which leads to the first line in \eqref{eqn:ub-m1}.
When $\theta < 1/2$, 
\eqref{pf:bound-m1-1} can be bounded by $\frac{2}{w} \B_\theta(m/2 - m\xi)$. By following essentially the same proof as above, one can obtain the second line in \eqref{eqn:ub-m1}.

If $\frac{1}{2m\xi} \leq \mu < \frac{\Lambda}{\sqrt{2m}}$,
using $(g/b)(m/2 + m\xi) \asymp w^{-1}$ again and $\bar \Phi(x) \leq \phi(x)/x$ by Lemma \ref{lem:gaussian-bound}, we have
\begin{align*}
\frac{2}{w} \bar \B_\theta(m/2 + m\xi) 
& \leq 
\frac{2 g(m/2 + m\xi)}{b(m/2 + m\xi)} 
\phi \left(\frac{2\sqrt{m}|\xi - \mu|}{\sqrt{1-4\mu^2}}\right)
\times
\frac{\sqrt{1-4\mu^2}}{2\sqrt{m|\xi - \mu|}}  \\
& \leq \frac{\sqrt{2m}}{1+m} \sqrt{\frac{1-4\xi^2}{|\xi - \mu|}} 
e^{-\frac{2m(\xi - \mu)^2}{1-4\mu^2} + mT(1/2 + \xi, 1/2)} \\
& \leq \sqrt{\frac{2(1-4\mu^2)}{m|\xi - \mu|}} 
e^{-\frac{2m(\xi - \mu)^2}{1-4\mu^2} + 2m\xi^2 + o(1)}  \\
& \lesssim
 \sqrt{\frac{2(1-4\mu^2)}{m|\xi - \mu|}} 
e^{\frac{4m}{1-4\mu^2}(\mu\xi - \mu^2/2)}.
\end{align*}
Similarly, we obtain 
$$
\frac{1}{w} \bar \Phi \left(2\sqrt{m}|\xi - \mu|\right)
\left(
\frac{|\xi - \mu|}{\mu \sqrt{1-4\xi^2}} 
\right)
\lesssim \frac{1}{(m+1)\mu}e^{\frac{4m}{1-4\mu^2}(\mu\xi - \mu^2/2)}.
$$
By combining the preceding two upper bounds and then using the inequality that $1 - 4\mu^2 \leq 1$, we obtain \eqref{eqn:ub-m1-2}.

Last, if $0 < \mu < \frac{1}{2m\xi}$, then
\begin{align*}
(I) & \leq \frac{2}{w} \bar\B_\theta(m/2 + m\xi)
\leq\frac{2}{w} e^{- mT(1/2 + \xi, 1/2 + \mu)}
\leq\frac{2}{w} e^{- \frac{2m (\mu - \xi)^2}{1-4\mu^2} + C_4 m\xi^4} \\
& 
\leq \frac{2}{w}e^{- \frac{2m (\mu - \xi)^2}{1-4\mu^2}}
\leq \frac{2}{w}e^{\frac{- 2m\xi^2 + 2}{1-4\mu^2}}
\leq \frac{2e^{\frac{2}{1-4\mu^2}}}{\sqrt{m}w^{C_\mu}},
\end{align*}
where $C_\mu = \frac{4\mu^2}{1-4\mu^2}$. We used the fact that $2m\xi^2 \asymp \log (1/w) + \log(\sqrt{m})$ to obtain the last inequality in the last display.
From \eqref{pf:m1-ub-7}, we also have
\begin{align}
(II) 
& \leq 
\frac{2e^{C_1 m\xi^4} }{m+1} 
\sum_{m \nu \leq \tilde u \leq m\xi}
e^{\frac{4\mu\tilde u - 2m\mu^2}{1-4\mu^2}}
\leq 
\frac{2 e^{\frac{4m\mu \nu}{1-4\mu^2} + C_1 m\xi^4} }{m+1} \sum_{0\leq \tilde u\leq m\zeta}
e^{\frac{4\mu\tilde u}{1-4\mu^2}} 
\nonumber \\
& \lesssim \frac{2m\zeta}{m+1} e^{\frac{4m\mu (\nu + \zeta)}{1-4\mu^2}}
\leq 2\zeta e^{\frac{2}{1-4\mu^2}},
\label{pf:m1-ub-8}
\end{align}
as $\mu < 1/(2m\xi)$.
The second inequality in the last display is obtained by using $\xi - \nu \leq \zeta$.
Thus, \eqref{eqn:ub-m1-3} is obtained by combining the upper bounds for $(I)$ and $(II)$.

To prove \eqref{eqn:ub-m1-4}, 
let's write $\bar{\mathbf{B}}_\theta(m/2 + m\xi) = P_\theta(U \geq m/2 + m\xi)$ for $U \sim \text{Bin}(m,  \theta)$ and denote $\tilde U = U - m/2$, $Z_m = \frac{2\tilde U - 2m\mu}{\sqrt{1 - 4\mu^2}}$, and
$W_m$ as the standard Brownian motion, 
then
\begin{align*}
P_\theta(U \geq m/2 + m\xi)
& = \mathbb{E}_\theta\Bigg(
\mathbbm{1}\left\{Z_m \geq \frac{2 m\xi - 2m\mu}{\sqrt{1 - 4\mu^2}}\right\}
\mathbbm{1}\left\{|Z_m - W_m| \geq \log m + x\right\}
 \\
& \quad \quad +
\mathbbm{1}\left\{Z_m \geq \frac{2 m\xi - 2m\mu}{\sqrt{1 - 4\mu^2}}\right\}
\mathbbm{1}\left\{|Z_m - W_m| < \log m + x\right\}
\Bigg)\\
& = (c) + (d).
\end{align*}
By applying the KMT approximation theorem in Lemma \ref{lem:KMT} (w.l.o.g we choose $C = 1$ in the lemma),
$(c) \leq P_\theta(|Z_m - W_m| \geq \log m + x) \lesssim e^{-K x}$ for some positive constant $K$.
Also, we have
\begin{align*}
(d) & \leq P \left(W_m \geq \frac{2m\xi - 2m\mu}{\sqrt{1-4\mu^2}} - (\log m + x) \right)  \\
& =  P \left(\frac{W_m}{\sqrt{m}} \geq \frac{2\sqrt{m}(\xi - \mu)}{\sqrt{1-4\mu^2}} 
- \frac{\log m + x}{\sqrt{m}} \right) \\
& = \bar \Phi\left(
\frac{2\sqrt{m}(\xi - \mu)}{\sqrt{1-4\mu^2}} - \frac{\log m + x}{\sqrt{m}}
\right).
\end{align*}
Using the lower bound in $(iii)$ in Lemma \ref{lem:gaussian-bound-carter} and
letting $\delta = (\log m + x)^2/(2m)$ and $z = \frac{2\sqrt{m}(\xi - \mu)}{\sqrt{1-4\mu^2}} - \delta$,
we obtain 
$\bar \Phi(z) \leq \bar \Phi(z + \delta) \exp(\rho(z) \delta + \delta^2/2)$,
where $\rho(z) = \phi(z)/\bar\Phi(z)$ and by Lemma \ref{lem:gaussian-bound}, $\rho(z) \leq z + 1/z$. 
By plugging-in the expressions for $z$ and $\delta$, the last display can be bounded by 
$$
\bar \Phi\left(\frac{2\sqrt{m}(\xi - \mu)}{\sqrt{1-4\mu^2}} \right)
e^{
\frac{2(\xi - \mu) (\log m + x)}{\sqrt{1-4\mu^2}} 
- \frac{(\log m + x)^2}{2m} 
+ \frac{(\log m + x)^2}{m} \left(\frac{2\sqrt{m}(\xi - \mu)}{\sqrt{1-4\mu^2}} -  \frac{(\log m + x)^2}{m}\right)^{-1}
}.
$$
Choosing $x = \log m$, 
then $\frac{\log m + x}{\sqrt{m}} = \frac{2\log m}{\sqrt{m}}=  o(1)$ and $\left(\frac{2\sqrt{m}(\xi - \mu)}{\sqrt{1-4\mu^2}} -  \frac{(\log m + x)^2}{m}\right)^{-1} \to 0$ as $m \to \infty$, thus
\begin{align*}
(c) + (d) 
& \lesssim e^{-K\log m} + 
\bar \Phi\left(
\frac{2\sqrt{m}(\xi - \mu)}{\sqrt{1-4\mu^2}}
\right)  
\exp\left(
\frac{4(\xi - \mu) \log m}{\sqrt{1-4\mu^2}} 
-o(1)
\right).
\end{align*}
From the last display, one can see that the assumption $\xi \log m \to 0$ is needed when $\xi > \mu$, as otherwise, the second term in the last display will diverge.
The result immediately follows by using $e^{-K\log m} \to 0$ and $\xi \log m \to 0$.
\end{proof}

\begin{lemma} (Lower bound for $m_1(\theta, w)$)
\label{lem:lb-m1}
For $m_1(\theta, w)$ given in \eqref{m1} with $\theta \in (0, 1)$,
let $\xi = |\xi_n|$ and $\nu = |\nu_n|$ for $\xi_n$ and $\nu_n$ in \eqref{eqn:xi} and \eqref{eqn:nu} respectively,
suppose $m \xi^4 \to 0$ as $m \to \infty$, 
then for any $w \leq w_\circ \in (0, 1)$ and a fixed $\mu_0 < 1/2$ such that $\mu = \theta - 1/2$, $\mu = \mu_0 \geq \frac{\Lambda}{\sqrt{2m}}$, and some positive constant $\Lambda$,
\begin{align}
\label{eqn:lb-m1}
m_1(\theta, w) \gtrsim 
\frac{1}{w}
\left(
\bar \B_\theta(m/2 + m\xi) +
\bar \Phi\left( \frac{2\sqrt{m}(\xi - \mu)}{\sqrt{1-4\mu^2}}\right)  T'_m(\mu, \xi)\right), 
\end{align}
where for $D_\mu = e^{\frac{4\mu}{1-4\mu^2}} - 1$ and $\xi := \xi(w)$,
$$
T'_m(\mu, \xi) = \frac{|\xi - \mu|}{D_\mu \sqrt{1-4\mu^2}}.
$$
Moreover, for a sufficiently large $m$, we have
\begin{align}
\label{eqn:lb-m2}
m_1(\theta, w) \gtrsim 
\begin{cases}
\frac{1}{w}
\left(
\bar \Phi\left( \frac{2\sqrt{m}(\xi - \mu)}{\sqrt{1-4\mu^2}}\right)  (1+ T'_m(\mu, \xi))
\right), \ & \ \text{if} \ \xi > \mu, \\
\frac{1}{w}
\left(
\bar \Phi\left( \frac{\sqrt{m}(\xi - \mu)}{\sqrt{1/2 + \mu}} \right)+ 
\bar \Phi\left( \frac{2\sqrt{m}(\xi - \mu)}{\sqrt{1-4\mu^2}} \right) T'_m(\mu, \xi)
\right) \ & \ \text{if} \ \xi \leq \mu.
\end{cases}
\end{align}
\end{lemma}

\begin{proof}
We write $m_1(\theta, w)$ as a summation of three terms given by
\begin{equation}
\begin{split}
\label{pf:m1-lb-1}
m_1(\theta, w) & = 
\underbrace{\sum_{|u - m/2| \geq m\xi} \beta(u, w) b_\theta(u)}_{(I)} + 
\underbrace{\sum_{m\nu < |u - m/2| < m\xi} \beta(u, w) b_\theta(u)}_{(II)} \\
& \quad + \underbrace{\sum_{|u - m/2| \leq m\nu} \beta(u, w) b_\theta(u)}_{(III)} 
\end{split}
\end{equation}
Since $1+w\beta(u) \leq 2w\beta(u)$ if $|u - m/2| \geq m\xi$, 
\begin{align*}
(I) \geq \frac{1}{2w} \sum_{|u - m/2| \geq m\xi} b_{\theta}(u) 
& = \frac{1}{2w} \left(\bar{\mathbf{B}}_{\theta} (m/2 + m\xi) + {\B}_\theta(m/2 - m\xi + 1)\right) 
\nonumber \\
& \geq \frac{1}{2w} \bar{\mathbf{B}}_{\theta} (m/2 + m\xi).
\end{align*}

Next, by Lemma \ref{lem:beta},
$0 \leq w\beta(u) \leq 1$ 
for $u \in [m/2 + m\nu, m/2 + m \xi]$,
and by Lemma \ref{lem:bound-beta},
$\beta(m/2) > -1$, we have
\begin{align}
(II) 
& = \sum_{m\nu < |u - m/2| < m\xi} 
\beta(u, m) b_{\theta}(u) 
\geq 
\sum_{m\nu < |u - m/2| < m\xi} 
\frac{\beta(u)}{1 + w\beta(m/2)} b_\theta(u) 
\nonumber \\
& \geq
\sum_{m\nu < |u - m/2| < m\xi}  \frac{g(u) b_\theta(u)}{(1-w) b(u)} -
 \frac{1}{1-w} \sum_{m\nu < |u - m/2| < m\xi}  b_\theta(u)\\
&  = (a) + (b).
\label{pf:m1-lb-2}
\end{align}
We first derive a lower bound for $(a)$.
Since $g(u) = (1+m)^{-1}$,
\begin{align}
(a) = \frac{1}{(1+m)(1-w)} \sum_{m\nu < |u - m/2| < m\xi} \frac{b_\theta(u)}{b(u)}.
\label{pf:m1-lb-3}
\end{align}
Let $T(a, p) = a\log (a/p) + (1-a) \log ((1-a)/(1-p))$ and $\tilde u = u - m/2$, then the ratio 
\begin{align}
\label{pf:m1-lb-3-0}
b_\theta(u)/b(u)
= \exp\left(m T(1/2+\tilde u/m, 1/2) - mT(1/2+\tilde u/m, 1/2 + \mu)\right).
\end{align}
By $(2)$ in Lemma \ref{lem:entropy-bound}, 
$$
mT(1/2 + \tilde u/m, 1/2) \geq 2 m (\tilde u/m)^2.
$$
If $\tilde u/m > \mu$, then by $(3)$ in Lemma \ref{lem:entropy-bound},
$$
mT(1/2+\tilde u/m, 1/2 + \mu)
\leq 
\frac{2m (\mu - \tilde u/m)^2}{1-4\mu^2} + 
\frac{16 (\tilde u/m - \mu)^3 \tilde u}{3(1-4(\tilde u/m)^2)^2}.
$$
Since $\nu < \tilde u/m < \xi$, the second term in the last display is $Km \xi^4$ for some constant $K > 16/3$.
If $0 < \tilde u/m < \mu$, then by $(3)$ in Lemma \ref{lem:entropy-bound},  
$$
mT(1/2+\tilde u/m, 1/2 + \mu) \leq \frac{2m (\mu - \tilde u/m)^2}{1-4\mu^2}.
$$
Therefore, by combining the relevant bounds, we obtain 
$$
\eqref{pf:m1-lb-3-0} \geq \exp\left( 2 m (\tilde u/m)^2 -  \frac{2m (\mu - \tilde u/m)^2}{1-4\mu^2} - Km\xi^4 \right).
$$
Note that the first two terms can be understood as the ratio between two Gaussian distributions, in which
their means and variables correspond to those of the two binomial distributions $b_\theta(u)$ and $b(u)$ respectively. 
By the lower bound in the last display and noticing that $(1-w)^{-1} \geq 1/2$ for a sufficiently large $n$, we have 
\begin{align*}
\eqref{pf:m1-lb-3} & \geq 
\frac{1}{2(m+1)} 
\sum_{m\nu < \tilde u <m\xi} 
e^{\frac{4m}{1-4\mu^2}\left(\frac{\mu \tilde u}{m} - \frac{\mu^2}{2}\right) - \frac{8 \mu^2\tilde u^2}{m(1-4\mu^2)} - Km\xi^4}\\
& \geq 
\frac{1}{2(m+1)} 
e^{- \frac{2m\mu^2}{1-4\mu^2} - \frac{8m\mu^2 \xi^2}{1-4\mu^2} - Km\xi^4}
\sum_{m\nu < \tilde u <m\xi} 
e^{\frac{4\mu\tilde u}{1-4\mu^2}} \\
& = 
\frac{1}{2(m+1)} 
e^{- \frac{2m\mu^2}{1-4\mu^2} - \frac{8m\mu^2 \xi^2}{1-4\mu^2} - Km\xi^4}
\times 
\frac{e^{\frac{4m\mu \xi}{1-4\mu^2}} - e^{\frac{4m\mu \nu}{1-4\mu^2}}}
{e^{\frac{4\mu}{1-4\mu^2}} - 1}
\\
& = 
\frac{1}{2 D_\mu (m+1)} 
e^{- \frac{2m\mu^2}{1-4\mu^2} - \frac{8m\mu^2 \xi^2}{1-4\mu^2} - Km\xi^4}
\left(
e^{\frac{4m\mu \xi}{1-4\mu^2}}
- e^{\frac{4m\mu \nu}{1-4\mu^2}}
\right) \\
& \geq 
\frac{1}{4 D_\mu (m+1)} 
e^{- \frac{2m\mu^2 - 4m\mu \xi }{1-4\mu^2} - \frac{8m\mu^2 \xi^2}{1-4\mu^2} - Km\xi^4},
\end{align*}
for $D_\mu = e^{\frac{4\mu}{1-4\mu^2}} - 1$.
Using $g(u) = (1+m)^{-1}$ and $(g/b)(m/2 + m\xi) \asymp 1/w$, 
by the assumption $m\xi^4 \to 0$, then for a sufficiently large $m$, $e^{-Km\xi^4} \geq 1/2$, then
the last line in the last display can be written as
\begin{align}
& \frac{b(m/2 + m\xi)}{2w D_\mu}
\frac{\phi \left(\frac{2\sqrt{m}(\xi - \mu)}{\sqrt{1-4\mu^2}} \right)}
{\phi \left(\frac{2\sqrt{m}\xi}{\sqrt{1-4\mu^2}} \right)}
e^{-\frac{8m\mu^2\xi^2}{1-4\mu^2}} 
\nonumber \\
& \quad \geq 
\frac{\sqrt{m}|\xi - \mu|}{w D_\mu \sqrt{1-4\mu^2}}
\bar \Phi
\left(
\frac{2\sqrt{m}(\xi - \mu)}{\sqrt{1-4\mu^2}}
\right)
\frac{b(m/2 + m \xi)}{\phi\left( \frac{2\sqrt{m} \xi}{\sqrt{1-4\mu^2}}\right)} 
e^{-\frac{8m\mu^2\xi^2}{1-4\mu^2}}, 
\label{pf:m1-lb-4}
\end{align}
which we used the inequality $\phi(x) \geq \bar\Phi(x) x$ for any $x > 0$ by Lemma \ref{lem:gaussian-bound}. 
Applying Lemma \ref{lem:bin-coef-bound} and then $(2)$ in Lemma \ref{lem:entropy-bound}, we then obtain
\begin{align*}
\frac{b(m/2 + m \xi)}{\phi\left( \frac{2\sqrt{m} \xi}{\sqrt{1-4\mu^2}}\right)} 
e^{-\frac{8m\mu^2\xi^2}{1-4\mu^2}}
\gtrsim \frac{e^{-mT(1/2 + \xi, 1/2) + 2m\xi^2} }{\sqrt{m(1-4\xi^2)}}
\geq \frac{1}{\sqrt{m(1-4\xi^2)}}\geq \frac{1}{\sqrt{m}}.
\end{align*}
Thus, $(a) \gtrsim \bar \Phi
\left(
\frac{2\sqrt{m}(\xi - \mu)}{\sqrt{1-4\mu^2}}
\right)
\frac{|\xi - \mu|}{w D_\mu \sqrt{1-4\mu^2}}$.

Next, we will bound $(b)$ in \eqref{pf:m1-lb-2} and $(III)$ in \eqref{pf:m1-lb-1} together, as they are both negative.
Using $0 > \beta(u) > -1$ for $|u - m/2| < m\nu_n$, then 
\begin{align*}
|(b) + (III)| \leq \frac{1}{1 + w\beta(m/2)} \sum_{|u - m/2| \leq m\nu} b_\theta(u) 
\leq \frac{1}{1-w} \sum_{|u - m/2| \leq m\nu} b_\theta(u)
\leq \frac{1}{1-w}.
\end{align*}
Last, by combining lower bounds for $(I)$, $(II)$, and $(III)$ and note that $w/(1-w) \to 0$ as $n \to \infty$, we thus obtain \eqref{eqn:lb-m1}.

To prove \eqref{eqn:lb-m2}, 
it is sufficient to obtain a lower bound for $\bar{\mathbf{B}}_{\theta} (m/2 + m\xi)$. 
If $\mu \geq \xi$, then by Lemma \ref{lem:binom-tail-slud}, we immediately have
$$
\bar{\mathbf{B}}_{\theta} (m/2 + m\xi)
\geq \bar \Phi\left(
\frac{\sqrt{m}(\xi - \mu)}{\sqrt{1/2 + \mu}}
\right).
$$
If $0 < \mu < \xi$, then by Lemma \ref{lem:binom-tail-mckay}, since $m \xi^2 \to \infty$, as $m \to \infty$, 
let $\sigma = \sqrt{m(1-4\mu^2)}/2$ and $Y(z) = \bar \Phi(z)/\phi(z)$ with $z = \frac{2 \sqrt{m}(\xi - \mu)}{\sqrt{1-4\mu^2}}$, and denote $b_\theta(m-1; k) = P(X = k)$,  where $X \sim \text{Bin}(m-1, \theta)$, 
we have
\begin{align*}
\bar{\mathbf{B}}_{\theta} (m/2 + m\xi)
= \sigma \bar \Phi\left(\frac{2\sqrt{m}(\xi - \mu)}{\sqrt{1 - 4\mu^2}}\right) 
\frac{b_\theta(m-1; m\xi+m/2 - 1)}{\phi \left(\frac{2\sqrt{m}(\xi - \mu)}{\sqrt{1 - 4\mu^2}}\right)} 
e^{\frac{|\xi - \mu|}{m}}.
\end{align*}
By Lemma \ref{lem:bin-coef-bound} and then $(2)$ in Lemma \ref{lem:entropy-bound}, we obtain
\begin{align*}
\frac{b_\theta(m-1; m\xi+m/2 - 1)}{\phi \left(\frac{2\sqrt{m}(\xi - \mu)}{\sqrt{1 - 4\mu^2}}\right)} 
& = \frac{1/2 + \xi}{1/2 + \mu} \times \frac{b_\theta(m\xi+m/2)}{\phi \left(\frac{2\sqrt{m}(\xi - \mu)}{\sqrt{1 - 4\mu^2}}\right)} \\
& > \frac{2}{\sqrt{m(1-4\xi^2)}} 
e^{-mT(1/2 + \xi, 1/2+\mu) + \frac{2m(\xi - \mu)^2}{1-4\mu^2}} \\
& \geq \frac{2}{\sqrt{m(1-4\xi^2)}} 
e^{- K m\xi^4},
\end{align*}
for some $K > 16/3$.
Therefore, we obtain
$$
\bar{\mathbf{B}}_{\theta} (m/2 + m\xi)
\geq \frac{2 \sigma}{\sqrt{m(1-4\xi^2)}} 
e^{- K m\xi^4 + \frac{|\xi - \mu|}{m}} 
\bar \Phi\left(\frac{2\sqrt{m}(\xi - \mu)}{\sqrt{1 - 4\mu^2}}\right) 
\gtrsim \bar \Phi\left(\frac{2\sqrt{m}(\xi - \mu)}{\sqrt{1 - 4\mu^2}}\right)
$$
by plugging-in the expression of $\sigma$ and using $|\xi - \mu|/m \to 0$ and the assumption $m\xi^4 \to 0$.
\end{proof}

Results in Lemmas \ref{lem:ub-m1} and \ref{lem:lb-m1} lead to the following corollary:
\begin{corollary}
\label{cor:m1}
For $m_1(\theta, w)$ given in \eqref{m1}, $\theta \in (0, 1)$, and $\xi = |\xi_n|$ and $\nu = |\nu_n|$ for $\xi_n$ and $\nu_n$ in \eqref{eqn:xi} and \eqref{eqn:nu} respectively,
suppose $m \xi^4 \to 0$, $w \to 0$, as $m \to \infty$,
then for any $1/2 > \mu  = \mu_0 \geq \frac{\Lambda}{\sqrt{2m}}$, $\mu = \theta - 1/2$ and $\Lambda > 0$ and some positive constant $C$, 
\begin{align*}
m_1(\theta, w) & \gtrsim
\begin{cases}
\frac{1}{w}
\left(
\bar \B_\theta(m/2 + m\xi) 
+ \bar \Phi\left(\frac{2\sqrt{m}(\xi - \mu)}{\sqrt{1-4\mu^2}}\right) T_m'(\mu, \xi)
\right) & \ \text{if} \ \mu \geq 0, \\
\frac{1}{w}
\left({\B}_\theta(m/2 - m\xi + 1) 
+ \bar \Phi\left(\frac{2\sqrt{m}(\xi - |\mu|)}{\sqrt{1-4\mu^2}}\right) T_m'(\mu, \xi)
\right) & \ \text{if} \ \mu < 0, 
\end{cases}
\\
m_1(\theta, w) & \lesssim
\begin{cases}
\frac{1}{w}
\left( \bar \B_\theta(m/2 + m\xi) 
+ \bar \Phi\left(2\sqrt{m}(\xi - \mu) \right)T_m(\mu, \xi)\right)
 & \ \text{if} \ \mu \geq 0, \\
\frac{1}{w}
\left( {\B}_\theta(m/2 - m\xi + 1) 
+ \bar \Phi\left(2\sqrt{m}(\xi - |\mu|) \right) T_m(\mu, \xi)\right)
& \ \text{if} \ \mu < 0.
\end{cases}
\end{align*}
\end{corollary}

\begin{remark}
\label{rmk:diff-Gaussian-binomial}
We compare the bounds obtained for $\tilde m(w)$ and $m_1(\theta, w)$ here to those in \citetalias{cast20MT} for the Gaussian sequence model. 
In our model, establishing bounds for $\tilde m_0(w)$ are relatively easier than that in the Gaussian sequence model, as $g(x)$ is a constant. However, bounding $m_1(\theta, w)$ presents some  challenges due to 1) the need to establish precise bounds for non-centered binomial distributions with $\theta_0 \neq 1/2$;
and 2) the necessity to control the ratio between two distributions---one arising from the null hypothesis and the other from the alternative hypothesis. In the Gaussian sequence model, the primary difference lies in the means of these distributions. In our model, both the means and variances are different. Consequently, the difference in the variance between the two distributions leads to a gap between the upper and lower bounds for $m_1(\theta, w)$ in Corollary \ref{cor:m1}. However, in the Gaussian sequence model, the upper and lower bounds for $m_1(\theta ,w)$ are of the same order, as the variances of distributions under the null and the alternative hypotheses are the same.
\end{remark}

\begin{lemma}
\label{lem:lb-for-m1}
For $m_1(\theta, w)$ given in \eqref{m1}, $\theta \in (0, 1)$, and $\xi = |\xi_n|$ given in \eqref{eqn:xi},
suppose $m \xi^4 \to 0$, $w \to 0$, as $m \to \infty$,
then for any $1/2 > \mu  \geq (1 + \rho) \xi(w)$ with any $w \leq w_0 \in (0, 1)$ and $\rho > 0$, there exists a $\varepsilon \in (0, 1)$ such that $m_1(\theta, w) \geq {(1-\varepsilon)}/{w}$.
\end{lemma}

\begin{proof}
Let $a = 1 + \rho/2$, for $w$ is small enough, we can write 
\begin{align*}
wm_1(\mu, w) 
& =\left\{ 
\sum_{\tilde u = - am\xi}^{am\xi}  + 
\sum_{|\tilde u| > am\xi}
\right\}
\frac{w\beta(\tilde u + m/2)}{1 + w \beta(\tilde u + m/2)}
b_\theta (\tilde u + m/2) \\
& \geq \sum_{|\tilde u| > am\xi}
\frac{w\beta(\tilde u + m/2)}{1 + w \beta(\tilde u + m/2)} b_\theta (\tilde u + m/2)
- \sum_{\tilde u = - am\xi}^{am\xi}  b_\theta (\tilde u + m/2) \\
& \geq 
\frac{w\beta(a m \xi + m/2)}{1 + w \beta(a m \xi + m/2)} b_\theta (am\xi + m/2)
- {\B}_{\theta}(am\xi + 1).
\end{align*}
Since $\mu \geq (1+\rho)\xi$, we have $b_\theta (am\xi + m/2)  \to 1$ when $w \to 0$.  
If $w \beta(am\xi + m/2 ) \to \infty$ (which we will soon show $w \beta(am\xi + m/2 ) \to \infty$), then the first term in the last display $\to 1$. 
Let $\varepsilon = {\B}_{\theta}(am\xi + 1) \in (0, 1)$, then the proof is completed. 
What remains is to show $w \beta(am\xi + m/2 ) \to \infty$.
Since $1/w = \beta(m\xi + m/2 )$, $\beta(u) \asymp (g/\phi)(u)$ and $g(u) = (1+m)^{-1}$, 
using Lemma \ref{lem:bin-coef-bound} and then $(2)$ in Lemma \ref{lem:entropy-bound}, we obtain
$$
\frac{\beta(am\xi + m/2 )}{\beta(m\xi + m/2 )}
\asymp \frac{b(m\xi + m/2)}{b(am\xi + m/2)}
\gtrsim e^{m T(1/2 + a\xi, 1/2) - mT(1/2 + \xi, 1/2)}
\geq e^{2(a-1)m\xi^2 + o(1)}.
$$
Since $a > 1$, the lower bound in the last display goes to $\infty$ as $m\to \infty$. 
\end{proof}


\subsection{Upper bound for $m_2(\theta, w)$}
\begin{lemma}
\label{lem:ub-m2}
Consider $m_2(\theta, w)$ as in \eqref{m2} and let $\xi = |\xi_n|$ and $\nu = |\nu_n|$ for $\xi_n$ and $\nu_n$ given in \eqref{eqn:xi} and \eqref{eqn:nu} respectively, 
suppose 
$m\xi^4 = o(1)$, 
$w \asymp s_n/n$ and $s_n = n^{v_1}$ and $m = n^{v_2}$ for $v_1\in (0,1)$,
then for any $(0, 1) \ni \theta \neq 1/2$,
\begin{align}
\label{eqn:ub-m2}
m_2(\theta, w) \leq 
\begin{cases}
\frac{2}{w^2} \left(
{\bar \B}_\theta(m/2 + m\xi) + 
4\sqrt{2} \bar \Phi\left(2\sqrt{m} (\xi - |\mu|)\right)
\right) & \ \text{if} \ \theta \geq 1/2,\\
\frac{2}{w^2} \left(
{\B}_\theta(m/2 - m\xi + 1) + 
4\sqrt{2} \bar \Phi\left(2\sqrt{m} (\xi - |\mu|)\right)
\right)  & \ \text{if} \  \theta < 1/2.
\end{cases}
\end{align}
\end{lemma}

\begin{proof}
We only consider the case when $\theta \geq 1/2$. The proof of the case $\theta < 1/2$ is similar and thus omitted. 
We split $m_2(\theta, w)$ into three parts as follows:
\begin{equation}
\begin{split}
\label{pf:m2-1}
m_2(\theta, w) 
& = \sum_{|u - m/2| \geq m\xi} \beta(u, w)^2 b_\theta(u) 
+ \sum_{m\nu < |u - m/2| < m\xi} \beta(u, w)^2 b_\theta(u)   \\
& \quad + \sum_{|u - m/2| \leq m\nu} \beta(u, w)^2 b_\theta(u) \\
& \quad = (a) + (b) + (c).
\end{split}
\end{equation}
As $\beta(u) \geq 1/w$ on $\{u: |u - m/2| \geq m\xi\}$,
$$
(a) \leq \frac{1}{w^2} \sum_{|u - m/2| \geq m\xi} b_\theta(u)
\leq \frac{2}{w^2} \bar \B_\theta(m/2 + m\xi).
$$
Next, since $\beta(u) \leq 0$ and $\beta(u)^2 \leq 1$, we immediately obtain 
$$
(c) 
\leq  \sum_{|u-m/2| \leq m\nu} \frac{\beta(u)^2 b_\theta(u)}{(1-w\beta(m/2))^2}
\leq \frac{1}{(1-w)^2} \sum_{|u-m/2| \leq m\nu} b_\theta(u) <  \frac{1}{(1-w)^2}.
$$
Last, using that $0 < w\beta(u) < 1$ for $\{u:|u - m/2| \in (m \nu, m\xi)\}$, 
$$
(b) \leq \frac{2}{w(1-w)^2}  \sum_{m\nu < |u-m/2| < m\xi} \beta(u) b_\theta(u)
< \frac{8}{w}  
\sum_{m\nu < |u - m/2| < m\xi} 
\frac{g(u)b_\theta(u)}{b(u)},
$$
which we used $1 - w > 1/2$ for a sufficiently large $n$.
To bound the upper bound in the last display, one can use the same argument as in the proof of Lemma \ref{lem:ub-m2}, then
the summation part in the last display can be bounded by \eqref{pf:m1-ub-4-1}, and thus leads to
\begin{align}
\label{pf:m2-2}
(b) \leq \frac{8}{w}
\sum_{\tilde u = m\nu}^{m\xi} g(u)
\frac{b_\theta(u)}{b(u)}
\leq 
\frac{8|\xi - \mu|}{w^2 \mu \sqrt{m(1-4\xi^2)}} 
\bar\Phi\left(2\sqrt{m}(\xi - \mu)\right).
\end{align}
If $\mu \geq \xi/K$ for some fixed $K > 0$, then 
\begin{align}
\label{pf:m2-3}
(b) \leq \frac{8\sqrt{2}K}{w^2 \sqrt{m}}\bar\Phi\left(2\sqrt{m}(\xi - \mu)\right),
\end{align}
which we used $1 - 4\xi^2 > 1/2$, as $\xi < 1/\sqrt{8}$ for a sufficiently large $m$.

If $\frac{1}{2m\xi} \leq \mu \leq \xi/K$, then the ratio
$\frac{|\xi - \mu|}{\mu}$ is at most $2m\xi^2$,
and then, \eqref{pf:m2-2} is bounded by 
\begin{align*}
\frac{16m\xi^2}{w^2\sqrt{m(1-4\xi^2)}} \bar\Phi\left({2\sqrt{m}(\xi - \mu)}\right)
& \leq \frac{8\sqrt{2}}{w^2} \bar\Phi\left({2\sqrt{m}(\xi - \mu)}\right),
\end{align*}
as $2m\xi^2 \sim \log (n/s_n) + \log(\sqrt{m}) \sim \log n \ll \sqrt{m}$.

Last, if $0 < \mu < \frac{1}{2m\xi}$, then from \eqref{pf:m1-ub-8} in the proof of Lemma \ref{lem:ub-m1}, we obtain 
$$
\sum_{m\nu \leq |u - m/2| \leq m\xi} \beta(u) b_\theta(u) \lesssim 2\zeta e^{\frac{2}{1-4\mu^2}} \to 0, \ \text{as} \ \zeta \to 0.
$$
By combining the three cases considered above and then combining the bounds we obtained for $(a)$, $(b)$ and $(c)$, the upper bound for $m_2(\theta, w)$ follows by using $w \to 0$ as $n \to \infty$.
\end{proof}

\begin{corollary}
\label{cor:m1-m2}
Consider $m_1(\theta, w)$ and $m_2(\theta, w)$ as in \eqref{m1} and \eqref{m2} respectively, 
let $\xi = |\xi_n|$ and $\nu = |\nu_n|$ for $\xi_n$ and $\nu_n$ given in \eqref{eqn:xi} and \eqref{eqn:nu} respectively, 
suppose 
$m\xi^4 = o(1)$, 
$w \asymp s_n/n$ and $s_n = n^{v_1}$ and $m = n^{v_2}$ for $v_1\in (0,1)$,
then for any $(0, 1) \ni \theta \neq 1/2$,
$$
m_2(\theta, w) \lesssim \frac{m_1(\theta, w)}{w}.
$$
\end{corollary}

\subsection{Controlling $m_1(\theta ,w)$ on the set containing relatively small signals}
\label{sec:control-small-signals}

Consider the following set:
\begin{align}
\label{J0}
\mathcal{J}_0 := \mathcal{J}(\theta_0, w, K) = \{1 \leq j \leq n: |\theta_{0,j} - 1/2| \geq \xi(w)/K\},
\end{align}
which is a subset of $\mathcal{S}_0 = \{1\leq j\leq n: \theta_{0,j} \neq 1/2\}$. 
Define the following two quantities:
$$
\mathcal{M}^{\mathcal{S}_0}(w) = \sum_{j\in \mathcal{S}_0} m_1(\theta_{0, j}, w),
\quad 
\mathcal{M}^{\mathcal{J}_0}(w, K) = \sum_{j\in \mathcal{J}_0} m_1(\theta_{0, j}, w).
$$
In the next lemma, we will bound the difference between the above two quantities. 
This bound is essential to obtain the uniform FDR control result for our multiple testing procedures given in Section \ref{sec:unif-fdr}. It is also used to derive the concentration bound for $\hat w$ in Section \ref{sec:bound-w-hat}.

\begin{lemma}
\label{lem:m1-cntl-small-signals}
Consider the set $\mathcal{J}_0$ given in \eqref{J0},
suppose $s_n \leq n^{v_1}$ and $m \geq n^{v_2}$ for some fixed $v_1 \in (0, 1)$ and $v_2 \geq \log \log n/\log n$, then there exists a constant $D > 0$ depending on $v_1, v_2$ and some fixed constants $\Lambda$ and a constant $K > \left(1 - \sqrt{\frac{v_1}{1 + v_1 + v_2/2}}\right)^{-1}$ such that for a sufficiently large $n$, we have
$$
\sup_{\theta_0 \in l_0[s_n]} 
\sup_{w \in [1/n, 1/\log n]} 
|\mathcal{M}^{\mathcal{S}_0}(w) - \mathcal{M}^{\mathcal{J}_0}(w, K)| 
\lesssim  
n^{1-D}.
$$
\end{lemma}

\begin{proof}
Allow us to slightly abuse the notation to introduce $\mathcal{J}_0^c$ such that $\mathcal{J}_0^c \cup \mathcal{J}_0 = \mathcal{S}_0$, where
$$
\mathcal{J}_0^c = \{1 \leq j \leq n: 0 < |\theta_{0,j} - 1/2| < \xi(w)/K\}.
$$
Let $\mu_{0,j} = \theta_{0,j} - 1/2$, one can write
\begin{align}
\sum_{j \in \mathcal{J}_0^c} m_1(\theta_{0, j}, w) 
& = 
\left\{
\sum_{0 < |\mu_{0,j}| < \frac{1}{2m\xi}}
+ 
\sum_{\frac{1}{2m\xi} \leq |\mu_{0,j}| \leq \frac{\Lambda}{\sqrt{2m}}}
+ 
\sum_{\frac{\Lambda}{\sqrt{2m}} < |\mu_{0,j}| \leq \frac{\xi(w)}{K}}
\right\} m_1(\theta_{0,j}, w) 
\label{pf:m1-cntl-small-signals-1}\\
& = (I) + (II) + (III).
\nonumber 
\end{align}

First, we bound $(I)$.
By Lemma \ref{lem:ub-m1},
using the fact that $|\mathcal{J}_0^c| \leq |\mathcal{S}_0| \leq s_n$ and $C_{\mu_{0,j}} \to 0$ for any $\mu_{0,j} < \frac{1}{2m\xi}$, as $C_{\mu_{0,j}} = \frac{4\mu_{0,j}^2}{1-4\mu_{0,j}^2} \to 0$, 
let $\tilde C = \max_j C_{\mu_{0,j}}$ and $C_1 = \max_{j} \exp(\frac{2}{1-4\mu_{0,j}^2})$,  
we thus obtain 
\begin{align*}
(I) \lesssim 2C_1 s_n \left(\zeta(w) + \frac{w^{\tilde C}}{\sqrt{m}}\right)
& \leq 2 C_1 \left(n^{v_1 - v_2/2} \sqrt{\log n} + n^{ v_1 - v_2/2 - \tilde C (\log\log n/\log n)}\right)\\
& \leq 4C_1 n^{v_1 - v_2/2 + \log \log n/(2\log n)},
\end{align*}
as $s_n = n^{v_1}$, $\zeta(w) = \sqrt{-\log w/(2m)}$, $1/n \leq w \leq 1/\log n$, and $m = n^{v_2}$.

Next, we bound $(II)$. By Lemma \ref{lem:ub-m1} again and using that $(2m\xi)^{-1} \leq \mu_{0,j} \leq \Lambda/\sqrt{2m}$ and $\sqrt{2m}\xi(w) \sim \sqrt{\log (1/w) + \log \sqrt{m}}$, we obtain
\begin{align*}
(II) & \lesssim s_n \max_{j} \left( \frac{1}{\sqrt{m |\xi(w) - \mu_{0,j}|}} + \frac{1}{(m+1)\mu_{0,j}} \right) e^{-2m(|\mu_{0,j}| - \xi)^2 + 2m\xi^2}\\
& \leq 4s_n \xi(w) e^{2\Lambda \sqrt{2m}\xi(w)} 
\leq \sqrt{2 \log n}C(v_1, v_2) n^{v_1 - v_2/2}  e^{\Lambda C(v_1, v_2) \sqrt{\log n}} \\
& = \sqrt{2} C(v_1, v_2) n^{v_1 - v_2/2 + \Lambda C(v_1, v_2)/\sqrt{\log n} - \log\log n/(2\log n)} \\
& \leq \sqrt{2}C(v_1,v_2) n^{1-(1-v_1 + v_2/2 -\Lambda C(v_1,v_2)/\sqrt{\log n})},
\end{align*}
where $C(v_1, v_2) = 2\sqrt{1 - v_1 + v_2/2}$.
As $n \to \infty$, $\Lambda C(v_1,v_2)/\sqrt{\log n} \to 0$.

Last, 
for $\Lambda/\sqrt{2m} < \mu_{0,j} \leq \xi(w)/K$, 
$T_m(\mu_{0,j}, m) \leq (1 - K^{-1}) \sqrt{2m} \xi(w)/\Lambda$.
Using that $s_n /w \leq n^{1+v_1}$, we obtain
\begin{align*}
(III) 
& \leq 2n \max_j \left(\bar \B_{\theta_0} (m/2 + m\xi(w)) + 
\bar \Phi \left(2\sqrt{m}(\xi(w) - |\mu_{0,j}|) \right) T_m(\mu_{0,j}, m) \right)\\
& \leq 2n \max_j \left(
e^{-mT(1/2 + \xi, 1/2 + |\mu_{0,j}|)} +
2 \Lambda^{-1} (1-K^{-1}) \sqrt{2m} \xi(w) e^{- 2m(\xi - |\mu_{0,j}|)^2 }\right)\\
& \leq 4 C(\Lambda, K, v_1, v_2) n^{1+v_1}  \sqrt{\log n}  e^{-2m(\xi - |\mu_{0,j}|)^2}, 
\end{align*}
as $mT(1/2 + \xi, 1/2 + |\mu_{0,j}|) \leq 2m(\xi - |\mu_{0,j}|)^2 + 6m \xi^4 \to 2m(\xi - |\mu_{0,j}|)^2$ by (3) in Lemma \ref{lem:entropy-bound}, as $m \xi^4 \to 0$ by assumption. 
Let $C_2 = C(\Lambda, K, v_1, v_2)$, then the upper bound in the last display can be bounded by 
$$
4 C_2 n^{1+v_1} \sqrt{\log n} e^{- 2m(1-1/K)^2 \xi^2} \leq 4C_2 n^{(1 + v_1)(1- (1-1/K)^2) + \log\log n/(2\log n) - v_2(1-1/K)^2/2}.
$$

Now we combine the above upper bounds for $(I)$, $(II)$, and $(III)$, then for a sufficiently large $n$,
\begin{align*}
\eqref{pf:m1-cntl-small-signals-1}
\lesssim 
n^{1-(1-v_1 + v_2/2)} + n^{1 + v_1 - (1+v_1+v_2/2)(1-1/K)^2}.
\end{align*}
Choosing $D = \min\{1 - v_1 +v_2/2, (1+v_1+v_2/2)(1-1/K)^2 - v_1\}$, 
if $K > \left(1 - \sqrt{\frac{v_1}{1 + v_1 + v_2/2}}\right)^{-1}$, then $D > 0$, providing that $v_1$ is bounded away from $1$ (which is true as we assume $w \leq 1/\log n$). Thus, we obtain $\eqref{pf:m1-cntl-small-signals-1} \lesssim n^{1-D}$.
\end{proof}

\begin{lemma}
\label{lem:m1-ratio}
Consider the set $\mathcal{J}_0$ given in \eqref{J0},
suppose $s_n \leq n^{v_1}$ and $m \geq n^{v_2}$ for some fixed $v_1 \in (0, 1)$ and $v_2 \geq \log\log n/\log n$, 
then for a sufficiently large $K >  A > 1$ and any $w \in [n^{-1}, (\log n)^{-1}]$, if $n$ is sufficiently large, then
$$
\mathcal{M}^{\mathcal{J}_0}(w/A, K) \geq K\mathcal{M}^{\mathcal{J}_0}(w, K)
$$
\end{lemma}

\begin{proof}
Recall that $ \mathcal{M}^{\mathcal{J}_0}(w, K)  = \sum_{j \in \mathcal{J}_0} m_1(\theta_{0,j}, w)$.
Using the lower bound of $m_1(\cdot, w)$ in Lemma \ref{lem:lb-m1}, we obtain
\begin{align}
& \mathcal{M}^{\mathcal{J}_0}(w/A, K) 
\nonumber \\
& \quad \geq 
\frac{A}{w} 
\sum_{j \in \mathcal{J}_0} 
\left(
\bar \B_{\theta_{0,j}} (m/2 + m\xi(w/A))
+ 
\bar \Phi \left( \frac{2\sqrt{m}(\xi(w/A) - |\mu_{0,j}|)}{\sqrt{1-4\mu^2_{0,j}}} \right) \
T_m'(\mu_{0,j}, w/A) \right) 
\nonumber \\
& \quad = (a) + (b).
\label{pf:m1-ratio}
\end{align}
We shall derive a lower bound for $(a)$ and $(b)$ respectively. 
For $(a)$, using that $\xi(w) < \xi(w/A)$ as long as $A > 1$ and $|\theta_{0,j} - 1/2| = |\mu_{0,j}| \geq \xi(w)/K$ for each $j \in \mathcal{J}_0$, we have
\begin{align}
\label{pf:m1-ratio-1}
\bar \B_{\theta_{0,j}} (m/2 + m\xi(w/A))
= \bar \B_{\theta_{0,j}} (m/2 + m\xi(w)) - 
\sum_{|\tilde u| = m\xi(w)}^{m\xi(w/A)} b_{\theta_{0,j}}(m/2 + |\tilde u|).
\end{align}
By plugging-in the expression of the density function of a binomial distribution, the second term in the last display can be written as 
$$
\sum_{|\tilde u| = m\xi(w)}^{m\xi(w/A)} b_{\theta_{0,j}}(m/2 + |\tilde u|)
= 
\sum_{|\tilde u| = m\xi(w)}^{m\xi(w/A)}
{m \choose {m/2 + |\tilde u|}} 
\theta_{0,j}^{m/2 + |\tilde u|} (1-\theta_{0,j})^{m/2 - |\tilde u|}.
$$
By Lemma \ref{lem:bin-coef-bound}, the last display equals to
$$
\sum_{|\tilde u| = m\xi(w)}^{m\xi(w/A)}
\frac{\sqrt{2}}{\sqrt{\pi m(1-4(\tilde u/m)^2)}} 
e^{- mT(1/2 + |\tilde u|/m, 1/2) + mT(1/2 + |\tilde u|/m, 1/2 + |\mu_{0,j}|) + o(1)}. 
$$
By Lemma \ref{lem:entropy-bound} and using $|\mu_{0,j}| \geq \xi(w)/K$,
the last display is bounded by
$$
\sum_{|\tilde u| = m\xi(w)}^{m\xi(w/A)}
\frac{\sqrt{2} e^{- 2 m (|\tilde u|/m)^2 + Cm (\tilde u/m)^4}}{\sqrt{\pi m(1-4(\tilde u/m)^2)}} 
\leq
\frac{\sqrt{2}e^{-2m \xi^2(w) - C m\xi^4(w)}}{\sqrt{\pi m(1-4 \xi^2(w))}}
\leq \frac{C' e^{-2m\xi^2(w)}}{\sqrt{m}} ,
$$
where $C > 16/3$ is a fixed constant and $C' = \frac{\sqrt{2} e^{-o(1)}}{\sqrt{\pi m (1-4\xi^2(w/A))}}$ as $m \xi^4(w) = o(1)$ by assumption.
Therefore, 
$$
(a) \geq \frac{1}{w}\sum_{j \in \mathcal{J}_0} \bar \B_{\theta_{0,j}} (m/2 + m\xi(w/A)) -
\frac{e^{-2m\xi^2(w)}}{\sqrt{m} w}.
$$
Since $2m\xi^2(w) \sim \log (1/w) + \log (m)/2$, the second term in the last display is $O(1/m) = o(1)$.

Next, we derive a lower bound for $(b)$ in \eqref{pf:m1-ratio}. 
By Lemma \ref{lem:lb-T}, as $A > 1$, we have $T_m(\mu_{0,j}, w/A) \geq \frac{1}{2}T_m(\mu_{0,j}, w)$, and thus 
$$
T_m'(\mu_{0,j}, w/A) = \frac{C_{\mu_{0,j}}}{\mu_{0,j}} \sqrt{\frac{1-4\mu_{0,j}^2}{1-4\xi^2}} T_m(\mu_{0,j}, w/A)
\geq \frac{C_{\mu_{0,j}}}{2\mu_{0,j}} \sqrt{\frac{1-4\mu_{0,j}^2}{1-4\xi^2}} T_m(\mu_{0,j}, w),
$$
where $C_{\mu} \leq \exp(\frac{4\mu}{1-4\mu^2}) - 1$.
For a sufficiently large $m$, $1-4\xi^2 \leq 3/4$, let $K_{\mu_{0,j}} = \frac{C_{\mu_{0,j}}}{\mu_{0,j}} \sqrt{\frac{1-4\mu_{0,j}^2}3}$, then the last display implies
$T_m'(\mu_{0,j}, w/A) \geq K_{\mu_{0,j}} T_m(\mu_{0,j}, w)$.
In addition, by Lemma \ref{lem:lb-Hmu}, 
$H_{\mu_{0,j}}(w/A)  \geq A^{1/(4K_0)} H_{\mu_{0,j}}(w)$,
where
$H_\mu = \frac{1}{w} \bar \Phi \left(\frac{2\sqrt{m}(\xi(w) - |\mu|}{\sqrt{1-4\mu^2}} \right)$ 
and a fixed $K_0 \geq \sqrt{2}/{4}$.
Therefore, 
\begin{align*}
(b) & \geq \frac{A^{1/(4K_0)}}{w} 
\sum_{j \in \mathcal{J}_0}  K_{\mu_{0,j}}
\bar \Phi \left( \frac{2\sqrt{m}(\xi(w) - |\mu_{0,j}|)}{\sqrt{1-4\mu^2_{0,j}}} \right) 
T_m(\mu_{0,j}, w)\\
& \geq \frac{A^{1/(4K_0)}}{2w}\underline K_{\mu}
\sum_{j \in \mathcal{J}_0} 
\bar \Phi \left( \frac{2\sqrt{m}(\xi(w) - |\mu_{0,j}|)}{\sqrt{1-4\mu^2_{0,j}}} \right),
\end{align*}
where $\underline K_{\mu}  = \min_j K_{\mu_{0,j}}$.
By combining the lower bounds of $(a)$ and $(b)$,
the result follows by letting $K = A \vee \frac{\underline K_{\mu}A^{1+1/(4K_0)}}{2}$.
\end{proof}

\begin{lemma} 
\label{lem:lb-T}
Consider $T_m(\mu, \xi(w))$ in \eqref{eqn:Tm}, for any $w \in (0, 1)$ and $\mu_0 > \mu \geq \xi(w)/K_0$  $\mu_0 < 1/2$, there exists a $w_0 = w_0(K_0, z)$ such that for all $w \leq w_0$, $z > 1$, and $\mu \geq \xi(w)/K_0$, we have
\begin{align*}
& T_m(\mu, \xi(w/z)) \geq \frac{T_m(\mu, \xi(w))}{2}.
\end{align*}
\end{lemma}

\begin{proof}
By the definition of $T_m(\mu, \xi(w))$, we have
\begin{align*}
\frac{|\xi(w/z) - \mu|}{\mu \sqrt{m(1-4\xi^2(w/z))}}
& \geq 
\frac{|\xi(w) - \mu|}{\mu \sqrt{m(1-4\xi^2(w/z))}}
- 
\frac{|\xi(w) - \xi(w/z)|}{\mu \sqrt{m(1-4\xi^2(w/z))}}\\
& \geq 
\frac{|\xi(w) - \mu|}{\mu \sqrt{m(1-4\xi^2(w))}}
-
\frac{|\xi(w) - \xi(w/z)|}{\mu \sqrt{m(1-4\xi^2(w/z))}},
\end{align*}   
as $\xi(w/z) \geq \xi(w)$ for any $z > 1$.
Since $\xi(u) \sim \sqrt{\frac{1}{2m} (\log u^{-1} + \log \sqrt{m})}$, 
$\xi(w/z) \sim \sqrt{\frac{\log z}{2m} + \xi^2(w)}$.
Using that $\mu \geq \xi(w)/K_0$,
the second term in the last line is bounded by 
$$
\frac{K_0 \left(\sqrt{\frac{\log z}{2m} + \xi^2(w)} - \xi(w)\right)}{\xi(w) \sqrt{m(1-4\xi^2(w))}}
\leq 
\frac{K_0 \left(\sqrt{\frac{\log z}{\log (1/w)}} \right)}{\sqrt{m(1-4\xi^2(w))}} \to 0, \ \text{as} \ m \to \infty.
$$
Thus, for a sufficiently large $m$,
$\displaystyle
T_m(\mu, \xi(w/z)) \geq  \frac{|\xi(w) - \mu|}{2\mu \sqrt{m(1-4\xi^2(w))}} 
= \frac{T_m(\mu, \xi(w))}{2}.
$
\end{proof}

\begin{lemma} 
\label{lem:lb-Hmu}
Consider the function 
\begin{align}
\label{eqn:Hmu}
H_\mu(w) = \frac{1}{w} \bar \Phi \left(\frac{2\sqrt{m}(\xi(w) - \mu)}{\sqrt{1-4\mu^2}} \right),
\end{align}
for any $w \in (0, 1)$ and $1/2 > \mu_0 > \mu \geq \xi(w)/K_0$, then there exists a $w_0 = w_0(K_0, z)$,
$K_0 \geq \sqrt{2}/4$ and $z > 1$, such that for any $w \leq w_0$, 
\begin{align*}
H_\mu(w/z) \geq  z^{1/(4K_0)} H_\mu(w).
\end{align*}
\end{lemma}
\begin{proof}
The proof is inspired by the proof of Lemma 19 of \citetalias{cast20MT} with substantial modifications are made due to $\xi(w)$ here is different from it in their model. 
Let $\Upsilon(u) = \log H_\mu(e^{-\mu})$, then the goal is to show the following inequality:
$$
\Upsilon(\log (z/w)) - \Upsilon(\log (1/w)) \geq \frac{1}{2K_0} 
\left(\log( z/w) - \log (1/w)\right).
$$
By the mean-value theorem, it is then sufficient to show that $\Upsilon'(u) \geq 1/(2K_0)$ for any $u \in [\log 1/w, \log z/w]$.
Note that 
$$
\xi'(w) = - \frac{1}{mw^2 \beta'(m/2 + m\xi(w))}.
$$
Thus, we have
\begin{align}
\label{pf:Hmu-1}
\Upsilon'(u) = 
1 - \frac{2 e^u }{\beta'(m/2 + m\xi(e^{-u})) \sqrt{m(1-4\mu^2)}} 
\frac{\phi}{\bar \Phi}\left(\frac{2\sqrt{m} (\xi(e^{-u}) - \mu)}{\sqrt{1-4\mu^2}}\right).
\end{align}
In addition,
$$
\beta'(x) = (\beta(x) + 1) \left( \Psi(x+1) - \Psi(m-x+1) \right) 
:= (\beta(x) + 1)Q(x),
$$
where $\Psi(x+1) = \frac{d}{dx} \log \Gamma(x+1)$, $\Gamma(\cdot)$ is the gamma function,
and using that $\beta(m/2 + m\xi(w)) = 1/w$ by Lemma \ref{lem:beta}, we thus have
$$
\beta'(m/2 + m\xi(e^{-u})) = (\beta(m/2 + m\xi(e^{-u})) + 1)Q(m/2 + m\xi(e^{-u}))
= Q(m/2 + m\xi(e^{-u})) (e^{-u} + 1).
$$
By plugging-in the above expression into \eqref{pf:Hmu-1} and 
let $C_m(\mu) = \frac{2}{\sqrt{m(1-4\mu^2)}}$, one obtains
\begin{align}
\label{pf:Hmu-2}
\Upsilon'(u) = 1 - \frac{C_m(\mu) e^u}{(1 + e^{u})Q(m/2 + m\xi(e^{-u}))} 
\frac{\phi}{\bar \Phi}\left(\frac{2\sqrt{m} (\xi(e^{-u}) - \mu)}{\sqrt{1-4\mu^2}}\right).
\end{align}
One needs to further bound the function $Q(\cdot)$.
Using the mean-value theorem again, then there exist $\xi^\star \in [-\xi, \xi]$ such that 
$$
Q(m/2 + m\xi(x)) = 
\Psi(m/2 - m\xi(x) + 1) - \Psi(m/2 + m\xi(x) +1) = 2m\xi(x) \Psi'(m/2 + m\xi^\star(x)+ 1).
$$
By Stirling's approximation, $\Gamma(x + 1) \sim \sqrt{2\pi} e^{(x + 1/2)\log x - x}$ for a sufficiently large $x$. We thus have
$$
\Psi(x+1) \sim \log x + \frac{1}{2x}, \ \text{and} \ \Psi'(x+1) \sim \frac{1}{x} - \frac{1}{2x^2}.
$$
Therefore, there exists a sufficiently large $u$ such that
$Q(m/2 + m\xi(e^{-u})) \sim 4\xi(e^{-u})$.
By plugging this bound into \eqref{pf:Hmu-2}, we arrive at
\begin{align}
\label{pf:lb-Hmu-1}
\Upsilon'(u) = 1 - \frac{C_m(\mu) e^u}{4\xi(e^{-u})(1 + e^{u})} 
\frac{\phi}{\bar \Phi}\left(\frac{2\sqrt{m} (\xi(e^{-u}) - \mu)}{\sqrt{1-4\mu^2}}\right).
\end{align}
Since the map $u \to e^u(1+e^u)^{-1}$ has limit $1$ as $u, m \to \infty$,
for large enough $u, m$, $e^u(1+e^u)^{-1} \leq 1 + \epsilon$ for some $\epsilon > 0$ to be specify later. 
Applying the lower bound in Lemma \ref{lem:gaussian-bound}, 
if $\mu < \xi(e^{-u}) - 1$, then
\begin{align}
\frac{C_m(\mu)}{4\xi(e^{-u})} 
\frac{\phi}{\bar \Phi}\left(\frac{2\sqrt{m} (\xi(e^{-u}) - \mu)}{\sqrt{1-4\mu^2}}\right)
& \leq \frac{C_m(\mu)}{4\xi(e^{-u})} 
\frac{
1 + \frac{4 m (\xi(e^{-u}) - \mu)^2}{1-4\mu^2}
}{
\frac{2\sqrt{m} (\xi(e^{-u}) - \mu)}{\sqrt{1-4\mu^2}}
} 
\nonumber \\
& = 
\frac{1}{4m\xi(e^{-u})(\xi(e^{-u}) - \mu)} 
+ \frac{\xi(e^{-u}) - \mu}{\xi(e^{-u})(1-4\mu^2)}.
\nonumber
\end{align}
The first term in the last display $\to 0$ as $m\to \infty$. 
By the assumption $\xi(w) -1 > \mu \geq \xi(w)/K_0$ for a sufficiently large $K_0$, the second term in the last display is bounded by 
$$
\frac{\xi(e^{-u}) - \mu}{(1-4\xi^2(e^{-u}))\xi(e^{-u})}
\leq \frac{1 + \epsilon}{1-4\xi^2(w)}\left(1 - \frac{1}{K_0}\right).
$$
When $\mu \geq \xi(e^{-u}) - 1$, 
\begin{align}
\label{pf:lb-Hmu-3}
\frac{C_m(\mu)}{4\xi(e^{-u})} 
\frac{\phi}{\bar \Phi}\left(\frac{2\sqrt{m} (\xi(e^{-u}) - \mu)}{\sqrt{1-4\mu^2}}\right)
\leq 
\frac{\phi(0)}{2\sqrt{m} \xi(e^{-u}) \sqrt{1-4\mu^2}} \left(\bar \Phi \left(\frac{2\sqrt{m}}{\sqrt{1-4\mu^2}}\right)\right)^{-1}.
\end{align}
By Lemma \ref{lem:gaussian-bound} again, 
$$
\bar \Phi \left(\frac{2\sqrt{m}}{\sqrt{1-4\mu^2}}\right)
\geq \frac{2\sqrt{m}}{\sqrt{1-4\mu^2}} \left(1 + \frac{4m}{1-4\mu^2}\right)^{-1} 
\phi\left(
\frac{2\sqrt{m}}{\sqrt{1-4\mu^2}}
\right),
$$
and then,
$$
\eqref{pf:lb-Hmu-3} \geq 
\frac{\phi(0)}{2\sqrt{2\pi} m\xi(w)} \exp\left(-\frac{2m}{1-4\mu^2}\right) \to 0,
\ \text{as} \ m \to \infty,
$$
providing that $\mu \leq \mu_0$ is bounded away from $1/2$.
Now we combine the upper bound for each case (either $\mu \geq \xi(e^{-u}) - 1$ or $\mu \leq \xi(e^{-u}) - 1$) to obtain
\begin{align*}
1 - \Upsilon'(u) & \leq  (1+\epsilon)\left(1 + \frac{4\xi^2(w)}{1-4\xi^2(w)} - \frac{1}{K_0}\right).
\end{align*}
For a sufficiently large $m$, since $\xi(w) \to 0$, 
if choosing 
$\epsilon^{-1} = 4K_0- 2$ for an $K_0 \geq \sqrt{2}/4$, then
$$
1 - \Upsilon'(u) \leq  (1+\epsilon)\left(1 - \frac{1}{2K_0}\right)
= 1- \frac{1}{4K_0},
$$
which implies $\Upsilon'(u)  \geq \frac{1}{4K_0}$. The proof is thus completed. 
\end{proof}

\section{Analyzing $\beta(u)$}
\label{sec:beta}
\begin{lemma}
\label{lem:monotone-beta}
$\beta(x) = (g/b)(x) - 1$ is non-decreasing on $x \in [m/2, m]$
and non-increasing on $x\in [0, m/2)$.
\end{lemma}
\begin{proof}
By plugging-in the expressions of $g(x)$ and $b(x)$, we obtain 
$$
d \beta(x) /du = d(g/b)(x)/dx = \frac{2^m}{m+1} \cdot \frac{d {m \choose x}}{dx}
=
\frac{2^m d (\Gamma(x+1) \Gamma(m-x+1)) / dx}{(m+1) \Gamma(m+1)}.
$$
We need to show $d \beta(x) / dx \geq 0$. 
By calculation, 
$$
\frac{d\Gamma(x+1) \Gamma(m-x+1)}{dx} 
= \Gamma(x+1) \Gamma(m-x+1) \left[ \Gamma'(x+1) - \Gamma'(m-x+1) \right],
$$
where $\Gamma'(\cdot)$ is the first derivative of $\Gamma(\cdot)$.
The last display is non-negative for any $x \in [m/2, m]$ as $\Gamma'(x+1) \geq \Gamma'(m-x+1)$ for $x \in [m/2, m]$. Note that
$\Gamma'(x+1)$ is a monotone increasing function, $\Gamma'(m - x + 1)$ is a monotone decreasing function, and $\Gamma'(x+1) = \Gamma'(m-x+1)$ if and only if $x = m/2$.
Thus, $d(g/b)(x)/dx \geq 0$, which implies that $\beta(x)$ is non-decreasing for any $x \in [m/2, m]$. 
In particular, when $x > m/2$, $\beta(x)$ is a strictly increasing function.
Since $b(x)$ is symmetric at $m/2$, we immediately have that $\beta(x)$ is non-increasing on $[0, m/2)$.
\end{proof}

\begin{lemma}
\label{lem:beta}
Define $\beta(u) = (g/b)(u) - 1$, let $\xi$ be the solution of $\beta(m/2 + m\xi) = 1/w$ for $w \in (0, 1)$, 
then there exists a fixed $\xi_\circ \in (0, 1/2)$ such that for any $|\xi| \leq \xi_\circ$,
\begin{align}
\label{ub-xi}
2m\xi^2 & \leq \log (1 + 1/w) + \log \left( \frac{\sqrt{2}(1+m)}{\sqrt{\pi m(1-4{\xi_\circ}^2)}} \right) + \frac{1}{12m},\\
2m\xi^2 & \geq \left[\log (1+1/w) + \log \left(\frac{\sqrt{2}(1+m)}{\sqrt{\pi m}}\right)\right]\left(1 + \frac{8\xi_\circ^2}{3(1-4\xi_\circ^2)^2}\right)^{-1}.
\label{lb-xi}
\end{align}
If $\xi = \xi_n$ and $m\xi_n^4 \to 0$, as $m\to \infty$, then
\begin{align}
\label{eqn:xi}
|\xi_n| \sim \sqrt{\frac{\log(1+1/w) + \log \left(\frac{\sqrt{2} (m + 1)}{\sqrt{\pi m}}\right)}{2m}}.
\end{align}
\end{lemma}

\begin{proof}
By the definition of $\beta(\cdot)$, $\beta(m/2 + m\xi)) = 1/w$ implies
$(g/b)(m/2 + m\xi) = 1 + 1/w$,
By plugging-in $g(x) = (1+m)^{-1}$ and $b(x) = \text{Bin}(m, 1/2)$ and then taking the logarithm of both sides, we obtain
\begin{align}
\label{pf-beta-1}
-\log (1+1/w) = \log (1+m) + \log {m \choose {m/2 + m\xi}} - 2\log m.
\end{align}
Without loss of generality, we assume $0 \leq \xi < 1/2$, as the binomial coefficient is symmetric at $m/2$, then by Lemma \ref{lem:bin-coef-bound}, we have
\begin{align}
\label{pf-beta-2}
{\frac{\sqrt{2} e^{-m T(1/2 + \xi, 1/2) + 2\log m} }{\sqrt{\pi m(1- 4 \xi^2)}}}
\leq 
{m \choose m/2 + m\xi}
\leq 
\frac{\sqrt{2} e^{-m T(1/2 + \xi, 1/2) + 2\log m + \omega(\xi)}}{\sqrt{\pi m(1- 4 \xi^2)}}, 
\end{align}
where $\omega(\xi) \leq (12m)^{-1}$.
By $(2)$ in Lemma \ref{lem:entropy-bound}, one can further bound
\begin{align}
\label{pf-beta-3}
2 \xi^2 \leq T(1/2 + \xi, 1/2) \leq 2 \xi^2 + \frac{16 \xi^4}{3(1- 4\xi^2)^2}.
\end{align}
By combining \eqref{pf-beta-2} and \eqref{pf-beta-3} and using that $\xi \leq \xi_\circ$, \eqref{pf-beta-1} implies 
$$
2m\xi^2 \leq \log (1+1/w) + \log \left(\frac{\sqrt{2}(1+m)}{\sqrt{\pi m (1-4\xi_\circ^2)}}\right) + \frac{1}{12m}.
$$
and
$$
2m\xi^2\left(1 + \frac{8\xi^2}{3(1-4\xi^2)^2}\right) \geq \log (1+1/w) + \log \left(\frac{\sqrt{2}(1+m)}{\sqrt{\pi m}}\right).
$$ 
Since $\frac{\xi^2}{(1-4\xi^2)^2} \leq \frac{\xi_\circ^2}{(1-4\xi_\circ^2)}$, we obtain the lower bound.

If $m\xi_n^4 = o(1)$, then the second term in the upper bound of \eqref{pf-beta-3} is $o(1)$, which implies $T(1/2 + \xi_n, 1/2) \sim 2\xi_n^2$. 
Also, we have $1-4\xi_n^2 \sim 1$. Therefore,
$$
2m\xi_n^2 \sim \log(1+1/w) + \log \left( \frac{\sqrt{2}(1+m)}{\sqrt{\pi m}} \right).
$$
We thus proved \eqref{eqn:xi}.
\end{proof}

\begin{lemma}
\label{lem:lb-beta}
Let $\beta(u) = (g/b)(u) - 1$ where $g(u) = 1/(m+1)$ and $b(u) = \text{Bin}(u; m, 1/2)$, then 
$$
\frac{\sqrt{\pi m}}{\sqrt{2}(m+1)(1+(12m)^{-1})} - 1 
\leq \beta\left(\frac{m}{2}\right) 
\leq 
\frac{\sqrt{\pi m}}{\sqrt{2}(m+1)} - 1.
$$
If $m \to \infty$, then $\beta(m/2) + 1\sim \frac{\sqrt{\pi m}}{\sqrt{2}(m+1)} $.
\end{lemma}
\begin{proof}
Since $b(m/2) = {m \choose m/2} 2^{-m}$, Lemmas \ref{lem:bin-coef-bound} and \ref{lem:entropy-bound} imply
$\sqrt{2/(\pi m)} b(m/2) \leq \sqrt{2/(\pi m)}e^{\omega(0)}$, where $\omega(s)$ is given in Lemma \ref{lem:bin-coef-bound}. In fact, by examing the proof of Lemma \ref{lem:bin-coef-bound}, $\omega(0) \leq (12m)^{-1}$, which implies $e^{\omega(0)} \leq (1 + (12m)^{-1}) \to 1$ as $m \to \infty$.
Last, using both the upper and lower bounds for $b(m/2)$, we obtain the bounds for $\beta(m/2)$.
\end{proof}

The sharp bounds we have derived for the binomial coefficient in Lemma \ref{lem:bin-coef-bound} enable us to establish tight bounds for the solution when $\beta(x) = 0$. These bounds are presented in the following lemma.

\begin{lemma}
\label{lem:solve-beta}
Define $\beta(x) = (g/b)(x) - 1$, let $\nu_n$ be the solution of $\beta(m/2 + m\nu_n) = 0$, then
there exists a fixed $\nu_\circ \in (0, 1/2)$ such that for $|\nu_n| < \nu_\circ$,
\begin{align}
\label{ub-nu}
2m\nu_n^2 & \leq  \log \left( \frac{\sqrt{2}(1+m)}{\sqrt{\pi m(1-4{\nu_\circ}^2)}} \right) + \frac{1}{12m},\\
2m\nu_n^2 & \geq \left[\log \left(\frac{\sqrt{2}(1+m)}{\sqrt{\pi m}}\right)\right]\left(1 + \frac{8\nu_\circ^2}{3(1-4\nu_\circ^2)^2}\right)^{-1} .
\label{lb-nu}
\end{align}
In particular, if $m\nu_n^4 \to 0$ as $m \to \infty$, then
\begin{align}
\label{eqn:nu}
|\nu_n| \sim \sqrt{\frac{1}{2m} \log \left(\frac{\sqrt{2} (m + 1)}{\sqrt{\pi m}}\right)}.
\end{align}
\end{lemma}
\begin{proof}
The proof is essentially the same as that of Lemma \ref{lem:beta}. One only needs to replace $-\log (1+1/w)$, $\xi_n$, and $\xi_\circ$ in the proof of Lemma \ref{lem:beta} with $0$, $\nu_n$, and $\nu_\circ$ respectively.
\end{proof}

\begin{figure}[!h]
\centering
    \subfigure[]{\includegraphics[width=0.32\textwidth]{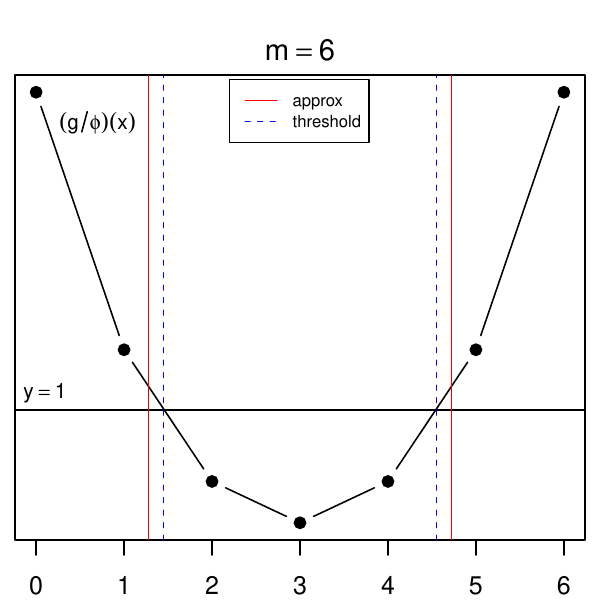}} 
    \subfigure[]{\includegraphics[width=0.32\textwidth]{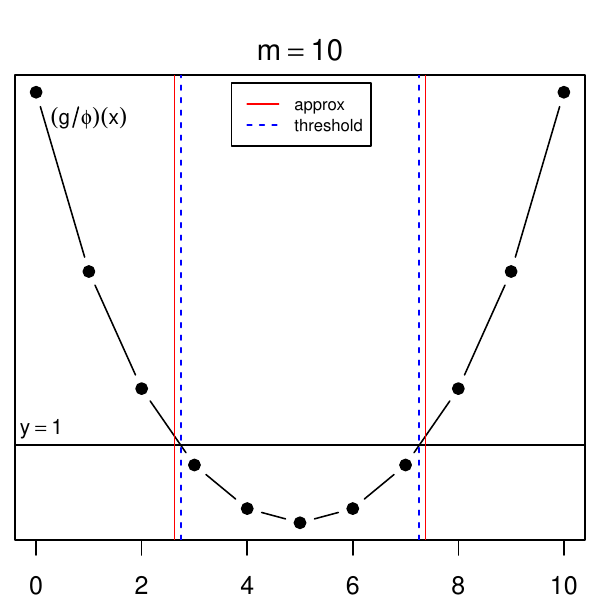}} 
    \subfigure[]{\includegraphics[width=0.32\textwidth]{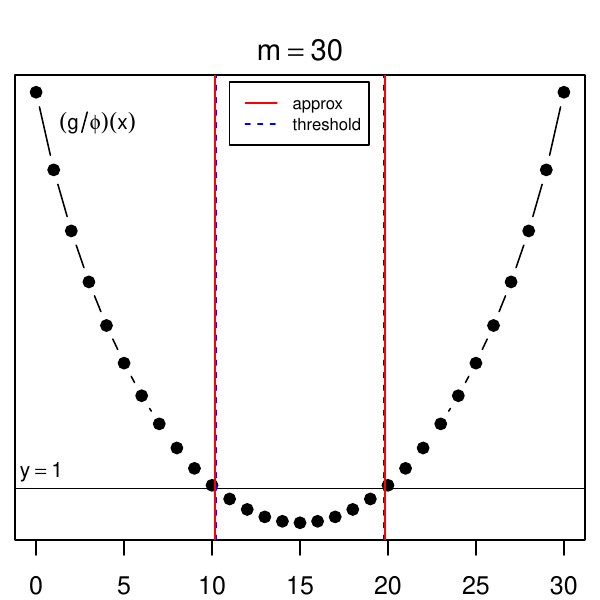}}
	\caption{Plot of the function $(g/b)(x)$ when (a) $m = 6$, (b) $m = 10$, and (c) $m = 30$. The blue dashed line indicates the exact value of $x$ when $(g/b)(x) = 1$ and the red solid line its approximated value, $\nu_n$, in Lemma \ref{lem:solve-beta}. }
	\label{fig:beta}
\end{figure}

In Figure \ref{fig:beta}, we plot the relation between $(g/\phi)(x)$ and $1$. 
Linear interpolation between points are used as the binomial distribution is discrete. The blue dash line indicates the threshold, which is the intersection between the two functions, and the red solid line represents its approximated value using \eqref{eqn:nu}. 
We choose three different values for $m$,  $6, 10$, and $30$, and plot their results in (a)--(c) accordingly.
We found that $\nu_n$ is already close to the threshold when $m = 6$. As $m$ increases, the two functions become closer. When $m = 30$, they almost overlap.

\begin{lemma}
\label{lem:relation-xi-zeta-nu}
Define $\beta(x) = (g/b)(x) - 1$, let $\nu_n$ be the solution of $\beta(m/2 + m\nu_n) = 0$ and $\xi_n(w)$ be the solution of $\beta(m/2+m\xi_n) = 1/w$, for $\zeta_n(w)$ given in \eqref{eqn:zeta}, if $m\xi_n^4 \to 0$ as $m \to \infty$, then for any $w \in (0, 1)$,
$$
\xi_n^2(w) \sim \nu_n^2 + \zeta_n^2(w).
$$
\end{lemma}

\begin{proof}
The result follows by directly plugging-in the bound we obtained for $\xi_n(w)$ and $\nu_n(w)$ in Lemmas \ref{lem:beta} and \ref{lem:solve-beta} respectively.
\end{proof}

\begin{lemma}
\label{lem:bound-beta}
Given $\beta(x) = (g/b)(x) - 1$ for any $x \in [m/2 - m\nu_n, m/2 + m\nu_n]$, where
$\nu_n$ is given in Lemma \ref{lem:solve-beta}, then for a sufficiently large $m$, $-1 < \beta(x) \leq 0$.
\end{lemma}
\begin{proof}
By Lemma \ref{lem:monotone-beta}, $\beta(x)$ is a monotone increasing function on $[m/2 , m]$.
Also, by Lemma \ref{lem:solve-beta}, $\beta(m/2 + m\nu_n) = 0$. 
Therefore, for any $x \in [m/2 - m\nu_n, m/2 + m\nu_n]$, $\beta(x) \leq 0$, which gives the upper bound.
Also, $\beta(x) \geq \beta(m/2)$. 
From Lemma \ref{lem:bin-coef-bound}, we have $\sqrt{\frac{2}{\pi m}} \leq \phi(m/2) \leq \sqrt{\frac{2}{\pi m}}e^{\omega(\nu_n)}$, where $\omega(\nu_n) \to 0$ as $m \to \infty$ by Lemma \ref{lem:bin-coef-bound}. 
Therefore,
\begin{align*}
0 < \frac{\sqrt{\pi m}}{(m+1)\sqrt{2}} e^{- \omega(\nu_n)} \leq 
\beta(m/2) + 1 \leq \frac{\sqrt{\pi m}}{(m+1)\sqrt{2}},
\end{align*}
which implies $\beta(m/2) > -1$.
\end{proof}


\section{Proof of Proposition \ref{prop:lowerbound}}
\label{sec:pf-lowerbound}

Recall that the multiple testing risk is given by 
\begin{align}
\label{pf:risk-lb-0}
\mathfrak{R}(\theta_0, b) 
= \mathbb{E}_0(\FDP(\theta_0, \T)+ \FNP(\theta_0, \T))
= \int P_{\theta_0}(\FDP(\theta_0, \T)+ \FNP(\theta_0, \T)\geq q) dq.
\end{align}
The proof consists two steps: first, we show that for any $\epsilon \in (0, 1)$ and $s_n > 0$, 
there exists $d_{n, \epsilon}$, where
$m d_{n,\epsilon}^2 = \bar{\mathbf{B}}^{-1}((1+\epsilon^{-1})s_n/(n-s_n)) - \bar{\mathbf{B}}^{-1}(\epsilon/4)$, such that
\begin{align}
\label{pf:risk-lb-1}
\sup_{\T \in \mathcal{T}} \sup_{\theta_0 \in \Theta_0^{-}[s_n; d_{n, \epsilon}]} 
(P_{\theta_0} (\FDP(\theta_0, \T)+ \FNP(\theta_0, \T)\leq 1-\epsilon) )
\leq 3e^{-s_n\epsilon/6},
\end{align}
where, for any $d_n \geq 0$ (possibly $d_n \to 0$ as $n \to \infty$),
$$
\Theta_0^{-} [s_n; d_{n}] = \{\theta_0 \in l_0[s_n]: |\theta_{0,j} - 1/2| \leq d_n, |\mathcal{S}_{\theta_0}| = s_0\}.
$$

To prove \eqref{pf:risk-lb-1}, we first obtain a lower bound for $\FDP(\theta_0, \T)$. 
Let $\delta$ and $\tau$ be arbitrary positive numbers and $\max_j |\theta_{0,j} - 1/2| \leq d_{n, \epsilon}$, then 
\begin{align}
\FDP(\theta_0, \T) 
& \geq \frac{
s_n^{-1} \sum_{j=1}^n \mathbbm{1} 
\left\{\theta_{0,j} = 1/2, |X_j - m/2| > m\tau \right\}
}{
1 + s_n^{-1} \sum_{j=1}^n \mathbbm{1} 
\left\{\theta_{0,j} = 1/2, |X_j - m/2| > m\tau  \right\}
} \nonumber \\
& \geq 
1 - \left(
\frac{1}{s_n} \sum_{j=1}^n \mathbbm{1} \left\{
\theta_{0,j} = 1/2, |X_j - m/2| > m\tau
\right\}
\right)^{-1}.
\label{eqn:pf-lowerbound-1}
\end{align} 
Let $\mathcal{A} = \{\tau: \tau \leq d_{n, \epsilon} + \delta\}$, then 
\begin{align}
\FDP(\theta_0, \T) & \geq \FDP(\theta_0, \T)\mathbbm{1}_{\mathcal{A}} 
\nonumber \\
& \geq  
1 - \Big(
\frac{1}{s_n} \sum_{j=1}^n \mathbbm\{\theta_{0,j} = 1/2, |X_j - m/2| > m(d_{n, \epsilon} + \delta) \}
\Big)^{-1}
\nonumber \\
& \geq 1 - \max \Bigg\{\Big(
\frac{1}{s_n} \sum_{j=1}^n \mathbbm\{\theta_{0,j} = 1/2, X_j > m/2 + m(d_{n, \epsilon} + \delta) \}
\Big)^{-1}, 
\label{fdp-lb-1}\\
& \qquad \qquad \quad \quad \Big(
\frac{1}{s_n} \sum_{j=1}^n \mathbbm\{\theta_{0,j} = 1/2, X_j < m/2 - m(d_{n, \epsilon} + \delta) \}
\Big)^{-1} \Bigg\}.
\label{fdp-lb-2}
\end{align}
Next, we obtain a lower bound for $\FNP(\theta_0, \T)$ as follows:
let's write
\begin{align*}
\FNP(\theta_0, \T)
& = 
\frac{1}{s_n}
\sum_{i=1}^n \{\theta_{0,j} \neq 1/2\}(1-\T_j)\\
& = \frac{1}{s_n} \sum_{j=1}^n \mathbbm{1}\{
\theta_{0,j} \neq 1/2, -m\tau < X_j - m/2 < m\tau
\}\\
& = \frac{1}{s_n} \sum_{j=1}^n \mathbbm{1}\{
\theta_{0,j} \neq 1/2, -m \tau - m(\theta_{0,j} - 1/2) < X_j - m\theta_{0,j}  < m \tau - m(\theta_{0,j} - 1/2)\}\\
& \geq 
\frac{1}{s_n} \sum_{j=1}^n 
\mathbbm{1} \{\theta_{0,j} \neq 1/2, - m(\tau - d_{n, \epsilon}) < X_j - m \theta_{0,j} < m(\tau - d_{n, \epsilon})\}.
\end{align*}
For $\mathcal{A}^c = \{\tau: \tau > d_{n, \epsilon} + \delta\}$, we have
\begin{align}
\label{fnp-lb}
\FNP(\theta_0, \T)\geq \FNP(\theta_0, \T)\mathbbm{1}_{\mathcal{A}^c} \geq \frac{1}{s_n} \sum_{j=1}^n \mathbbm{1}
\left\{\theta_{0,j} \neq 1/2, \left|X_j - m\theta_{0,j}\right| < m\delta \right\}.
\end{align}
By combining the lower bounds in \eqref{fdp-lb-1}--\eqref{fnp-lb}, one obtains
\begin{align*}
& \FDP(\theta_0,\T) + \FNP(\theta_0, \T)\\
& \quad \geq 
\min \Bigg\{\frac{1}{s_n} \sum_{j=1}^n \mathbbm{1} \{\theta_{0,j} \neq 1/2, |X_j - m\theta_{0,j}| < m\delta\}, \\
& \quad \quad \quad \quad \quad \Big(1 - \frac{1}{s_n} \sum_{j=1}^n \mathbbm\{\theta_{0,j} = 1/2, X_j > m/2 + m(d_{n, \epsilon} + \delta) \} \Big)^{-1}, \\
& \quad \quad \quad \quad \quad \Big(1 - \frac{1}{s_n} \sum_{j=1}^n \mathbbm\{\theta_{0,j} = 1/2, X_j < m/2 - m(d_{n, \epsilon} + \delta) \}
\Big)^{-1} \Bigg\}.
\end{align*}
The last display implies that for any $\epsilon \in (0, 1)$, we have
\begin{align}
& P_{\theta_0} \left(\FDP(\theta_0, \T)+ \FNP (\theta_0, \T) \leq  1- \epsilon\right) \nonumber \\
& \quad \leq P_{\theta_0} 
\left(
\frac{1}{s_n} \sum_{j: \theta_{0,j}\neq 1/2} \mathbbm{1} \{|X_j - m\theta_{0,j}| < m\delta\} \leq 1-\epsilon
\right) 
\label{fdrfnrlb-1}\\
& \quad \quad + 2 P_{\theta_0} 
\left(
\frac{1}{s_n} \sum_{j: \theta_{0,j} =1/2} \mathbbm{1} \{X_j - m/2 \geq m(d_{n, \epsilon} + \delta)\} \leq \frac{1}{\epsilon}
\right).
\label{fdrfnrlb-2}
\end{align}
We choose $\epsilon = 4\bar \B(m\delta)$ with $\bar \B(m\delta) = P(|X_j - m\theta_{0,j}| \geq m\delta)$,
\eqref{fdrfnrlb-1} can be bounded by
\begin{align}
& P_{\theta_0} \left(\frac{1}{s_n} \sum_{j: \theta_{0,j} \neq 1/2} \mathbbm{1}\{ | X_j - m \theta_{0,j} | < m \delta\} \leq 1-\epsilon \right) \nonumber \\
& \quad = P_{\theta_0} \left(\sum_{j: \theta_{0,j} \neq 1/2} \mathbbm{1}\{ | X_j - m \theta_{0,j} | \geq m \delta\} \geq s_n\epsilon \right) \nonumber \\
& \quad = P_{\theta_0} \left(\sum_{j: \theta_{0,j} \neq 1/2} \left(\mathbbm{1}\{ | X_j - m \theta_{0,j} | \geq m \delta\} - 2 \bar \B(m\delta) \right) \geq s_n\epsilon/2 \right).
\label{fdrfnrlb-3}
\end{align}
By applying the Bernstein's inequality in Lemma \ref{bernstein-inq}, which we let
$A = s_n\epsilon/2$, $M \leq 1$, and 
$V = \sum_{j:\theta_{0,j} \neq 1/2} \text{Var}\left(\mathbbm{1}\{|X_j - m\theta_{0,j}| > m\delta\}\right)
\leq 2 s_n \bar \B(m\delta) = A$, 
\eqref{fdrfnrlb-3} is bounded by
$$\exp\left(- \frac{s_n^2 \epsilon^2}{8(V + s_n\epsilon/6)}\right)
\leq \exp\left(- \frac{s_n^2 \epsilon^2}{4(s_n\epsilon + s_n\epsilon/3)}\right)
= \exp\left(- \frac{3s_n \epsilon}{16}\right),
$$ 
The upper bound for \eqref{fdrfnrlb-2} can be obtained in a similar way:
first, subtracting $\bar{\mathbf{B}}(md_{n, \epsilon}+m\delta)$ on both sides and obtain
{\small $$
2 P_{\theta_0}\left(
\sum_{j: \theta_{0,j} = 1/2} \left( \mathbbm{1} \{X_j - m/2 \geq m(d_{n, \epsilon} + \delta) \} - \bar{\mathbf{B}}(md_{n, \epsilon} + m\delta) \right)
\leq \frac{s_n}{\epsilon} - (n-s_n) \bar{\mathbf{B}}(md_{n, \epsilon}+m\delta)
\right) .
$$}

Next, we choose $d_{n, \epsilon}$ such that 
$(n - s_n) \bar{\mathbf{B}}(md_{n, \epsilon}+m\delta) = s_n/\epsilon + s_n$, 
which implies $\delta = \bar{\mathbf{B}}^{-1}(\epsilon/4)/m$ due to
$\epsilon = 4\bar{\mathbf{B}}(m\delta)$ 
and 
$d_{n, \epsilon} = \frac{1}{m} \left(\bar{\mathbf{B}}^{-1} 
\left(\frac{s_n}{n-s_n}(1+\epsilon^{-1}) \right) - \bar{\mathbf{B}}^{-1}(\epsilon/4)\right).$
By applying the Bernstein's inequality again, the last display is bounded by 
\begin{align*}
& 2 P_{\theta_0} 
\left(\sum_{j: \theta_{0,j} = 1/2} 
\left(\mathbbm{1} \{X_j - m/2 \geq m(a_n + \delta) \} - \bar{\mathbf{B}}(a_n + \delta) \right)
\leq - s_n
\right) \\
& \quad \leq 2 \exp\left(- \frac{s_n^2}{2(2s_n/\epsilon + s_n/3)} \right)
= 2 e^{- 3s_n \epsilon/14}.
\end{align*}
We now combine the lower bound for both \eqref{fdrfnrlb-1} and \eqref{fdrfnrlb-2} to get
$$
\sup_{\T \in \mathcal{T}} \sup_{\theta_0 \in \Theta_0^-[s_n, d_{n, \epsilon}]} P_{\theta_0} (\FDP(\theta_0, \T)+ \FNP(\theta_0, \T)\leq 1- \epsilon) \leq e^{-3s_n\epsilon/16} + 2e^{-3s_n\epsilon/14} \leq 3e^{-s_n\epsilon/6}.
$$
We thus verified \eqref{pf:risk-lb-1}.

In the second step, we derive a lower bound for \eqref{pf:risk-lb-0}.
By taking the integral with respect to $\epsilon \geq 1/t_n$ for some $t_n \to \infty$, we obtain
$$
\inf_{\T \in \mathcal{T}} \inf_{\theta_0 \in \Theta_0^-[s_n, d_{n, \epsilon}]} 
P_{\theta_0} (\FDP(\theta_0, \T)+ \FNP(\theta_0, \T)> 1- \epsilon) 
> 1- 3e^{-s_n\epsilon/6}.
$$
Let $q = 1-\epsilon$ and $\theta_0 \in \Theta_0^-[s_n;b_n]$ with
\begin{align}
b_n =\frac{1}{m} \left( \bar{\mathbf{B}}^{-1}\left((t_n + 1) \frac{s_n}{n-s_n}\right) - \bar{\mathbf{B}}^{-1} \left(\frac{1}{4t_n}\right)\right),
\label{fdrfnr-4}
\end{align}
then
\begin{align*}
& \inf_{\T \in \mathcal{T}} \inf_{\theta_0 \in \Theta_0^-[s_n; b_n]} 
\mathfrak{R}(\theta_0, \T) \\
& \quad = 
\inf_{\T \in \mathcal{T}} \inf_{\theta_0 \in \Theta_0^-[s_n; b_n]} 
\int P_{\theta_0}(\FDP(\theta_0, \T)+ \FNP(\theta_0, \T)> 1-\epsilon) d(1-\epsilon) \\
& \quad \geq 
\inf_{\T \in \mathcal{T}} \inf_{\theta_0 \in \Theta_0^-[s_n; b_n]} 
\int (1 - 3e^{-s_n \epsilon/6}) d(1-\epsilon)\\
& \quad 
\geq 
\inf_{\T \in \mathcal{T}} \inf_{\theta_0 \in \Theta_0^-[s_n; b_n]} 
\int_0^{1-1/t_n} (1 - 3e^{-s_n (1-y)/6}) dy \\
& \quad 
\geq 1 - 1/t_n - 18/s_n.
\end{align*}
By invoking Lemma \ref{lem:inv-bin-bound}, since $s_n/n \to 0$, $m\gg \log^2n$,
and $m\varepsilon^4 = \log^2(2\log(n/s_n))/m \leq 4(1-v_1)^2\log^2n/m \to 0$ as 
$\varepsilon \sim \sqrt{2\log(n/s_n)/m} = \sqrt{2(1-v_1)\log n/m}$,
\eqref{fdrfnr-4} implies 
\begin{align}
\label{fdrfnr-5}
b_n \sim \sqrt{\frac{\log(n/s_n - 1) - \log (1+t_n)}{2m}} - \sqrt{\frac{\log 4t_n}{2m}}.
\end{align}
Choosing $t_n = \log (n/s_n) \to \infty$ as $n \to \infty$, 
then $1 - t_n^{-1} - 18s_n^{-1} \to 1$ and $b_n \sim \sqrt{\log (n/s_n)/(2m)}$.
By combining the above results, then for $a < 1$, we have
\begin{align*}
& \liminf_{n \to \infty} \inf_{\T \in \mathcal{T}} \sup_{\theta_0 \in \Theta_0[s_n, a]}
(\FDR(\theta_0, \T) + \FNR(\theta_0, \T)) \\
& \quad \geq 
\liminf_{n \to \infty} \inf_{\T \in \mathcal{T}} \inf_{\theta_0 \in \Theta_0^-[s_n; b_n]}
(\FDR(\theta_0, \T) + \FNR(\theta_0, \T)) \geq 1.
\end{align*}
We thus complete the proof.


\section{Useful lemmas for binomial distributions}
\label{sec:bound-binom}

\begin{lemma}
\label{lem:bin-coef-bound}
Let ${m \choose m/2 + ms}$ be the binomial coefficient for any $s \in [0, 1/2)$
and $T(a, p) = a\log (a/p) + (1-a) \log ((1-a)/(1-p))$ for $a, p \in (0, 1)$, then 
\begin{align}
\label{binomial-bound}
{m \choose m/2 + ms}
=
\frac{\sqrt{2} e^{-m T(1/2 + s, 1/2) + m\log 2 + \omega(s)}}{\sqrt{\pi m(1- 4 s^2)}}, 
\end{align}
where $\omega(s) = a_1 - a_2(s) - a_3(s)$, $(12m+1)^{-1} \leq a_1 \leq (12m)^{-1}$, 
$(6m+12ms + 1)^{-1} \leq a_2(s) \leq (6m+12ms)^{-1}$,
and 
$(6m-12ms + 1)^{-1} \leq a_3(s) \leq (6m-12ms)^{-1}$.

In particular, if $ms^2 \to 0$, then 
\begin{align}
\label{binomial-bound-2}
\log {m \choose m/2 + ms} \sim 
- \frac{1}{2}\log\left(\frac{\pi m}{2}\right) + m\log 2.
\end{align}
\end{lemma}
\begin{proof}
We write
$$
{m \choose {m/2 + ms}} = \frac{m!}{(m/2 + ms)!(m/2 - ms)!},
$$
then \eqref{binomial-bound} is obtained by directly applying the Sterling approximation:
$$
n! = \sqrt{2\pi} \exp \left((n + 1/2)\log n - n + a(n) \right),
$$
for any $n \in \mathbb{Z}^+$,
where $(12n+1)^{-1} \leq a(n) \leq (12n)^{-1}$. 
The result in \eqref{binomial-bound-2} can be proved by using the fifth point of $(d)$ in Lemma \ref{lem:entropy-bound}, i.e., $T(1/2+s, 1/2) \sim 2ms^2$ when $s = o(1)$, and then using that $a_1, a_2(s), a_3(s) \to 0$ as $m\to \infty$.
\end{proof}

\begin{lemma} 
\label{lem:binom-tail}
Let $X \sim \text{Bin}(m, \theta)$ and
$\bar{\mathbf{B}}_\theta(k) = \sum_{k' = k}^{m} b_\theta(k' = k)$, 
if $k = ma$ for any $1 > a > \theta \geq 1/2$, then 
\begin{align}
\label{eqn:binom-tail}
\frac{ e^{-mT(a, \theta)} }{\sqrt{2\pi ma (1-a)}}
\leq \bar{\mathbf{B}}_\theta(ma) 
\leq \frac{a(1-\theta) e^{-mT(a, \theta) + 1/(12m)} }{(a-\theta) \sqrt{2\pi ma (1-a)}},
\end{align}
where $T(a, \theta) = a \log \frac{a}{\theta} + (1 - a) \log \frac{1- a}{1-\theta}$.
Moreover, 
$\bar{\mathbf{B}}_\theta(ma) \leq e^{-mT(a, \theta)}$ and 
$- \frac{1}{m} \log \bar{\mathbf{B}}_\theta(ma) \sim T(a, \theta)$ as $m \to \infty$.
\end{lemma}

\begin{proof}
The lower bound in \eqref{eqn:binom-tail} is obtained by first noting that $\bar{\mathbf{B}}_\theta(ma) \geq b_\theta(x = ma) \geq {m \choose ma} \theta^{ma}(1-\theta)^{m - ma}$ and then using the lower bound of the binomial coefficient in \eqref{binomial-bound}.
The upper bound in \eqref{eqn:binom-tail} is obtained by first noting that $\bar{\mathbf{B}}_\theta(ma) =\frac{ \bar{\mathbf{B}}_\theta(ma)}{b(ma)} b(ma)$, and then invoking Lemma \ref{lem:mills-ratio-bin} and using the upper bound of the binomial coefficient given in \eqref{binomial-bound}.
The next inequality is the Chernoff bound for the binomial cdf. 
By taking logarithm of both upper and lower bounds in \eqref{eqn:binom-tail} and noting that the remaining terms in either side are of a smaller order of $m T(a, \theta)$, one obtains the last inequality.
\end{proof}

\begin{lemma} \citep{carter04}
\label{lem:binom-tail-carter}
Let $X \sim \text{Bin}(m, 1/2)$ and $\bar{\mathbf{B}}(k) = P(X \geq k)$, if $m \geq 28$, define 
$$
\gamma(\varepsilon) = \frac{(1+\varepsilon)\log (1+\varepsilon) + (1-\varepsilon)\log(1-\varepsilon) - \varepsilon^2}{2\varepsilon^4}
=
\sum_{r = 0}^\infty \frac{\varepsilon^{2r}}{(2r+3)(2r+4)},
$$
which is an increasing function. Define $\varepsilon = \frac{2K-M}{M}$ where $K = k -1$ and $M = m-1$, then 
there exists a $\lambda_k \in [(12k+1)^{-1}, (12k)^{-1}]$ and an $r_k \in [-C\log M/M, C/M]$ for some positive constant $C$ such that 
\begin{align}
\label{eqn:binom-tail}
\bar{\mathbf{B}}(k) = P(X \geq k) = \bar \Phi(\varepsilon \sqrt{M}) e^{A_m(\varepsilon)},
\end{align}
where $A_m(\varepsilon) = -M \varepsilon^4 \gamma(\varepsilon) - \log (1-\varepsilon^2)/2 - \lambda_{m-k} + r_k$
for all $\varepsilon$ corresponding to the range $m/2 < k \leq m-1$.
\end{lemma}

\begin{lemma} \citep{mckay89}
\label{lem:binom-tail-mckay}
Let $X \sim \text{Bin}(m, \theta)$, where
$0 < \theta < 1$, $m \geq 1$, and $m\theta \leq k \leq m$.
Define $z = (k - m\theta)/\sigma$, $\sigma = \sqrt{m\theta(1-\theta)}$, then 
\begin{align}
\label{eqn:binom-tail-mckay}
\bar \B_\theta(k) = \sigma \text{Bin}(k-1; m-1, \theta) Y(z) \exp(E_\theta(k, m)/\sigma),
\end{align}
where $\text{Bin}(k-1; m-1, \theta)$ is the binomial distribution at $k -1$ with parameters $m -1$ and $\theta$, $Y(z) = \bar\Phi(z)/\phi(z)$, and 
$0 \leq E_\theta(k, m) \leq \min\left\{\sqrt{\pi/8}, 1/z\right\}.$
\end{lemma}

\begin{lemma} \citep{slud77}
\label{lem:binom-tail-slud}
Let $\bar \B_\theta(k) = \sum_{q = k}^{m} b_\theta(q)$,
if $k \leq m\theta$, then 
$$
\bar \B_\theta(k) \geq 1 - \Phi 
\left(\frac{k - m\theta}{\sqrt{m\theta}} \right).
$$
\end{lemma}
\begin{lemma} \citep{diaconis91}
\label{lem:mills-ratio-bin}
Let $b_\theta(k) = \text{Bin}(m, \theta)$ and $\bar \B_\theta(k) =  \sum_{q = k}^{m} b_\theta(q)$, then for any $k > m\theta$, $\theta \in (0, 1)$, and $m \geq 1$, 
$$
\frac{k}{m} \leq \frac{\bar \B_\theta(k)}{b_\theta(k)} \leq \frac{k(1-\theta)}{k - m\theta}.
$$
\end{lemma}

\begin{lemma}
\label{lem:inv-bin-bound}
Let $\bar{\mathbf{B}}(m/2 + mx) = \sum_{q = m/2 + mx}^{m} b(q)$, 
define $M = m-1$, $K = k -1= m/2+mx-1$, and $\varepsilon = 2K/M - 1$,
if $\varepsilon \leq 0.957$ and $m \geq 28$, then 
for any $y \in (0, 1/2)$,
$$
\frac{m+1}{2} + \frac{M}{2} \underline{\digamma}_m(y)
\leq \bar \B^{-1}(y) 
\leq \frac{m+1}{2} + \frac{M}{2} \overline{\digamma}_m(y),
$$
where $- C\log M/M \leq r_k \leq C/M$ for some fixed constant $C$ and 
$\lambda_{m-k} \in [\frac{1}{12m+1}, \frac{1}{12m}]$, 
\begin{align*}
& \underline{\digamma}_m(y) =
\sqrt{
\left(
\left\{\left(\frac{2\log(1/y) - \log\log(1/y) + 2r_k - 2\lambda_{m-k} - \log(16\pi)}{2M} + \frac{1}{2}\right)_+\right\}^{1/2}
- \frac{1}{\sqrt{2}}
\right)_+},
\end{align*}
and 
\begin{align*}
& \overline{\digamma}_m(y)
= 
\sqrt{\frac{2\log(1/y) - 2\lambda_{m-k} + 2r_k}{M-2}}.
\end{align*}
In particular, if $m \to \infty$ and $m \varepsilon^4 \to 0$, then 
$$
\frac{m}{2} + \sqrt{\frac{m}{2}\left(\log(1/y) - \log\sqrt{\log(1/y)} - \log(4\sqrt{\pi})\right)_+}
\leq \bar \B^{-1}(y) 
\leq \frac{m}{2} + \sqrt{\frac{m\log(1/y)}{2}};
$$
furthermore, if $y \to 0$, then 
$\bar \B^{-1}(y)  \sim m/2 + \sqrt{m\log(1/y)/2}$.
\end{lemma}

\begin{proof}
We start with introducing the upper and lower bounds for the inverse of the standard Gaussian cdf in Lemma 36 of \citetalias{cast20MT}: denote $h = \bar \Phi(z)$, for any $h \in (0, 1/2)$, the upper tail probability of the standard Gaussian is 
\begin{align}
\label{pf:inv-bin-1}
\{(2\log(1/h) - \log\log(1/h) - \log(16\pi))_+\}^{1/2}
\leq \bar\Phi^{-1}(h) 
\leq \{2\log (1/h)\}^{1/2}.
\end{align}
Let $M = m-1$ and $\varepsilon = (2mx - 1)/M$, from 
\eqref{eqn:binom-tail}, we have
\begin{align}
\label{pf:inv-bin-2}
y = \bar \B(m/2 + mx) = \bar\Phi(\varepsilon \sqrt{M}) e^{A_m(\varepsilon)},
\end{align}
where $A_m(\varepsilon)$ is given in Lemma \ref{lem:binom-tail-carter}.
We then can obtain the bounds by relating them to the bounds of the inverse of the Gaussian cdf. 

{\it Upper bound.} 
By combining \eqref{pf:inv-bin-2} with the upper bound in \eqref{pf:inv-bin-1}, we obtain
$$
{\varepsilon^2}{M} \leq 2\log(1/y) + 2A_m(\varepsilon).
$$
By plugging-in the expression of $A_m(\varepsilon)$, we have
\begin{align*}
2\log(1/y) 
& \geq 
\varepsilon^2M + 2M\varepsilon^4\gamma(\varepsilon) 
+ \log (1-\varepsilon^2) + 2\lambda_{m-k} - 2r_k \\
& \geq 
\varepsilon^2M - 2\varepsilon^2 + 2\lambda_{m-k} - 2r_k,
\end{align*}
as $\gamma(\varepsilon) \geq 0$ and $\log (1-\varepsilon^2) > -2\varepsilon^2$ for $\varepsilon \in (0, 1/2)$.
The last display implies 
$$
\varepsilon \leq \sqrt{\frac{2\log(1/y) - 2\lambda_{m-k} + 2r_k}{M-2}}.
$$
Since $\bar \B^{-1}(y) = \frac{m}{2} + \frac{\varepsilon M + 1}{2}$, the last display implies the upper bound
$$
\bar \B^{-1}(y)
\leq \frac{m+1}{2} + \frac{M}{2} \sqrt{\frac{2\log(1/y) - 2\lambda_{m-k} + 2r_k}{M-2}}.
$$

{\it Lower bound.} By the lower bound given in \eqref{pf:inv-bin-1} and \eqref{pf:inv-bin-2}, we have
$$
\varepsilon^2 M \geq 2\log(1/y) + 2A_m(\varepsilon) - \log(\log(1/y) + A_m(\varepsilon)) - \log(16\pi),
$$
which implies 
\begin{align}
2\log(1/y) 
& \leq 
\varepsilon^2 M - 2A_m(\varepsilon)
+ \log (\log(1/y) + A_m(\varepsilon)) + \log (16\pi) 
\nonumber \\
& \leq
\varepsilon^2 M - 2A_m(\varepsilon)
+ \log \log(1/y) + \log (16\pi) 
\nonumber\\
& 
= 
\varepsilon^2 M + 2M\varepsilon^4 \gamma(\varepsilon) + \log(1-\varepsilon^2) + 2\lambda_{m-k} - 2r_k
+ \log \log(1/y) + \log (16\pi)
\nonumber\\
& 
\leq 
\varepsilon^2 M + 2M\varepsilon^4 - \varepsilon^2 + 2\lambda_{m-k} - 2r_k
+ \log \log(1/y) + \log (16\pi)
\label{pf:inv-bin-3-0} \\
& \leq 
2M\left(\varepsilon^2 + 1/\sqrt{2}\right)^2 - M + 2\lambda_{m-k} - 2r_k
+ \log \log(1/y) + \log (16\pi),
\label{pf:inv-bin-3}
\end{align}
where we used $A_m(\varepsilon) \leq 0$ to obtain the second inequality in the last display
and $\log(1-x) \leq -x$ as long as $1-x > 0$ and $\gamma(\varepsilon) \leq 1$ to obtain the third inequality. 
Note that
\begin{align*}
\gamma(\varepsilon) 
& = \sum_{r = 0}^\infty \frac{\varepsilon^{2r}}{(2r+3)(2r+4)}
\leq \frac{1}{12} \sum_{r = 0}^\infty \varepsilon^{2r}
= \frac{1}{12(1-\varepsilon^2)} \leq 1,
\end{align*}
as long as $\varepsilon < \sqrt{1- 1/12} \approx 0.957$.
The lower bound in \eqref{pf:inv-bin-3} implies that
$$
\left(\varepsilon^2 + 1\sqrt{2}\right)^2 
\geq \frac{1}{2M} \left(
2\log(1/y) - \log\log(1/y) + M + 2r_k - 2\lambda_{m-k} - \log(16\pi)
\right).
$$
By taking the square root on both sides and subtracting $1/\sqrt{2}$ in the preceding display, the lower bound for $\bar \B^{-1}(y)$ follows by plugging the lower bound of $\varepsilon$ into $(m + 1 + \varepsilon M)/2 = \bar \B^{-1}(y)$.

The second inequality in the lemma can be proved as follows: for a sufficiently large $M$, $m \approx M$, $r_k = o(1)$,
and $\lambda_{m-k} = o(1)$, the expression of the upper bound for $\bar \B^{-1}(y)$ is then of the same order as $m/2 + \sqrt{m\log(1/y)/2}$.
Since $m\varepsilon^4 = o(1)$, the upper bound \eqref{pf:inv-bin-3-0} implies
\begin{align}
2\log (1/y) \leq M\varepsilon^2 + \log\log(1/y) + \log (16\pi) + o(1),
\end{align}
which implies 
$\varepsilon^2 \geq 2\log(1/y) - \log\log(1/y) - \log(16\pi)$ for a sufficiently large $m$.
Using $\bar \B^{-1}(y) =(m + 1 + \varepsilon M)/2 \approx m/2 + m\varepsilon/2$ for a sufficiently large $m$, we obtain the desired lower bound.
If $y \to 0$, then $\log(1/y) \gg \log\log(1/y)/2 + \log(4\sqrt{\pi})$, which leads to the last inequality.
\end{proof}

\begin{lemma} 
\label{lem:entropy-bound}
For $a \geq 1/2$ and $p \geq 1/2$,
let
\begin{align}
T(a, p) = a \log \left( \frac{a}{p} \right) + (1-a) \log \left(\frac{1-a}{1-p}\right),
\label{eqn:T}
\end{align}
and define $h_p(\epsilon) = T(p + \epsilon, p)$,
then 
\begin{enumerate}
\item[(a)] $h_p(\epsilon)$ is a monotone increasing function on $\epsilon \in (0, 1-p)$
and a monotone decreasing function on $\epsilon \in (-p, 0)$;
\item[(b)] $h_p(\epsilon)$ is continuous and nonnegative; it achieves the global minimum at $\epsilon = 0$. 
\item[(c)] $\epsilon^3 h_p'''(\epsilon)/6$ is positive if $\epsilon \in [0, 1-p)$ or $\epsilon \in (-p, 1/2-p)$; it is negative if $\epsilon \in (1/2 - p, 0)$.
When $p = 1/2$, $\epsilon^3 h_p'''(\epsilon)/6 \geq 0$ for any $\epsilon \in (-1/2, 1/2)$.
\item[(d)] There exists $\epsilon^\star$ such that $\epsilon^\star \in [0, \epsilon]$ if $\epsilon \geq 0$ 
or $\epsilon^\star \in [\epsilon, 0]$ if $\epsilon < 0$, 
\begin{align}
\label{eqn:hp-bound}
h_p(\epsilon) = \frac{\epsilon^2}{2p(1-p)} + \frac{\epsilon^3 (2p+2\epsilon^\star - 1)}{6(p+\epsilon^\star)^2(1-p-\epsilon^\star)^2}.
\end{align}
In particular, we have the following:
\begin{enumerate}
\item[(1)] if $\epsilon = 0$, then $h_p(\epsilon) = 0$;
\item[(2)] if $0 < \epsilon < 1-p$, then,
$$
\frac{\epsilon^2}{2p(1-p)} 
\leq h_p(\epsilon)
\leq \frac{\epsilon^2}{2p(1-p)}
+ \frac{8 \epsilon^3 (2p + 2\epsilon - 1)}{3 (1-4(p + \epsilon - 1/2)^2)^2};
$$ 

\item[(3)] if $1/2- p < \epsilon < 0$, then 
$$
\frac{\epsilon^2}{2p(1-p)} + \frac{8 \epsilon^3 (2p -1)}{3 (1-4(p + \epsilon - 1/2)^2)^2}
\leq h_p(\epsilon)
\leq \frac{\epsilon^2}{2p(1-p)};
$$

\item[(4)] if $-p < \epsilon < 1/2-p$, then
$$
\frac{\epsilon^2}{2p(1-p)} 
\leq h_p(\epsilon)
\leq \frac{\epsilon^2}{2p(1-p)}
+ \frac{8 \epsilon^3 (2p - 1)}{3 (1-4(p - 1/2)^2)^2};
$$ 

\item[(5)] If $p = 1/2$ and $\epsilon = o(1)$, then $h_p(\epsilon) \sim 2\epsilon^2$ for any $\epsilon \in (-1/2, 1/2)$.
\end{enumerate}
\end{enumerate}
\end{lemma}

\begin{proof}
The following results are useful for our proof:
\begin{align*}
h_p(\epsilon) & = (p + \epsilon) \log \left(\frac{p + \epsilon}{p}\right) 
+ (1 - p - \epsilon) \log \left(\frac{1 - p- \epsilon}{1-p}\right), \\
h_p'(\epsilon) & = \log (1 + \epsilon p^{-1}) - \log (1 - {\epsilon}(1 - p)^{-1}), \\
h_p''(\epsilon) & = (p + \epsilon)^{-1} + (1-p-\epsilon)^{-1}, \\
h_p'''(\epsilon) & = \frac{2p + 2\epsilon - 1}{(p + \epsilon)^2(1- p - \epsilon)^2}.
\end{align*}
First, let us verify (a)--(c).
(a) is easy to verify as $h_p'(\epsilon) > 0$ if $\epsilon > 0$ and  $h_p'(\epsilon) < 0$ if $\epsilon < 0$. 
For (b), the proof of $h_p(\epsilon)$ for $\epsilon \in (-p, 1-p)$ is continuous is trivial and is omitted. since $h''(\epsilon) > 0$, $h_p(\epsilon)$ is a convex function; also, $h_p(\epsilon)$ achieves the global minimum at $\epsilon =  0$ and $h(0) = 0$, so, $h_p(\epsilon)$ is nonnegative. 

Next, we prove (d). By applying the Taylor's theorem up to the third term together with the mean-value theorem, we obtain 
\begin{align}
\label{taylor-hp}
h_p(\epsilon) = \frac{\epsilon^2}{2p(1-p)} + \frac{\epsilon^3 h_p'''(\epsilon^\star)}{6},
\end{align}
for an $\epsilon^\star$ between $0$ and $\epsilon$.

Last, we prove (1)--(5).  (1) is trivial.
(2) and (3) can be verified by plugging-in the expression for $h'''_p(\epsilon^\star)$ and noting that $ h'''(\epsilon^\star) > 0$ if $\epsilon \in (0, 1-p)$ and $\epsilon^3 h'''(\epsilon^\star) < 0$ if $\epsilon \in (1/2 - p, 0)$ respectively. 
(4) can be proved in a similar way but noticing that $h'''_p(\epsilon) < 0$ but $\epsilon h'''_p(\epsilon) > 0$ if $\epsilon \in (-p, 1/2-p)$.
The last result can be verified easily by plugging $p = 1/2$ into \eqref{taylor-hp} and then using $\epsilon = o(1)$, then $\epsilon^3h'''_p(\epsilon^\star) = o(12\epsilon^2)$ for any $\epsilon \in (-1/2, 1/2)$.
\end{proof}

\begin{lemma}
\label{lem:lb-bin-cdf}
Let $X \sim \text{Bin}(m, \theta)$ and $\bar \B_\theta(\cdot)$ be one minus of its cdf, for $\xi := \xi(w)$ in \eqref{eqn:xi}, for any $w \in (0, 1)$, if $m/2 < m \theta \leq m/2 + m\xi \leq m$ and $m\xi^4 \to 0$ as $m \to \infty$, then 
$$
\bar \B_{\theta}(m/2 + m\xi) \geq 
\frac{1}{2}\sqrt{\frac{1-2(\theta - 1/2)}{1-2\xi}}
\bar \Phi \left(
\frac{2\sqrt{m}(\xi - (\theta - 1/2))}{\sqrt{1-4(\theta - 1/2)^2}}.
\right)
$$
\end{lemma}
\begin{proof}
Denote $\mu = \theta - 1/2$ and let $\sigma = \sqrt{m (1-4\mu^2)}/2$, $z = \frac{2\sqrt{m}(\xi - \mu)}{\sqrt{1-4\mu^2}}$, and $Y(z) = \bar\Phi(z) /\phi(z)$, then by Lemma \ref{lem:binom-tail-mckay}, 
\begin{align}
\bar \B_\theta(m/2 + m\xi)
& \geq \sigma \text{Bin}(m/2 + m\xi - 1; m -1, \mu + 1/2) Y(z) 
\nonumber \\
& = \frac{(1/2 + \xi) \sqrt{m(1-4\mu^2)}}{2(1/2 + \mu)}\frac{b_{\theta}(m/2 + m\xi)}{\phi(z)} \bar\Phi(z)
\label{pf:lb-bin-cdf-1}
\end{align}
By Lemma \ref{lem:bin-coef-bound}, the ratio
$$
\frac{b_{\theta}(m/2 + m\xi)}{\phi(z)}
\geq \frac{2}{\sqrt{m(1-4\xi^2)}} e^{-mT(1/2 + \xi, 1/2+\mu) + z^2/2}.
$$
Since $\mu < \xi$, by (3) in Lemma \ref{lem:entropy-bound} and the assumption $m\xi^4 \to 0$, 
the last display can be further bounded below by 
${2}(1-o(1))/{\sqrt{m(1-4\xi^2)}} .$
By plugging-in the above lower bound, 
for a sufficiently large $m$,
\eqref{pf:lb-bin-cdf-1} can be bounded from below by 
$$
(1-o(1))\frac{(1/2 + \xi) \sqrt{1-4\mu^2}}{(1/2 + \mu) \sqrt{1-4\xi^2}} \bar \Phi(z)
\geq \frac{1}{2}\sqrt{\frac{1-2\mu}{1-2\xi}} \bar \Phi(z).
$$
The result follows by plugging-in the expression of $z$.
\end{proof}

\begin{lemma}
\label{lem:ratio-bin-cdf}
Let $X \sim \text{Bin}(m, \theta)$ and $\bar \B_\theta(\cdot)$ be its upper tail probability, for positive $a_1, a_2 \leq \frac{1}{m}$ such that $|2ma_1^2 - 2ma_2^2| \leq 1/4$ and $1/2 \leq \theta < 1$, if $m \to \infty$, then there exists a $C > 0$ depending on $\theta, a_1, a_2$ such that
$$
\frac{\bar \B_\theta(m/2 + ma_1)}{\bar \B_\theta(m/2 + ma_2)}
\geq C \exp\left(- \frac{2m |a_1^2 - a_2^2|}{1-4(\theta - 1/2)^2}\right).
$$
\end{lemma}
\begin{proof}
If $a_2 \geq a_1$, the result is trivial. Let's focus on the case $a_1 > a_2$. 
Denote $b_\theta(x; m-1) = \text{Bin}(x, m-1, \theta)$ and $b_\theta(x) = b_\theta(x; m)$.
Let $\mu = \theta - 1/2$,
if $0 \leq a_2 < a_1 \leq \mu$, then
$$
\frac{\bar \B_\theta(m/2 + ma_1)}{\bar \B_\theta(m/2 + ma_2)} \geq 1/2 > 1/4 \geq \frac{1}{4} e^{- 2m |a_1^2 - a_2^2|},
$$
as $1/2 < \bar \B_\theta(m/2 + ma_1) < \bar \B_\theta(m/2 + ma_2) < 1$.

If $0 \leq a_2 \leq \mu \leq a_1$, 
then $\bar \B_\theta(m/2 + ma_2) \geq 1/2$. 
By Lemma \ref{lem:lb-bin-cdf}, we have
$$
\bar \B_\theta(m/2 + ma_1)
\geq 
\frac{1}{2} \sqrt{\frac{1 -2\mu}{1-2a_1}}
\bar\Phi
\left(
\frac{2\sqrt{m}(a_1 - \mu)}{\sqrt{1-4\mu^2}}
\right)
\geq \frac{1}{2} \sqrt{\frac{1 -2\mu}{1-2a_1}}
\bar\Phi
\left(
\frac{2\sqrt{m}(a_1 - a_2)}{\sqrt{1-4\mu^2}}
\right).
$$
Using that $2m|a_1^2 - a_2^2| \leq \frac{1}{4}$ and $\sqrt{\frac{1 -2\mu}{1-2a_1}} \geq \frac{1}{\sqrt{2}}$ and Lemma \ref{lem:gaussian-bound}, we have
\begin{align*}
\bar\Phi
\left(
\frac{2\sqrt{m}(a_1 - a_2)}{\sqrt{1-4\mu^2}}
\right) 
& \geq \frac{\frac{2\sqrt{m}(a_1 - a_2)}{\sqrt{1-4\mu^2}}}{1 + \frac{4{m}(a_1 - a_2)^2}{1-4\mu^2}}
\phi
\left(
\frac{2\sqrt{m}(a_1 - a_2)}{\sqrt{1-4\mu^2}}
\right)  \\
& \geq 
\min\left\{
\frac{1}{2 \sqrt{2(1-4\mu^2})},
\sqrt{\frac{1-4\mu^2}{2}}
\right\}
\frac{1}{\sqrt{2\pi}}
\exp\left(- \frac{2m(a_1^2 - a_2^2)}{1-4\mu^2}\right)\\
& = C_1 \exp\left(- \frac{2m(a_1^2 - a_2^2)}{1-4\mu^2}\right).
\end{align*}
Last, if $0 \leq \mu < a_2 < a_1$, by invoking Lemma \ref{lem:binom-tail-mckay}, 
let $\sigma = \sqrt{m(1-4\mu^2)}/2$ and $z_i = (ma_i - m\mu)/\sigma$ for $i = 1, 2$, then
\begin{align}
\label{pf:ratio-bin-cdf-1}
\frac{\bar \B_\theta(m/2 + ma_1)}{\bar \B_\theta(m/2 + ma_2)}
= 
\frac{b_\theta(m/2 + ma_1 -1; m-1) Y(z_1)}{b_\theta(m/2 + ma_2 -1; m-1) Y(z_2)}
\exp(A_m),
\end{align}
where 
$Y(z) = \bar \Phi(z)/\phi(z)$ and
$A_m = (E_\theta(m/2 + ma_1, m) - E_\theta(m/2 + ma_2, m))/\sigma$.
By Lemma \ref{lem:gaussian-bound}, we have
$\frac{\bar \Phi(z_1)}{\bar \Phi(z_2)}
\geq 
\frac{z_1 z_2}{1 + z_1^2}
\frac{\phi(z_1)}{\phi(z_2)}.$
By plugging-in expressions of $z_1$ and $z_2$ and using that $z_2^2 > 1$,
we obtain 
$$
\frac{\bar \Phi(z_1)}{\bar \Phi(z_2)}
\geq 
\frac{a_2-\mu}{a_1-\mu}\exp\left(- \frac{2{m}(a_1^2 - a_2^2)}{1-4\mu^2}\right).
$$
Next, by Lemma \ref{lem:bin-coef-bound} and then $(d)$ in Lemma \ref{lem:entropy-bound}, as $m\to \infty$, $m \approx m-1$ and $m/2 - 1 \approx m/2$, then
\begin{align*}
\frac{b_\theta(m/2 + ma_1 -1; m-1)/\phi(z_1)}{b_\theta(m/2 + ma_2 -1; m-1)/\phi(z_2)} 
& =
\sqrt{\frac{1 - 4a_2^2}{1-4a_1^2}}
e^{- (m - 1) (T(1/2 + a_1, \theta) - T(1/2 + a_2, \theta)) - z_1^2/2 + z_2^2/2}\\
& \geq
\sqrt{\frac{1 - 4a_2^2}{1-4a_1^2}}
e^{-m (a_1^2 - a_2^2) K},
\end{align*}
where $T(a, p) = a\log(a/p) + (1-a)\log((1-a)/(1-p))$
and $K = \min\{K(\epsilon_1^\star), K(\epsilon_2^\star)\}$,
$K(\epsilon_i^\star) = \frac{8(\mu + \epsilon_i^\star)}{3 (1/2 + \mu + \epsilon_i^\star)^2 (1/2 - \mu- \epsilon_i^\star)^2} > 0$.
Since $2m|a_1^2 - a_2^2| \leq 1/4$ by assumption, $e^{- m(a_1^2 - a_2^2)K} > e^{-K/8}$. Thus, 
$$
\frac{b_\theta(m/2 + ma_1 -1; m-1)/\phi(z_1)}{b_\theta(m/2 + ma_2 -1; m-1)/\phi(z_2)} 
\geq \sqrt{\frac{1 - 4a_2^2}{1-4a_1^2}} \exp\left(- \frac{K}{8}\right).
$$
Moreover, as $m \to \infty$, $E_\theta(m/2 + ma_1, m) = \frac{1}{2\sqrt{m} a_1^2} \to 0$.
By combining the relevant bounds above, we obtain
$$
\eqref{pf:ratio-bin-cdf-1} \geq
\frac{(a_2-\mu)\sqrt{1-4a_2^2)} e^{-K/8}}{(a_1-\mu)(1-4a_1^2)}
\exp\left(- \frac{2{m}(a_1^2 - a_2^2)}{1-4\mu^2}\right)
\geq C_2 \exp\left(- \frac{2{m}(a_1^2 - a_2^2)}{1-4\mu^2}\right).
$$
The proof is completed by taking $C = \min\{C_1, C_2, 1/4\}$.
\end{proof}

\section{Auxiliary lemmas}
\label{sec:aux-lemma}

\begin{lemma}
\label{bound-Omega0}
Consider the event $\Omega_n = \{\# \{j \in \mathcal{S}_0, |X_j - {m}/{2} | > bm \zeta_n \} \geq s_n - K_n\}$ for $X_j \sim \text{Bin}(m, p_j)$, $p_j > 1/2$, $\zeta_n$ is given in \eqref{eqn:zeta}, $s_n = |\mathcal{S}_0|$, and $K_n = o(s_n)$, if $s_n \ll (1 + 1/w)^{-(a-b)^2/24}$,
then $P(\Omega_n^c) = o(1)$. 
\end{lemma}

\begin{proof}
By definition, 
\begin{align*}
\Omega_n^c & = \{ \# \{j \in \mathcal{S}_0: |X_j - m/2 | > b m \zeta_n\} < s_n - K_n\}\\
& = \{ \# \{j \in \mathcal{S}_0: |X_j - {m}/{2} | \leq bm \zeta_n\} > K_n\}
\end{align*}
Thus, $P(\Omega_n^c) = P(\text{Bin}(s_n, h_n) > K_n)$,
where 
\begin{align*}
h_n & = P(|X_j - m/2| \leq b m\zeta_n) \\
& = P(|X_j - mp_j  + mp_j - m/2| \leq bm \zeta_n)  \\
& \leq P(|mp_j - m/2| - |X_j - mp_j| \leq bm\zeta_n/2) \\
& = P(|X_j - mp_j| > |mp_j - m/2| - b m \zeta_n/2) \\
& \leq P(|X_j - mp_j| > (a-b)m \zeta_n/2),
\end{align*}
which we used the inequality $|a +b| \geq |a| - |b|$ to obtain the first inequality.
Since $\mathbb{E}(X_j) = mp_j$, by applying the Chernoff bound:
$
P(|x - \mu| > \eta \mu) \leq 2e^{-\eta^2 \mu/3}
$
for $0 < \eta < 1$, $\zeta_n^2 \sim \frac{1}{2m} \log(1+w^{-1})$, and
choosing $\eta = (a-b) \zeta_n/(2p_j)$, then
$$
h_n \leq 
2 e^{-m (a-b)^2 \zeta_n^2/(12p_j)}
\leq
2e^{-(a-b)^2 \log(1+1/w)/24}
= 2 
(1+1/w)^{- (a-b)^2/24} := \tilde h_n,
$$
which $\to$ 0 as $m \to \infty$.
We thus obtain
$$
P(A_n^c) = P(\text{Bin}(s_n, h_n) > K_n) 
\leq P(\text{Bin}(s_n, \tilde h_n) > K_n).
$$
We use the Bernstein's inequality (see Lemma \ref{bernstein-inq}) to control the probability in the upper bound of last display. 
Let $Z_i \sim \text{Bern}(\tilde h_n)$, $1 \leq i \leq s_n$ as $s_n$ independent Bernoulli variables,
we choose $A = K_n \geq 2s_n \tilde h_n$ and $\sum_{j \in \mathcal{S}_0} \text{Var}(Z_i) = s_n \tilde h_n (1-\tilde h_n) \leq s_n \tilde h_n =V$ and $M = 1$, 
then
$$
P(\text{Bin}(s_n, \tilde h_n) > K_n) = P\left(\sum_{i=1}^{s_n} Z_i > K_n\right)
\leq \exp\left(- \frac{K_n^2}{2 s_n \tilde h_n + 2K_n/3}\right)
\leq \exp\left(-\frac{6}{5}s_n\tilde h_n\right) \to 0,
$$
as $s_n\tilde h_n \to 0$ for a large enough $n$.
\end{proof}

\begin{lemma} (Bernstein's inequality)
\label{bernstein-inq}
Let $W_i$, $1 \leq i \leq n$, be centered independent variables with $|W_i| \leq D$ and $\sum_{i=1}^n \text{Var}(W_i) \leq V$, then for any $A \geq 0$, 
\begin{align*}
& P\left(
\sum_{i=1}^n W_i \geq A
\right)
\leq \exp \left(
- \frac{A^2}{2(V + DA/3)}
\right), \\
& P\left(
\sum_{i=1}^n W_i \leq -A
\right)
\leq \exp \left(
- \frac{A^2}{2(V + DA/3)}
\right).
\end{align*}
\end{lemma}

\begin{lemma} (KMT approximation theorem \citep{kmt75})
\label{lem:KMT}
Let $\epsilon_1, \dots, \epsilon_n$ be i.i.d. random variables with $\mathbb{E}(\epsilon_1) = 0$ and $\mathbb{E}(\epsilon_1^2) = 1$, and $\mathbb{E}e^{\theta \epsilon_1} < \infty$ for some $\theta > 0$. 
For each $k$, let $S_k = \sum_{i=1}^k \epsilon_i$. 
Then for any $n$, it is possible to construct a version of $(S_k)_{0 \leq k \leq n}$ and a standard Brownian motion $(W_k)_{0 \leq k \leq n}$ on the same probability space such that for all $x \geq 0$, 
\begin{align}
\label{eqn:KMT}
P \left(\max_{k\leq n}|S_k - W_k| \geq C\log n + x\right) \leq K_1 e^{- K_2 x},
\end{align}
for some positive constants $C_1, K_1, K_2$ do not depend on $n$. 
\end{lemma}

\begin{lemma}
\label{lem:gaussian-bound}
For any $x>0$, let $\phi(\cdot)$ and $\Phi(\cdot)$ be the pdf and cdf of the standard normal distribution respectively. Denote $\bar \Phi(\cdot) = 1-\Phi(\cdot)$, then for any $x > 0$,
\begin{align*}
\frac{x \phi(x)}{1+x^2} < & \bar \Phi(x) < \frac{\phi(x)}{x},
\end{align*}
In particular, for any $x \geq 1$, $\bar\Phi(x) \geq \frac{\phi(x)}{2x}$ and, if $x \to \infty$, $\bar \Phi(x) \sim \frac{\phi(x)}{x}$.
If $x \to 0$ is small, we also have
\begin{align*}
\frac{1}{\sqrt{2\pi}} e^{-x^2/2} < \bar \Phi(x) < \frac{1}{2} e^{-x^2/2}.
\end{align*}
\end{lemma}

\begin{lemma} \citep{carter04}
\label{lem:gaussian-bound-carter}
Let $\phi(x)$ and $\Phi(x)$ be the pdf and the cdf of the standard normal distribution respectively,
define $\bar \Phi(x) = 1-\Phi(x)$,
$$
\rho(x) = \phi(x)/\bar\Phi(x), \quad 
r(x) = \rho(x) - x,
$$
then for $x\in \mathbb{R}$ and $\delta \geq 0$, the relation between $\bar \Phi(x + \delta)$ and $\bar \Phi(x)$ satisfies the following inequality:
\begin{itemize}
\item[(i)] 
$e^{-\delta \rho(x + \delta)} \leq \bar \Phi(x + \delta)/\bar \Phi(x) \leq e^{-\delta \rho(x)}$,
\item[(ii)] 
$e^{-\delta r(x)} \leq e^{x\delta + \delta^2/2} \bar \Phi(x + \delta)/\bar \Phi(x) \leq e^{-\delta r(x +\delta)}$,
\item[(iii)]
$e^{-\rho(x) \delta - \delta^2/2} \leq \bar \Phi(x + \delta)/\bar \Phi(x) \leq e^{-x \delta - \delta^2/2}$.
\end{itemize}
\end{lemma}

\begin{lemma} (Lemma 40 of \citetalias{cast20MT})
\label{lem:boundFDR}
For $m \geq 1$ and $p_1,\dots, p_m \in (0, 1)$, consider $U = \sum_{i=1}^m B_i$, where $B_i \sim \text{Ber}(p_i)$, $1\leq i\leq m$, are independent. For any nonnegative variable $T$ independent of $U$, we have
$$
\mathbb{E}\left(
\frac{T}{T+U} \mathbbm{1}_{\{T >0\}}
\right)
\leq \exp\left(-\mathbb{E} U\right) + \frac{12 \mathbb{E}T}{\mathbb{E} U}.
$$
\end{lemma}

\section{Additional numerical experiments}
\label{sec:ad-sim}

\begin{figure}[h!]
\centering
	  \subfigure[$m = 85$, $s_n/n = 0.001$]{\includegraphics[width=0.33\textwidth]{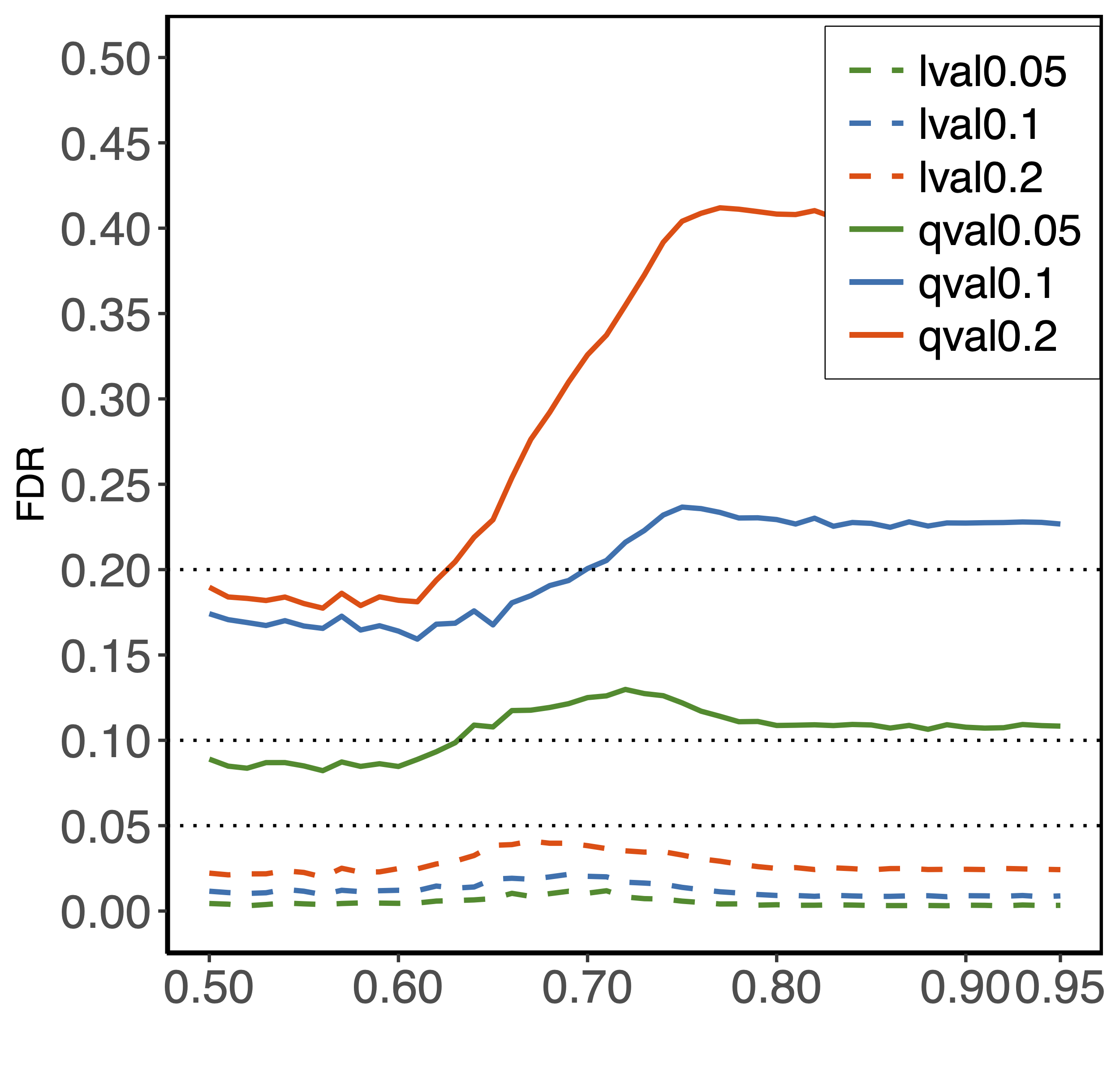}}%
     \subfigure[$m = 85$, $s_n/n = 0.1$]{\includegraphics[width=0.33\textwidth]{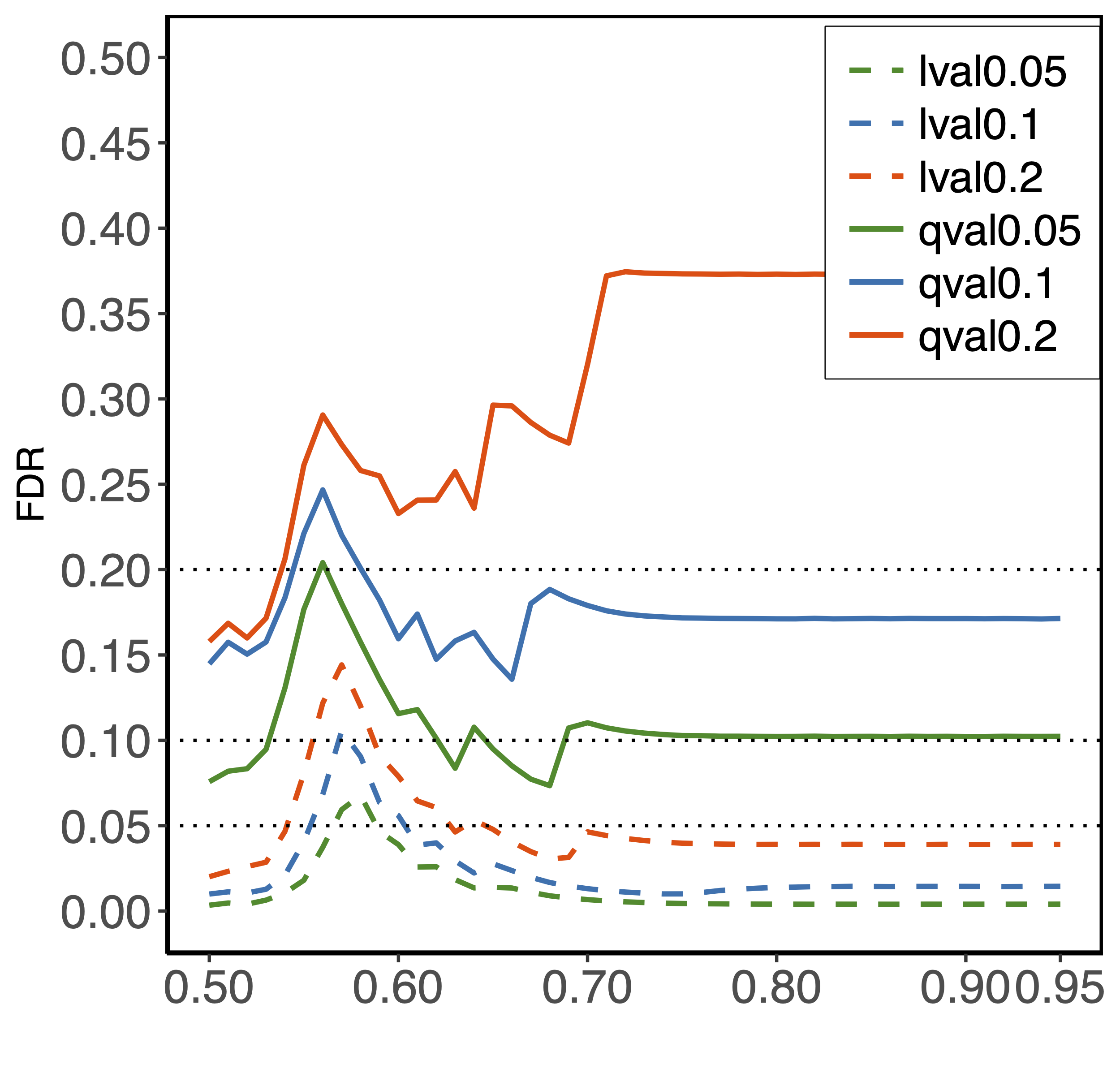}}%
    \subfigure[$m = 85$, $s_n/n = 0.5$]{\includegraphics[width=0.33\textwidth]{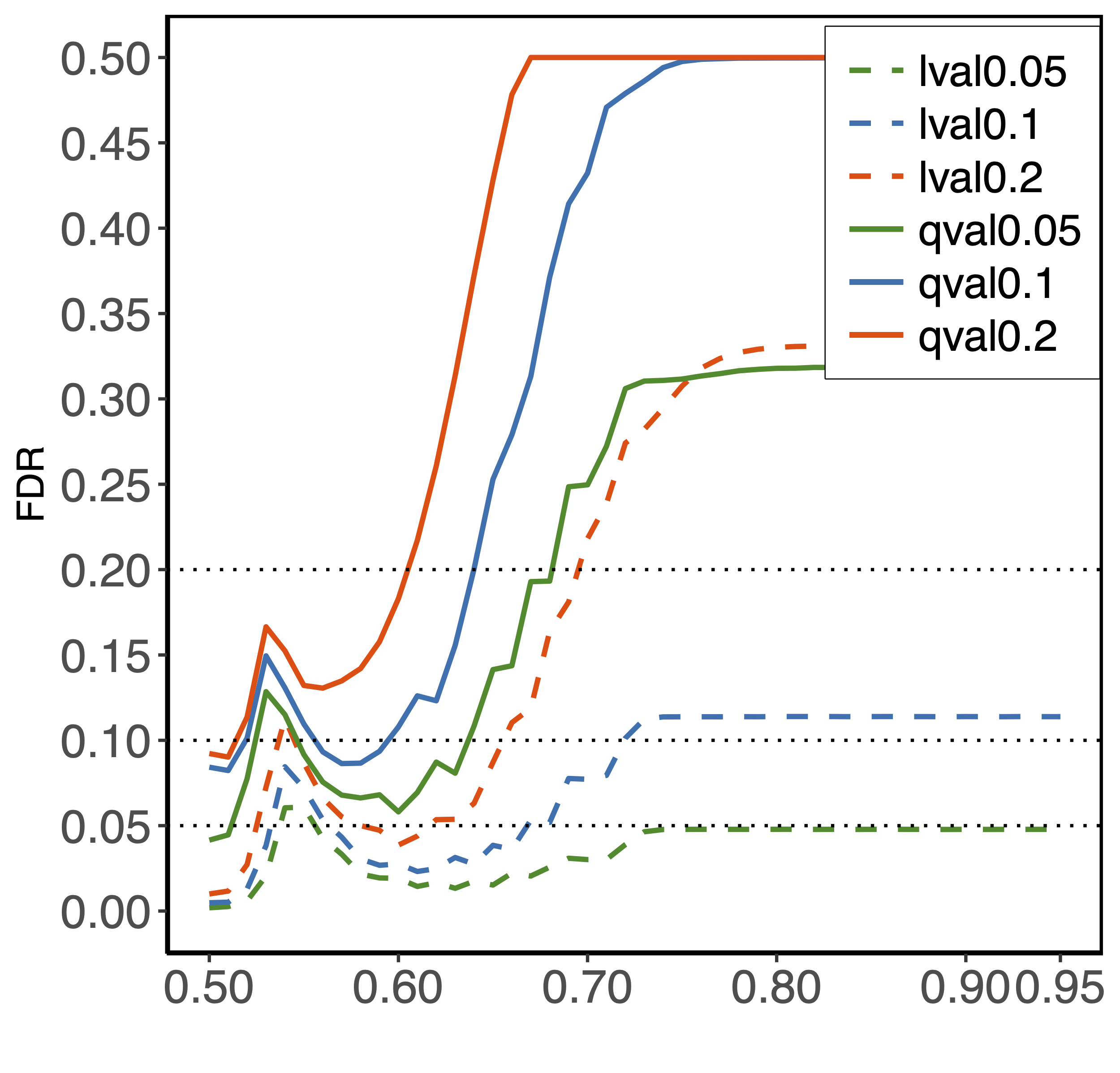}} 
    \subfigure[$m = 200$, $s_n/n = 0.001$]{\includegraphics[width=0.33\textwidth]{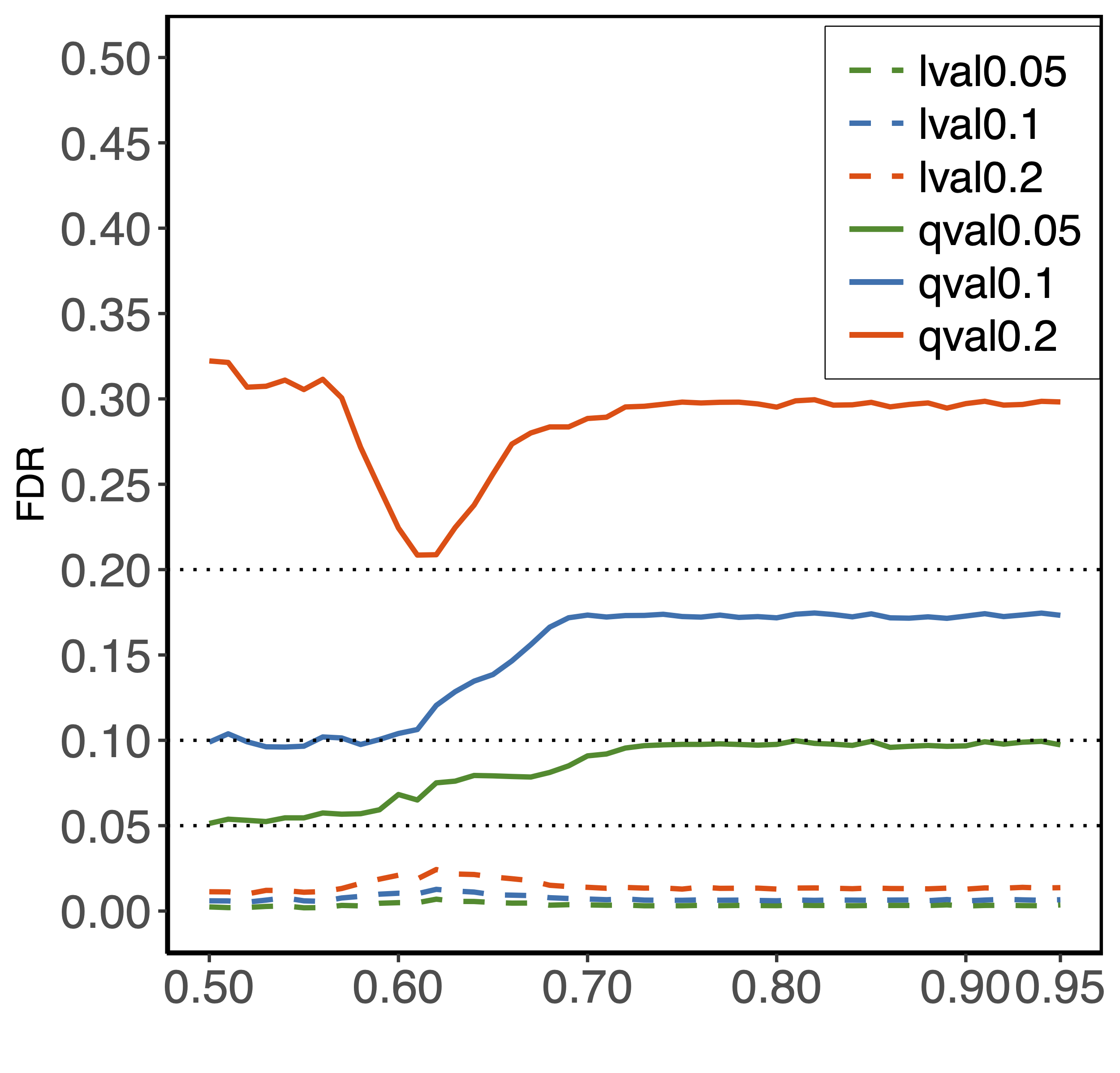}}%
     \subfigure[$m = 200$, $s_n/n = 0.1$]{\includegraphics[width=0.33\textwidth]{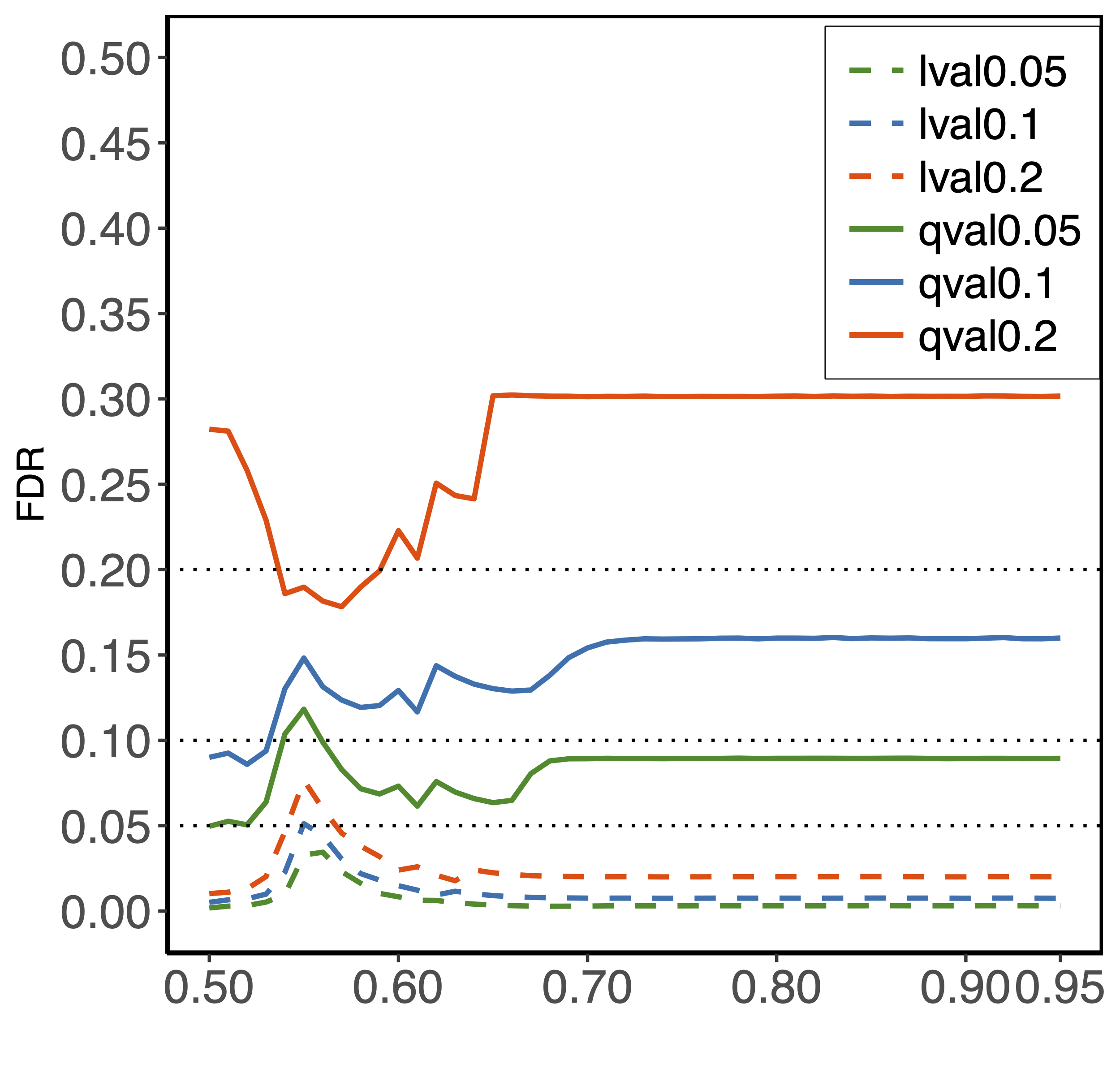}}%
    \subfigure[$m = 200$, $s_n/n = 0.5$]{\includegraphics[width=0.33\textwidth]{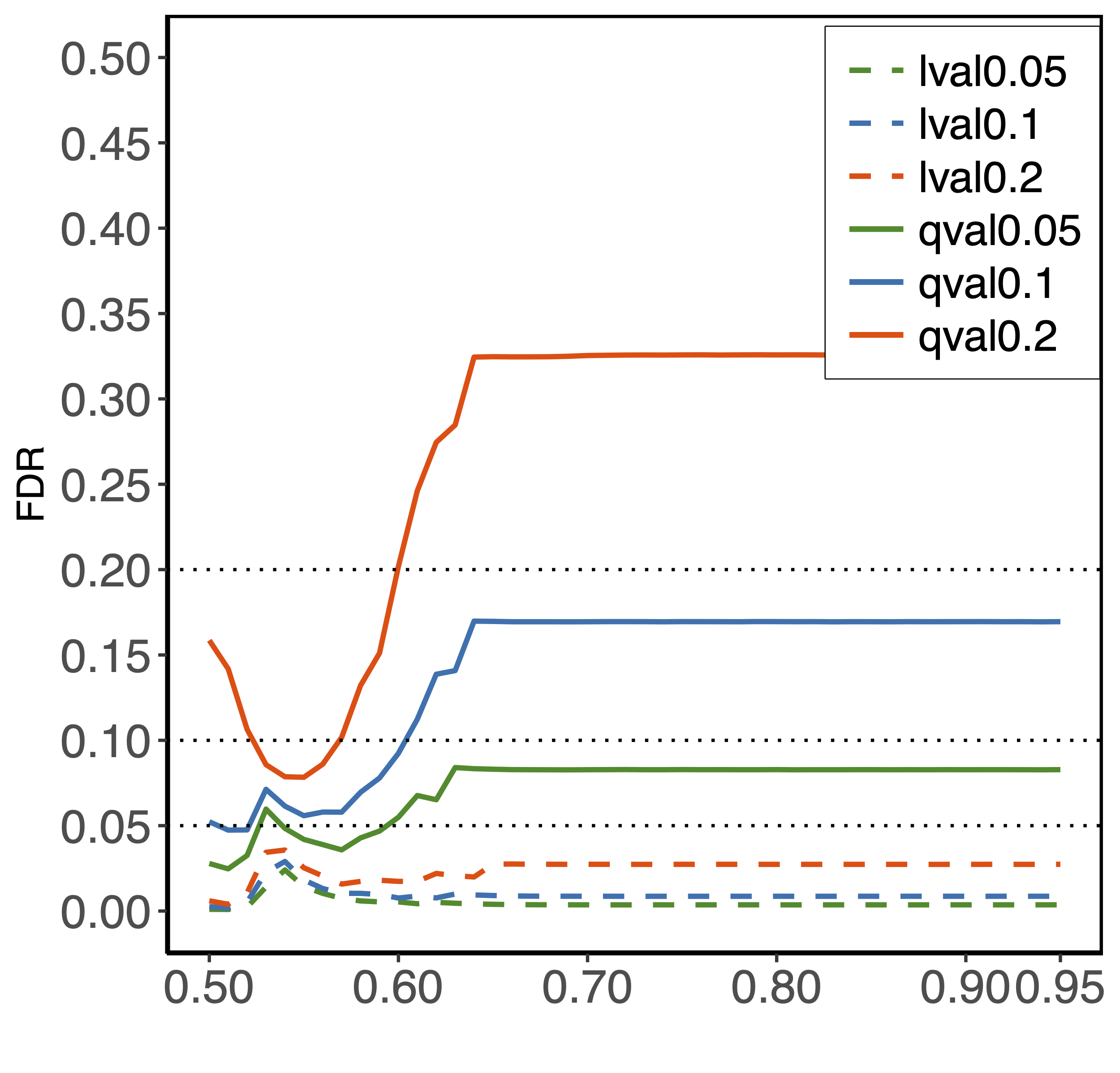}} 
    \subfigure[$m = 1,000$, $s_n/n = 0.001$]{\includegraphics[width=0.33\textwidth]{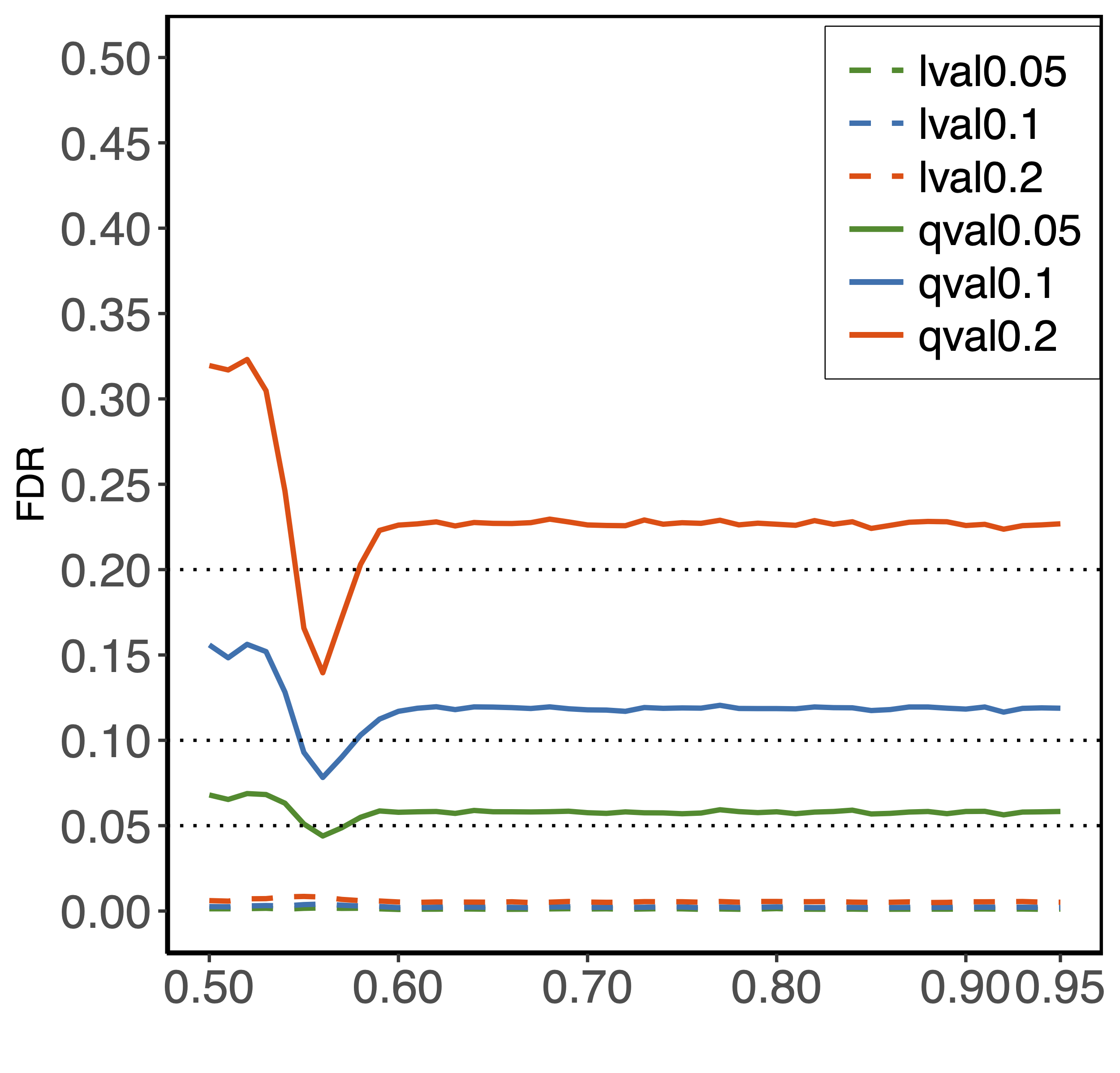}}%
     \subfigure[$m = 1,000$, $s_n/n = 0.1$]{\includegraphics[width=0.33\textwidth]{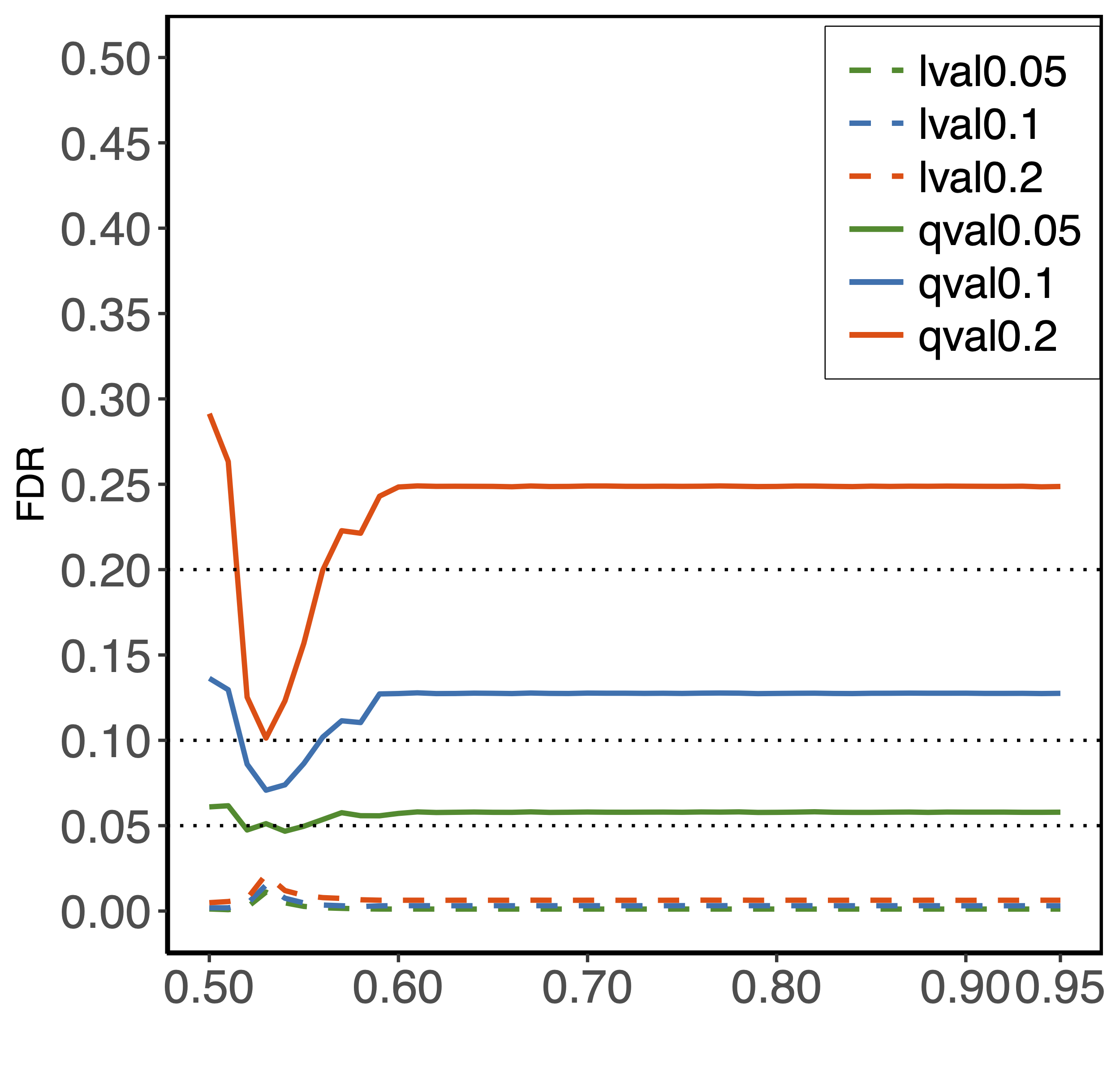}}%
    \subfigure[$m = 1,000$, $s_n/n = 0.5$]{\includegraphics[width=0.33\textwidth]{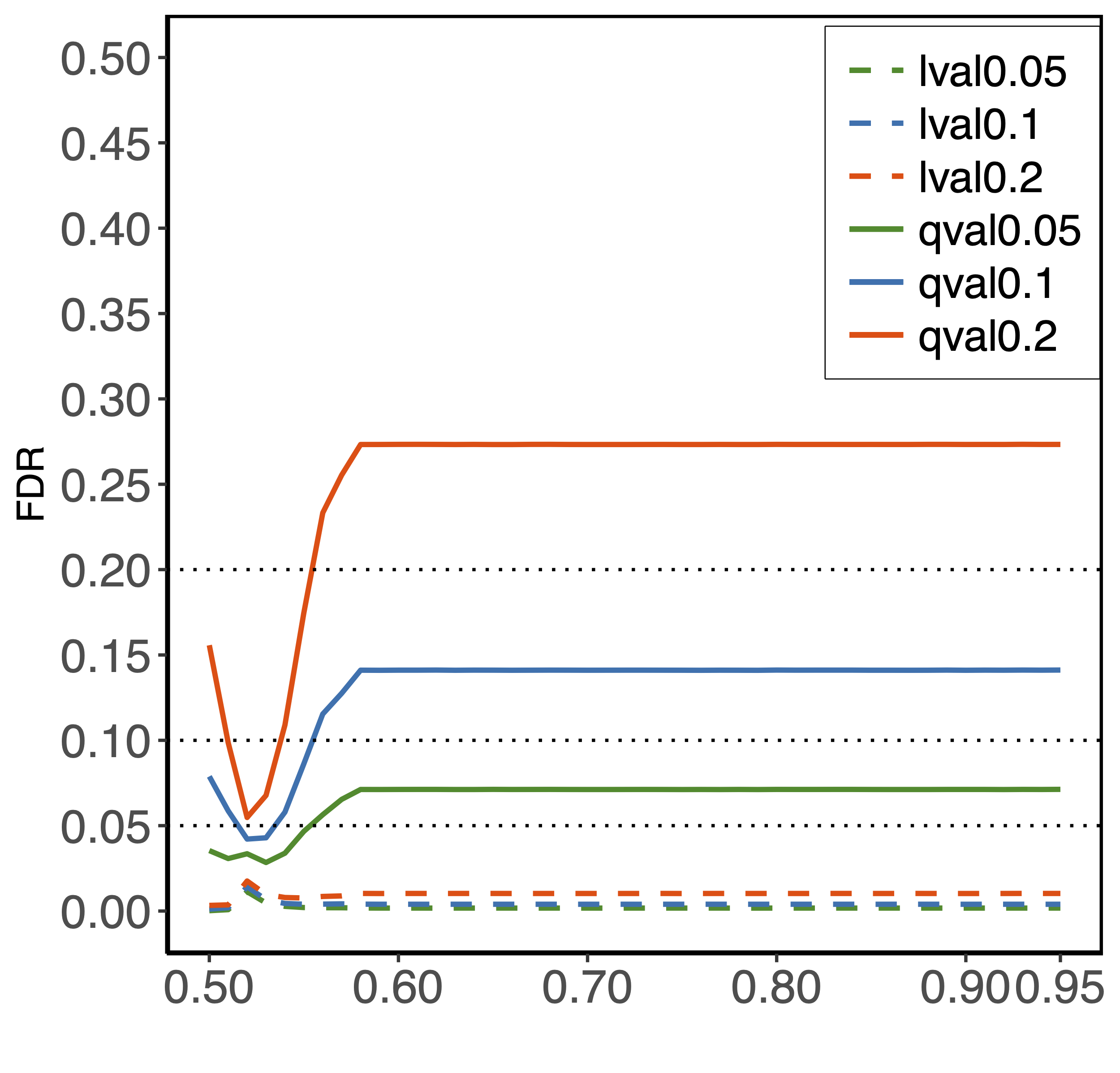}} 
	\caption{The estimated FDR of the $\ell$-value (dash) and the $q$-value (solid) procedures at $t = 0.05$ (blue), $t = 0.1$ (green), and $t = 0.2$ (red) with $m = (\log n)^2, 200$, and $1000$ and $s_n/n = 0.001, 0.1$, and $0.5$ respectively when choosing $\gamma \sim \text{Beta}(5, 5)$.}
	\label{fig:compare-q-l-beta5}
\end{figure}

\begin{figure}[h!]
\centering
	  \subfigure[$m = 85$, $s_n/n = 0.001$]{\includegraphics[width=0.33\textwidth]{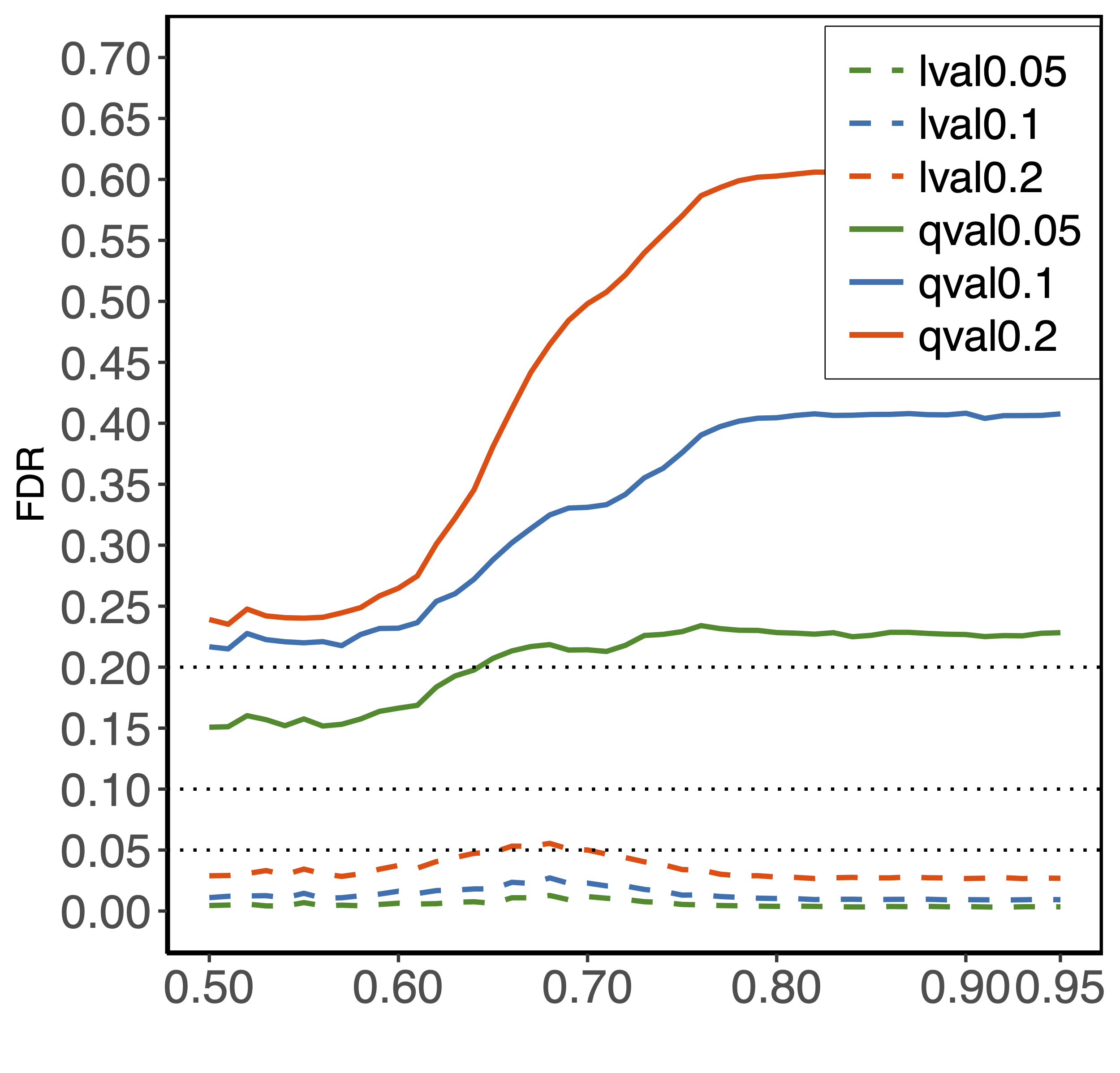}}%
     \subfigure[$m = 85$, $s_n/n = 0.1$]{\includegraphics[width=0.33\textwidth]{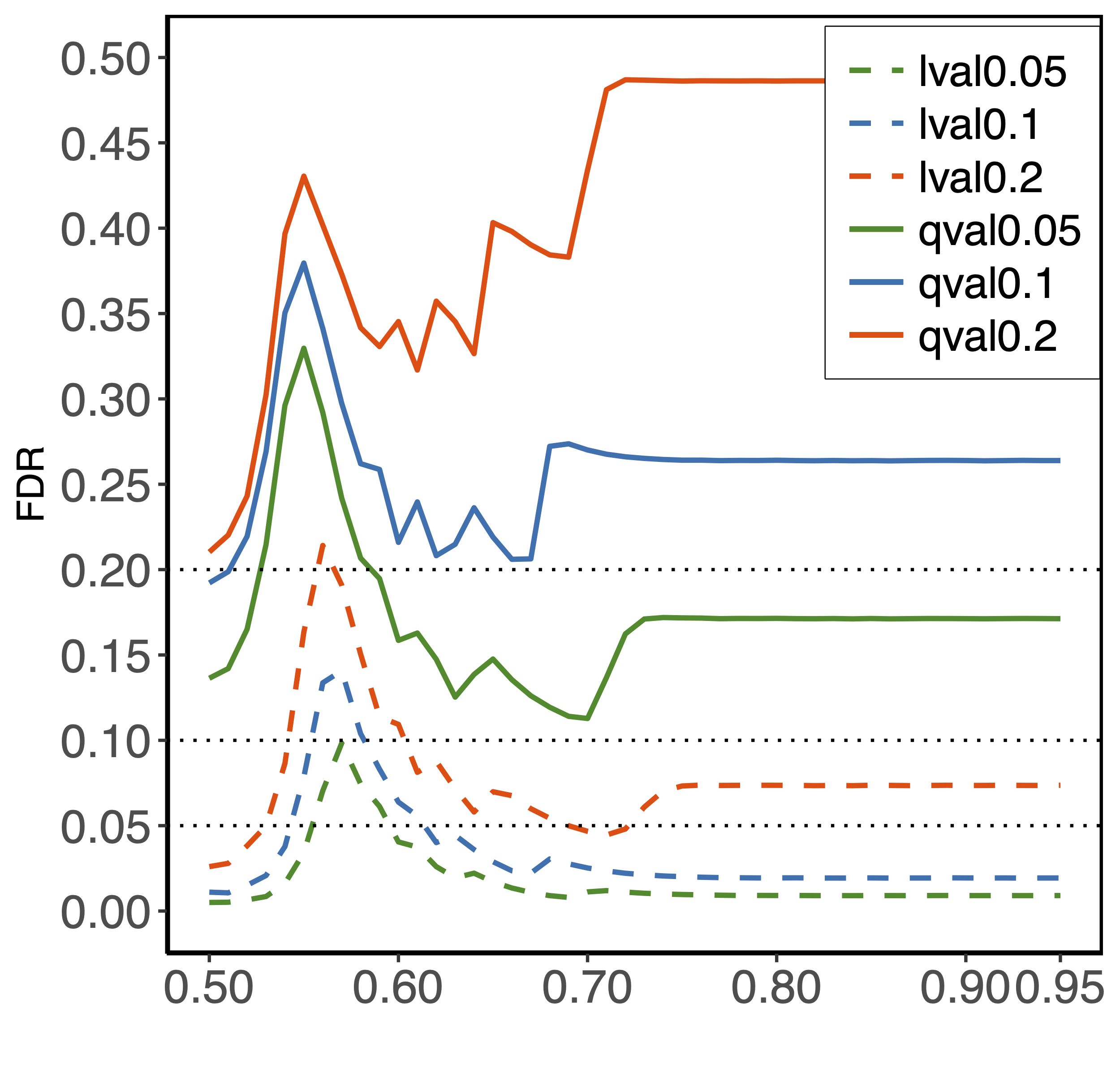}}%
    \subfigure[$m = 85$, $s_n/n = 0.5$]{\includegraphics[width=0.33\textwidth]{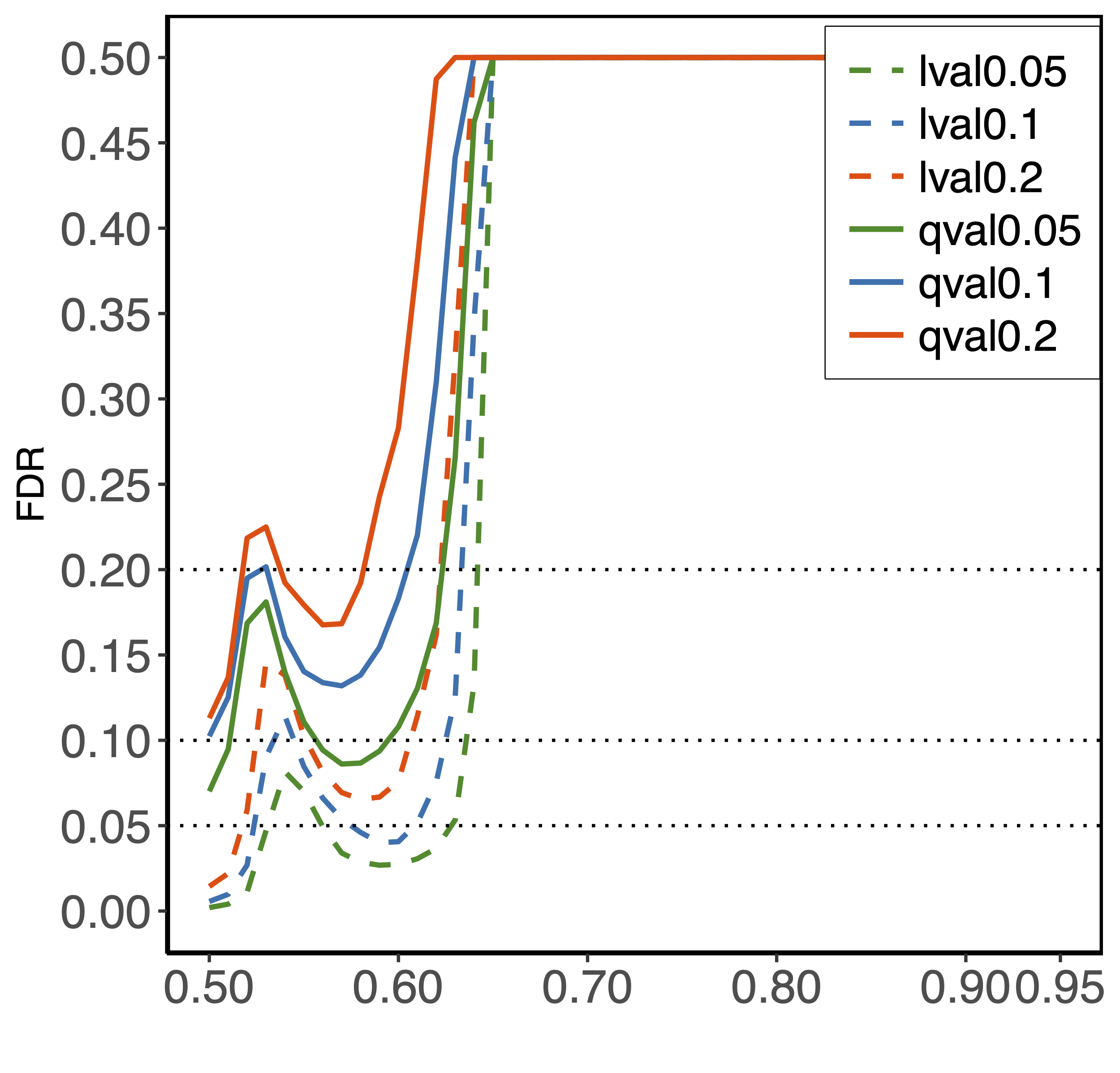}} 
    \subfigure[$m = 200$, $s_n/n = 0.001$]{\includegraphics[width=0.33\textwidth]{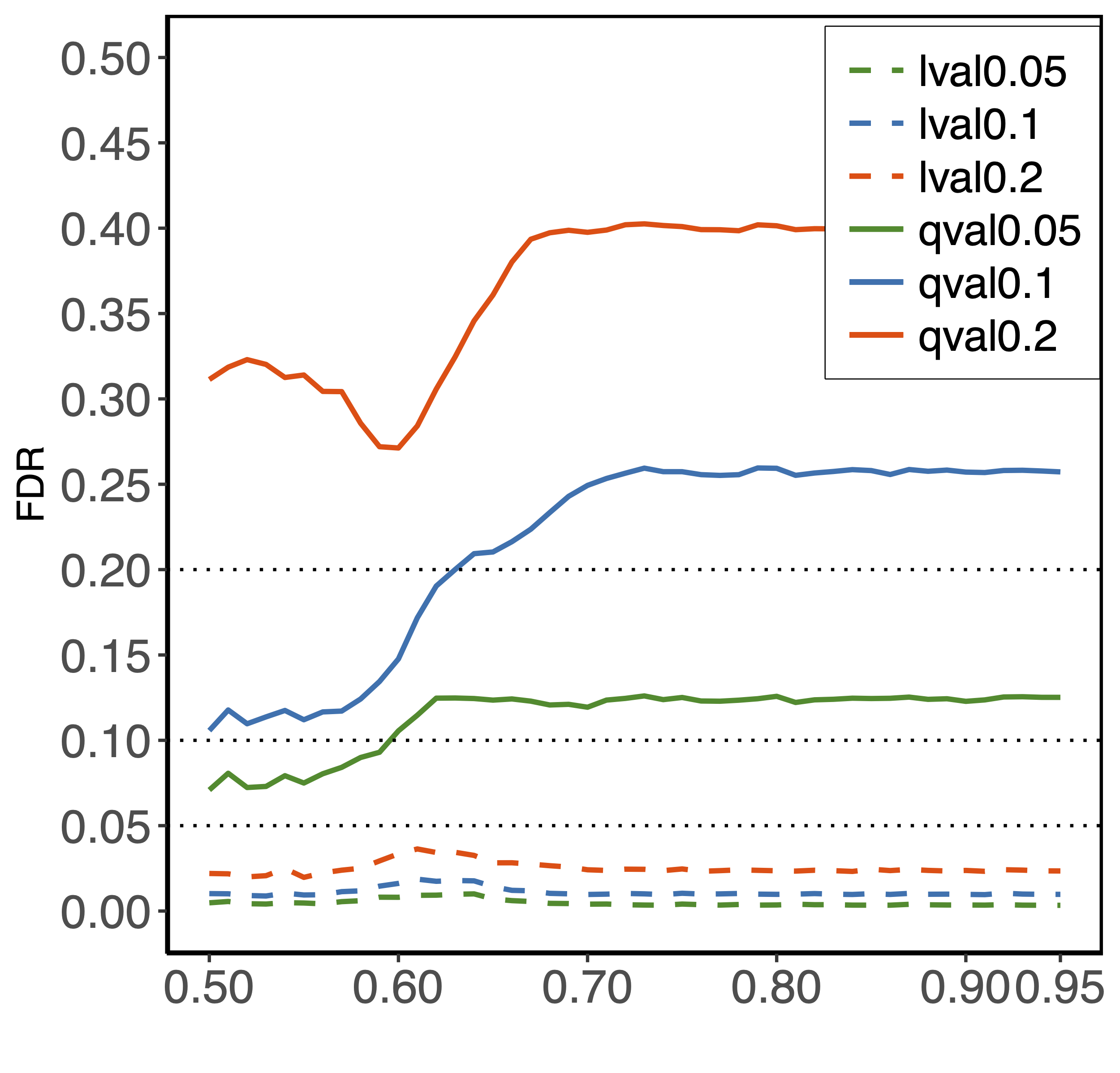}}%
     \subfigure[$m = 200$, $s_n/n = 0.1$]{\includegraphics[width=0.33\textwidth]{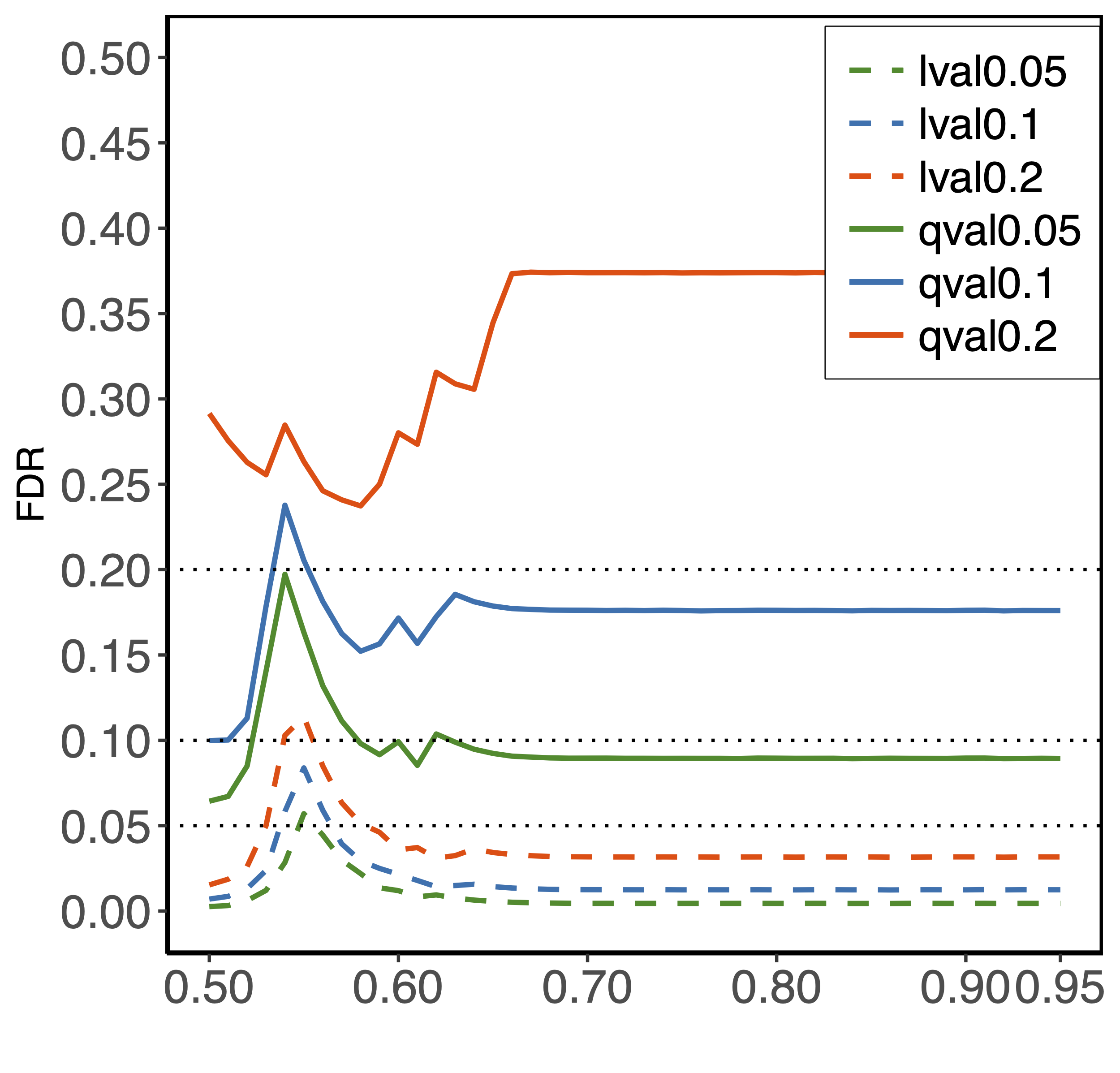}}%
    \subfigure[$m = 200$, $s_n/n = 0.5$]{\includegraphics[width=0.33\textwidth]{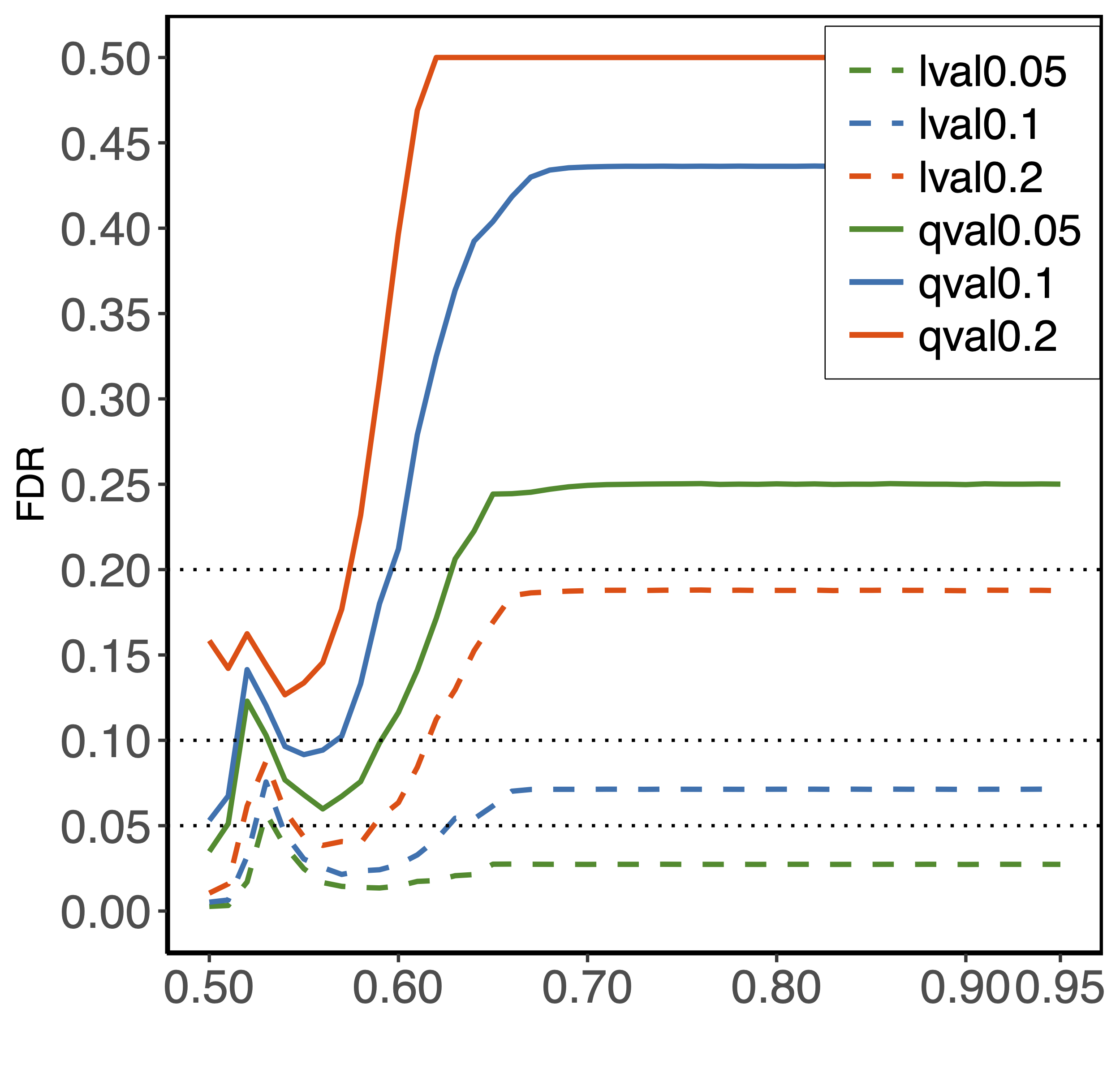}} 
    \subfigure[$m = 1,000$, $s_n/n = 0.001$]{\includegraphics[width=0.33\textwidth]{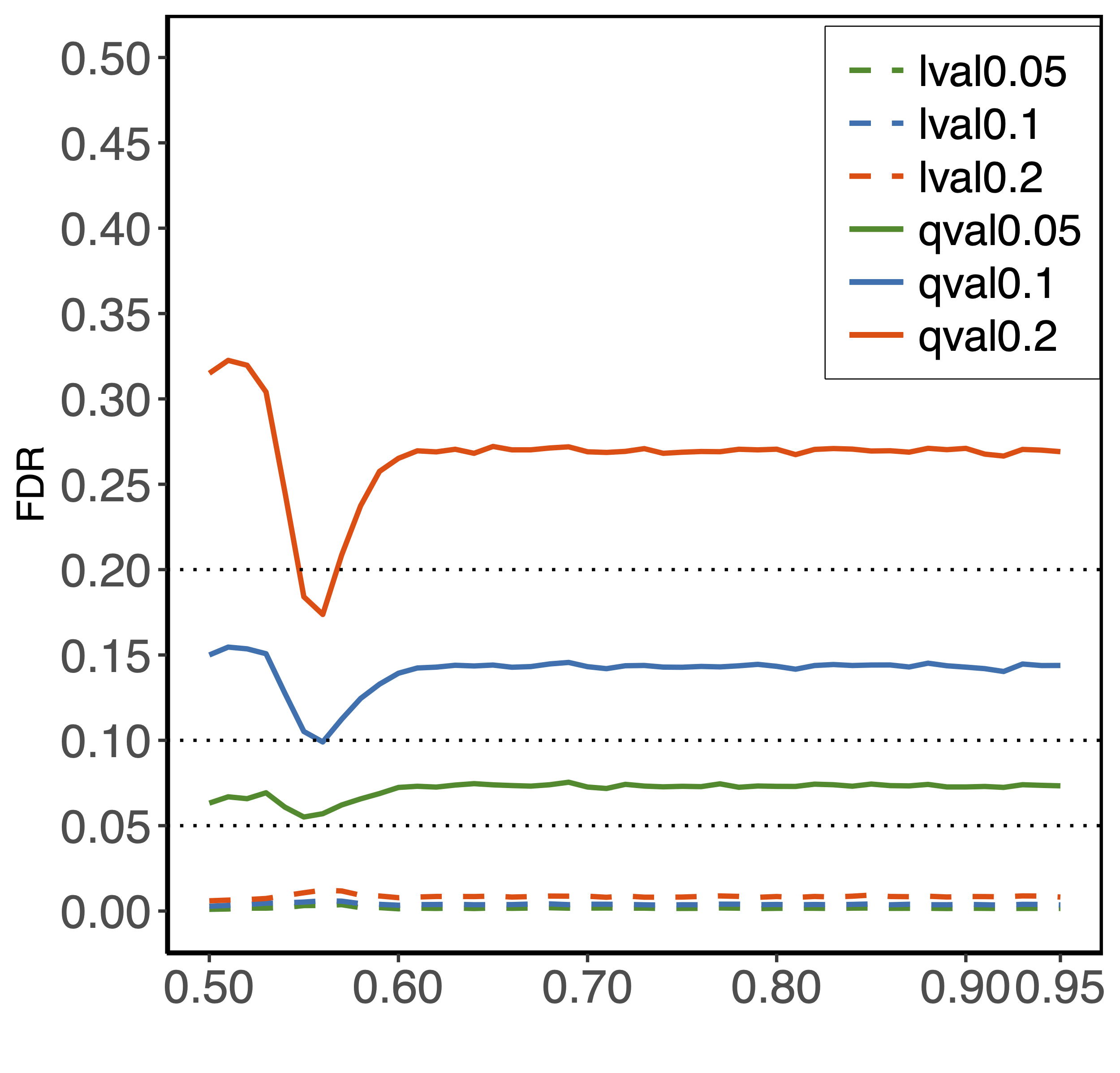}}%
     \subfigure[$m = 1,000$, $s_n/n = 0.1$]{\includegraphics[width=0.33\textwidth]{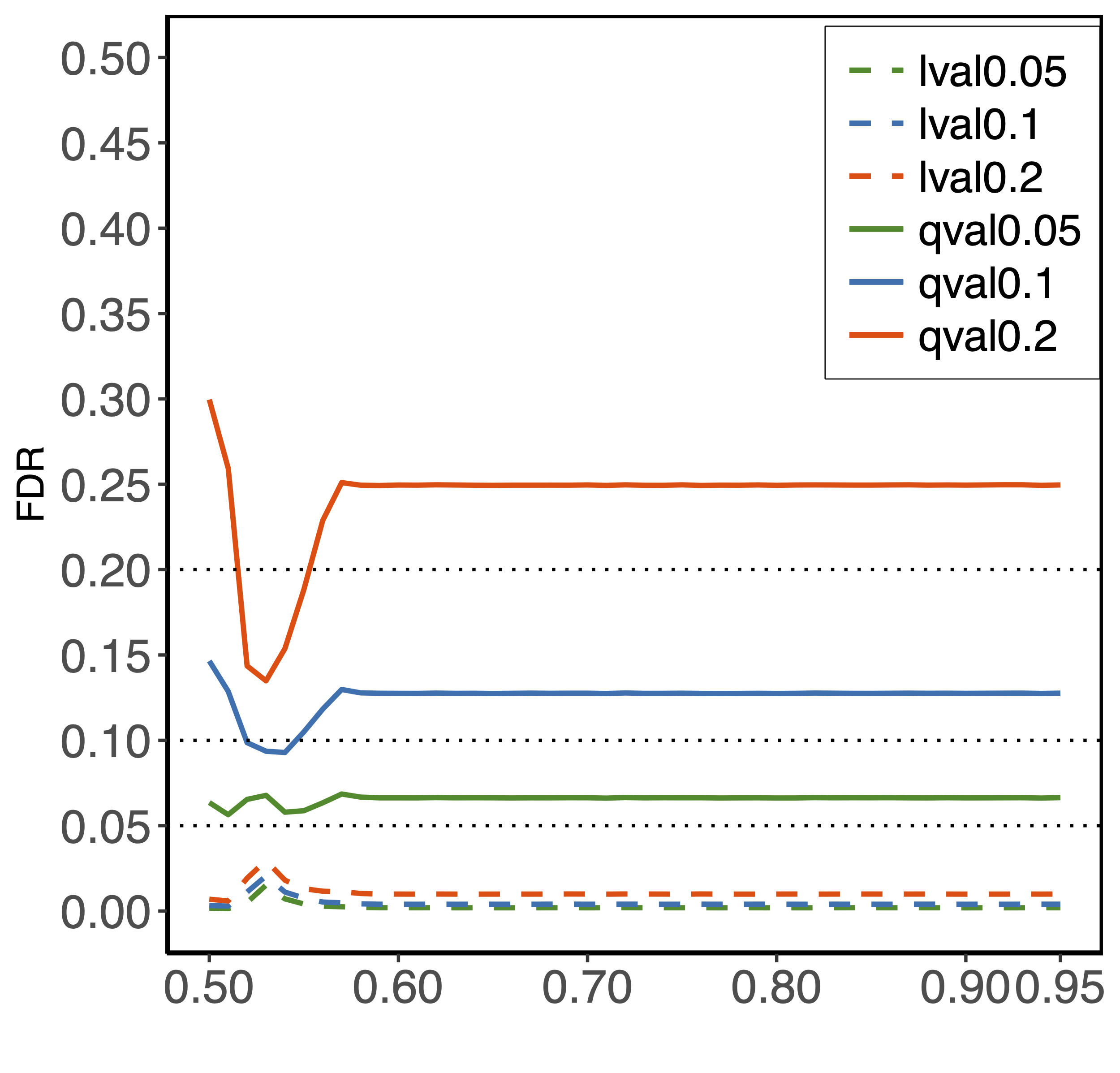}}%
    \subfigure[$m = 1,000$, $s_n/n = 0.5$]{\includegraphics[width=0.33\textwidth]{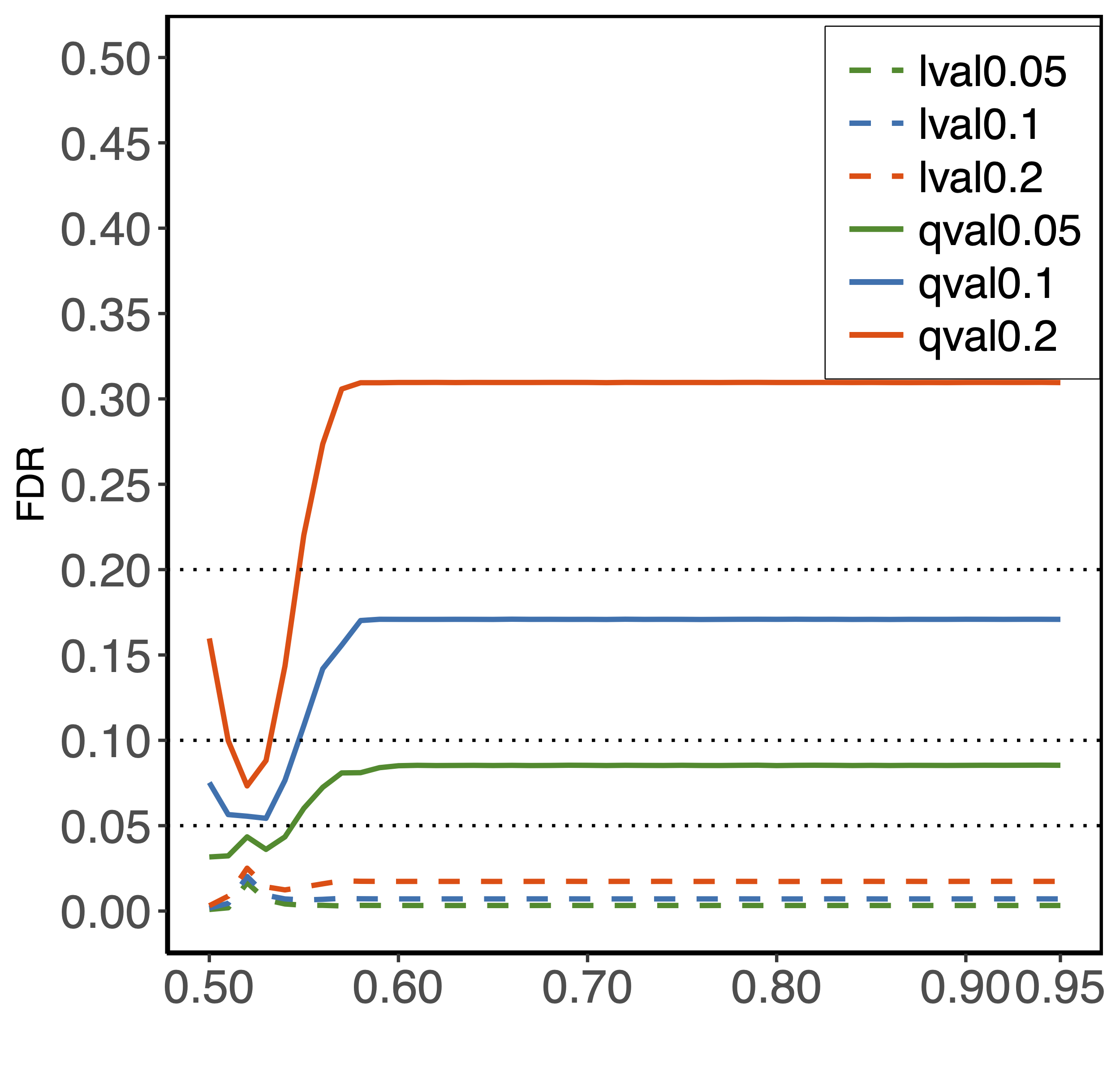}} 
	\caption{The estimated FDR of the $\ell$-value (dash) and the $q$-value (solid) procedures at $t = 0.05$ (blue), $t = 0.1$ (green), and $t = 0.2$ (red) with $m = (\log n)^2, 200$, and $1000$ and $s_n/n = 0.001, 0.1$, and $0.5$ respectively when choosing $\gamma \sim \text{Beta}(10, 10)$.}
	\label{fig:compare-q-l-beta10}
\end{figure}

In this section, we conduct two additional simulation studies to verify our conjecture that selecting a different parameter value for the beta prior $\text{Beta}(\alpha, \alpha)$ will not address the issue for the $\ell$-value procedure using the uniform prior $\text{Unif}[0, 1]$. 
The data generation process remains the same as it in Section \ref{sec:sim-1}. Here, we consider two different priors $\text{Beta}(5, 5)$ and $\text{Beta}(10, 10)$ for $\gamma$. Simulation results are presented in Figures \ref{fig:compare-q-l-beta5} and \ref{fig:compare-q-l-beta10} respectively. 
Upon comparing each subplot in the these two figures, we indeed observe that the FDR of $\ell$-value procedure does not increase. In both figures, for all nine scenarios, the FDR values are consistently small, close to zero. An exception is found in (c) and (f) when $s_n/n = 0.5$ and $m$ is relatively small, where the FDR tends to be very large. However, as $m$ increases, as seen in (i), the FDR returns to values close to zero, way below the target level. On the other hand, the $q$-value procedure did not perform as well as in the case when using the uniform prior either. Changing the prior appears to cause the $q$-value procedure to estimate the FDR significantly higher than the target level in all nine scenarios. 
In sum, we provide two examples that demonstrate that, in general, choosing $\alpha > 1$ in the $\text{Beta}(\alpha, \alpha)$ prior does not lead to an improvement in FDR control.

\bibliographystyle{chicago}
\bibliography{citation.bib}
\stopcontents

\end{document}